\documentclass[11pt,a4paper]{article}
\usepackage[a4paper,margin=25mm,headheight=14pt,headsep=8mm,footskip=12mm]{geometry}
\usepackage[T1]{fontenc}
\usepackage{lmodern}
\usepackage{microtype}
\usepackage{authblk}
\usepackage{amsmath,amssymb,amsfonts,mathtools}
\usepackage{amsthm}
\usepackage[round,authoryear]{natbib}
\usepackage{bm,mathrsfs}
\usepackage{graphicx,booktabs,xcolor,tcolorbox,multicol,enumitem}
\usepackage{algorithm}
\usepackage{algpseudocode}
\usepackage{float}
\usepackage{needspace}
\usepackage{placeins}
\usepackage{flafter}
\usepackage{subcaption}
\usepackage{multirow}
\usepackage{xurl}
\usepackage[
    colorlinks=false,
    pdfborder={0 0 1},
    linkbordercolor={0 1 0},
    citebordercolor={0 1 0},
    urlbordercolor={0 1 0}
]{hyperref}

\definecolor{abstractbg}{RGB}{238,238,238}
\definecolor{headerrule}{RGB}{180,180,180}

\AtBeginDocument{%
  \setlength{\floatsep}{5pt plus 2pt minus 1pt}%
  \setlength{\textfloatsep}{7pt plus 2pt minus 1pt}%
  \setlength{\intextsep}{6pt plus 2pt minus 1pt}%
  \addtolength{\abovedisplayskip}{-1.5pt}%
  \addtolength{\belowdisplayskip}{-1.5pt}%
  \setlength{\emergencystretch}{2em}%
  \renewcommand{\topfraction}{0.95}%
  \renewcommand{\bottomfraction}{0.90}%
  \renewcommand{\textfraction}{0.05}%
  \renewcommand{\floatpagefraction}{0.75}%
  \setcounter{topnumber}{3}%
  \setcounter{bottomnumber}{2}%
  \setcounter{totalnumber}{4}%
  \captionsetup{skip=3pt}%
  \raggedbottom
}

\theoremstyle{plain}
\newtheorem{theorem}{Theorem}[section]
\newtheorem{lemma}[theorem]{Lemma}
\newtheorem{proposition}[theorem]{Proposition}
\newtheorem{corollary}[theorem]{Corollary}
\theoremstyle{definition}
\newtheorem{definition}[theorem]{Definition}
\newtheorem{assumption}[theorem]{Assumption}
\theoremstyle{definition}
\newtheorem{remark}[theorem]{Remark}

\newcommand{\R}{\mathbb{R}}

\newcommand{\M}{\mathcal{M}}

\newcommand{\dH}{d_{H}}

\hypersetup{
  pdftitle={Beyond Curvature Bias: Higher-Order Manifold Fitting via Multiscale Cancellation},
  pdfauthor={Zhijun Zhu, Haoyu Yang, Ruoyu Xiong, Bingjie Li, Zhigang Yao},
  pdfsubject={Manifold fitting and geometric bias cancellation},
  bookmarksnumbered=true,
  pdflinkmargin=2.5pt
}

\title{Beyond Curvature Bias: Higher-Order Manifold Fitting via Multiscale Cancellation}

\renewcommand{\Affilfont}{\small}

\makeatletter
\renewcommand{\AB@affilsepx}{%
  \protect\\[4pt]\protect\Affilfont}
\makeatother

\author[1]{Zhijun Zhu}
\author[2]{Haoyu Yang}
\author[3]{Ruoyu Xiong}
\author[4]{Bingjie Li}
\author[5]{Zhigang Yao}

\makeatletter
\g@addto@macro\AB@authlist{%
  \textsuperscript{,}%
  \thanks{Corresponding author:
    \texttt{zhigang.yao@nus.edu.sg}.}%
}
\makeatother

\affil[1]{School of Mathematical Sciences,
  Fudan University, Shanghai, China}

\affil[2]{Department of Mathematics,
  National University of Singapore, Singapore}

\affil[3]{School of Mathematics and Physics,
  Xi'an Jiaotong-Liverpool University, Suzhou, China}

\affil[4]{Department of Epidemiology and Biostatistics,
  School of Public Health, Tongji Medical College,
  Huazhong University of Science and Technology,
  Wuhan, China}

\affil[5]{Department of Statistics \& Data Science,
  National University of Singapore, Singapore}

\date{}
\begin{document}
\maketitle
\begin{abstract}
Local manifold-fitting procedures smooth noisy observations
 over neighborhoods of scale $h$, inducing a systematic
normal bias of order $h^2+\sigma^2$ that persists at the
population level.
We develop geometric bias cancellation (GBC),
a multiscale framework for higher-order manifold
recovery under ambient Gaussian noise.
For procedures admitting a
normal-bias expansion $c(h,\sigma)b(x)+O(h^4)$,
we construct cross-scale correspondence using estimated
normal fibers and show that second-order normal accuracy
suffices for fourth-order correspondence accuracy.
Canceling the leading geometric bias yields $O(h^4)$
population Hausdorff accuracy while preserving an
embedded manifold diffeomorphic to the underlying manifold.
We establish the required expansions and finite-sample
recovery guarantees for successive tangent-space projection
and Gaussian density ridges. For STSP, geometric bias cancellation also achieves an
asymptotically vanishing Hausdorff error ratio
relative to oracle bandwidth-tuned single-scale fitting.
For the local contraction estimator
of \citet{yao2023manifold}, we identify an obstruction to
scalar bias cancellation and develop a directional cancellation
achieving $O(\sigma^4)$ population Hausdorff accuracy with
uniform positive reach, as well as fourth-order finite-sample
accuracy when $n\gtrsim\sigma^{-(d+6)}\log(1/\sigma)$.
Simulations confirm the predicted improvement
from second- to fourth-order population bias.
\end{abstract}

\section{Introduction}
\label{sec:intro}

Let $\M\subset\R^D$ be an unknown compact
$d$-dimensional embedded submanifold, and
suppose we observe
\begin{equation}
\label{eq:model}
Y_i=X_i+\xi_i,\qquad i=1,\ldots,n,
\end{equation}
where $X_1,\ldots, X_n$ are i.i.d. with
distribution $\omega$ supported on $\M$, and
$\xi_1,\ldots,\xi_n$ are i.i.d. with distribution
$\phi_\sigma$,
independent of the latent sample.
The goal of manifold fitting is to construct
a $d$-dimensional embedded manifold
$\widehat{\M}\subset\R^D$ that approximates
$\M$ in Hausdorff distance.

\subsection{The Bias--Variance Barrier in Manifold Fitting}
\label{intro::barrier}

Many local manifold-fitting procedures depend on a fitting scale,
such as a neighborhood radius or kernel bandwidth.
Larger scales stabilize estimation but introduce systematic bias,
whereas smaller scales reduce bias at the cost of increased
sampling variability \citep{aamari2019nonasymptotic}.
However, local averaging over a curved manifold induces a
population-level normal bias that persists even with infinitely
many observations. Increasing the sample size cannot remove
this bias at a fixed scale, while bandwidth tuning only
moves along the bias--variance trade-off.
Such geometric bias is well documented in principal curves
and density ridges
\citep{hastie1989principal,tibshirani1992principal,
chen2015asymptotic}.
In noisy manifold fitting, the leading bias can depend
jointly on the fitting scale, noise level, and local
curvature, as illustrated by the local contraction
estimator of \citet[Section~4.3]{yao2023manifold}.
Following the terminology of \citet{yao2026review},
we refer to this method as the YY estimator.
For example, the population normal error
of a Gaussian density ridge admits the expansion
\[
u_{h,\sigma}^{\mathrm{DR}}(x)
=
\frac{d}{2}(h^2+\sigma^2)H(x)+O(h^4),
\qquad \sigma\lesssim h,
\]
where $H(x)$ is the mean-curvature vector.

\subsection{Geometric Bias Cancellation: The Multiscale Principle}
\label{intro::multi-bias}
In nonparametric statistics, structured deterministic bias
can often be canceled rather than controlled solely through
bandwidth selection. When the leading bias separates into
scale-dependent and spatial components, combining estimators
across bandwidths can eliminate the leading term.
This principle underlies higher-order kernels
\citep{bartlett1963statistical}, Tukey's twicing,
the generalized jackknife of \citet{schucany1977improvement},
multiplicative bias correction \citep{jones1995simple},
and bias reduction in local polynomial smoothing
\citep{fan1993local,ruppert1994multivariate}.

Applying these bias-correction ideas to manifold fitting
under ambient Gaussian noise introduces additional
geometric difficulties. Unlike classical regression,
where the covariate coordinates are prescribed, manifold
fitting requires local tangent spaces and coordinate
systems to be inferred from noisy observations.
Classical errors-in-variables methods
\citep{fan1993nonparametric,delaigle2009design}
do not directly address this uncertainty in the
underlying geometry.
In particular, while higher-order local polynomial
methods can reduce approximation bias
\citep{fan1993local,ruppert1994multivariate},
\citet{aizenbud2025localgeo} showed that errors in
tangent-space estimation can introduce additional
systematic bias in manifold fitting. Thus, increasing
the polynomial degree alone does not guarantee
higher-order manifold recovery or preservation of
the embedded-manifold structure.

We develop geometric bias cancellation (GBC),
which exploits the scale dependence of the
population normal bias to achieve higher-order
manifold recovery through multiscale cancellation.
\citet{xia2026stsp} introduced successive tangent-space
projection (STSP) and its multiscale extension, MS-STSP,
using squared-bandwidth extrapolation to reduce
curvature-induced bias. Under Gaussian noise, their
analysis establishes second-order localization of the
population fixed-point manifold, together with
finite-sample localization guarantees.
However, their analysis does not establish higher-order recovery under Gaussian noise or guarantee that the extrapolated estimator remains a smooth embedded manifold diffeomorphic to the underlying manifold. 
GBC addresses these limitations by canceling the
actual noise-dependent leading normal bias, rather
than relying on squared-bandwidth extrapolation,
and by constructing cross-scale correspondence
through estimated normal fibers. This yields
fourth-order Hausdorff recovery under Gaussian noise
while preserving a smooth embedded manifold
diffeomorphic to the underlying manifold.

We now describe the basic scalar mechanism underlying our method. Let
\(\mathcal M_h\) denote the population fitted manifold at scale \(h\). For
each \(x\in\mathcal M\), suppose that \(\mathcal M_h\) can be represented
locally as a normal graph over \(\mathcal M\), and let
\(X_h(x)\in \mathcal M_h\cap (x+N_x\mathcal M)\)
denote the corresponding point above \(x\). We assume that the normal-direction
error admits the factorized expansion
\begin{equation}
\label{eq:factorization}
X_h(x)-x
=
c(h,\sigma)b(x)+O(h^4),
\qquad x\in\mathcal M,
\end{equation}
where \(b:\mathcal M\to\mathbb R^D\) is a scale-independent geometric bias
field and \(c(h,\sigma)\) is a scalar coefficient satisfying
\(c(h,\sigma)=O(h^2+\sigma^2)\).
The key feature of \eqref{eq:factorization} is that
the leading bias shares a common geometric field
across scales, while its coefficient retains the joint
dependence on $h$ and $\sigma$ needed for cancellation.

Let \(h_1<\cdots<h_K\) be fitting scales and choose weights
\(w_1,\ldots,w_K\) satisfying
\begin{equation}
\label{eq:weights}
\sum_{k=1}^{K} w_k = 1,
\qquad
\sum_{k=1}^{K} w_k c(h_k,\sigma)=0.
\end{equation}
If the fitted manifolds are matched through their common base point
\(x\in\mathcal M\), define
\(\widetilde X(x) = \sum_{k=1}^{K} w_k X_{h_k}(x)\).
Then \eqref{eq:factorization} and \eqref{eq:weights} give
\[
\widetilde X(x)-x
=
O\!\left(
\sum_{k=1}^{K}|w_k|h_k^4
\right).
\]

\begin{figure}[!htbp]
    \centering
    \includegraphics[width=\linewidth]{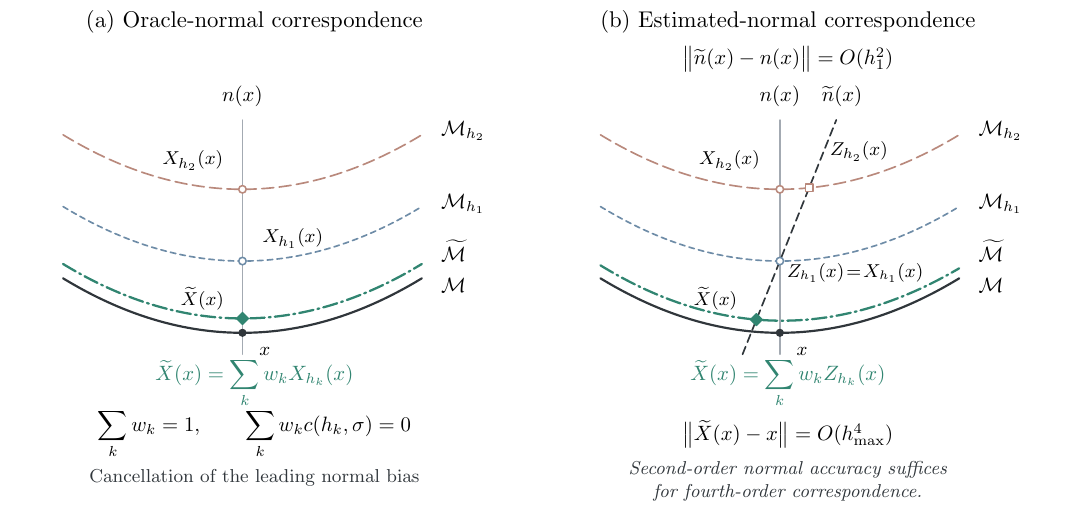}
\caption{
Geometric intuition behind geometric bias cancellation (GBC).
\textbf{(a) Oracle-normal correspondence.}
If the true normal directions were known, fitted manifolds
at different scales could be matched along common normal
fibers. Their leading geometric biases share the same
spatial structure but differ in their scale dependence.
A suitable weighted combination therefore cancels the
leading second-order bias, leaving a fourth-order error.
\textbf{(b) Estimated-normal correspondence.}
Since the true normal directions are unknown in practice,
GBC instead constructs cross-scale correspondence using
normal directions estimated from a reference fitted manifold.
The key insight is that second-order accuracy of these
estimated normals is sufficient to achieve fourth-order
correspondence accuracy. Consequently, the leading
geometric bias can still be canceled, yielding fourth-order
manifold recovery without oracle geometric information.
}
    \label{fig:multiscale-workflow}
\end{figure}

Thus the leading second-order normal bias is canceled,
leaving a fourth-order population remainder,
as illustrated in the left panel of
Figure~\ref{fig:multiscale-workflow}.
Hausdorff accuracy alone does not ensure geometric regularity.
Positive reach guarantees a tubular neighborhood on which the
nearest-point projection is uniquely defined and stable
\citep{federer1959curvature,niyogi2008finding}.
However, a small Hausdorff error and a smooth embedded image do not
by themselves provide a uniform positive lower bound on reach.
We therefore distinguish preservation of the embedded-manifold
structure from quantitative reach guarantees.
The general framework in this paper establishes the former, while the last example
also provides the latter through additional second-order estimates.

The construction above assumes oracle correspondence
along the true normal fibers of $\mathcal M$.
Since these normal spaces are unknown, GBC instead
constructs cross-scale correspondence using normal
fibers estimated from the reference manifold
$\mathcal M_{h_1}$, as illustrated in the right panel of
Figure~\ref{fig:multiscale-workflow}.
The key geometric insight is that second-order
normal estimation is sufficient to preserve
fourth-order bias cancellation.
Indeed, an $O(h_1^2)$ error in the estimated normal
projector, combined with an $O(h_k^2)$ displacement
between fitted manifolds, induces only an
$O(h_1^2h_k^2)$ correspondence error.
For comparable scales $h_k\asymp h$ and bounded
cancellation weights, this error is of fourth order.
Consequently, the leading geometric bias can still
be canceled, yielding
\[
\sup_{x\in\mathcal M}
\|\widetilde X(x)-x\|=O(h^4).
\]
Thus, GBC achieves fourth-order manifold recovery
through estimated-normal correspondence without
oracle knowledge of the underlying geometry.

\begin{figure}[!htbp]
\centering
\includegraphics[width=\linewidth]{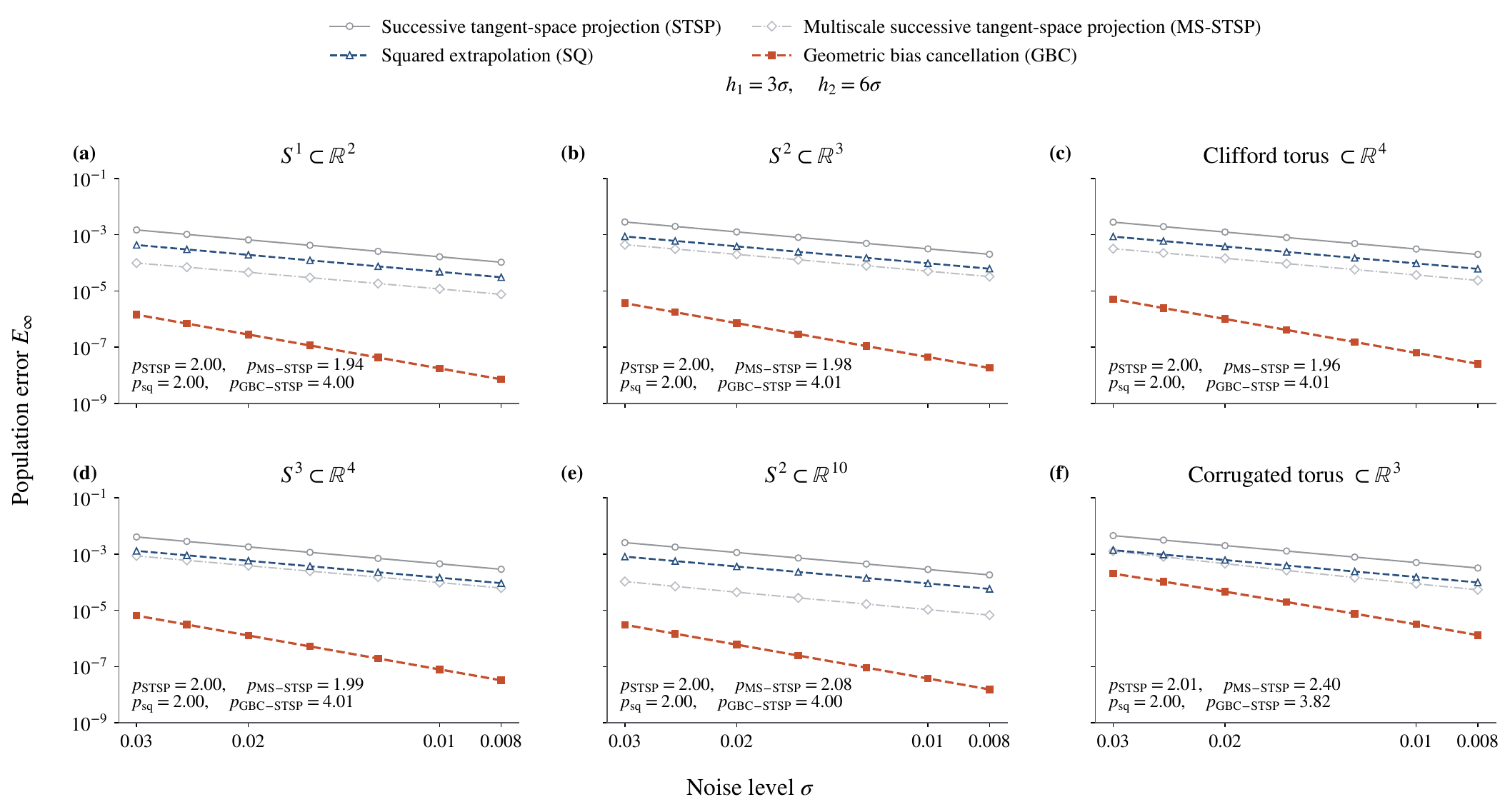}
\caption{
Population maximum errors for single-scale STSP,
naive squared-bandwidth extrapolation,
the proposed coefficient-based extrapolation,
and the original MS-STSP of Xia, Li and Yao (2026)
across six manifolds.
The first three methods use the same smooth-kernel
STSP fits; the two extrapolation methods also share
the same estimated-normal correspondences and differ
only in their weights.
Naive extrapolation uses $\sum_k w_k h_k^2=0$,
whereas our method uses
$\sum_k w_k c_{\mathrm{STSP}}(h_k,\sigma)=0$.
Original MS-STSP uses hard-ball local moments and
serves as a separate benchmark.
The fitted log--log slopes are reported in the
legend order for each panel.
}
\label{fig:intro_stsp_population}
\end{figure}

The role of the noise-dependent bias coefficient is
illustrated in Figure~\ref{fig:intro_stsp_population}.
Across six manifolds, conventional squared-bandwidth
extrapolation retains approximately second-order
population error, even when applied to the same
smooth-kernel STSP fits and estimated-normal
correspondences as our method.
In contrast, extrapolation using the correct
noise-dependent coefficient achieves approximately
fourth-order convergence.
The original MS-STSP construction of
\citet{xia2026stsp}, included as a separate
benchmark, also exhibits predominantly
second-order behavior in these experiments.
These results highlight that multiscale combination
alone is insufficient: higher-order recovery requires
cancellation of the actual leading geometric bias.

The scalar formulation also has a structural limitation.
It applies when the leading normal bias shares a common
scale dependence across directions. For the smooth-kernel YY estimator studied in Section~\ref{sec:yao_contraction}, the two leading curvature components
in higher codimension have scale-dependent coefficients
$\sigma^2(1+\alpha_k)$ and $\sigma^2\alpha_k$.
Canceling both components with scalar weights would require
\(\sum_k w_k(1+\alpha_k)=0, \sum_k w_k\alpha_k=0\).
Subtracting these conditions gives $\sum_k w_k=0$,
contradicting the affine condition $\sum_k w_k=1$.
Thus scalar extrapolation cannot, in general, remove
both leading components, even with exact correspondence.
This motivates the directional extension developed in
Section~\ref{sec:yao_contraction}, which assigns separate
weights to different normal components and achieves
fourth-order population Hausdorff accuracy under the
stated conditions.

Table~\ref{tab:accuracy_sample_size} contrasts published
single-scale guarantees with our fourth-order multiscale
results. For example, \citet{yao2023manifold} achieve
$O(\sigma^2\log(1/\sigma))$ accuracy with
$n\asymp\sigma^{-(d+3)}$, whereas our GBC-YY achieves $O(\sigma^4)$ accuracy under
$n\gtrsim W^2\sigma^{-(d+6)}\log(1/\sigma)$
and the stated preliminary-manifold conditions.
These results concern different accuracy targets
rather than a direct sample-efficiency comparison.

\Needspace{2.5\baselineskip}
\subsection{Main contributions}
Our main contribution is geometric bias cancellation
(GBC), a framework for higher-order manifold
recovery under Gaussian noise.
Beyond removing the leading geometric bias,
GBC addresses the cross-scale correspondence
and stochastic stability required to achieve
higher-order accuracy in finite samples.
Our contributions are as follows.

\emph{Geometric correspondence and debiasing.}
We establish a general higher-order recovery
principle for manifold estimators admitting
factorized normal-bias expansions.
Since fitted manifolds at different scales do not
share an observable parameterization, we develop
estimated-normal-fiber correspondence and prove
that second-order normal accuracy suffices for
fourth-order correspondence accuracy
(Proposition~\ref{thm:normal_guided_manifold_recovery}).
The resulting extrapolated manifold achieves
fourth-order population Hausdorff accuracy
while preserving its embedded-manifold structure.

\emph{Higher-order manifold recovery and statistical guarantees.}
We derive explicit noise-dependent normal-bias
expansions for STSP and Gaussian density ridges,
and establish finite-sample recovery guarantees.
For comparable scales $h\asymp\sigma$, the
GBC estimators attain
\[
d_H(\widehat{\mathcal M},\mathcal M)
=
O_p\left(
h^4+h\sqrt{\frac{\log n}{nh^d}}
\right),
\]
with fourth-order accuracy under the sufficient
sample-size condition
$n\gtrsim\sigma^{-(d+6)}\log(1/\sigma)$.
More importantly, for STSP, we prove that
the GBC estimator achieves an
asymptotically vanishing Hausdorff error ratio
relative to the oracle bandwidth-tuned
single-scale estimator over $h\asymp\sigma$
(Corollary~\ref{cor:stsp_oracle_separation}).
This separation holds under the weaker condition
$n\sigma^{d+2}/\log n\to\infty$,
demonstrating a statistical advantage of
bias cancellation beyond bandwidth selection.

\emph{Directional geometric bias cancellation.}
For the YY estimator,
we derive a directional decomposition of the leading
normal bias and identify an obstruction to its
cancellation by scalar extrapolation. This motivates
a componentwise multiscale construction that achieves
fourth-order population Hausdorff accuracy while
retaining a uniform positive-reach guarantee
(Propositions~\ref{thm:ysl_common_raw_expansion} and~\ref{thm:ysl_common_population}). We further establish finite-sample refinement bounds
for the proposed GBC-YY estimator.

\begin{table}[!htbp]
\centering
\caption{Single-scale benchmarks and GBC recovery guarantees.}
\label{tab:accuracy_sample_size}

\footnotesize
\setlength{\tabcolsep}{3pt}
\renewcommand{\arraystretch}{1.2}

\begin{tabular*}{\linewidth}
{@{\extracolsep{\fill}}llcc@{}}
\toprule
Estimator & Construction & Sample condition & Error bound \\
\midrule

\multirow{2}{*}{STSP}
& Single-scale
& $n\gtrsim\sigma^{-(d+2)}L_\sigma$
& $O(\sigma^2)^\dagger$ \\

& Geometric bias cancellation
& $n\gtrsim\sigma^{-(d+6)}L_\sigma$
& $O(\sigma^4)$ \\

\addlinespace[5pt]

\multirow{2}{*}{DR}
& Single-scale
& $n\gtrsim\sigma^{-(D+8)}L_\sigma^{-(D+6)/2}$
& $O_p(\sigma^2L_\sigma)^\ddagger$ \\

& Geometric bias cancellation
& $n\gtrsim\sigma^{-(d+6)}L_\sigma$
& $O(\sigma^4)$ \\

\addlinespace[5pt]

\multirow{2}{*}{YY}
& Single-scale
& $n\asymp\sigma^{-(d+3)}$
& $O(\sigma^2L_\sigma)$ \\

& Geometric bias cancellation
& $n\gtrsim W^2\sigma^{-(d+6)}L_\sigma$
& $O(\sigma^4)$ \\

\bottomrule
\end{tabular*}

\par\smallskip
\begin{minipage}{\linewidth}
\scriptsize
Single-scale benchmarks:
\citet{xia2026stsp}, STSP;
\citet{genovese2014nonparametric}, DR;
\citet{yao2023manifold}, YY.
Geometric bias cancellation denotes the constructions
developed here; the YY construction uses directional weights.

Here $L_\sigma=\log(1/\sigma)$ and
$r_n=(\log n/n)^{2/(D+8)}$.

$^\dagger$ The cited STSP result provides finite-sample
localization, not Hausdorff recovery.

$^\ddagger$ The DR sample condition is inferred from
$O_p(r_n+\sigma^2L_\sigma)$, not a uniform joint
$(n,\sigma)$ guarantee.

Our bounds hold with high probability, assuming selected
empirical branches for STSP and DR or an admissible pilot
for YY. Global initialization is not addressed.
The YY condition excludes pilot cost; $W$ is the
directional weight amplification factor.

The bounds concern different accuracy targets and
do not establish a direct sample-efficiency comparison.
\end{minipage}
\end{table}

\subsection{Related Work}
\label{sec:state}

Requiring the estimator $\widehat{\M}$ itself to be a
$d$-dimensional manifold distinguishes manifold fitting
from related problems surveyed by \citet{yao2026review}.
Manifold embedding methods, including Isomap, locally linear
embedding, Laplacian eigenmaps, local tangent space alignment,
diffusion maps, and uniform manifold approximation and projection
\citep{tenenbaum2000global,roweis2000nonlinear,belkin2003laplacian,%
zhang2004principal,coifman2006diffusion,mcinnes2018umap},
construct low-dimensional coordinates designed to preserve
selected geometric or statistical structure.
Their primary target is such a representation, rather than
an explicit reconstruction of the manifold as a geometric
object in the ambient space $\R^D$.
Manifold denoising methods such as
\citet{wang2010manifold} primarily return corrected
data points, whereas the moving least-squares
projection approach of \citet{sober2020manifold}
constructs a smooth approximating manifold.
Principal curves and principal flows
\citep{hastie1989principal,panaretos2014principal}
summarize dominant variation, while principal
boundaries \citep{yao2020principal} address
classification on manifolds.
More recently, \citet{su2025principal} developed
principal nested submanifolds for nonlinear data
decomposition, constructing a hierarchy of
nested submanifolds.
Our focus, by contrast, is the geometric accuracy
of a fitted manifold as an estimator of the
underlying manifold, measured by
$\dH(\widehat{\M},\M)$.
This allows systematic geometric displacement
of the estimator to be defined, quantified and,
as we show, removed.

The modern formulation descends from the geometric Whitney problem
\citep{fefferman2020reconstruction}: given $\mathcal{A}\subset\R^{D}$, when
does a smooth $d$-dimensional submanifold approximate $\mathcal{A}$ well in
both distance and smoothness? Its statistical counterpart was posed by \citet{fefferman2016testing} as a
question of model adequacy, leading naturally to statistical formulations
such as the noisy observation model in \eqref{eq:model}.

Early reconstruction methods based on Delaunay triangulation
\citep{lee1980two,cheng2005manifold,boissonnat2009manifold} required noiseless
dense samples and returned piecewise linear objects.
\citet{genovese2012minimax, genovese2012manifold} established statistical
limits for manifold estimation under several noise models. Under bounded
perpendicular noise, the minimax Hausdorff rate depends on the intrinsic
dimension, while under additive Gaussian noise the problem becomes a
singular deconvolution problem and the corresponding minimax bounds depend
only logarithmically on the sample size. These results highlight the severe
statistical difficulty created by Gaussian contamination, but they do not
remove the need to track the joint roles of sample size, noise level, and
fitting scale. In the setting studied
here, the stochastic approximation error depends on both the sample size and
the fitting scale, while \(h\) and \(\sigma\) jointly determine the leading
population bias. Our aim is not full deconvolution at a fixed noise level,
but finite-order geometric bias correction in a joint
small-noise and small-bandwidth regime. For Gaussian
density ridges, this has a simple interpretation:
extrapolating matched ridge locations to zero total
smoothing variance; see Remark~\ref{rem:dr_deconvolution}.

To obtain computationally tractable geometric estimators under noise,
\citet{genovese2014nonparametric} studied density ridges as a proxy, with
asymptotic theory by \citet{chen2015asymptotic} and
\citet{qiao2016theoretical} and computable characterizations by
\citet{ozertem2011finding}; the resulting rates depend on \(D\) rather than
\(d\), and the ridge need not be regular.

\begingroup\emergencystretch=1em
Geometric regularity of the output was secured by the sequence
\citet{mohammed2017learning}, \citet{fefferman2018fitting} and
\citet{yao2025unbounded}, all of which define $\widehat{\M}$ as the zero or
ridge set of a smooth surrogate built from local principal component analysis,
and all of which obtain reach bounded below by $c\tau$.
\citet{mohammed2017learning} treat the noiseless case, with accuracy
$O(r^{2})$. \citet{fefferman2018fitting} admit ambient Gaussian noise and, with
$r=O(\sigma^{1/2})$, obtain
\(\dH(\widehat{\M},\M)=O(\sigma)\).
Their guarantee is established under the sample-size restriction
\(N\leq \exp(D)\), and therefore does not directly yield an
\(N\to\infty\) consistency statement at fixed ambient dimension.
\citet{yao2025unbounded} remove this restriction and replace the successive
projection construction by the bias function
\(
f(z)=\widehat{\Pi}_{z}^{\perp}
\Bigl[z-\sum_{i\in I_{z}}\alpha_{i}(z){y}_{i}\Bigr].
\)
Their method attains the same order of Hausdorff accuracy and a quantitative
reach guarantee with probability
\(1-c\exp(-Cr^{d+2}N)\), while allowing arbitrarily large sample sizes and
extending the analysis to unbounded noise.\par
\endgroup

In a related development, \citet{li2026curvature} developed a curvature-driven
manifold-fitting method under isotropic Gaussian noise. Using a
normalized local kernel estimator, they obtained the expansion
\(F(z)=\pi_{\M}(z) +\frac{d}{2}\sigma^2 H(\pi_{\M}(z))+O(\sigma^3)\),
and established $O(\sigma^2)$ Hausdorff accuracy together with
a uniform positive-reach guarantee.
This result makes explicit the role of mean curvature in the
leading geometric bias.

\subsection{Outline}
\label{sec:outline}

The remainder of the paper is organized as follows.
Section~\ref{Multiscale Bias Cancellation} develops the general theory of geometric bias cancellation,
including estimated-normal correspondence and finite-sample analysis.
Section~\ref{sec:examples} validates the general framework for successive
tangent-space projection (STSP) and Gaussian density-ridge
fitting (DR), with the complete DR analysis provided in
Appendix~\ref{sec:gaussian_density_ridge}.
Section~\ref{sec:yao_contraction} develops a directional multiscale refinement for the YY estimator.
Section~\ref{sec:simulation} presents numerical experiments.
Section~\ref{sec:conclusion} concludes with a discussion of limitations and future directions.
Technical proofs are given in Appendices~\ref{app:general_ms}--\ref{app:ysl23_full_pipeline}.

\section{General Theory for Geometric Bias Cancellation}
\label{Multiscale Bias Cancellation}
Throughout the paper, we assume isotropic Gaussian
noise with density $\phi_\sigma$ and covariance
matrix $\sigma^2 I_D$.
For each fitting scale $h$,
let $X_h:\mathcal M\to\mathbb R^D$ denote the
population embedding map associated with the fitting
procedure, and define the corresponding fitted
manifold by
\(
\mathcal M_h:=X_h(\mathcal M).
\)
We use $\widehat X_h$ and
$\widehat{\mathcal M}_h:=\widehat X_h(\mathcal M)$
for their empirical counterparts, respectively.
For multiscale extrapolation, $\widetilde X$ and
$\widetilde{\mathcal M}:=\widetilde X(\mathcal M)$
denote the population extrapolation map and its
image, whereas $\widehat X$ and
$\widehat{\mathcal M}:=\widehat X(\mathcal M)$
denote the corresponding empirical quantities.
Thus, hats distinguish empirical quantities from
their population counterparts, while tildes on
fitted maps and manifolds indicate multiscale
extrapolation.

\begin{assumption}[Geometry, sampling, and noise]
\label{ass:geometry-observation}
The following conditions are assumed throughout the paper.

\begin{enumerate}[label=\textup{(A\arabic*)},leftmargin=3em]

\item
The manifold \(\mathcal M\subset\mathbb R^D\) is a compact,
boundaryless, \(d\)-dimensional \(C^r\) embedded manifold, where
\(1\le d<D\) and \(r\ge4\). The second fundamental form and the
derivatives of curvature terms involved in the analysis are uniformly
bounded. Moreover, the latent variable \(X\) is uniformly distributed
on \(\mathcal M\).

\item
The manifold satisfies the positive reach condition
\(\operatorname{reach}(\mathcal M)\ge\tau>0\).
Consequently, for any fixed \(0<\tau<\operatorname{reach}(\mathcal M)\), the nearest-point
projection
\(\pi_{\mathcal M}:\operatorname{Tub}_{\tau}(\mathcal M)\rightarrow \mathcal M\)
is uniquely defined and smooth.

\item
The bandwidth $h$ and noise level $\sigma$ satisfy
\(0<\sigma\lesssim h\ll\tau\).

\end{enumerate}

\end{assumption}

Uniform sampling is imposed for the estimator-specific verification results. The geometric cancellation argument itself does not directly use this restriction, provided the required bias expansion and regularity conditions hold. A numerical study of GBC-STSP under nonuniform
sampling is provided in Appendix~\ref{supp:nonuniform_sampling}.

\begin{definition}[Normal-direction error]
Let \(\mathcal M_h\) be a fitted manifold with bandwidth \(h\).
Suppose that, for each \(x\in\mathcal M\), there exists a unique point
\(x_h\in \mathcal M_h\) satisfying
\(
\pi_{\mathcal M}(x_h)=x.
\)
Equivalently,
\(x_h=x+u_h(x), u_h(x)\in N_x\mathcal M\).
We call \(u_h(x)\) the normal-direction error at \(x\). The global
normal-direction error is defined by
\(\|u_h\|_{\infty} := \sup_{x\in\mathcal M}\|u_h(x)\|\).
\end{definition}
\begin{assumption}
\label{ass:normal-section}
For sufficiently small $h$, we assume that the fitted manifold admits a unique $C^1$
normal graph representation over the underlying manifold.
\end{assumption}

\begin{assumption}
\label{assump:normal_error_t}
We assume that the normal-direction error admits the expansion
\[
u_h(x)=c(h,\sigma)b(x)+O(h^4),
\]
where \(c(h,\sigma)=O(h^2+\sigma^2)\) and
\(b:\mathcal M\to\mathbb R^D\) is a bounded, scale-independent
field, typically the mean-curvature vector.
The \(O(h^4)\) remainder is uniform over \(x\in\mathcal M\)
and the admissible scale regime. We assume \(\sigma\lesssim h\),
so \(\|u_h\|_\infty=O(h^2)\).
\end{assumption}
\begin{remark}
Retaining the joint dependence of $c(h,\sigma)$ on bandwidth
and noise allows multiscale weights to cancel the leading
geometric bias. This factorization is estimator-specific
and is verified for STSP in Section~\ref{sec:examples}
and Gaussian density ridges in Appendix~\ref{sec:gaussian_density_ridge}.
\end{remark}
\begin{remark}
We adopt the factorization \(c(h,\sigma)b(x)\) because this form is common
in manifold fitting. More generally, one may have
\(u_h(x)=c(h)b(\sigma,x)+O(h^4)\),
provided that \(b(\sigma,\cdot)\) is common across the fitting scales.
The leading bias may also contain multiple terms,
\(u_h(x)=\sum_{j=1}^J c_j(h,\sigma)b_j(x)+O(h^4)\),
as in the YY estimator example in Section~\ref{sec:yao_contraction}.
Although the cancellation constraints and weights must be
adapted to the specific bias structure, the underlying
principle and the main steps of the geometric and
statistical analysis remain broadly applicable.
\end{remark}

\begin{remark}[A Common Approach to Obtaining the Normal-Bias Expansion]
\label{rem:verification_of_bias_expansion}
We present a common approach to deriving the normal-bias expansion,
which applies to STSP in Section~\ref{sec:examples} and Gaussian
density-ridge fitting in Appendix~\ref{sec:gaussian_density_ridge}.
In many manifold fitting procedures, the normal-bias expansion can be derived from the population fitting equation using the implicit-function theorem.
Here, \(\mathcal F_s(x,\beta)=0\) denotes the suitably
normalized population fitting equation in local normal
coordinates \(y=x+s^2\beta\), where
\(\beta\in N_x\mathcal M\). Its solution \(\beta_s(x)\)
determines the normal displacement
\(u_s(x)=s^2\beta_s(x)\).
For STSP, taking \(s=h\), the equation is
\[
\mathcal F_h^{\mathrm{STSP}}(x,\beta)
=
h^{-2}P_{N_x\mathcal M}
R_h^{\mathrm{STSP}}(x+h^2\beta)=0.
\]
For Gaussian density ridges, taking
\(s=\rho_h=(h^2+\sigma^2)^{1/2}\), the equation is
\[
F_{\rho_h}(x,\beta)
=
P_{N_x\mathcal M}
Q_{h,\sigma}^{\mathrm{DR}}(x+\rho_h^2\beta)
\nabla\log p_{h,\sigma}(x+\rho_h^2\beta)=0.
\]

To solve these equations, we first expand them as
\(\mathcal F_s(x,\beta) = \mathcal F_0(x,\beta)+O_{C^2}(s^2)\).
Suppose the leading equation
\(\mathcal F_0(x,\beta)=0\) has a solution
\(\beta_0(x)\), with
\(D_\beta\mathcal F_0(x,\beta_0(x))\)
uniformly invertible. By the implicit-function theorem,
there exists a unique nearby solution \(\beta_s(x)\).
Taylor expansion around \(\beta_0(x)\) gives
\[
0=
D_\beta\mathcal F_0(x,\beta_0)
(\beta_s-\beta_0)
+O\!\left(\|\beta_s-\beta_0\|^2+s^2\right).
\]
Consequently, under the corresponding uniform
smoothness conditions,
\[
\beta_s(x)=\beta_0(x)+O_{C^2}(s^2),
\qquad
u_s(x)=s^2\beta_0(x)+O_{C^2}(s^4).
\]
The detailed derivations of these normal-bias expansions are provided in
Appendices~\ref{sec:app_stsp} and~\ref{sec:app_density_ridge}.

\end{remark}

\subsection{Geometric bias cancellation with Known Normal Directions}
\label{known_direction}
We first consider the idealized setting of known normal directions,
where fitted manifolds at different scales are matched through their
common projection onto $\mathcal M$.
Although manifolds cannot be linearly combined intrinsically,
their corresponding embedding maps can be combined pointwise
in the ambient space.

We construct the multiscale
extrapolated estimator as
\(\widetilde X(x) = \sum_{i=1}^{K} w_i X_{h_i}(x)\),
where the weights satisfy$\sum_{i=1}^{K}w_i=1$ and $\sum_{i=1}^{K}w_i c(h_i,\sigma) =0$.
Under Assumption~\ref{assump:normal_error_t}, the corresponding
normal-direction error is given by
\[
\widetilde u(x)
=
\sum_{i=1}^{K}w_i u_{h_i}(x)
=
\sum_{i=1}^{K}w_i c(h_i,\sigma)b(x)
+
O\left(\sum_{i=1}^{K}|w_i| h_i^4\right)
=
O\left(\sum_{i=1}^{K}|w_i| h_i^4\right).
\]

\begin{proposition}[Exact-correspondence extrapolation]
\label{thm:manifold_recovery_1}
Under Assumptions~\ref{ass:geometry-observation},
\ref{ass:normal-section} and \ref{assump:normal_error_t},
the multiscale extrapolated map satisfies
\[
d_H(\widetilde{\mathcal M},\mathcal M)
\leq
\sup_{x\in\mathcal M}
\bigl\|
\widetilde X(x)-x
\bigr\|
\leq
C_R\sum_{i=1}^{K}|w_i|h_i^4.
\]
where \(C_R>0\) denotes the uniform constant in the remainder bound
\(
\sup_{x\in\mathcal M}
\bigl\|
u_h(x)-c(h,\sigma)b(x)
\bigr\|
\leq
C_R h^4
\)
for all sufficiently small \(h\).
If
\(
C_R\sum_{i=1}^{K}|w_i|h_i^4<\tau,
\)
then \(\widetilde X\) is a \(C^1\) embedding. Consequently,
\(\widetilde{\mathcal M}
=\widetilde X(\mathcal M)\)
is a \(d\)-dimensional \(C^1\) embedded manifold that is
\(C^1\)-diffeomorphic to \(\mathcal M\).
\end{proposition}

\subsection{Geometric Bias Cancellation via Estimated-Normal-Fiber Correspondence}
\label{sec:normal_guided_extrapolation}
Exact normal correspondence provides the benchmark, but it depends on the unknown manifold. We now construct correspondence from the fitted manifolds themselves and show that its error is of the same order as the residual bias after cancellation. The
key observation is that, under suitable normal-estimation accuracy, the
resulting correspondence error does not affect the bias-cancellation rate.

For each
\(x\in\mathcal M\), define the reference point
\(
y_1(x)
:=
X_{h_1}(x)
=
x+u_{h_1}(x)
\in\mathcal M_{h_1}.
\)
By Assumption~\ref{assump:normal_error_t},
\(
u_{h_1}(x)
=
c(h_1,\sigma)b(x)+O(h_1^4).
\)
Rather than prescribing a particular normal-space estimation procedure,
we impose the following general accuracy condition.

\begin{assumption}[Accuracy and regularity of estimated-normal correspondence]
\label{ass:estimated-normal-bundle}
Let \(h_1<h_2<\cdots<h_K\).
In codimension one, assume that
\[
\sup_{x\in\mathcal M}
\|\widetilde n(x)-n(x)\|
\leq C_Nh_1^2,
\]
where \(n(x)\) and \(\widetilde n(x)\) are compatibly chosen local
unit normals. For codimension greater than one, let
\(\widetilde N(x)\subset\mathbb R^D\) be the estimated
\((D-d)\)-dimensional normal space at \(y_1(x)\), and replace the
preceding condition by
\[
\sup_{x\in\mathcal M}
\left\|P_{\widetilde N(x)}-P_{N_x\mathcal M}\right\|_{\mathrm{op}}
\leq C_Nh_1^2.
\]
In codimension one, write
\(\widetilde N(x)=\operatorname{span}\{\widetilde n(x)\}\).
In any codimension, regularity of the estimated normal field refers
to the projector field \(x\mapsto P_{\widetilde N(x)}\).
Orthonormal frames are used only locally; no global trivialization
of the normal bundle is assumed.
\end{assumption}

The following lemma demonstrates how this assumption can be met.
\begin{lemma}[Population PCA accuracy]
\label{lem:pca_normal_accuracy}
Under Assumption~\ref{ass:geometry-observation}, suppose additionally that
\(\mathcal M\) is \(C^6\) and \(\lambda=\sigma/h\in\Lambda\), where
\(\Lambda\Subset(0,\infty)\) is fixed and compact.
Use \(K(t)=(1-t)_+^6\) and
\(K_h(y,z)=K(\|z-y\|^2/h^2)\), and define the weighted mean and
centered covariance by
\[
\mu_h(y)=\frac{\mathbb E[Y K_h(y,Y)]}{\mathbb E[K_h(y,Y)]},
\qquad
\Sigma_h(y)=\frac{\mathbb E[(Y-\mu_h(y))(Y-\mu_h(y))^\top K_h(y,Y)]}
{\mathbb E[K_h(y,Y)]}.
\]
Let \(P_h(y)\) be the orthogonal projector onto the leading
\(d\)-dimensional eigenspace of \(\Sigma_h(y)\), and set
\(Q_h(y)=I_D-P_h(y)\). For
\(y=x+h^2\beta,\quad x\in\mathcal M,\quad \beta\in N_x\mathcal M\),
with \(\beta\) uniformly bounded,
\[
\Sigma_h(y)
=
h^2\left\{
\tau_{\sigma/h}P_{T_x\mathcal M}
+
a_{\sigma/h}P_{N_x\mathcal M}
\right\}
+
O_{C^2}(h^4),
\quad
\inf_{\lambda\in\Lambda}
(\tau_\lambda-a_\lambda)>0.
\]
The constants \(\tau_\lambda=\mathbb E_{\nu_\lambda}[v_{T,1}^2]\) and
\(a_\lambda=\mathbb E_{\nu_\lambda}[v_{N,1}^2]\) are the coordinate
second moments under a probability measure \(\nu_\lambda\) on \(B_1(0)\).
Its density is proportional to
\(K(\|v\|^2)\exp(-\|v_N\|^2/(2\lambda^2))\), where
\(v_T=P_{T_x\mathcal M}v\) and \(v_N=P_{N_x\mathcal M}v\).
Consequently,
\[
P_h(y)
=
P_{T_x\mathcal M}+O_{C^2}(h^2),
\qquad
Q_h(y)
=
P_{N_x\mathcal M}+O_{C^2}(h^2).
\]
The bounds are uniform in \(x\), \(\lambda\in\Lambda\), and bounded
\(\beta\); \(C^2\) derivatives are taken in local normal-bundle
coordinates \((x,\beta)\), with \(h,\lambda\) held fixed.
\end{lemma}

Define the reference point
\(
y_1(x)
:=
X_{h_1}(x)
=
x+u_{h_1}(x)
\in\mathcal M_{h_1}.
\)
In codimension one, define the affine normal line through
\(y_1(x)\) by
\(\widetilde L_{h_1}(x;s) = y_1(x) + s\,\widetilde n(x), \quad s\in\mathbb R\).
In arbitrary codimension, use the frame-independent affine normal fiber
\(\widetilde{\mathcal F}_{h_1}(x) :=y_1(x)+\widetilde N(x)\).
\begin{assumption}[Tangent-space stability of fitted manifolds]
\label{ass:tangent-stability}
There exists a constant \(C_T>0\), independent of \(h_k\),
such that
\[
\sup_{x\in\mathcal M}
\left\|
P_{T_{X_{h_k}(x)}\mathcal M_{h_k}}
-
P_{T_x\mathcal M}
\right\|_{\mathrm{op}}
\leq
C_T h_k^2
\]
for all sufficiently small \(h_k\).
\end{assumption}

\begin{lemma}[Well-posedness of estimated-normal correspondence]
\label{lem:wellposed_correspondence}
\leavevmode\par\nobreak

Suppose Assumptions~\ref{ass:geometry-observation},
\ref{ass:normal-section}, \ref{assump:normal_error_t},
\ref{ass:estimated-normal-bundle} and \ref{ass:tangent-stability} hold.
Then, for sufficiently small
\(
h_{\max}:=\max_{1\leq k\leq K}h_k,
\)
and for every \(k=1,\ldots,K\), the affine normal fiber
\(\widetilde{\mathcal F}_{h_1}(x)\) has a unique local intersection
with \(\mathcal M_{h_k}\) in a sufficiently
small tubular neighborhood of \(x\). We denote this intersection by
\(
Z_{h_k}(x)\in\mathcal M_{h_k}.
\)

Moreover, if the estimated normal-projector field
\(x\mapsto P_{\widetilde N(x)}\) is \(C^1\), then the resulting
correspondence map
\(
Z_{h_k}:\mathcal M\longrightarrow\mathcal M_{h_k}
\)
is \(C^1\).
\end{lemma}
The collection
\(
\{Z_{h_k}(x)\}_{k=1}^{K}
\)
therefore defines the estimated-normal cross-scale correspondence
associated with the reference point \(y_1(x)\). When \(D-d=1\), the estimated normal space is one-dimensional and the
affine normal fiber reduces to
\[
\widetilde{\mathcal F}_{h_1}(x)
=
\left\{
y_1(x)+s\widetilde n(x):
s\in\mathbb R
\right\}.
\]
Here \(\widetilde n\) may be chosen locally. For higher codimension,
the construction uses the entire estimated normal space
\(\widetilde N(x)\); the intersection point does not depend on a
choice of local orthonormal frame.

For each corresponded point \(Z_{h_k}(x)\in\mathcal M_{h_k}\), let
\(
x_k
=
\pi_{\mathcal M}\bigl(Z_{h_k}(x)\bigr)
\in\mathcal M
\)
denote its nearest-point projection onto the underlying manifold.
Since \(Z_{h_k}(x)\in\mathcal M_{h_k}\), it can be written as
\(Z_{h_k}(x) = x_k+u_{h_k}(x_k)\).
By Assumption~\ref{assump:normal_error_t},
\(u_{h_k}(x_k) = c(h_k,\sigma)b(x_k)+O(h_k^4)\).
The correspondence error and the variation of the normal displacement
are controlled by Lemmas~\ref{lem:normal_guided_correspondence}
and~\ref{lem:normal_field_stability} in
Appendix~\ref{app:general_ms}.

Next, using the estimated-normal correspondence
\(\{Z_{h_k}(x)\}_{k=1}^{K}\), we construct the GBC estimator
\[
\widetilde X(x)
=
\sum_{k=1}^{K} w_k Z_{h_k}(x),
\]
where the weights satisfy
\(
\sum_{k=1}^{K} w_k=1,
\sum_{k=1}^{K} w_k c(h_k,\sigma)=0 .
\)
By Assumption~\ref{assump:normal_error_t} and the cancellation constraint,
\(\|\sum_{k=1}^K w_k u_{h_k}(x)\|\leq C\sum_{k=1}^K|w_k|h_k^4\).
Lemmas~\ref{lem:normal_guided_correspondence}
and~\ref{lem:normal_field_stability} give
\(\|x_k-x\|\leq Ch_1^2h_k^2\) and
\(\|u_{h_k}(x_k)-u_{h_k}(x)\|\leq Ch_1^2h_k^4\), respectively.
Using \(h_1\leq h_k\), these bounds yield the fourth-order recovery
estimate stated below. In particular, uniformly bounded weights and
comparable scales give \(\widetilde X(x)-x=O(h_{\max}^4)\).

\begin{proposition}
\label{thm:normal_guided_manifold_recovery}

Suppose Assumptions~\ref{ass:geometry-observation},
\ref{ass:normal-section},
\ref{assump:normal_error_t},
\ref{ass:estimated-normal-bundle}, and
\ref{ass:tangent-stability} hold. Let
\(Z_{h_k}(x)\in\mathcal M_{h_k}\) be the estimated-normal correspondence
points defined in Lemma~\ref{lem:wellposed_correspondence}, and define
\(\widetilde X(x) = \sum_{k=1}^{K}w_k Z_{h_k}(x), \quad x\in\mathcal M\),
where the weights satisfy
\(\sum_{k=1}^{K}w_k=1, \quad \sum_{k=1}^{K}w_kc(h_k,\sigma)=0\).

Then, for sufficiently small
\(h_{\max}:=\max_k h_k\),
\[
\sup_{x\in\mathcal M}
\|\widetilde X(x)-x\|
\leq
C
\sum_{k=1}^{K}|w_k|h_k^4,
\]
where \(C>0\) is independent of the scales \(h_k\). Consequently, for
\(
\widetilde{\mathcal M}
=
\widetilde X(\mathcal M),
\)
we have
\(
d_H(\widetilde{\mathcal M},\mathcal M)
\leq
C
\sum_{k=1}^{K}|w_k|h_k^4.
\)

\begingroup\sloppy
Moreover, suppose that \(K\) is fixed and that the extrapolation weights
satisfy
\(
\sup_{h_{\max}}\sum_{k=1}^{K}|w_k|\leq W_0<\infty.
\)
Suppose also that the estimated normal-projector field
\(x\mapsto P_{\widetilde N(x)}\) is \(C^1\) with a uniformly bounded
\(C^1\) norm. Then, for sufficiently small
\(h_{\max}\), the estimated-normal GBC map
\(\widetilde X\) is a \(C^1\) embedding. Consequently,
\(
\widetilde{\mathcal M}
=
\widetilde X(\mathcal M)
\)
is a \(d\)-dimensional \(C^1\) embedded manifold that is
\(C^1\)-diffeomorphic to \(\mathcal M\).
\par\endgroup
\end{proposition}
\begin{remark}
We focus on removing the leading second-order normal bias
and establishing fourth-order accuracy. Further cancellation
would require more detailed estimator-specific expansions,
additional cancellation constraints, and, for estimated-normal
extrapolation, more accurate cross-scale correspondence,
without necessarily improving finite-sample accuracy.
\end{remark}
\begin{remark}
The population results above establish positional accuracy and,
under their respective regularity conditions, a \(C^1\) embedded
image that is diffeomorphic to \(\mathcal M\).
They do not assert a uniform positive lower bound on the reach
of the corrected image.
\end{remark}
\begin{remark}
\label{rem:reference_normal_space}
A natural choice of the estimated normal space is the normal space of
the reference fitted manifold itself,
\(
\widetilde N(x)
=
N_{X_{h_1}(x)}\mathcal M_{h_1}.
\)
Under Assumption~\ref{ass:tangent-stability},
\[
\sup_{x\in\mathcal M}
\left\|
P_{\widetilde N(x)}
-
P_{N_x\mathcal M}
\right\|_{\mathrm{op}}
\leq
C_T h_1^2,
\]
since the orthogonal normal projector is the complement of the tangent
projector.
Moreover, if \(\mathcal M_{h_1}\) is \(C^2\) and its second fundamental
form is uniformly bounded for sufficiently small \(h_1\), then the normal
projector field
\(
x\longmapsto P_{\widetilde N(x)}
\)
is \(C^1\) with a uniform bound
\(
\sup_{x\in\mathcal M}
\left\|
D P_{\widetilde N(x)}
\right\|_{\mathrm{op}}
\leq C_N',
\)
where \(C_N'\) is independent of \(h_1\). Hence the estimated normal
bundle satisfies both the required \(O(h_1^2)\) approximation accuracy
and the \(C^1\) regularity condition used in the estimated-normal
correspondence.

Alternatively, the normal space may be estimated locally by PCA. PCA has the practical advantage of relatively low computational cost, especially when an explicit tangent or normal representation of the fitted manifold is unavailable. 
For the smooth-kernel PCA construction used here,
the required second-order normal-projector accuracy
follows from kernel-moment expansions and a uniform
population eigengap. This accuracy is sufficient to
preserve the fourth-order correspondence error.
To further guarantee that the population GBC
estimator is a $C^1$ embedded manifold, the estimated
normal-projector field must also be $C^1$ with a
uniformly bounded derivative. These properties are
verified for STSP and DR in Appendices~\ref{sec:app_stsp} and~\ref{sec:app_density_ridge}.
\textit{}
\end{remark}

\begin{remark}[Computational interpretation of normal-fiber matching]
\label{rem:enf_computational_view}
The estimated-normal-fiber (ENF) construction replaces
cross-scale point matching by a local transverse intersection.
This avoids global nearest-point selection and admits
an a posteriori numerical-error certificate based on the
residual of the local intersection equation. The following
result makes the geometric conditioning explicit; the
corresponding implementable residual equations for STSP and
Gaussian density ridges are analyzed in Sections~\ref{sec:stsp_validation}--\ref{sec:dr_validation}
and Appendices~\ref{app:stsp_finite_sample_proof} and~\ref{app:dr_finite_sample_proof}.
\end{remark}

\begin{proposition}[Uniform local conditioning of ENF correspondence]
\label{prop:enf_local_conditioning}
Suppose Assumptions~\ref{ass:geometry-observation}, \ref{ass:normal-section}, \ref{assump:normal_error_t}, \ref{ass:estimated-normal-bundle}, and~\ref{ass:tangent-stability} hold.
Let $q=D-d$, and let $h_1\le h_k$ tend to zero. For each
$x\in\mathcal M$, choose a local $C^2$ chart
$\psi_x:U_x\subset\mathbb R^d\to\mathcal M$ with
$\psi_x(0)=x$ and $D\psi_x(0)$ an isometry, and let
$E_x\in\mathbb R^{D\times q}$ be an orthonormal basis of
the estimated normal space $\widetilde N(x)$. Define
\[
  \mathcal F_{k,x}(u,t)
  :=X_{h_k}(\psi_x(u))-X_{h_1}(x)-E_xt,
  \qquad (u,t)\in\mathbb R^d\times\mathbb R^q.
\]
For every fixed sufficiently large $C_{\mathrm{loc}}$,
there exist orthogonal matrices $O_x\in\mathbb R^{D\times D}$
and a constant $C$, uniform in $x$ and the admissible scales,
such that
\begin{equation}\label{eq:enf_geometric_jacobian}
 \sup_{\|u\|+\|t\|\le C_{\mathrm{loc}}h_k^2}
 \|D_{(u,t)}\mathcal F_{k,x}(u,t)-O_x\|_{\mathrm{op}}
 \le C(h_k^2+h_1^2)=:\delta_k.
\end{equation}
Consequently, for $\delta_k<1$,
\begin{equation}\label{eq:enf_geometric_cond}
 1-\delta_k\le\sigma_{\min}(D\mathcal F_{k,x})
 \le\sigma_{\max}(D\mathcal F_{k,x})\le1+\delta_k,
 \qquad
 \kappa_2(D\mathcal F_{k,x})
 \le\frac{1+\delta_k}{1-\delta_k}=1+O(h_k^2+h_1^2).
\end{equation}
If $(u_k^*,t_k^*)$ denotes the local intersection provided
by Lemma~\ref{lem:wellposed_correspondence} and $(\bar u,\bar t)$ belongs to the same
convex local search region, then
\begin{equation}\label{eq:enf_geometric_residual}
 \|(\bar u,\bar t)-(u_k^*,t_k^*)\|
 \le\frac{\|\mathcal F_{k,x}(\bar u,\bar t)\|}{1-\delta_k}.
\end{equation}
The same conclusions hold for empirical branches and
aligned empirical normal frames on any event where their
$C^1$ branch errors and normal-frame errors are bounded by
$e_{n,k}$ and $s_{n,1}$, respectively, with $\delta_k$
replaced by $C(h_k^2+h_1^2+e_{n,k}+s_{n,1})$, provided
the corresponding local empirical intersection exists.
\end{proposition}

\begin{remark}[Comparison with closest-point matching]
\label{rem:enf_versus_cp}
Continuous closest-point (CP) correspondence can also
achieve fourth-order localization. However, uniqueness
and smoothness of the global CP map require an additional
uniform curvature condition in the general framework;
without it, CP can be multivalued even when the local estimated normal-fiber(ENF)
intersection is unique and well conditioned. The precise CP
statement, proof, and a counterexample are collected in
Proposition~\ref{prop:cp_supplement} and
Appendix~\ref{app:general_enf_proofs}. 
\end{remark}

\subsubsection*{Finite-sample consistency of GBC with estimated-normal correspondence}
\label{sec:finite_sample_normal_guided}
In practice,
both the fitted manifolds and the normal directions are estimated from
finite observations. We therefore study how these empirical errors
propagate through the estimated-normal correspondence and multiscale
extrapolation.
Define the population and empirical reference points by
\(
y_1(x)=X_{h_1}(x),
\widehat y_1(x)=\widehat X_{h_1}(x).
\)

\begin{assumption}[Uniform empirical approximation of single-scale estimators]
\label{assump:normal_guided_empirical}
There exist deterministic quantities \(r_n(h_k)\to0\) and
\(\delta_n^{X}\to0\) such that
\[
\mathbb P\left(
\bigcap_{k=1}^{K}
\left\{
\sup_{x\in\mathcal M}
\left\|
\widehat X_{h_k}(x)-X_{h_k}(x)
\right\|
\leq r_n(h_k)
\right\}
\right)
\geq
1-\delta_n^{X}.
\]
\end{assumption}

\begin{assumption}[Finite-sample normal estimation accuracy]
\label{assump:finite_normal_estimation}
In codimension one, let \(\widehat n(x)\) denote the finite-sample
estimate of the population normal direction \(\widetilde n(x)\) at
the reference estimator, with compatible local orientations.
There exist deterministic quantities
\(s_n(h_1)\to0\) and \(\delta_n^{N}\to0\) such that
\[
\mathbb P\left(
\sup_{x\in\mathcal M}
\left\|
\widehat n(x)-\widetilde n(x)
\right\|
\leq
s_n(h_1)
\right)
\geq
1-\delta_n^{N}.
\]
For codimension greater than one, let \(\widehat N(x)\) denote the
finite-sample estimate of \(\widetilde N(x)\), and interpret the norm
above as
\(\|P_{\widehat N(x)}-P_{\widetilde N(x)}\|_{\mathrm{op}}\).
In codimension one, write
\(\widehat N(x)=\operatorname{span}\{\widehat n(x)\}\).
\end{assumption}

Using the empirical reference point and normal estimate, define the
empirical normal trajectory in codimension one by
\(\widehat L_{h_1}(x;s) = \widehat y_1(x)+s\widehat n(x)\).
In arbitrary codimension, use the affine normal fiber
\(\widehat{\mathcal F}_{h_1}(x) :=\widehat y_1(x)+\widehat N(x)\).
Let
\(
\widehat Z_{h_k}(x)\in\widehat{\mathcal M}_{h_k}
\)
denote the resulting local estimated-normal correspondence point, and let
\(Z_{h_k}(x)\in\mathcal M_{h_k}\) denote its population counterpart.

\begin{lemma}[Finite-sample estimated-normal correspondence error]
\label{lem:finite_normal_correspondence}
\leavevmode\par\nobreak
Suppose Assumptions~\ref{assump:normal_guided_empirical} and
\ref{assump:finite_normal_estimation} hold, together with the local
regularity conditions of Lemma~\ref{lem:wellposed_correspondence}.
Then, provided that the empirical affine normal fiber admits a local
intersection, the following bound holds on the intersection of the two
high-probability events.
\[
\sup_{x\in\mathcal M}
\left\|
\widehat Z_{h_k}(x)-Z_{h_k}(x)
\right\|
\leq
C
\left[
r_n(h_1)
+
r_n(h_k)
+
s_n(h_1)h_k^2
\right],
\qquad
k=1,\ldots,K.
\]
\end{lemma}

We now define the finite-sample GBC map by
\[
\widehat X(x)
=
\sum_{k=1}^{K}
w_k\widehat Z_{h_k}(x),
\qquad
x\in\mathcal M,
\]
where
\(
\sum_{k=1}^{K}w_k=1,
\sum_{k=1}^{K}w_kc(h_k,\sigma)=0.
\)
The corresponding finite-sample GBC estimator is
\(
\widehat{\mathcal M}
=
\widehat X(\mathcal M).
\)

\begin{theorem}[Finite-sample error decomposition for GBC]
\label{thm:finite_normal_guided_decomposition}
Suppose Assumptions~\ref{ass:geometry-observation},
\ref{ass:normal-section},
\ref{assump:normal_error_t},
\ref{ass:estimated-normal-bundle},
\ref{ass:tangent-stability},
\ref{assump:normal_guided_empirical}, and
\ref{assump:finite_normal_estimation} hold. Let the extrapolation
weights satisfy
\(\sum_{k=1}^{K}w_k=1, \quad \sum_{k=1}^{K}w_kc(h_k,\sigma)=0\).

Let \(E_X\) and \(E_N\) denote the high-probability events in
Assumptions~\ref{assump:normal_guided_empirical} and
\ref{assump:finite_normal_estimation}, respectively. Then
\(\mathbb P(E_X\cap E_N) \geq 1-\delta_n^{X}-\delta_n^{N}\).
On \(E_X\cap E_N\), whenever the empirical estimated-normal correspondence
points are locally well defined,
\[
\begin{aligned}
d_H(\widehat{\mathcal M},\mathcal M)
&\leq \sup_{x\in\mathcal M}\|\widehat X(x)-x\|\\
&\leq C\sum_{k=1}^{K}|w_k|h_k^4
+C\sum_{k=1}^{K}|w_k|
\bigl[r_n(h_1)+r_n(h_k)+s_n(h_1)h_k^2\bigr].
\end{aligned}
\]
\end{theorem}

The first term is the residual population bias after Geometric bias cancellation. The second and third terms arise from finite-sample
approximation of the reference and fitted manifolds, while the last
term quantifies the additional correspondence error caused by estimating
the normal direction.

\begingroup\emergencystretch=1em
Consequently, for fixed \(K\) and bounded extrapolation
weights, if
\(
\max_{1\leq k\leq K}r_n(h_k)\to0
\) and
\(
s_n(h_1)\to0,
\)
then the finite-sample GBC estimator converges to its
population GBC counterpart. Moreover, if
\(
s_n(h_1)=O(h_1^2)
\)
and the scales are of the same order, then
\(
s_n(h_1)h_k^2=O(h_{\max}^4),
\)
so normal estimation does not change the fourth-order population
bias-cancellation rate.\par
\endgroup

\begin{remark}[Choice of extrapolation weights]
The extrapolation weights affect both bias cancellation and finite-sample
error amplification. Let
\(
a_{n,k}
=
h_k^4+r_n(h_1)+r_n(h_k)+s_n(h_1)h_k^2.
\)
By Theorem~\ref{thm:finite_normal_guided_decomposition},
\(
d_H(\widehat{\mathcal M},\mathcal M)
\lesssim
\sum_{k=1}^{K}|w_k|a_{n,k}.
\)
Hence a natural choice is
\(\min_{w\in\mathbb R^K} \sum_{k=1}^{K}a_{n,k}|w_k|\)
subject to
\(\sum_{k=1}^{K}w_k=1, \sum_{k=1}^{K}w_kc(h_k,\sigma)=0\).
Thus, while the constraints remove the leading bias, excessively large
signed weights should be avoided because they amplify both higher-order
bias and finite-sample errors.
\end{remark}

\begin{proposition}[Error propagation from an estimated noise level]
\label{prop:noise_level_perturbation}
Suppose the conditions of
Theorem~\ref{thm:finite_normal_guided_decomposition}, other than its
true-noise cancellation constraint, hold at the true noise level
\(\sigma\) and fixed deterministic scales \(h_1<\cdots<h_K\). Let \(\widehat{\sigma}>0\) be an estimate of the noise level, and let the
spatially constant, possibly data-dependent weights
\(\{\widehat w_k\}_{k=1}^{K}\) satisfy
\[
\sum_k\widehat w_k=1,\qquad
\sum_k\widehat w_k c(h_k,\widehat\sigma)=0,\qquad
\sum_k|\widehat w_k|\leq W_0<\infty.
\]
The noise estimate enters only the weights; the single-scale estimators
and correspondences are unchanged. Define
\[
\widehat X_{\widehat\sigma}(x)
=\sum_k\widehat w_k\widehat Z_{h_k}(x),\qquad
\widehat{\mathcal M}_{\widehat\sigma}
=\widehat X_{\widehat\sigma}(\mathcal M),\qquad
\Delta_{n,c}
=\sum_k\widehat w_k\{c(h_k,\sigma)-c(h_k,\widehat\sigma)\}.
\]
The uncanceled leading field is \(\Delta_{n,c}b(x)\). On \(E_X\cap E_N\),
whenever the local empirical correspondences are well defined,
\[
\begin{aligned}
d_H(\widehat{\mathcal M}_{\widehat\sigma},\mathcal M)
&\leq\sup_{x\in\mathcal M}
\|\widehat X_{\widehat\sigma}(x)-x\|\\
&\leq C\sum_k|\widehat w_k|
\{h_k^4+r_n(h_1)+r_n(h_k)+s_n(h_1)h_k^2\}
+\|b\|_\infty|\Delta_{n,c}|.
\end{aligned}
\]
Suppose \(E_\sigma=\{|\widehat\sigma-\sigma|\leq\rho_n\}\) has probability
at least \(1-\delta_n^\sigma\), where \(\rho_n\) is deterministic.
If \(c(h_k,\cdot)\) is continuously differentiable on a deterministic
interval \(I_n\subset(0,\infty)\) containing \(\sigma\) and all values of
\(\widehat\sigma\) on \(E_\sigma\), set
\(L_{n,k}=\sup_{s\in I_n}|\partial_s c(h_k,s)|<\infty\).
Then, on \(E_\sigma\),
\[
|\Delta_{n,c}|\leq\rho_n\sum_k|\widehat w_k|L_{n,k},
\qquad
\mathbb P(E_X\cap E_N\cap E_\sigma)
\geq1-\delta_n^X-\delta_n^N-\delta_n^\sigma.
\]
The added term preserves a target error rate when it is no larger
than that rate; consistency of \(\widehat\sigma\) alone does not ensure
fourth-order bias accuracy.
\end{proposition}

\begin{algorithm}[!htbp]
\caption{Finite-Sample GBC with Estimated-Normal-Fiber Correspondence}
\label{alg:normal_guided_ms}
\begin{algorithmic}[1]

\Require Data $\mathcal Y=\{y_i\}_{i=1}^n$,
fitting procedure $\mathcal A$ with residuals
$\widehat R_{h_k}^{\mathcal A}$, scales $h_1<\cdots<h_K$,
selected branch $\widehat{\mathcal M}_{h_1}$,
queries $\mathcal Z_N\subset\widehat{\mathcal M}_{h_1}$,
$C_{\mathrm{loc}}>0$, tolerances $\varepsilon_{n,k}\geq0$,
validity thresholds, and weights $\{w_k\}_{k=1}^K$
satisfying
$\sum_k w_k=1$ and $\sum_k w_kc(h_k,\sigma)=0$.

\Ensure Output cloud $\widehat{\mathcal P}_N^{\mathrm{num}}$
and retained query set $S$.

\State Construct $\widehat R_{h_k}^{\mathcal A}$ for $k=2,\ldots,K$;
set $S\gets\varnothing$.

\For{each $z\in\mathcal Z_N$}
    \State Estimate $\widehat N(z)$ and an orthonormal
    basis $\widehat E(z)\in\mathbb R^{D\times(D-d)}$.
    \If{reference validity checks fail}
        \State \textbf{continue}.
    \EndIf
    \State Set $\widehat t_1(z)=0$.
    \For{$k=2,\ldots,K$}
        \State Find $\widehat t_k(z)\in\mathbb R^{D-d}$
        satisfying
        \Statex \hspace{\algorithmicindent}
        $\|\widehat t_k(z)\|\leq C_{\mathrm{loc}}h_k^2,$
        $\left\|\widehat E(z)^\top
        \widehat R_{h_k}^{\mathcal A}
        (z+\widehat E(z)\widehat t_k(z))\right\|
        \leq\varepsilon_{n,k}.$
        \If{the solve or validity checks fail}
            \State \textbf{break}.
        \EndIf
    \EndFor
    \If{all local solves succeed}
        \State $\displaystyle
        \widehat F^{\mathrm{num}}(z)
        =z+\widehat E(z)\sum_{k=1}^K
        w_k\widehat t_k(z)$;
        $S\gets S\cup\{z\}$.
    \EndIf
\EndFor

\State \Return
$\widehat{\mathcal P}_N^{\mathrm{num}}
=\{\widehat F^{\mathrm{num}}(z):z\in S\}$ and $S$.

\end{algorithmic}
\end{algorithm}

\noindent\textbf{Reference initialization and finite queries.}
Our focus is higher-order refinement through geometric bias cancellation rather than global manifold initialization.
Algorithm~\ref{alg:normal_guided_ms} takes a selected
reference branch $\widehat{\mathcal M}_{h_1}$ and a finite
query set
$\mathcal Z_N=\{z_1,\ldots,z_N\}
\subset\widehat{\mathcal M}_{h_1}$.
For STSP and Gaussian density ridges, we establish the
existence and regularity of the relevant empirical branches,
but do not address their global identification or complete
sampling from arbitrary noisy initializations.

Each local correspondence search is restricted to
$\|t\|\leq C_{\mathrm{loc}}h_k^2$ and must satisfy
the prescribed residual tolerance and estimator-specific
validity checks.
Let $S\subseteq\mathcal Z_N$ be the set of successfully
retained queries. When $S\neq\varnothing$, define
\(
\delta_S=
\sup_{z\in\widehat{\mathcal M}_{h_1}}
\min_{z_j\in S}\|z-z_j\|.
\)
The numerical covering bound depends on $\delta_S$,
rather than the fill distance of the original query set.
An empty retained set is reported as a failed run.
The formal estimator is a continuous manifold obtained
by exact correspondence and extrapolation, whereas
the algorithm returns a finite numerical sample of it.

\section{Validation of the GBC Framework}
\label{sec:examples}
We now examine the geometric bias cancellation principle in two
representative manifold-fitting procedures.  These two examples verify
the factorized normal expansion required by the general theory in
Section~\ref{Multiscale Bias Cancellation} through two structurally different constructions.  We begin with successive tangent-space projection (STSP),
where the leading curvature bias arises from local means
and PCA-based tangent estimation, and the full
estimated-normal finite-sample framework of Section~\ref{Multiscale Bias Cancellation}
can be instantiated. We then consider Gaussian
density-ridge fitting (DR) in Section~\ref{sec:dr_validation}, where the
same factorized bias structure follows from the local
ridge equation, providing a second verification of the
general framework. The complete DR analysis is given
in Appendix~\ref{sec:gaussian_density_ridge}.

\subsection{Successive tangent-space projection}
\label{sec:stsp_validation}
We consider a smooth-kernel version of successive
tangent-space projection (STSP), originally introduced
by \citet{xia2026stsp}.
Unlike the original hard-radius construction, we use
a compactly supported smooth kernel to define the
local mean and covariance.

Throughout this subsection, the superscript
$\mathrm{STSP}$ identifies quantities associated
with successive tangent-space projection.
The subscript $h$ denotes the fitting bandwidth,
and $\lambda_h=\sigma/h$ denotes the normalized
noise level. Population quantities are written
without hats, while hats denote their empirical
counterparts, following the convention in
Section~\ref{Multiscale Bias Cancellation}.

Our goal is to show that its population fixed-point manifold admits the
factorized normal expansion
\(u_h^{\mathrm{STSP}}(x) = c_{\mathrm{STSP}}(h,\sigma)H(x) + O_{C^2}(h^4)\),
thereby verifying the structural assumptions of Section~\ref{Multiscale Bias Cancellation}.  
Here we fix the radial kernel
\(
K(t)=(1-t)_+^6, t\geq0,
\)
and put \(K_h(y,z)=K(\|z-y\|^2/h^2)\).
For \(y\in\mathbb R^D\), define the population kernel mass, local mean,
and covariance by
\[
m_{0,h}(y)=\mathbb E[K_h(y,Y)],\qquad
\mu_{h}(y)=\frac{\mathbb E[Y K_h(y,Y)]}{m_{0,h}(y)},
\]
\[
\Sigma_{h}(y)
=
\frac{\mathbb E\!\left[
(Y-\mu_{h}(y))(Y-\mu_{h}(y))^\top K_h(y,Y)
\right]}{m_{0,h}(y)}.
\]
We suppress the kernel subscript within this subsection.
\[
P_{h}(y)=P_{d}(\Sigma_{h}(y)),\qquad Q_{h}(y)=I-P_{h}(y),
\]
\[
T_{h}^{\mathrm{STSP}}(y)
=
\mu_{h}(y)+P_{h}(y)\{y-\mu_{h}(y)\},
\qquad
R_{h}^{\mathrm{STSP}}(y)
=
Q_{h}(y)\{y-\mu_{h}(y)\}.
\]
Thus \(y\) is a population STSP fixed point if and only if
\(R_{h}^{\mathrm{STSP}}(y)=0\).
Whenever the \(d\)-th and \((d+1)\)-st eigenvalues of
\(\Sigma_{h}(y)\) are separated, let \(P_{h}(y)\) denote the orthogonal
projector onto its leading \(d\)-dimensional eigenspace.

\begin{assumption}[STSP population scaling]
\label{ass:stsp-pop}
The manifold \(\mathcal M\) is \(C^6\).
The normalized noise level satisfies
\(
\lambda_h
=
\frac{\sigma}{h}
\in
\Lambda_{\mathrm{STSP}},
\)
where
\(
\Lambda_{\mathrm{STSP}}\Subset(0,\infty)
\)
is compact. In a multiscale comparison, the same noise level \(\sigma\)
is used at every bandwidth, while \(\lambda_h=\sigma/h\) varies with
\(h\). All statements below are uniform over admissible pairs
\((h,\sigma)\) for which \(\lambda_h\in\Lambda_{\mathrm{STSP}}\).
Moreover, \(h/\tau\) is sufficiently small.
\end{assumption}

For \(\lambda>0\), define the probability measure on
\(B_1(0)\subset\mathbb R^D\) by
\[
d\nu_{\lambda}(v)
=
\frac{
K(\|v\|^2)
\exp\!\left(
-\frac{\|v_N\|^2}{2\lambda^2}
\right)\,dv
}{
\displaystyle
\int_{B_1(0)}
K(\|w\|^2)\exp\!\left(
-\frac{\|w_N\|^2}{2\lambda^2}
\right)\,dw
},
\]
where
\(
v=v_T+v_N,
v_T\in T_x\mathcal M,
v_N\in N_x\mathcal M.
\)
By rotational invariance, the quantities
\(
a_{\lambda}
=
\mathbb E_{\lambda}[v_{N,1}^2],
c_{\lambda}
=
\mathbb E_{\lambda}[v_{N,1}^2v_{T,1}^2],
\tau_{\lambda}
=
\mathbb E_{\lambda}[v_{T,1}^2]
\)
are independent of the choice of tangent and normal coordinates.
These are the moments of the fixed smooth kernel, not those of a hard-ball
neighborhood. In particular, the extrapolation weights below use
\(c_{\mathrm{STSP}}\).

We next analyze the fixed-point equation in true normal coordinates. For
\(x\in\mathcal M\) and \(\beta\in N_x\mathcal M\), define the scaled
normal residual
\(\mathcal F_h^{\mathrm{STSP}}(x,\beta) = h^{-2} P_{N_x\mathcal M} R_{h}^{\mathrm{STSP}} \bigl(x+h^2\beta\bigr)\).

\begin{lemma}[Scaled normal equation and nondegeneracy]
\label{lem:stsp_scaled_normal_equation}
Suppose Assumptions~\ref{ass:geometry-observation} and
\ref{ass:stsp-pop} hold. Then, uniformly for bounded
\(\beta\in N_x\mathcal M\),
\[
\mathcal F_h^{\mathrm{STSP}}(x,\beta)
=
\frac1{\lambda_h^2}
\left[
a_{\lambda_h}\beta
-
\frac d2
\left(
c_{\lambda_h}
+
\lambda_h^2a_{\lambda_h}
\right)
H(x)
\right]
+
O_{C^2}(h^2).
\]
For each \(\lambda\in\Lambda_{\mathrm{STSP}}\), the corresponding
reference normal equation has the unique solution
\[
\beta_0(x,\lambda)
=
\frac d2
\left(
\lambda^2+
\frac{c_{\lambda}}{a_{\lambda}}
\right)
H(x).
\]
Furthermore,
\[
D_\beta
\mathcal F_0^{\mathrm{STSP}}
\bigl(
x,\beta_0(x,\lambda);\lambda
\bigr)
=
\frac{a_{\lambda}}{\lambda^2}
I_{N_x\mathcal M},
\]
and therefore
\[
\inf_{\substack{x\in\mathcal M\\
\lambda\in\Lambda_{\mathrm{STSP}}}}
\sigma_{\min}
\left[
D_\beta
\mathcal F_0^{\mathrm{STSP}}
\bigl(
x,\beta_0(x,\lambda);\lambda
\bigr)
\right]
>0.
\]
\end{lemma}

The \(C^2\) convergence and uniform nondegeneracy of the scaled normal
equation imply the existence of a unique \(C^2\) local fixed-point branch.
For later use, we record explicitly the uniform parity consequence that
removes any cubic normal remainder while retaining the full dependence on
\(\lambda_h=\sigma/h\).

\begin{lemma}[Absence of a cubic normal remainder]
\label{lem:stsp_cubic_cancellation}
Suppose Assumptions~\ref{ass:geometry-observation} and
\ref{ass:stsp-pop} hold. Let the same noise level \(\sigma\) be used while
\(h\) varies over admissible bandwidths satisfying
\(\lambda_h=\sigma/h\in\Lambda_{\mathrm{STSP}}\). If a \(C^2\)
population STSP fixed-point section is written in normal coordinates as
\(y_h(x)=x+h^2\beta_h(x), \quad \sup_h\|\beta_h\|_{C^2}<\infty\),
then
\[
\beta_h(x)
=
\frac d2
\left(
\lambda_h^2
+
\frac{c_{\lambda_h}}{a_{\lambda_h}}
\right)
H(x)
+
O_{C^2}(h^2).
\]
Equivalently,
\[
y_h(x)
=
x
+\frac d2h^2
\left(
\lambda_h^2
+
\frac{c_{\lambda_h}}{a_{\lambda_h}}
\right)
H(x)
+
O_{C^2}(h^4).
\]
Thus, after retaining the full scale dependence through \(\lambda_h\),
there is no additional cubic normal remainder.
\end{lemma}

We can now characterize the entire local population fixed-point set.

\begin{proposition}[Population STSP fixed-point manifold and bias expansion]
\label{thm:stsp_bias_expansion}
Suppose Assumptions~\ref{ass:geometry-observation} and
\ref{ass:stsp-pop} hold. Then there exist constants
\(R<\infty\) and \(h_0>0\) such that, for \(0<h<h_0\), the population
STSP fixed points in
\(\left\{ x+v: x\in\mathcal M,\; v\in N_x\mathcal M,\; \|v\|\le Rh^2 \right\}\)
form a unique \(C^2\) normal graph
\(\mathcal M_h^{\mathrm{STSP}} = \left\{ X_h^{\mathrm{STSP}}(x): x\in\mathcal M \right\}\),
where
\(X_h^{\mathrm{STSP}}(x) = x+u_h^{\mathrm{STSP}}(x)\).
Moreover,
\[
u_h^{\mathrm{STSP}}(x)
=
c_{\mathrm{STSP}}(h,\sigma)H(x)
+
O_{C^2}(h^4),
\]
uniformly over \(x\in\mathcal M\), where
\[
c_{\mathrm{STSP}}(h,\sigma)
=
\frac d2h^2
\left[
\left(\frac{\sigma}{h}\right)^2
+
\frac{c_{\sigma/h}}{a_{\sigma/h}}
\right].
\]

In addition,
\(\mathcal M_h^{\mathrm{STSP}}\) is a \(C^2\) embedded manifold
\(C^2\)-diffeomorphic to \(\mathcal M\).
\end{proposition}

Proposition~\ref{thm:stsp_bias_expansion} verifies the factorized
normal-bias expansion with
$b(x)=H(x)$ and $c(h,\sigma)=c_{\mathrm{STSP}}(h,\sigma)$.
For scales $h_1<\cdots<h_K$, choose weights satisfying
\[
\sum_{k=1}^K w_k=1,
\qquad
\sum_{k=1}^K w_k c_{\mathrm{STSP}}(h_k,\sigma)=0.
\]
Using $\mathcal M_{h_1}^{\mathrm{STSP}}$ as the reference manifold,
we define the estimated normal space by
\(
\widetilde N^{\mathrm{STSP}}(x)
=
\operatorname{Ran}Q_{h_1}
\bigl(X_{h_1}^{\mathrm{STSP}}(x)\bigr).
\)
Let $Z_{h_k}^{\mathrm{STSP}}(x)$ denote the corresponding
normal-fiber intersections with
$\mathcal M_{h_k}^{\mathrm{STSP}}$, and define
\[
\widetilde X^{\mathrm{STSP}}(x)
=
\sum_{k=1}^K w_k Z_{h_k}^{\mathrm{STSP}}(x),
\qquad
\widetilde{\mathcal M}^{\mathrm{STSP}}
=
\widetilde X^{\mathrm{STSP}}(\mathcal M).
\]

Proposition~\ref{thm:stsp_bias_expansion} gives
$\|u_{h_1}^{\mathrm{STSP}}\|_{C^2}=O(h_1^2)$,
while the PCA projector expansion in the proof of
Lemma~\ref{lem:stsp_scaled_normal_equation} yields
\(
\bigl\|
P_{\widetilde N^{\mathrm{STSP}}}
-
P_{N\mathcal M}
\bigr\|_{C^1}
\le Ch_1^2.
\)
Thus, Proposition~\ref{thm:normal_guided_manifold_recovery}
applies and gives
\[
d_H(\widetilde{\mathcal M}^{\mathrm{STSP}},\mathcal M)
\le
\sup_{x\in\mathcal M}
\|\widetilde X^{\mathrm{STSP}}(x)-x\|
\le C\sum_{k=1}^K |w_k|h_k^4.
\]
For fixed $K$, uniformly bounded weights, and sufficiently small
$h_{\max}$, the corrected map is a $C^1$ embedding whose
image is $C^1$-diffeomorphic to $\mathcal M$.
In particular, for comparable scales,
$d_H(\widetilde{\mathcal M}^{\mathrm{STSP}},\mathcal M)
=O(h_{\max}^4)$.

\begin{proposition}[Failure of squared-bandwidth extrapolation
under Gaussian noise]
\label{prop:stsp-naive}
Suppose Assumptions~\ref{ass:geometry-observation} and~\ref{ass:stsp-pop} hold, with
$\{1/3,1/6\}\subset\Lambda_{\mathrm{STSP}}$.
Let $h_1=3\sigma$ and $h_2=6\sigma$, and let
$Z_{h_k}(x)$, $k=1,2$, be the population STSP
correspondence points obtained either from the exact
normal fibers or from the estimated-normal construction
in Section~\ref{sec:normal_guided_extrapolation}, using the reference PCA normal field
at scale $h_1$.

Consider the squared-bandwidth extrapolation
\[
\widetilde X_{\mathrm{sq}}(x)
=\frac43 Z_{h_1}(x)-\frac13 Z_{h_2}(x),
\qquad
\widetilde{\mathcal M}_{\mathrm{sq}}
=\widetilde X_{\mathrm{sq}}(\mathcal M).
\]
Its weights satisfy
\(
\sum_{k=1}^2 w_k=1,
\sum_{k=1}^2 w_kh_k^2=0.
\)
Then, uniformly over $x\in\mathcal M$,
\[
\widetilde X_{\mathrm{sq}}(x)-x
=\Gamma_{d,D}\sigma^2H(x)+O(\sigma^4),
\]
where
\[
\Gamma_{d,D}
=\frac d2\left(
1+12\frac{c_{1/3}}{a_{1/3}}
 -12\frac{c_{1/6}}{a_{1/6}}
\right)
>\frac{d(d+2)}{2(d+14)}>0.
\]
Consequently,
\[
d_H(\widetilde{\mathcal M}_{\mathrm{sq}},\mathcal M)
=\Gamma_{d,D}\|H\|_\infty\sigma^2
+O(\sigma^4)
=\Theta(\sigma^2).
\]
\end{proposition}

\paragraph*{Finite-sample consistency of GBC-STSP}
We next instantiate the finite-sample estimated-normal theory for STSP.

For \(y\in\mathbb R^D\), define the empirical kernel mass, local mean,
and covariance by
\[
\widehat m_{0,h}(y)=\frac1n\sum_{i=1}^n K_h(y,Y_i),\qquad
\widehat\mu_{h}(y)
=
\frac{\sum_{i=1}^nY_iK_h(y,Y_i)}{\sum_{i=1}^nK_h(y,Y_i)},
\]
\[
\widehat\Sigma_{h}(y)
=
\frac{\sum_{i=1}^n
\bigl(Y_i-\widehat\mu_{h}(y)\bigr)
\bigl(Y_i-\widehat\mu_{h}(y)\bigr)^\top K_h(y,Y_i)}
{\sum_{i=1}^nK_h(y,Y_i)}.
\]
Let \(\widehat P_{h}(y)\) be the orthogonal projector onto the leading
\(d\)-dimensional eigenspace of \(\widehat\Sigma_{h}(y)\), and set
\(\widehat Q_{h}(y)=I-\widehat P_{h}(y)\).
The empirical STSP operator and residual are
\[
\widehat T_h^{\mathrm{STSP}}(y)
=\widehat\mu_h(y)+\widehat P_h(y)\bigl(y-\widehat\mu_h(y)\bigr),
\qquad
\widehat R_h^{\mathrm{STSP}}(y)
=\widehat Q_h(y)\bigl(y-\widehat\mu_h(y)\bigr).
\]

On every open region where the empirical mass and the empirical
\(d\)-th eigengap are positive, these moments, spectral projectors, and
residuals are at least \(C^2\). For a fixed sufficiently large \(R\), put
\(\mathcal T_h = \{x+v:x\in\mathcal M,\ v\in N_x\mathcal M,\ \|v\|\leq Rh^2\}\).
The empirical fitted set is the local zero set of the computable residual,
\[
\widehat{\mathcal M}_h^{\mathrm{STSP}}
=
\{y\in\mathcal T_h:\widehat R_{h}^{\mathrm{STSP}}(y)=0\}.
\]
The tube identifies the branch associated with \(\mathcal M\) for the
analysis; neither the true normal coordinates nor \(\mathcal M\) are
inputs to the residual calculation. For convenience, define
\(\eta_n(h)=\sqrt{\frac{\log n}{nh^d}}\).

\begin{proposition}[Empirical STSP branch existence and stability]
\label{lem:stsp_empirical_rate}
Suppose Assumptions~\ref{ass:geometry-observation} and
\ref{ass:stsp-pop} hold, and assume
\(
\eta_n(h)=o(h), \text{equivalently},
\frac{nh^{d+2}}{\log n}\longrightarrow\infty.
\)
Choose \(R\) large enough to contain the population branch in the interior
of \(\mathcal T_h\). For sufficiently small \(h\), with probability tending
to one, uniformly on \(\mathcal T_h\),
\[
\widehat m_{0,h}(y)\geq ch^d,\qquad
\lambda_d\!\left(\widehat\Sigma_{h}(y)\right)
-\lambda_{d+1}\!\left(\widehat\Sigma_{h}(y)\right)\geq ch^2.
\]
On the same event, the entire local zero set is a unique \(C^2\) normal
graph
\[
\widehat{\mathcal M}_h^{\mathrm{STSP}}
=\{\widehat X_h^{\mathrm{STSP}}(x):x\in\mathcal M\},\qquad
\widehat X_h^{\mathrm{STSP}}(x)=x+\widehat u_h^{\mathrm{STSP}}(x),\quad
\widehat u_h^{\mathrm{STSP}}(x)\in N_x\mathcal M,
\]
and \(\widehat X_h^{\mathrm{STSP}}\) is a \(C^2\) embedding. Moreover,
\[
\|\widehat u_h^{\mathrm{STSP}}-u_h^{\mathrm{STSP}}\|_\infty
=O_p\!\left(h\eta_n(h)\right),\qquad
\|\widehat u_h^{\mathrm{STSP}}-u_h^{\mathrm{STSP}}\|_{C^1}
=O_p\!\left(\eta_n(h)\right).
\]

Consequently, on the same high-probability event,
the uniform empirical approximation condition
holds with the deterministic bound
\(
r_n^{\mathrm{STSP}}(h)=Ch\eta_n(h).
\)
The empirical tangent spaces satisfy
\[
\sup_{x\in\mathcal M}
\left\|
P_{T_{\widehat X_h^{\mathrm{STSP}}(x)}
\widehat{\mathcal M}_h^{\mathrm{STSP}}}
-P_{T_x\mathcal M}
\right\|_{\mathrm{op}}
=O_p\!\left(h^2+\eta_n(h)\right).
\]
\end{proposition}

For the estimated-normal correspondence, estimate the normal space at the
reference scale \(h_1\) by kernel-weighted local PCA. Define the population
and empirical reference projector fields by
\(
\widetilde Q_1(x)
=Q_{h_1}\!\left(X_{h_1}^{\mathrm{STSP}}(x)\right),
\widehat Q_1(x)
=\widehat Q_{h_1}\!\left(\widehat X_{h_1}^{\mathrm{STSP}}(x)\right),
\)
and put
\(\widetilde N_{h_1}^{\mathrm{STSP}}(x)=\operatorname{Ran}\widetilde Q_1(x)
=\widetilde N^{\mathrm{STSP}}(x)\) and
\(\widehat N_{h_1}^{\mathrm{STSP}}(x)=\operatorname{Ran}\widehat Q_1(x)\).
The kernel-weighted population projector expansion gives
\(\|\widetilde Q_1-P_{N\mathcal M}\|_{C^1}\leq Ch_1^2\),
where \(P_{N\mathcal M}\) denotes the field \(x\mapsto P_{N_x\mathcal M}\).
The differentiated moment and spectral-projector estimates in the proof of
Proposition~\ref{lem:stsp_empirical_rate} also give
\[
\|\widehat Q_1-\widetilde Q_1\|_\infty
=O_p\!\left(\eta_n(h_1)\right),\qquad
\|\widehat Q_1-\widetilde Q_1\|_{C^1}
=O_p\!\left(\frac{\eta_n(h_1)}{h_1}\right).
\]
Thus the empirical reference field is \(C^1\), its
\(C^1\) norm is bounded with probability tending
to one under \(\eta_n(h_1)=o(h_1)\), and the
finite-sample normal-estimation condition holds
on the corresponding high-probability event with
the deterministic bound
\(s_n^{\mathrm{STSP}}(h_1)=C\eta_n(h_1)\).

Using
\(\widehat{\mathcal M}_{h_1}^{\mathrm{STSP}}\)
as the reference estimator and
\(\widehat N_{h_1}^{\mathrm{STSP}}\)
as the estimated normal space, let
\(
\widehat Z_{h_k}^{\mathrm{STSP}}(x)
\in
\widehat{\mathcal M}_{h_k}^{\mathrm{STSP}}
\)
denote the empirical estimated-normal correspondence point, and let
\(
Z_{h_k}^{\mathrm{STSP}}(x)
\in
\mathcal M_{h_k}^{\mathrm{STSP}}
\)
denote its population counterpart. For fixed \(K\), if
\(\max_k\eta_n(h_k)/h_k\to0\), the empirical branch and projector bounds
above imply that these local intersections exist uniquely, are transverse,
and depend \(C^1\)-smoothly on the reference point, with probability tending
to one. The argument is included in the proof of
Theorem~\ref{thm:finite_sample_stsp}. The finite-sample correspondence bound
then gives
\[
\sup_{x\in\mathcal M}
\left\|
\widehat Z_{h_k}^{\mathrm{STSP}}(x)
-
Z_{h_k}^{\mathrm{STSP}}(x)
\right\|
=
O_p
\left(
h_1\eta_n(h_1)
+
h_k\eta_n(h_k)
+
h_k^2\eta_n(h_1)
\right).
\]

Define the finite-sample GBC-STSP map by
\(
\widehat X^{\mathrm{STSP}}(x)
=
\sum_{k=1}^K
w_k\widehat Z_{h_k}^{\mathrm{STSP}}(x),
\)
and
\(
\widehat{\mathcal M}^{\mathrm{STSP}}
=
\widehat X^{\mathrm{STSP}}(\mathcal M).
\)
This formal estimator is the continuous image of the smooth empirical zero
branches, not a finite point cloud. On the event where these branches and
correspondences are defined, its reference-branch parameterization is
\[
\widehat F^{\mathrm{STSP}}
=\widehat X^{\mathrm{STSP}}\circ
 \bigl(\widehat X_{h_1}^{\mathrm{STSP}}\bigr)^{-1},
\qquad
\widehat{\mathcal M}^{\mathrm{STSP}}
=\widehat F^{\mathrm{STSP}}
 \bigl(\widehat{\mathcal M}_{h_1}^{\mathrm{STSP}}\bigr).
\]

\begin{theorem}[Finite-sample recovery of GBC-STSP]
\label{thm:finite_sample_stsp}
Suppose the conditions of Proposition~\ref{thm:stsp_bias_expansion} hold.
Let the number \(K\) of fitting scales be fixed, assume
\(\max_k\eta_n(h_k)/h_k\to0\), and let the
deterministic extrapolation weights satisfy
\(\sum_{k=1}^K|w_k|\leq W_0<\infty\) and
\(
\sum_{k=1}^K w_k=1,
\sum_{k=1}^K
w_kc_{\mathrm{STSP}}(h_k,\sigma)=0.
\)
Then
\[
d_H
\left(
\widehat{\mathcal M}^{\mathrm{STSP}},
\mathcal M
\right)
=
O_p
\left(
\sum_{k=1}^K
|w_k|
\left[
h_k^4
+
h_1\eta_n(h_1)
+
h_k\eta_n(h_k)
+
h_k^2\eta_n(h_1)
\right]
\right).
\]
For sufficiently small \(h_{\max}\), the empirical local branches and
estimated-normal correspondence points are well defined with probability
tending to one, and \(\widehat{\mathcal M}^{\mathrm{STSP}}\) is a
\(C^1\) embedded manifold \(C^1\)-diffeomorphic to \(\mathcal M\).

For the numerical implementation, suppose that the supplied reference
branch is the branch \(\widehat{\mathcal M}_{h_1}^{\mathrm{STSP}}\)
identified above. Let
\(\mathcal Z_N=\{z_1,\ldots,z_N\}\subset
\widehat{\mathcal M}_{h_1}^{\mathrm{STSP}}\)
be a nonempty set of retained query points, with fill distance
\(\delta_N =\sup_{z\in\widehat{\mathcal M}_{h_1}^{\mathrm{STSP}}} \min_{1\leq j\leq N}\|z-z_j\|\).
Let \(q=D-d\), and choose an orthonormal basis
\(\widehat E(z)\in\mathbb R^{D\times q}\) of the reference PCA normal
space at each query point. For \(k=2,\ldots,K\), define
\(g_{k,z}(t) =\widehat E(z)^\top\widehat R_{h_k}^{\mathrm{STSP}} \bigl(z+\widehat E(z)t\bigr)\).
Fix a sufficiently large constant \(C_{\mathrm{loc}}\), independent of
\(n\) and the scales. If numerical solutions \(\widehat t_k(z)\) satisfy
\[
\|\widehat t_k(z)\|\leq C_{\mathrm{loc}}h_k^2,\qquad
\|g_{k,z}(\widehat t_k(z))\|\leq\varepsilon_{n,k}
\]
for all \(z\in\mathcal Z_N\), for deterministic tolerances
\(\varepsilon_{n,k}\geq0\), put \(\widehat t_1(z)=0\),
\(\varepsilon_{n,1}=0\), and define the finite numerical output
\[
\widehat{\mathcal P}_N^{\mathrm{STSP,num}}
=\left\{
z_j+\widehat E(z_j)\sum_{k=1}^K w_k\widehat t_k(z_j):
1\leq j\leq N
\right\}.
\]
Then, on the same high-probability event,
\[
\begin{aligned}
d_H\!\left(\widehat{\mathcal P}_N^{\mathrm{STSP,num}},\mathcal M\right)
&\leq
d_H\!\left(\widehat{\mathcal M}^{\mathrm{STSP}},\mathcal M\right)
+C\sum_{k=1}^K|w_k|\varepsilon_{n,k}+C\delta_N.
\end{aligned}
\]
For comparable scales,
\(\max_k\varepsilon_{n,k}=O(h_{\max}^4)\) and
\(\delta_N=O(h_{\max}^4)\) keep the numerical and covering errors within
the fourth-order deterministic bound. The embedding conclusion above
concerns the continuous exact-intersection estimator, not the finite
point cloud; the residual tolerance alone does not imply regularity of
a numerical interpolation of that cloud.
\end{theorem}

\begin{corollary}[A sufficient sample size for fourth-order STSP accuracy]
\label{cor:stsp_fourth_order_sample_size}
Suppose the assumptions of Theorem~\ref{thm:finite_sample_stsp}, other
than its sample-size condition, hold. Let
\(h_k=a_kh\), where \(0<a_1<\cdots<a_K<\infty\) are fixed and
\(h\to0\) in the admissible noise--bandwidth regime. There exist constants
\(C_0,C_1>0\) such that, if
\(
 n h^{d+6}\geq C_0\log n,
\)
then, for sufficiently small \(h\) and large \(n\), with probability
at least \(1-Kn^{-2}\),
\[
 d_H\!\left(\widehat{\mathcal M}^{\mathrm{STSP}},\mathcal M\right)
 \leq C_1h^4,
\]
and \(\widehat{\mathcal M}^{\mathrm{STSP}}\) is a \(C^1\) embedded
manifold \(C^1\)-diffeomorphic to \(\mathcal M\).
In particular, when \(h\asymp\sigma\to0\), there exist constants
\(C_0',C_1'>0\) such that the sufficient sample-size condition
\[
 n\geq C_0'\sigma^{-(d+6)}\log(1/\sigma)
\]
yields
\(d_H(\widehat{\mathcal M}^{\mathrm{STSP}},\mathcal M)
\leq C_1'\sigma^4\) on the same event.
Under the numerical assumptions of
Theorem~\ref{thm:finite_sample_stsp}, if also
\(\max_k\varepsilon_{n,k}=O(h^4)\) and \(\delta_N=O(h^4)\), then
\(d_H\!\left(\widehat{\mathcal P}_N^{\mathrm{STSP,num}},\mathcal M\right) \leq Ch^4\)
on that event.
\end{corollary}

\begin{remark}
Under the preceding assumptions and $\eta_n(h)=o(h)$,
the smooth kernel and differentiated concentration bounds
ensure the existence and $C^1$ regularity of the empirical
branch and PCA-guided correspondence. Consequently, the
continuous GBC-STSP estimator is a $C^1$ embedded manifold
diffeomorphic to $\mathcal M$.
True normal coordinates are used only in the analysis;
numerical implementation requires initialization on the
branch associated with $\mathcal M$.
\end{remark}

\paragraph*{Estimated noise level}
For fixed deterministic bandwidths, select one point
\(\widehat y_1\in\widehat{\mathcal M}_{h_1}^{\mathrm{STSP}}\) on the
reference branch. Let \(\widehat\ell_1\geq\cdots\geq\widehat\ell_D\)
be the eigenvalues of \(\widehat\Sigma_{h_1}(\widehat y_1)\), and put
\(q=D-d\). The calibration function \(a(\lambda)=a_\lambda\) is
strictly increasing. Choose a fixed compact interval
\(\Lambda_+\Subset(0,\infty)\) containing \(\Lambda_{\mathrm{STSP}}\),
and write \(\Pi_I\) for projection onto an interval \(I\). Define
\(
\widehat a_1=\frac{1}{qh_1^2}\sum_{j=d+1}^D\widehat\ell_j,
\widehat\sigma
=h_1a^{-1}\!\left(\Pi_{a(\Lambda_+)}(\widehat a_1)\right).
\)
Under the conditions of Proposition~\ref{lem:stsp_empirical_rate},
Proposition~\ref{prop:stsp_noise_estimation} in Appendix~\ref{sec:app_stsp} gives
\(
|\widehat\sigma-\sigma|
=O_p\!\left(h_1^3+h_1\eta_n(h_1)\right).
\)
Use this single estimate in spatially constant plug-in weights satisfying
Proposition~\ref{prop:noise_level_perturbation}. For fixed comparable
scales and uniformly bounded weights, the error rate and embedding
conclusion of Theorem~\ref{thm:finite_sample_stsp} remain valid.
In particular, Corollary~\ref{cor:stsp_fourth_order_sample_size}
retains its fourth-order guarantee without a stronger sample-size
condition. 

The preceding results establish fourth-order recovery
at prescribed fitting scales. A further statistical
question is whether this improvement can be reproduced
by selecting a better bandwidth for the original
single-scale estimator.
The following corollary addresses this question
within the STSP family. Rather than comparing
against a fixed bandwidth, it uses an oracle that
minimizes the realized Hausdorff error over an
entire admissible bandwidth interval.

\begin{corollary}[Oracle bandwidth-tuned rate separation for STSP]
\label{cor:stsp_oracle_separation}
Suppose Assumptions~\ref{ass:geometry-observation} and~\ref{ass:stsp-pop} hold. Let
$\sigma=\sigma_n\to0$ and
\(
    \frac{n\sigma^{d+2}}{\log n}\longrightarrow\infty.
\)
Fix $0<a<b<\infty$ with
$[1/b,1/a]\subset\Lambda_{\mathrm{STSP}}$,
and write $\mathcal H_n=[a\sigma,b\sigma]$.
For each $h\in\mathcal H_n$, let
$\widehat{\M}_h^{\mathrm{STSP}}$ denote the local
empirical STSP branch of Proposition~\ref{lem:stsp_empirical_rate},
constructed using one common sample. Define
\[
    \gamma(t):=\frac d2
    \left(1+t^2\frac{c_{1/t}}{a_{1/t}}\right),
    \qquad
    \gamma_*:=\min_{t\in[a,b]}\gamma(t)>0,
    \qquad
    \varepsilon_n:=
    \sqrt{\frac{\log n}{n\sigma^d}}.
\]
Then the oracle single-scale error satisfies
\begin{equation}\label{eq:stsp_oracle_single_error}
\begin{aligned}
    \inf_{h\in\mathcal H_n}
    d_H(\widehat{\M}_h^{\mathrm{STSP}},\M)
    &=
    \gamma_*\|H\|_\infty\sigma^2
    +O_p(\sigma^4+\sigma\varepsilon_n)\\
    &=
    \gamma_*\|H\|_\infty\sigma^2
    \{1+o_p(1)\}.
\end{aligned}
\end{equation}

Choose any fixed $a<a_1<a_2<b$, put
$h_k=a_k\sigma$ for $k=1,2$, and let
\(
    w_1=
    \frac{\gamma(a_2)}
         {\gamma(a_2)-\gamma(a_1)},
    w_2=
    -\frac{\gamma(a_1)}
          {\gamma(a_2)-\gamma(a_1)}.
\)
The denominators are nonzero, and these
weights are independent of $n$ and satisfy
the cancellation constraints of Theorem~\ref{thm:finite_sample_stsp}.
Let $\widehat{\M}^{\mathrm{STSP}}$
be the corresponding
GBC-STSP estimator. Then
\begin{equation}\label{eq:stsp_oracle_rate_separation}
\begin{aligned}
    &d_H(\widehat{\M}^{\mathrm{STSP}},\M)
    =O_p(\sigma^4+\sigma\varepsilon_n),\\
    &
    \frac{
      d_H(\widehat{\M}^{\mathrm{STSP}},\M)
    }{
      \displaystyle\inf_{h\in\mathcal H_n}
      d_H(\widehat{\M}_h^{\mathrm{STSP}},\M)
    }
    =
    O_p\left(
      \sigma^2+
      \sqrt{\frac{\log n}{n\sigma^{d+2}}}
    \right)
    =o_p(1).
\end{aligned}
\end{equation}
\end{corollary}

Corollary~\ref{cor:stsp_oracle_separation} shows that the gain from multiscale cancellation
cannot be reproduced by bandwidth selection alone within the
specified STSP family. More importantly, it identifies two
distinct statistical regimes. Under
\(
\frac{n\sigma^{d+2}}{\log n}\longrightarrow\infty,
\)
the GBC-STSP estimator achieves an asymptotically vanishing
error ratio relative to the oracle bandwidth-tuned single-scale
estimator, without requiring fourth-order stochastic accuracy.
The stronger condition
\(
n\sigma^{d+6}\gtrsim\log n
\)
is sufficient for fourth-order Hausdorff recovery, as established
in Corollary~\ref{cor:stsp_fourth_order_sample_size}. Thus, the statistical advantage of geometric
bias cancellation emerges before the full fourth-order accuracy
becomes attainable. This comparison concerns STSP over the
admissible bandwidth range $h\asymp\sigma$ and does not imply
minimax optimality among general manifold-fitting procedures.

\subsection{Gaussian Density-Ridge Fitting}
\label{sec:dr_validation}

We next verify the geometric bias cancellation framework
for Gaussian density-ridge fitting (DR), a structurally
different procedure from STSP. Rather than relying on
local averaging and PCA-based projection, density ridges
are characterized by the gradient and Hessian of a
smoothed density. Under Gaussian observation noise,
the effective smoothing variance is
$\rho_h^2=h^2+\sigma^2$.

Throughout this subsection, the superscript
$\mathrm{DR}$ identifies density-ridge quantities.
We use $h$ for the kernel bandwidth and
$\rho_h=(h^2+\sigma^2)^{1/2}$ for the effective
smoothing scale. The population bias coefficient
is denoted by $c_{\mathrm{DR}}(h,\sigma)$,
with empirical quantities following the
hat convention of Section~\ref{Multiscale Bias Cancellation}.

Under Assumption~\ref{ass:geometry-observation} and
the additional condition that $\mathcal M$ is $C^6$,
the population ridge branch associated with
$\mathcal M$ admits a $C^2$ normal-graph representation
satisfying
\[
u_{h,\sigma}^{\mathrm{DR}}(x)
=
\frac d2(h^2+\sigma^2)H(x)
+O_{C^2}(h^4),
\]
uniformly over $x\in\mathcal M$. Thus, the factorized
normal-bias expansion holds with
\(
b(x)=H(x),
c_{\mathrm{DR}}(h,\sigma)
=\frac d2(h^2+\sigma^2).
\)

For admissible scales $h_1<\cdots<h_K$, choose
bounded extrapolation weights satisfying
\(
\sum_{k=1}^K w_k=1,
\sum_{k=1}^K w_k(h_k^2+\sigma^2)=0.
\)
Applying the estimated-normal correspondence of
Section~\ref{sec:normal_guided_extrapolation}
yields fourth-order population Hausdorff accuracy.
For fixed $K$, scales $h_k\asymp h$, and uniformly
bounded weights, the empirical GBC-DR estimator
satisfies
\[
d_H\bigl(
\widehat{\mathcal M}^{\mathrm{DR}},
\mathcal M
\bigr)
=
O_p\left(
h^4+h\sqrt{\frac{\log n}{nh^d}}
\right).
\]
In particular, $nh^{d+6}\ge C\log n$ for a sufficiently
large constant $C$ is sufficient for fourth-order
accuracy, and the empirical estimator is a
$C^1$ embedded manifold diffeomorphic to
$\mathcal M$ with high probability.

The density-ridge construction provides a second
verification of the geometric bias-cancellation
principle through a different population fitting
equation. Unlike STSP, whose leading bias coefficient
involves kernel-dependent moment ratios, DR admits
an explicit coefficient determined by the total
smoothing variance $h^2+\sigma^2$.
In particular, conventional squared-bandwidth
extrapolation with $\sum_k w_kh_k^2=0$
leaves the leading error
$\frac d2\sigma^2 H(x)$, whereas cancellation
based on $c_{\mathrm{DR}}(h,\sigma)$ removes
the entire second-order term.
The complete population expansion, finite-sample
analysis, and geometric guarantees for DR
are provided in Appendix~\ref{sec:gaussian_density_ridge}.

\paragraph*{Numerical conditioning of ENF correspondence}
For both STSP and DR, the implementable ENF
correspondence equations are locally well conditioned.
Let $z$ lie on the empirical reference branch
and let $E(z)\in\mathbb R^{D\times q}$ be an
orthonormal estimated-normal basis, where $q=D-d$.
Define
\(
g_{k,z}(t)
=
E(z)^\top\widehat R_{h_k}(z+E(z)t),
\eta_n(h)=\sqrt{\frac{\log n}{nh^d}}.
\)
The differentiated residual estimates in
Appendices~\ref{app:stsp_finite_sample_proof} and~\ref{app:dr_finite_sample_proof} yield,
on the corresponding high-probability events,
\[
\|D_tg_{k,z}(t)-\gamma_k I_q\|_{\mathrm{op}}
\le
C\{h_k+h_1^2+\eta_n(h_k)+\eta_n(h_1)\}
=o(1),
\]
uniformly over the prescribed local search regions,
where $\gamma_k=a_{\sigma/h_k}/(\sigma/h_k)^2$ for STSP
and $\gamma_k=h_k^2/(h_k^2+\sigma^2)$ for DR.
In the respective admissible regimes,
$\inf_k\gamma_k\ge\gamma_{\mathrm{ENF},*}>0$.
Consequently, the Jacobian condition number
converges to one, and numerical correspondence
errors satisfy
\[
\|\widehat t_k-t_k^*\|
\le
\frac{2}{\gamma_{\mathrm{ENF},*}}
\|g_{k,z}(\widehat t_k)\|.
\]
Thus, local numerical residuals directly certify
correspondence accuracy. 

\section{Directional Geometric Bias Cancellation for the YY Estimator}
\label{sec:yao_contraction}

We consider the YY estimator with smooth compactly
supported kernels replacing the original truncated
Gaussian-type weights. Its two-stage structure is
retained: a spherical local mean determines the
contraction axis, followed by a cylindrical local
average. Across scales, we keep the preliminary point
and axis fixed and vary only the transverse radius.

In higher codimension, the leading normal bias separates
into axial and complementary components with distinct
scale dependence, which cannot, in general, be canceled
simultaneously by scalar extrapolation. This motivates
directional GBC using separate weights
for the two components. Under the preliminary-manifold
conditions in Assumption~\ref{ass:ysl_full_population},
we establish population and finite-sample guarantees
for the resulting refinement.
Throughout, $1\leq d<D$, $q=D-d$, $m=q-1$,
and $L_\sigma=\log(1/\sigma)$. Here $q$ denotes the codimension and $m$ the
dimension of the complementary normal subspace.
We use $r_0$, $r_{1,k}$, and $r_2$ for the
first-stage radius, the $k$th transverse radius,
and the longitudinal radius, respectively.
These radii are specific to the YY construction and are distinct from the fitting
bandwidths used in Sections~\ref{Multiscale Bias Cancellation} and~\ref{sec:examples}.

\subsection{Higher-Order Normal Expansion}
\label{sec:ysl_full_single_scale}

\begin{assumption}[Population regime]
\label{ass:ysl_full_population}
\begin{enumerate}
\item 
\label{item:ysl_population_geometry}
The underlying manifold \(\mathcal M\) is connected and \(C^6\).
\item
The common preliminary manifold \(\mathcal M^{\mathrm{pre}}\) is a nonempty,
compact, connected, boundaryless, \(C^6\) embedded \(d\)-manifold satisfying
\(
\mathcal M^{\mathrm{pre}}
\subset\operatorname{Tub}_{C_0\sigma}(\mathcal M),
\operatorname{reach}(\mathcal M^{\mathrm{pre}})\geq\tau_*>0.
\)
The constants \(C_0,\tau_*\) do not depend on \(\sigma\). 
\item
\label{item:ysl_population_radii}
For a fixed integer \(p\geq4\), the first-stage and longitudinal radii are
\(
 r_0=c_0\sigma,
 r_2=C_2\sigma\sqrt{L_\sigma},
\)
where \(c_0>0\) is fixed and \(C_2\) is sufficiently large.  The transverse
radii are spatially constant and satisfy
\(
 c_-\sigma\leq r_1\leq c_+\sigma,
0<c_-<c_+<\infty.
\)
\item
For the spherical population mean defined below, the first-stage signal
satisfies the independent nondegeneracy condition
\(
\label{eq:ysl_common_signal}
 \inf_{z\in\mathcal M^{\mathrm{pre}}}
 \|F_{r_0}(z)-z\|\geq c_F\sigma,
 c_F>0.
\)
\end{enumerate}
\end{assumption}
\begin{remark}[Preliminary manifold and first-stage nondegeneracy]
\label{first-stage signal}
The preliminary manifold is an input to the refinement procedure.
An \(O(\sigma)\)-accurate pilot with the required smoothness and positive
reach is not sufficient by itself: the first-stage signal condition in
Assumption~\ref{ass:ysl_full_population} must also hold. This independent
condition controls the normalization of \(F_{r_0}(z)-z\) in
\eqref{eq:ysl_common_first_stage}.
The first-stage expansion in the proof of
Lemma~\ref{lem:ysl_common_first_stage_bounds} shows that this regime requires
\(c\sigma \leq d(z,\mathcal M)\leq C\sigma, z\in\mathcal M^{\mathrm{pre}}\),
uniformly for sufficiently small \(\sigma\). Thus, a preliminary manifold
that is uniformly \(o(\sigma)\)-close to \(\mathcal M\), without an
additional offset, is outside the scope of the present theorem.

This order-\(\sigma\) separation can be obtained, for example, by a
normal-offset construction. Suppose a smooth pilot manifold
\(\mathcal P_\sigma\) satisfies
\(
    d_H(\mathcal P_\sigma,\mathcal M)\leq B\sigma,
    \operatorname{reach}(\mathcal P_\sigma)\geq\tau_{\mathrm p}>0.
\)
For a smooth unit normal field
\(\nu_{\mathrm p}(p)\in N_p\mathcal P_\sigma\) and a fixed constant
\(A>B\), define
\(
    \mathcal M^{\mathrm{pre}}
    =
    \left\{
        p+A\sigma\nu_{\mathrm p}(p):
        p\in\mathcal P_\sigma
    \right\}.
\)
For sufficiently small \(\sigma\), this gives
\(
    (A-B)\sigma
    \leq
    d(z,\mathcal M)
    \leq
    (A+B)\sigma,
    z\in\mathcal M^{\mathrm{pre}},
\)
while the required smoothness and positive-reach properties are preserved
under the corresponding uniform regularity conditions. Together with the first-stage expansion, this suggests a sufficient
route to the nondegeneracy condition, provided the offset normal field
has uniformly controlled derivatives and the offset image retains a
uniform reach bound. Proposition~\ref{prop:yy-sample-split-pilot}
formalizes this conditional construction using an independent pilot
sample. The controlled YY experiments in Sections~\ref{sec:population-bias}--\ref{sec:exp-directional} use
preliminary manifolds constructed from the known underlying geometry.
\end{remark}
We use the compactly supported polynomial kernels
\(K_{0,r_0}\) and \(K_{\perp,r_1}\), together with a \(C^3\)
longitudinal plateau cutoff \(K_{\parallel,r_2}\).
Their definitions and regularity properties are given in
Appendix~\ref{sec:ysl_common_localization_proof}. Define
\begin{equation}
\label{eq:ysl_common_first_stage}
 F_{r_0}(z)
 =\frac{\mathbb E[Y K_{0,r_0}(\|Y-z\|)]}
         {\mathbb E[K_{0,r_0}(\|Y-z\|)]},
 \qquad
 U_0(z)=\frac{V_0(z)V_0(z)^\top}{\|V_0(z)\|^2},
 \quad V_0(z)=F_{r_0}(z)-z.
\end{equation}
The complete second-stage population mean is
\begin{equation}
\label{eq:ysl_full_population_G}
 \begin{split}
 \Psi_{r_1,z}(y)
 &=K_{\parallel,r_2}(\|U_0(z)(y-z)\|)
   K_{\perp,r_1}(\|(I_D-U_0(z))(y-z)\|),\\
 G_{r_1}(z)
 &=\frac{\mathbb E[Y\Psi_{r_1,z}(Y)]}
         {\mathbb E[\Psi_{r_1,z}(Y)]},
 \qquad
 \mathcal M_{r_1}^{\mathrm{YY}}=G_{r_1}(\mathcal M^{\mathrm{pre}}).
 \end{split}
\end{equation}

\begin{lemma}[Preliminary normal parameterization]
\label{lem:ysl_full_preliminary_bounds}
Under Assumption~\ref{ass:ysl_full_population}, for sufficiently small
\(\sigma\), the nearest-point projection from
\(\mathcal M^{\mathrm{pre}}\) to \(\mathcal M\) is a \(C^2\)
diffeomorphism.  Its inverse has the form
\(X_0(x)=x+u_0(x), u_0(x)\in N_x\mathcal M\),
where
\begin{equation}
\label{eq:ysl_common_preliminary_bounds}
 \|u_0\|_{C^0}\leq C\sigma,
 \qquad \|u_0\|_{C^1}\leq C\sigma^{1/2},
 \qquad \|u_0\|_{C^2}\leq C.
\end{equation}
In particular, \(X_0\) is uniformly bi-Lipschitz.
\end{lemma}

For the analysis, write
\(
Q=P_{N_x\mathcal M}, T=I_D-Q,
v_0=V_0(X_0(x)),
\)
and define
\(
\nu_0=\frac{Qv_0}{\|Qv_0\|},
P_0=\nu_0\nu_0^\top,
S_0=Q-P_0,
b_\perp=S_0u_0.
\)
We also write
\(\overline U_0(x)=U_0(X_0(x))\). 
Thus, $Q$ is the full normal projector,
$P_0$ is the rank-one projector onto the
first-stage normal axis, and $S_0$ is the
projector onto its orthogonal complement
within the normal space. The vector $b_\perp$
represents the complementary component of
the preliminary normal displacement and is
distinct from the generic bias field $b(x)$
used in Section~\ref{Multiscale Bias Cancellation}.

\begin{lemma}[First-stage axis and residual offset]
\label{lem:ysl_common_first_stage_bounds}
Under Assumption~\ref{ass:ysl_full_population}, uniformly for \(j=0,1,2\),
\begin{equation}
\label{eq:ysl_common_axis_bounds}
 \|P_0\|_{C^j}\leq C\sigma^{-j/2},
 \qquad
 \|b_\perp\|_{C^j}\leq C\sigma^{2-j/2},
 \qquad
 \|\overline U_0-P_0\|_{C^j}\leq C\sigma^{2-j/2}.
\end{equation}
Moreover, \(\|Tv_0\|_{C^j}\leq C\sigma^{3-j/2}\) and
\(\inf_x\|Qv_0(x)\|\geq c\sigma\).
\end{lemma}
The residual \(b_\perp\) must be retained, since it contributes at
the same order as the leading curvature bias and is canceled explicitly
across scales. 
For \(m\geq1\), let \(\alpha_m(\kappa)\) and \(\beta_m(\kappa)\)
denote the kernel moments defined in
Appendix~\ref{sec:ysl_common_moments}.

\begin{figure}[!htbp]
\centering
\includegraphics[width=0.72\linewidth]{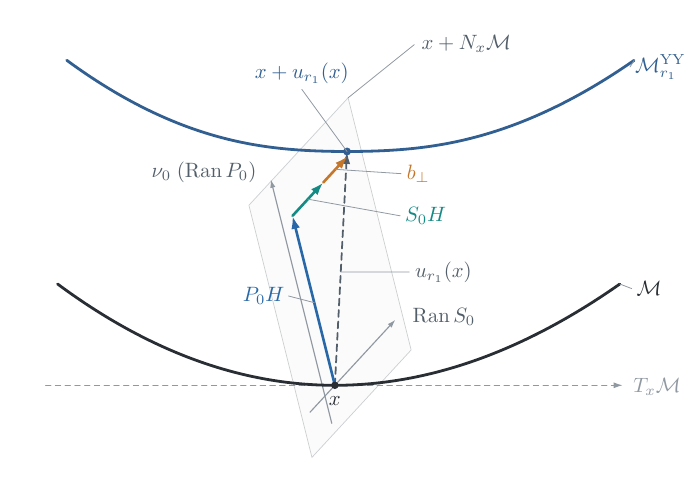}
\caption{Schematic decomposition of the leading normal
bias in the YY estimator into axial and complementary
components. The colored arrows represent the scale-weighted
contributions of $P_0H$, $S_0H$, and $b_\perp$, with
higher-order terms omitted.}
\label{fig:yy_normal_bias}
\end{figure}

\begin{proposition}[Raw single-scale expansion]
\label{thm:ysl_common_raw_expansion}
Suppose Assumption~\ref{ass:ysl_full_population} holds, and let
\(Y_{r_1}=G_{r_1}\circ X_0\).  If \(q>1\), then
\begin{equation}
\label{eq:ysl_common_raw_expansion}
 \begin{split}
 Y_{r_1}(x)-x
 ={}&t_{r_1}(x)
 +\frac d2\sigma^2\{1+\alpha_m(r_1/\sigma)\}P_0(x)H(x)\\
 &+\frac d2\sigma^2\alpha_m(r_1/\sigma)S_0(x)H(x)
 +\{1-\beta_m(r_1/\sigma)\}b_\perp(x)
 +R_{r_1}^N(x),
 \end{split}
\end{equation}
where \(H=d^{-1}\operatorname{tr}\mathrm{II}\),
\(t_{r_1}=T(Y_{r_1}-\operatorname{Id}_{\mathcal M})\), and
\(R_{r_1}^N(x)\in N_x\mathcal M\).  Uniformly over admissible radii,
\begin{equation}
\label{eq:ysl_common_raw_bounds}
\begin{aligned}
\|R_{r_1}^N\|_{C^0}&\leq C\sigma^4,\qquad
\|t_{r_1}\|_{C^j}\leq C\sigma^{3-j/2},\\
\|Q(Y_{r_1}-\operatorname{Id}_{\mathcal M})\|_{C^j}
&\leq C\sigma^{2-j/2},\qquad j=0,1,2.
\end{aligned}
\end{equation}
For \(q=1\), the same statement holds after deleting the terms involving
\(S_0\) and \(b_\perp\), with \(P_0=Q\) and \(\alpha_m=\alpha_0\).
\end{proposition}

Thus, unlike the scalar factorization in Section~\ref{Multiscale Bias Cancellation}, the leading normal bias
decomposes into the three common fields \(P_0H\), \(S_0H\), and \(b_\perp\),
with different dependence on the transverse scale \(r_1\). This motivates
the directional weighting scheme below.

\begin{corollary}[The completed single-scale image]
\label{thm:ysl_full_single_scale_section}
Under the assumptions of Proposition~\ref{thm:ysl_common_raw_expansion}, for
sufficiently small \(\sigma\), \(G_{r_1}|_{\mathcal M^{\mathrm{pre}}}\)
is a \(C^2\) embedding.  Its image is a unique normal graph
\(\mathcal M_{r_1}^{\mathrm{YY}} =\{x+u_{r_1}(x):x\in\mathcal M\}, \quad u_{r_1}(x)\in N_x\mathcal M\),
with
\begin{equation}
\label{eq:ysl_common_single_image_bounds}
 \|u_{r_1}\|_{C^j}\leq C\sigma^{2-j/2},\quad j=0,1,2,
 \qquad
 \operatorname{reach}(\mathcal M_{r_1}^{\mathrm{YY}})\geq c_*>0.
\end{equation}
The normal section \(u_{r_1}\) has the same normal leading terms as
\eqref{eq:ysl_common_raw_expansion}, with a \(C^0\) remainder bounded by
\(C\sigma^4\).
\end{corollary}

\subsection{Directional weights and common-start extrapolation}
\label{sec:ysl_apply_general_theory}

For \(q>1\), select a fixed number \(K\geq3\) of admissible transverse
radii and write
\(
 \kappa_k=r_{1,k}/\sigma,
 \alpha_k=\alpha_m(\kappa_k),
 \beta_k=\beta_m(\kappa_k).
\)
Choose spatially constant axial weights \(a_k\) and complementary weights
\(b_k\) satisfying
\begin{equation}
\label{eq:ysl_common_weights}
\begin{aligned}
&\sum_k a_k=1,\qquad \sum_k a_k(1+\alpha_k)=0,\\
&\sum_k b_k=1,\qquad \sum_k b_k\alpha_k=0,\qquad
  \sum_k b_k(1-\beta_k)=0,\\
&W:=\sum_k (|a_k|+|b_k|)\le W_0<\infty .
\end{aligned}
\end{equation}
Here $\alpha_k$ and $\beta_k$ are
scale-dependent kernel-moment coefficients,
whereas $a_k$ and $b_k$ denote the axial and
complementary cancellation weights,
respectively. The quantity $W$ measures the
total absolute weight amplification.

A single affine vector of scalar weights cannot separately cancel the two
curvature coefficients: their difference is \(\sigma^2\), independent of
\(r_1\).  The two directional weight vectors in
\eqref{eq:ysl_common_weights} remove this obstruction and additionally cancel
the exact field \(b_\perp\).  This is the componentwise version of the geometric bias cancellation principle, with the split supplied by the
YY estimator.

Define the implementable population map
\begin{equation}
\label{eq:ysl_common_algorithm}
 \begin{split}
 A_k(z)&=a_kU_0(z)+b_k(I_D-U_0(z)),\\
 \widetilde G(z)
 &=z+\sum_{k=1}^K A_k(z)\{G_{r_{1,k}}(z)-z\}
   =\sum_{k=1}^K A_k(z)G_{r_{1,k}}(z),\\
 \widetilde{\mathcal M}^{\mathrm{YY}}
 &=\widetilde G(\mathcal M^{\mathrm{pre}}).
 \end{split}
\end{equation}
The equality follows from \(\sum_k A_k(z)=I_D\).  At a given preliminary
point, all cylindrical means are computed from that same point with the same
axis.  The weights depend on \(d,q,p,\sigma\) and the selected radii, but not
on \(H\), \(P_0\), or \(b_\perp\).

\begin{proposition}[Fourth-order population refinement from an admissible
preliminary manifold]
\label{thm:ysl_common_population}
Suppose Assumption~\ref{ass:ysl_full_population} holds, \(q>1\), and the
weights satisfy \eqref{eq:ysl_common_weights}.
Put \(\widetilde{\Phi}=\widetilde G\circ X_0\).
Then for sufficiently small \(\sigma\), \(\widetilde{\Phi}\) is a
\(C^2\) embedding and
\(\widetilde\varphi=\pi_{\mathcal M}\circ\widetilde{\Phi}\)
is a \(C^2\) diffeomorphism of \(\mathcal M\). The normal
reparameterization
\(
 \widetilde X
 =\widetilde{\Phi}\circ\widetilde\varphi^{-1}
 =\operatorname{Id}_{\mathcal M}+\widetilde u
\)
satisfies
\begin{equation}
\label{eq:ysl_common_ms_image}
 \begin{split}
 \|\widetilde u\|_{C^0}&\leq C\sigma^4,
 \qquad \|\widetilde u\|_{C^2}\leq C\sigma,\\
 d_H(\widetilde{\mathcal M}^{\mathrm{YY}},\mathcal M)&\leq C\sigma^4,
 \qquad
 \operatorname{reach}(\widetilde{\mathcal M}^{\mathrm{YY}})\geq c_*>0.
 \end{split}
\end{equation}
The constants are uniform over the admissible scales and bounded weights.
\end{proposition}

The proposition asserts fourth-order accuracy of the output manifold, not
fourth-order pointwise accuracy at the original preliminary parameter.
The remaining \(O(\sigma^3)\) tangential displacement is absorbed by the
nearest-point reparameterization and contributes only \(O(\sigma^6)\) to
normal distance.  This reparameterization is an analytical device: the
algorithm is exactly \eqref{eq:ysl_common_algorithm}. The uniform reach bound follows from the second-order control
of the population map and the quantitative embedding estimates.

\begin{remark}[Hypersurface specialization]
\label{rem:ysl_common_hypersurface}
When \(q=1\), \(S_0=b_\perp=0\) and
\(\sigma^2\{1+\alpha_0(r_1/\sigma)\} =\sigma^2+r_1^2/(d+2p+2)\).
Two distinct transverse scales suffice.  Choose scalar weights satisfying
\(\sum_k a_k=1\) and \(\sum_k a_k(1+\alpha_0(\kappa_k))=0\), and set
\(b_k=a_k\).  Then \(A_k=a_kI_D\), and the same common-start argument gives
fourth-order Hausdorff accuracy.  No condition involving \(\beta_m\) is
imposed when the complementary normal block is absent.
\end{remark}

\subsection{Finite-sample approximation}
\label{sec:ysl_full_finite_sample}
Following \citet[Theorem~5.4]{yao2023manifold},
we analyze the finite-sample refinement conditional on
a fixed preliminary manifold $\mathcal M^{\mathrm{pre}}$
satisfying Assumption~\ref{ass:ysl_full_population}.
The results also apply to an independently constructed pilot
on the event that these conditions hold.
Define
\[
\label{eq:ysl_common_empirical_first_stage}
 \widehat F_{r_0}(z)
 =
 \frac{\sum_{i=1}^{n}Y_i
          K_{0,r_0}(\|Y_i-z\|)}
      {\sum_{i=1}^{n}K_{0,r_0}(\|Y_i-z\|)},
 \qquad
 \widehat U_0(z)
 =
 \frac{\widehat V_0(z)\widehat V_0(z)^\top}
      {\|\widehat V_0(z)\|^2},
 \quad
 \widehat V_0(z)=\widehat F_{r_0}(z)-z.
\]
For each selected transverse radius, let
\[
\label{eq:ysl_common_empirical_second_stage}
 \begin{split}
 \widehat\Psi_{k,z}(y)
 &=
 K_{\parallel,r_2}
 \bigl(\|\widehat U_0(z)(y-z)\|\bigr)
 K_{\perp,r_{1,k}}
 \bigl(\|(I_D-\widehat U_0(z))(y-z)\|\bigr),\\
 \widehat G_k(z)
 &=
 \frac{\sum_{i=1}^{n}Y_i
                  \widehat\Psi_{k,z}(Y_i)}
      {\sum_{i=1}^{n}\widehat\Psi_{k,z}(Y_i)}.
 \end{split}
\]
The notation $\widehat G_k$ denotes the
empirical counterpart of the population
second-stage map $G_{r_{1,k}}$, using the
same estimated first-stage axis across scales.
With
\(\widehat A_k(z) = a_k\widehat U_0(z)+b_k(I_D-\widehat U_0(z))\),
define the empirical common-start estimator
\[
\label{eq:ysl_common_empirical_map}
 \begin{split}
 \widehat G(z)
 &=
 z+\sum_{k=1}^K\widehat A_k(z)
       \{\widehat G_k(z)-z\}
 =
 \sum_{k=1}^K\widehat A_k(z)\widehat G_k(z),\\
  \widehat{\mathcal M}^{\mathrm{YY}}
  &=
  \widehat G(\mathcal M^{\mathrm{pre}}),
  \qquad
  \widehat{\Phi}
  =
  \widehat G\circ X_0.
 \end{split}
\]

For \(0<\delta<1/2\), let \(C_{\mathrm{ent}}\geq1\) be a fixed
constant depending only on \(d,D\) and the chosen kernel and projector-jet
families, and set
\begin{equation}
\label{eq:ysl_common_finite_rate}
\ell_{\delta}(\sigma)
=\log\!\left(\frac e\delta\right)
+C_{\mathrm{ent}}\log\!\left(\frac e\sigma\right),
\qquad
\zeta_{n,\delta}(\sigma)
=\sqrt{\frac{\ell_{\delta}(\sigma)}{n\sigma^d}}
+\sqrt{L_\sigma}\,
\frac{\ell_{\delta}(\sigma)}{n\sigma^d}.
\end{equation}
The first term is the variance contribution from the effective local sample
size \(n\sigma^d\).  The second is the Bernstein envelope contribution; the
factor \(\sqrt{L_\sigma}\) comes from the longitudinal support
\(r_2=C_2\sigma\sqrt{L_\sigma}\).  The noise level \(\sigma\) is a fixed
sufficiently small quantity throughout this finite-sample statement.

\begin{theorem}[Finite-sample refinement rate for GBC-YY]
\label{thm:ysl_common_finite_sample}
Suppose the population conditions of
Proposition~\ref{thm:ysl_common_population} hold and the deterministic weights
satisfy \eqref{eq:ysl_common_weights}.  Fix a sufficiently small
\(\sigma\) and \(0<\delta<1/2\).  There are constants \(C,c>0\),
independent of \(n,\delta\), and \(\sigma\) in the stated small-noise range,
such that, whenever
\( n\sigma^d\geq C\ell_{\delta}(\sigma),
 \zeta_{n,\delta}(\sigma)\leq c\sigma^{1/2},
\)
the following claims hold simultaneously with probability at least
\(1-\delta\):
\[
\label{eq:ysl_common_empirical_signal_mass}
 \inf_{z\in\mathcal M^{\mathrm{pre}}}
 \|\widehat F_{r_0}(z)-z\|
 \geq
 \frac{c_F}{2}\sigma.
\]
Moreover, all corresponding normalized empirical first- and second-stage
denominators are bounded below by \(c\sigma^d\), and, for \(j=0,1\),
\begin{equation}
\label{eq:ysl_common_explicit_single_scale_rates}
 \begin{split}
 \|\widehat U_0-U_0\|_{C^j(\mathcal M^{\mathrm{pre}})}
 &\leq
 C\sigma^{-j}\zeta_{n,\delta}(\sigma),\\
 \max_{1\leq k\leq K}
 \|\widehat G_k-G_{r_{1,k}}\|_{C^j(\mathcal M^{\mathrm{pre}})}
 &\leq
 C\sigma^{1-j}\zeta_{n,\delta}(\sigma).
 \end{split}
\end{equation}
Consequently,
\[
\label{eq:ysl_common_explicit_multiscale_rates}
 \|\widehat{\Phi}
       -\widetilde{\Phi}\|_{C^j(\mathcal M)}
 \leq
 CW\sigma^{1-j}\zeta_{n,\delta}(\sigma),
 \qquad
 j=0,1,
\]
and
\[
\label{eq:ysl_common_explicit_hausdorff_rate}
 d_H\!\left(
 \widehat{\mathcal M}^{\mathrm{YY}},\mathcal M
 \right)
 \leq
 C\sigma^4+CW\sigma\zeta_{n,\delta}(\sigma).
\]
For sufficiently small fixed \(\sigma\), the same event also makes
\(\widehat G|_{\mathcal M^{\mathrm{pre}}}\) a \(C^2\) embedding, because
its \(C^1\) distance from the population embedding is at most
\(CW\zeta_{n,\delta}(\sigma)\).

For the numerical implementation, let
\(\varnothing\neq\mathcal S\subseteq\mathcal Z_N
\subset\mathcal M^{\mathrm{pre}}\) be the retained query set, and define
\(\delta_{\mathcal S} =\sup_{z\in\mathcal M^{\mathrm{pre}}} \min_{z_j\in\mathcal S}\|z-z_j\|\).
If numerical evaluations satisfy
\(
 \max_{z\in\mathcal S}
 \|\widehat G^{\mathrm{num}}(z)-\widehat G(z)\|
 \leq\varepsilon_{\mathrm{num}},
 \varepsilon_{\mathrm{num}}\geq0,
\)
put
\(\widehat{\mathcal P}^{\mathrm{YY,num}}_N =\{\widehat G^{\mathrm{num}}(z):z\in\mathcal S\}\).
Then, on the same high-probability event,
\[
 \begin{aligned}
 d_H\!\left(\widehat{\mathcal P}^{\mathrm{YY,num}}_N,\mathcal M\right)
 &\leq d_H\!\left(\widehat{\mathcal M}^{\mathrm{YY}},\mathcal M\right)
       +\varepsilon_{\mathrm{num}}+C\delta_{\mathcal S}\\
 &\leq C\sigma^4+CW\sigma\zeta_{n,\delta}(\sigma)
       +\varepsilon_{\mathrm{num}}+C\delta_{\mathcal S}.
 \end{aligned}
\]
Here \(C\) is uniform under the stated bound on \(W\).
With exact evaluation, \(\varepsilon_{\mathrm{num}}=0\) and this is
\(\widehat{\mathcal P}^{\mathrm{YY}}_N\) returned by
Algorithm~\ref{alg:ysl23_ms}.  The fill distance uses the retained queries,
not the original set before failures; an empty retained set is reported
as a failed run.  The embedding conclusion concerns the continuous
estimator, not the finite point cloud.
\end{theorem}

\begin{algorithm}[!htbp]
\caption{Finite-Sample Directional GBC for the YY Estimator}
\label{alg:ysl23_ms}
\begin{algorithmic}[1]

\Require Observations $\mathcal Y=\{y_1,\ldots,y_n\}$,
query points
\(\mathcal Z_N=\{z_1,\ldots,z_N\}\subset\mathcal M^{\mathrm{pre}}\)
on a supplied preliminary manifold, radii
\(r_0,r_2,r_{1,1},\ldots,r_{1,K}\),
directional weights \(\{a_k,b_k\}_{k=1}^K\),
and positive mass and displacement thresholds.

\Ensure Finite output point cloud
\(\widehat{\mathcal P}^{\mathrm{YY}}_N\)
and successful query set \(\mathcal S\).

\State Initialize \(\mathcal S=\varnothing\).

\For{each \(z\in\mathcal Z_N\)}
    \State Check the first-stage kernel mass; if valid, compute
    \(\widehat V_0(z)=\widehat F_{r_0}(z)-z\).

    \If{the mass check fails or the displacement is too small}
        \State Record failure and \textbf{continue}.
    \EndIf

    \State Define the common contraction axis
    \(\widehat U_0(z) = \frac{\widehat V_0(z)\widehat V_0(z)^\top} {\|\widehat V_0(z)\|^2}\).

    \For{\(k=1,\ldots,K\)}
        \State Using the same center \(z\), axis
        \(\widehat U_0(z)\), and observations \(\mathcal Y\),
        check the second-stage kernel mass.

        \If{the mass check fails}
            \State Record failure and \textbf{break}.
        \EndIf

        \State Compute
        \(\widehat G_k(z)
        =\widehat G_{r_{1,k}}(z;\widehat U_0(z))\)
        and set
        \(\widehat A_k(z) = a_k\widehat U_0(z)+b_k(I_D-\widehat U_0(z))\).
    \EndFor

    \If{all scales succeed}
        \State Compute
        \(\widehat G(z) = z+\sum_{k=1}^K \widehat A_k(z)\{\widehat G_k(z)-z\}\),
        and add \(z\) to \(\mathcal S\).
    \EndIf
\EndFor

\State \Return
\(\widehat{\mathcal P}^{\mathrm{YY}}_N
=\{\widehat G(z):z\in\mathcal S\}\)
and \(\mathcal S\).

\end{algorithmic}
\end{algorithm}
\begin{corollary}[A sufficient sample size for fourth-order accuracy]
\label{cor:ysl_common_fourth_order_sample_size}
Under the conditions of
Theorem~\ref{thm:ysl_common_finite_sample}, there is a constant \(C\)
such that
\begin{equation}
\label{eq:ysl_common_fourth_order_sample_size}
 n\geq
 CW^2\sigma^{-(d+6)}\ell_{\delta}(\sigma)
\end{equation}
implies, for sufficiently small fixed \(\sigma\),
\(W\sigma\zeta_{n,\delta}(\sigma) \leq C\sigma^4\).
Hence, with probability at least \(1-\delta\),
\begin{equation}
\label{eq:ysl_common_fourth_order_finite_conclusion}
 d_H\!\left(
 \widehat{\mathcal M}^{\mathrm{YY}},\mathcal M
 \right)
 \leq
 C\sigma^4,
 \qquad
 \widehat{\mathcal M}^{\mathrm{YY}}
 \text{ is a }C^2\text{ embedded manifold}.
\end{equation}

On the same event, if the retained queries and numerical evaluations also
satisfy \(\delta_{\mathcal S}=O(\sigma^4)\) and
\(\varepsilon_{\mathrm{num}}=O(\sigma^4)\), then
\(d_H\!\left(\widehat{\mathcal P}^{\mathrm{YY,num}}_N,\mathcal M\right) \leq C\sigma^4\).
\end{corollary}
\begin{remark}
The \(C^2\) embedding conclusion describes the smoothness of the
empirical image, whereas the stochastic approximation bounds
control only its position and first derivatives.
These results therefore do not establish a deterministic positive
reach lower bound uniform in the sample size and the admissible
noise level.
In particular, the population reach bound is not transferred to
the empirical image by the \(C^1\) perturbation argument alone.
Such a uniform guarantee would require additional quantitative
second-order control.
\end{remark}
\begin{remark}[Relation to the YY estimator sample size]
\label{rem:ysl_common_ysl23_sample_size}
The proofs of
\citet[Theorems~3.5 and~4.5]{yao2023manifold} use a total sample size of
order \(\sigma^{-(d+3)}\), corresponding to an effective local sample size
of order \(\sigma^{-3}\).  Retaining the sample-size term in the same
weighted-average calculation, if
\(
 n\geq
 C\sigma^{-(d+3)}\ell_{\delta}(\sigma),
\)
then
\(
 \sigma\zeta_{n,\delta}(\sigma)
 \leq
 C\{\sigma^{5/2}+\sigma^4\sqrt{L_\sigma}\}.
\)
This fluctuation is below the original YY estimator approximation scale
\(\sigma^2L_\sigma\) for sufficiently small \(\sigma\).  Exposing the
fourth-order population improvement requires the stronger sample size
\eqref{eq:ysl_common_fourth_order_sample_size}; this is the statistical
price of resolving an \(O(\sigma^4)\), rather than an
\(O(\sigma^2L_\sigma)\), geometric bias.
\end{remark}

\begin{proposition}[Sample-split normal-offset initialization]
\label{prop:yy-sample-split-pilot}
Suppose the model and geometric conditions of Assumption~\ref{ass:geometry-observation}
and Assumption~\ref{ass:ysl_full_population}(\ref{item:ysl_population_geometry}) hold, with the radii specified in
Assumption~\ref{ass:ysl_full_population}(\ref{item:ysl_population_radii}). Let $\mathcal P_\sigma$ be a
pilot manifold constructed from an independent sample. Assume that,
with probability at least $1-\delta_0$, the following event
$\mathcal E_{\rm pil}$ occurs: $\mathcal P_\sigma$ is a compact,
connected, boundaryless $C^6$ embedded $d$-manifold, and
\(
 d_H(\mathcal P_\sigma,\mathcal M)\le B\sigma,
 \operatorname{reach}(\mathcal P_\sigma)\ge\tau_p>0.
\)
On this event, suppose that $\mathcal P_\sigma$ admits a $C^6$ unit
normal field $\nu_\sigma$ satisfying
$\|\nu_\sigma\|_{C^2(\mathcal P_\sigma)}\le C_\nu$.
The constants $B,\tau_p,C_\nu$ are independent of $\sigma$.
Fix $A>B$ and set
\[
 \mathcal M_{\rm pre}
 =\{p+A\sigma\nu_\sigma(p):p\in\mathcal P_\sigma\}.
\]
For sufficiently small $\sigma$, on $\mathcal E_{\rm pil}$,
$\mathcal M_{\rm pre}$ satisfies Assumption~\ref{ass:ysl_full_population}, including the
first-stage signal condition with a constant $c_F>0$ independent
of $\sigma$.

Use a second, independent sample of size $n_1$ to construct the
common-start GBC-YY refinement, with the radii and
uniformly bounded directional weights of Theorem~\ref{thm:ysl_common_finite_sample}. If
\[
 n_1\ge C W^2\sigma^{-(d+6)}\ell_{\delta_1}(\sigma),
 \qquad 0<\delta_1<1/2,
\]
then, for sufficiently small $\sigma$,
\[
 \mathbb P\!\left\{
 d_H(\widehat{\mathcal M}^{\rm YY},\mathcal M)
 \le C'\sigma^4,
 \ \widehat{\mathcal M}^{\rm YY}\text{ is a $C^2$ embedded manifold}
 \right\}\ge 1-\delta_0-\delta_1.
\]
The constants are uniform on $\mathcal E_{\rm pil}$.
\end{proposition}

The guarantee is conditional on a pilot satisfying the stated
geometric and normal-field conditions, with total sample size
$n_0+n_1$. The sample size $n_0$ required to construct such
a pilot is not quantified here.

\begingroup
\setlength{\floatsep}{5pt plus 2pt minus 1pt}
\setlength{\textfloatsep}{7pt plus 2pt minus 1pt}
\setlength{\intextsep}{6pt plus 2pt minus 1pt}
\captionsetup{font=small,skip=3pt}
\flushbottom
\renewcommand{\topfraction}{0.97}
\renewcommand{\bottomfraction}{0.92}
\renewcommand{\textfraction}{0.03}
\renewcommand{\floatpagefraction}{0.78}
\makeatletter
\setlength{\@fptop}{0pt}
\setlength{\@fpsep}{8pt plus 2pt minus 1pt}
\makeatother

\section{Experiments}
\label{sec:simulation}
This section presents numerical experiments for three
manifold-fitting procedures: successive tangent-space
projection (STSP), Gaussian density-ridge fitting (DR),
and the YY estimator implemented as in Section~\ref{sec:yao_contraction}.
We examine population bias reduction from second to
fourth order, finite-sample accuracy, estimated-normal
correspondence, and directional extrapolation in higher
codimension. Additional robustness studies are provided
in Appendix~\ref{app:simulation-experiments}.

The three methods use procedure-specific fitting scales
and stability conditions. Within each method, Single-scale and
GBC share the same kernel, query points, fitting procedure,
and observations in each finite-sample repetition.
Single-scale uses the specified reference scale, whereas GBC
combines scale-specific fits using bias-cancellation weights.

The finite-sample experiments follow the protocol of
\citet{yao2023manifold}, while population experiments use
expectation-level operators. Throughout, $n$ denotes the
number of observations and $N$ the number of query points
(written $N_q$ in the correspondence experiment).
Each successful query produces one fitted point.

\begin{enumerate}
    \item For each repetition, independently generate
    $\mathcal Y=\{Y_i=X_i+\xi_i\}_{i=1}^n$,
    where $X_i$ follows the specified distribution on
    $\mathcal M$ and
    $\xi_i\sim\mathcal N(0,\sigma^2I_D)$ independently.

    \item Apply the fitting procedure at each required scale.
    STSP uses smooth-kernel local means and PCA-based
    tangent projection; DR uses Gaussian density-ridge
    equations; and YY uses a common first-stage contraction
    axis and scale-dependent cylindrical means.

    \item Construct the GBC estimator using
    procedure-specific cancellation weights.
    STSP and DR use estimated-normal correspondence
    (Section~\ref{sec:normal_guided_extrapolation}),
    whereas YY uses common-start directional extrapolation
    (Section~\ref{sec:ysl_apply_general_theory}).
    All weights correspond to the implemented kernels.

    \item Let $\mathcal S\subseteq\{1,\ldots,N\}$ be the
    nonempty set of successful query indices, with outputs
    $\{\widehat z_j:j\in\mathcal S\}$.
    We report
    \[
    E_\infty=\max_{j\in\mathcal S}d(\widehat z_j,\mathcal M),
    \qquad
    E_{\mathrm{avg}}=
    \frac{1}{|\mathcal S|}
    \sum_{j\in\mathcal S}d(\widehat z_j,\mathcal M).
    \]
    Failures are recorded separately, with reliability
    criteria specified in the relevant experiments.
    These one-sided errors do not measure full Hausdorff
    distance; reverse-coverage results are provided in
    Appendix~\ref{app:bandwidth-robustness}.
\end{enumerate}

The main scale and sample-size experiments use six representative
manifolds:
$\mathbb S^1\subset\mathbb R^2$,
$\mathbb S^2\subset\mathbb R^3$,
the Clifford torus $\mathbb T_C^2\subset\mathbb R^4$,
$\mathbb S^3\subset\mathbb R^4$,
$\mathbb S^2\subset\mathbb R^{10}$,
and the corrugated torus
$\mathbb T_{\mathrm{corr}}^2\subset\mathbb R^3$.
The correspondence experiments additionally use the
nonconstant-curvature closed curve
\(
\gamma(t)=(1.5\cos t,\sin t,0.3\sin 2t)
\in\mathbb R^3.
\)

The directional ablation uses its embedding in \(\mathbb R^4\), obtained
by appending a zero coordinate. We also consider a corrugated torus in \(\mathbb R^3\) to assess performance
on a geometrically complex surface with spatially varying curvature. Together, these examples examine changes
in intrinsic dimension, ambient dimension, codimension, and local geometry.

\paragraph*{Initialization and scope of the YY experiments}
The YY experiments in Sections~\ref{sec:population-bias}--\ref{sec:exp-directional} use a controlled preliminary
manifold constructed from the known target geometry, in order to
isolate the refinement mechanism. Within each experiment, all
compared refinements share the same preliminary manifold, query
points, and observations.

\paragraph*{Implementation and code}
All experiments are implemented in Python using NumPy, SciPy,
and CuPy, with GPU acceleration, parallel execution, and
double-precision arithmetic.
Detailed numerical procedures and experiment configurations
are provided in Appendix~\ref{app:simulation-experiments}.
The source code, dependencies, and scripts for reproducing
all figures are available at
\url{https://github.com/zzjpp/Multiscale-Bias-Cancellation-for-Accurate-Manifold-Estimation}.

\subsection{Population Verification}
\label{sec:population-bias}

The population theory predicts that the single-scale error is second order,
whereas GBC removes the leading bias and leaves a
fourth-order remainder. We verify this prediction in two steps. We first
check the overall error scaling along the joint noise--bandwidth path used
in the theory. We then focus on the smallest noise levels and examine
whether the estimated orders remain stable as the experiment moves further
into the asymptotic regime.

\subsubsection{Scaling along the theoretical path}
\label{sec:population-joint-scaling}

We first examine the overall scaling behavior along the theoretical
noise--bandwidth path. Here, we illustrate this using STSP;
results for DR and YY, together with additional STSP results, are reported
in Appendix~\ref{app:simulation-experiments}.

For STSP, we vary $\sigma$ over seven values between
$0.008$ and $0.03$, with $h_1=3\sigma$ and
$h_2=6\sigma$. Thus, $\sigma/h_1=1/3$ and
$\sigma/h_2=1/6$ remain fixed, matching the
theoretical scaling regime.

To isolate the effect of the bias-cancellation
coefficient, we compare three estimators constructed
from the same population smooth-kernel STSP fits:
single-scale STSP at $h_2$, naive squared-bandwidth
extrapolation, and GBC-STSP.
Both extrapolation methods use identical
estimated-normal correspondences.
The naive method uses weights $(4/3,-1/3)$, satisfying
$\sum_k w_k h_k^2=0$, whereas our method uses
weights satisfying
$\sum_k w_k c_{\mathrm{STSP}}(h_k,\sigma)=0$.
We additionally include the original hard-ball
MS-STSP method of Xia, Li and Yao (2026) as a
separate benchmark.
For each manifold, we evaluate the population
maximum error $E_\infty$ and estimate its convergence
order by regressing $\log E_\infty$ on $\log\sigma$.

Figure~\ref{fig:intro_stsp_population} shows that naive squared-bandwidth
extrapolation retains approximately second-order
error, despite using the same fitted manifolds and
correspondences as our method.
For example, on the circle, the fitted slopes are
$2.00$ for single-scale STSP, $2.00$ for naive
extrapolation, and $4.00$ for our method.
The corresponding slopes on the Clifford torus
are $2.00$, $2.00$, and $4.01$.
The proposed method exhibits near-fourth-order
behavior on five of the six manifolds, with a slope
of $3.82$ on the corrugated torus over the tested
noise range.

These findings are consistent with
Proposition~\ref{prop:stsp-naive}: canceling the
nominal $h^2$ term leaves a nonzero second-order
bias under Gaussian noise.
The improvement therefore does not arise merely
from combining multiple scales or from the
correspondence construction, but from canceling
the correct noise-dependent geometric bias.

\subsubsection{Stability in the small-noise regime}
\label{sec:population-asymptotic-stability}

The preceding log--log regressions show the expected second- and
fourth-order scaling over the full noise grid. We now examine whether
these orders remain stable as the noise level decreases, rather than
relying solely on the fitted slopes.
To this end, we normalize the population maximum errors by their
predicted orders:
\(\frac{E_{\infty,\mathrm{Single}}}{\sigma^2}, \frac{E_{\infty,\mathrm{GBC}}}{\sigma^4}\).
We examine whether the normalized errors remain bounded
and approximately stable as $\sigma$ decreases. Such
stability is consistent with the predicted convergence
orders, while a nonzero limiting constant is not required
by the established upper bounds.

Figure~\ref{fig:e1-small-noise-stability} shows that, for all three
methods, the normalized errors remain nearly constant in the small-noise
regime, with only mild variation over the tested range. In particular,
the stability of \(E_{\infty,\mathrm{GBC}}/\sigma^4\) provides direct
numerical evidence that the fourth-order behavior persists at smaller
noise levels. Together with the full-grid slopes, these results support
the predicted reduction from second- to fourth-order population error
across three different manifold-fitting procedures.

\begin{figure}[!htbp]
    \centering
    \includegraphics[width=0.94\linewidth]
    {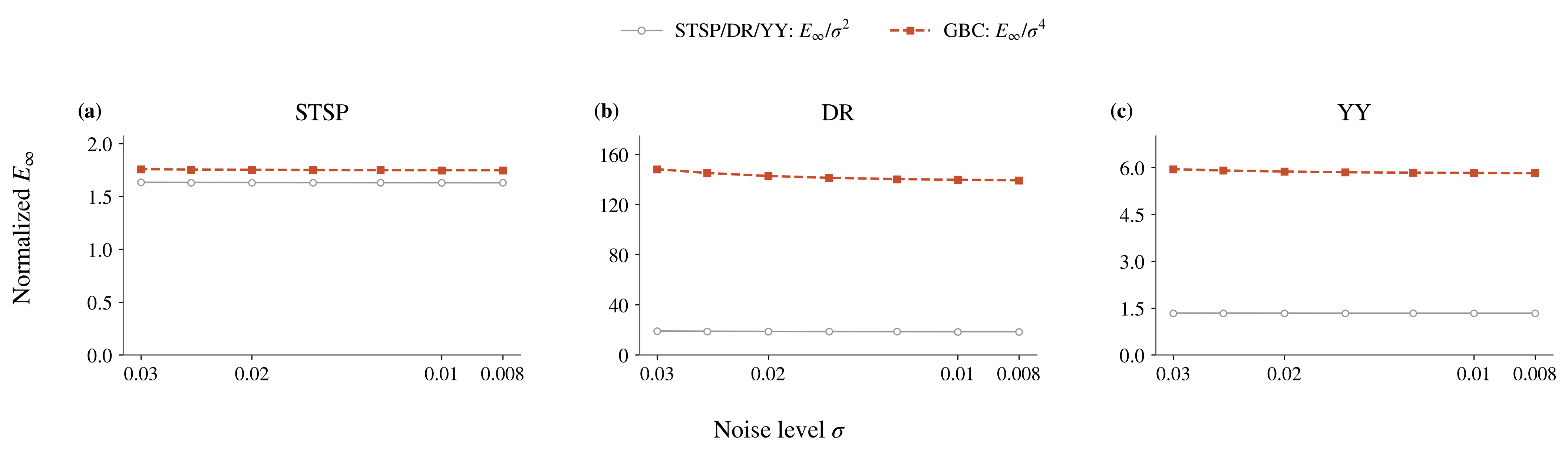}
    \caption{
    Normalized population maximum errors on the circle for STSP, DR,
    and YY (left to right). The curves show
    \(E_{\infty,\mathrm{Single}}/\sigma^2\) and
    \(E_{\infty,\mathrm{GBC}}/\sigma^4\) under the respective
    noise-dependent scale settings.
    Their nearly constant levels at small \(\sigma\) support the
    predicted second- and fourth-order error scaling.
    }
    \label{fig:e1-small-noise-stability}
\end{figure}
\subsection{Finite-Sample Behavior at Fixed Scales}
\label{sec:asymptotic}
We next examine whether increasing the sample size drives the empirical
multiscale error toward its $O(h^4)$ population approximation error, while the
single-scale estimator approaches its corresponding $O(h^2)$ limit. We fix the noise level at $\sigma=0.005$ and use the bandwidths
$h_1=0.1$ and $h_2=0.2$, while progressively increasing the sample size
$n$ to examine the finite-sample asymptotic behavior of the estimators.
For each manifold, we perform 20 paired repetitions. The single-scale
estimator uses $h_2$, while the GBC estimator combines the fits at
$h_1$ and $h_2$ using the corresponding bias-cancellation weights. Figure~\ref{fig:e2-finite-sample-main} reports maximum errors on all six
manifolds; average errors and the DR and YY results are given in Appendix~\ref{app:simulation-experiments}.

\begin{figure}[!htbp]
    \centering
    \includegraphics[width=0.98\linewidth]{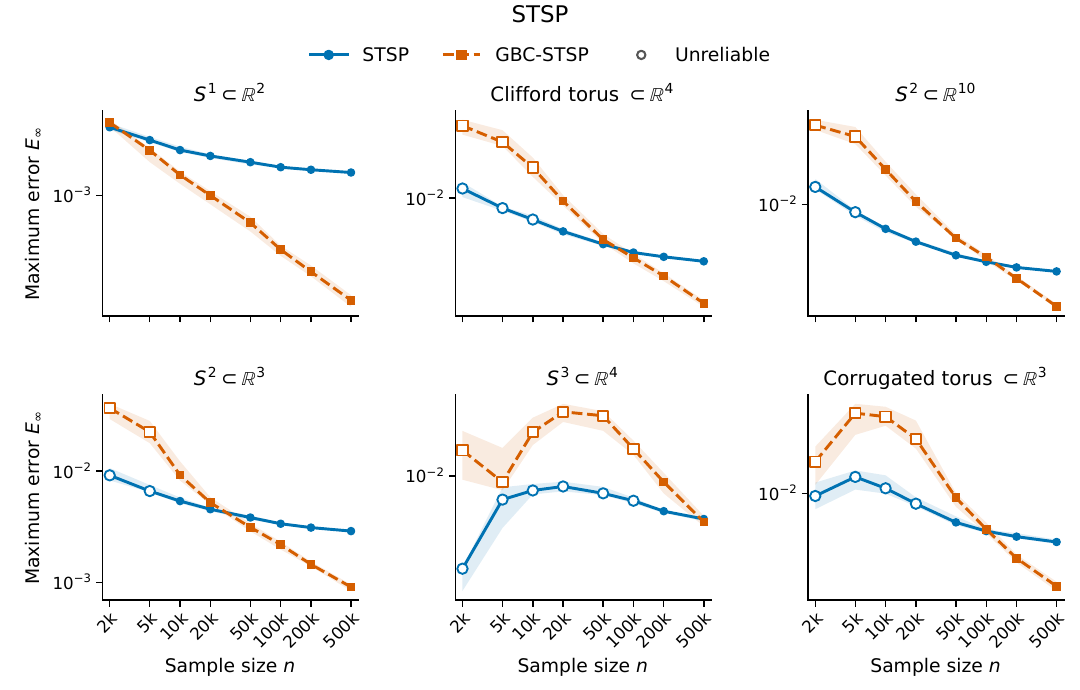}
    \caption{
    Finite-sample maximum error of single-scale STSP and GBC-STSP on the six manifolds shown,
    with \(\sigma=0.005\), \(h_1=0.1\), and \(h_2=0.2\).
    Curves show medians and shaded bands show interquartile ranges across
    20 repetitions as \(n\) increases from \(2{,}000\) to \(500{,}000\).
    Hollow markers indicate settings that fail the prescribed reliability
    criterion.
}
    \label{fig:e2-finite-sample-main}
\end{figure}

Figure~\ref{fig:e2-finite-sample-main} shows that
single-scale errors tend to plateau as the sample size
increases, whereas GBC cancellation achieves lower
errors over the reliable large-sample range. This behavior
is consistent with the removal of the leading geometric
bias, although sampling fluctuations remain visible.
For example, on the circle at $h_2=0.2$ and
$n=500{,}000$, the median average error decreases from
$1.343\times10^{-3}$ to $5.243\times10^{-5}$, a factor
of $25.6$ (Figure~\ref{fig:e2-stsp-finite-sample-appendix}).

\subsection{Bandwidth Robustness}
\label{sec:e3-bandwidth-robustness}

We next examine the sensitivity of STSP, DR, and YY to the choice of
fitting scale on six manifolds. We fix \(\sigma=0.015\) and vary the
scale parameter over the common grid in \([0.06,5]\).
For STSP and DR, we set \(h_1=h_2/2\) and vary \(h_2\).
For YY, we vary the central transverse radius \(r_1\), with
three scales \((0.4,1,1.6)r_1\).
The sample sizes for the six manifolds are
\(
n=(12000,24000,24000,60000,24000,36000),
\)
respectively. Each setting uses 20 repetitions, with 256 query
points for the first five manifolds and 512 for the corrugated
torus. Within each method, Single-scale and GBC use the same
observations and query points.

A configuration is considered reliable if the overall query-failure
rate is at most \(1\%\) and at least 18 of the 20 repetitions are
complete. Filled markers indicate reliable configurations, while
hollow markers indicate unreliable ones. Failed queries are not
imputed, and unreliable configurations are excluded from accuracy
comparisons.
We present the STSP results in
Figure~\ref{fig:e3-bandwidth-robustness-main} as a representative
example. The corresponding average errors, together with the
DR and YY results, are reported in Appendix~\ref{app:simulation-experiments}.

\begin{figure}[!htbp]
    \centering
    \includegraphics[width=0.98\linewidth]
    {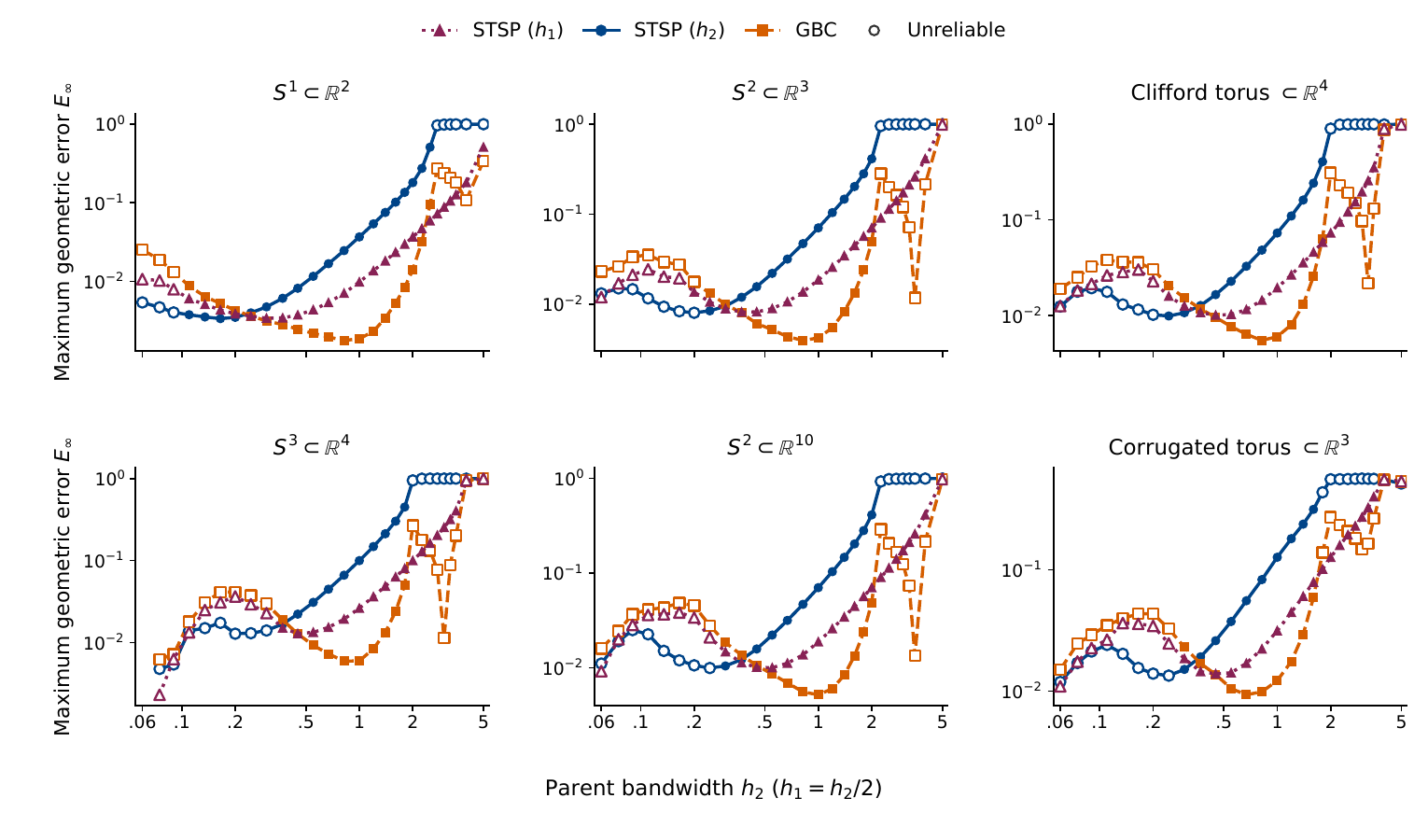}
    \caption{
    Bandwidth robustness of STSP on six manifolds, measured by
    maximum geometric error \(E_\infty\).
    Curves compare Single at \(h_1\), Single at \(h_2\), and GBC,
    with \(\sigma=0.015\), \(h_1=h_2/2\), and sample sizes
    \(n=(12000,24000,24000,60000,24000,36000)\)
    in panel order.
    Filled markers indicate reliable configurations, while
    hollow markers indicate those failing the reliability criterion.
    }
    \label{fig:e3-bandwidth-robustness-main}
\end{figure}

Figure~\ref{fig:e3-bandwidth-robustness-main} shows a broadly
U-shaped error pattern, reflecting the bias--variance tradeoff
in bandwidth selection. Multiscale cancellation reduces errors
at intermediate bandwidths, and on several manifolds, GBC
outperforms the best reliable single-scale baseline across
multiple tested bandwidths. This suggests improved robustness
to bandwidth selection, although the extent of the improvement
varies across fitting procedures and manifold geometries, as
shown by the additional results in Appendix~\ref{app:simulation-experiments}.

\subsection{Correspondence Accuracy and Computational Stability}
\label{sec:enf_cp_comparison}

We compare estimated-normal-fiber (ENF) correspondence with local
continuous closest-point (CP) matching, discrete nearest-neighbour (NN)
matching, and oracle-normal matching on the same empirical STSP fits.
The comparison separates the accuracy of a cross-scale match from its
consequences for geometric recovery and computation. Proposition~\ref{prop:enf_local_conditioning}
establishes uniform local conditioning for ENF. Proposition~\ref{prop:cp_supplement} shows that global continuous CP also attains
fourth-order localization, although its uniqueness and smoothness
may require additional uniform curvature control.

We use the nonconstant-curvature curve
\(
\gamma(\theta)
  =(1.5\cos\theta,\sin\theta,0.3\sin(2\theta))^\top
  \in\mathbb R^3,
\)
with $\sigma=0.06$, $h_1=0.12$, $h_2=0.24$,
$n=200{,}000$, and $N_q=32$ reference queries, over ten paired
repetitions. All methods share the observations, fitted reference
points, and extrapolation weights. Oracle matching uses true normal
fibers, while ENF uses normal spaces estimated by reference-scale
weighted PCA.

The implemented CP baseline constructs a continuous local chart of
the empirical target STSP branch using target-scale weighted PCA and
implicit residual solves. Given a common reference point
$z_{1,j}\in\widehat{\mathcal M}_{h_1}$, it numerically minimizes
$\|z-z_{1,j}\|^2$ over the selected local target sheet. Unlike
discrete NN, which selects among $M_{\mathrm{match}}$ separately
fitted target-scale candidates, local CP does not optimize over a
finite candidate library. Its search uses neither the true curve
parameterization nor oracle-normal information. The shared reference
queries nevertheless retain geometry-assisted initialization and
validation. Since the local search does not certify global optimality
over arbitrary additional components of the empirical STSP zero set,
we distinguish this \emph{local CP benchmark} from the global CP map
studied in Proposition~\ref{prop:cp_supplement}. Implementation and timing
details appear in Appendix~\ref{app:bandwidth-robustness}.

For method $A$, let $Z_{2,j}^{(A)}$ denote its target-scale match and
$Z_{2,j}^{\mathrm{Or}}$ the oracle-normal match on the same empirical
target branch. We report
\begin{equation}\label{eq:enf_cp_metrics}
E_{\infty}^{(A)}
=\max_{1\le j\le N_q}d\bigl(w_1z_{1,j}+w_2Z_{2,j}^{(A)},\mathcal M\bigr),
\qquad
E_{\mathrm{match}}^{(A)}
=\max_{1\le j\le N_q}\bigl\|Z_{2,j}^{(A)}-Z_{2,j}^{\mathrm{Or}}\bigr\|.
\end{equation}
The second quantity measures discrepancy from \emph{oracle-normal
matching on the empirical fit}, not error relative to an unobserved
population target point. In particular,
$E_{\mathrm{match}}^{\mathrm{Or}}=0$ by definition.

\begin{figure}[!htbp]
  \centering
  \includegraphics[width=0.99\linewidth]
  {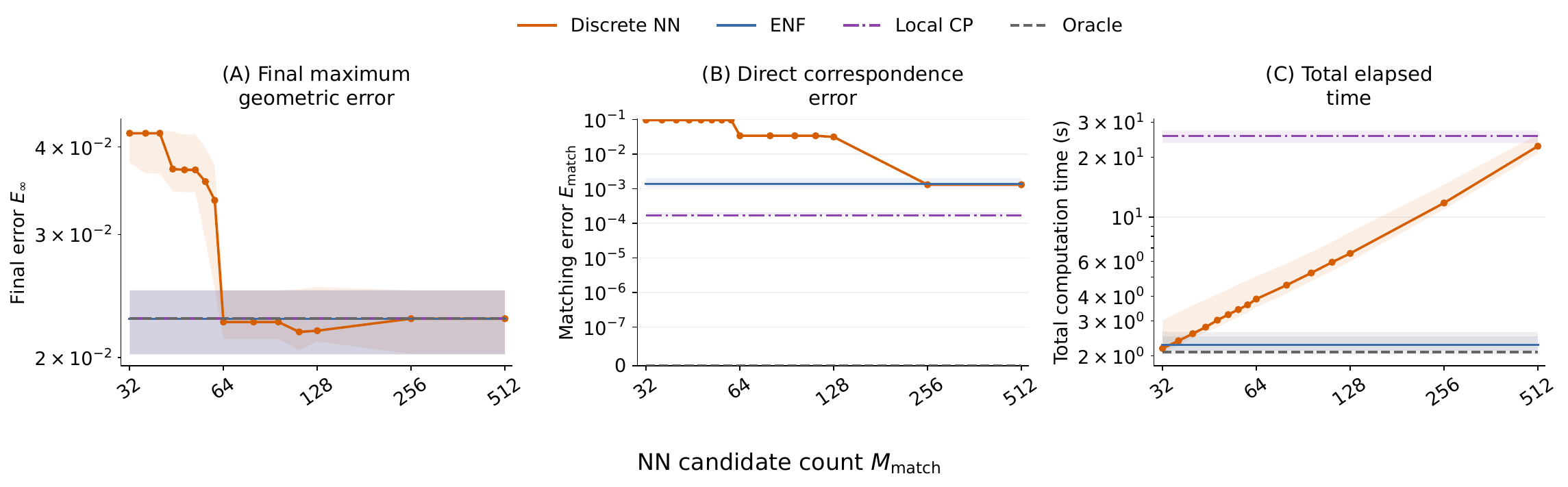}
  \caption{Correspondence accuracy and measured computational cost
  for empirical STSP on a nonconstant-curvature curve in
  $\mathbb R^3$ ($n=200{,}000$, $\sigma=0.06$, $h_1=0.12$,
  $h_2=0.24$, $N_q=32$; ten paired repetitions).
  (A) Final maximum point-to-manifold error after extrapolation.
  (B) Direct discrepancy relative to oracle-normal matching on the
  empirical target branch; the Oracle value is exactly zero.
  (C) Total measured elapsed time for the reported implementations.
  Bands represent interquartile
  ranges. ENF, local CP, and Oracle are independent of the NN
  candidate count. Local CP optimizes continuously on a constructed
  empirical target sheet, without a certificate of global
  closest-point optimality over the entire empirical zero set.}
  \label{fig:enf_cp_comparison}
\end{figure}

Figure~\ref{fig:enf_cp_comparison} distinguishes direct matching
precision from final geometric accuracy. Local CP has a smaller
discrepancy from oracle-normal matching than ENF, yet the two methods
have nearly identical final maximum geometric errors. Thus the
additional matching precision of local CP does not produce a
noticeable improvement in the final extrapolated reconstruction in
this experiment.

The reported total elapsed time for local CP is approximately
$25.6$ seconds per batch of $32$ queries (median over ten
repetitions), and the corresponding median CP correspondence-stage time is
$24.615$ seconds. The plotted ENF time is substantially lower,
whereas discrete NN becomes more expensive as its candidate set is
refined. These are implementation- and timing-protocol-specific
observations, not general complexity bounds. All $320$ queries are
accepted for each method and candidate-count setting; the experiment
therefore does not demonstrate a difference in observed failure rates.

The key geometric distinction is that ENF provides a unique,
locally well-conditioned transverse intersection, whereas
global CP correspondence can be multivalued without uniform
curvature control. Appendix~\ref{sec:supp_oscillatory_cp} illustrates this
distinction through an oscillatory normal-graph counterexample.

\subsection{Directional Ablation of the YY Estimator}
\label{sec:exp-directional}

We examine whether direction-dependent extrapolation is necessary to cancel the leading normal bias in higher codimension. We use the
manifold
\(
\mathcal M
=
\left\{
\bigl(
1.5\cos\theta,\,
\sin\theta,\,
0.3\sin(2\theta),\,
0
\bigr)
:\theta\in[0,2\pi)
\right\}
\subset\mathbb R^4,
\)
which has codimension three. All scales use the same preliminary point and
the same first-stage axis, so no additional cross-scale correspondence is
required. We set
\(
r_0=3\sigma,
r_2=12\sigma\sqrt{\log(1/\sigma)},
r_{1,k}=\kappa_k\sigma,
\)
with
\(
\kappa=(1.25,5,8).
\)
Let \(G_k(z)\) denote the scale-dependent YY update from the common
starting point \(z\), let \(U_0\) denote the projector onto the first-stage
axis, and write
\(
    \delta_k = G_k(z)-z.
\)
The single estimator uses only the middle transverse scale,
\(
    z_{\mathrm{Single}}
    =
    z+\delta_2.
\)
For the multiscale constructions, let \(a_k\) and \(b_k\) denote the axial and
complementary weights from Section~\ref{sec:yao_contraction}.
The first scalar control applies the axial weights to the entire displacement,
\(
 z_{\mathrm{Scalar-A}}
    =
    z+\sum_k a_k\delta_k,
\)
whereas the second applies the complementary weights to the entire
displacement,
\(
     z_{\mathrm{Scalar-C}}
    =
    z+\sum_k b_k\delta_k.
\)
The proposed directional estimator instead applies the two weight systems to
their corresponding components,
\(
     z_{\mathrm{Directional}}
    =
    z+
    \sum_k
    \left[
        a_kU_0
        +
        b_k(I-U_0)
    \right]\delta_k.
\)

\begin{table}[!htbp]
\centering
\caption{Fitted log--log slopes of the population errors
with respect to $\sigma$.}
\label{tab:e5-population-slopes}

\footnotesize
\setlength{\tabcolsep}{3pt}
\renewcommand{\arraystretch}{1.2}

\begin{tabular*}{\linewidth}
{@{\extracolsep{\fill}}lccc@{}}
\toprule
Method
& \shortstack{Total Maximum\\error}
& \shortstack{Axial normal\\error}
& \shortstack{Complementary\\normal error} \\
\midrule

\textsc{Single}
& 2.001 & 1.995 & 2.013 \\

\textsc{Scalar-A} (axial weights)
& 2.038 & 3.993 & 2.038 \\

\textsc{Scalar-C} (complementary weights)
& 1.992 & 1.992 & 4.005 \\

\shortstack[l]{\textsc{Directional}\\(componentwise weights)}
& \textbf{4.002} & \textbf{3.993} & \textbf{4.005} \\

\bottomrule
\end{tabular*}

\par\smallskip
\begin{minipage}{\linewidth}
\scriptsize
\textsc{Scalar-A} applies the axial weights to the entire
displacement, while \textsc{Scalar-C} applies the complementary
weights to the entire displacement.
\textsc{Directional} applies separate weight vectors along
the estimated first-stage axis and its ambient orthogonal
complement.
Slopes near $2$ indicate an uncanceled second-order term,
whereas slopes near $4$ indicate fourth-order accuracy.
\end{minipage}
\end{table}

Table~\ref{tab:e5-population-slopes} shows the resulting population behavior.
\textsc{Single} retains second-order bias in both normal components.
\textsc{Scalar-A} removes the axial second-order term but leaves the
complementary component, whereas \textsc{Scalar-C} does the reverse.
Consequently, both scalar controls remain second order in total geometric
error. Only \textsc{Directional} cancels both components simultaneously,
giving a geometric-error slope of \(4.002\), in agreement with the predicted
\(O(\sigma^4)\) population rate.

We next examine the finite-sample behavior of the directional construction.
The noise level is fixed at $\sigma=0.05$, and the three normalized
transverse radii are
\(
    (r_{1,1}, r_{1,2}, r_{1,3})/\sigma=(5,6,12).
\)
The single-scale estimator uses the smallest radius, $r_1=5\sigma$, whereas
\textsc{Scalar-A}, \textsc{Scalar-C}, and \textsc{Directional} combine all
three scales using their corresponding extrapolation weights. We consider
\(
    n\in\{240{,}000,480{,}000,960{,}000\},
\)
with 20 independent repetitions and 256 fixed query points per repetition.

\begin{figure}[!htbp]
    \centering
    \includegraphics[width=0.98\linewidth]{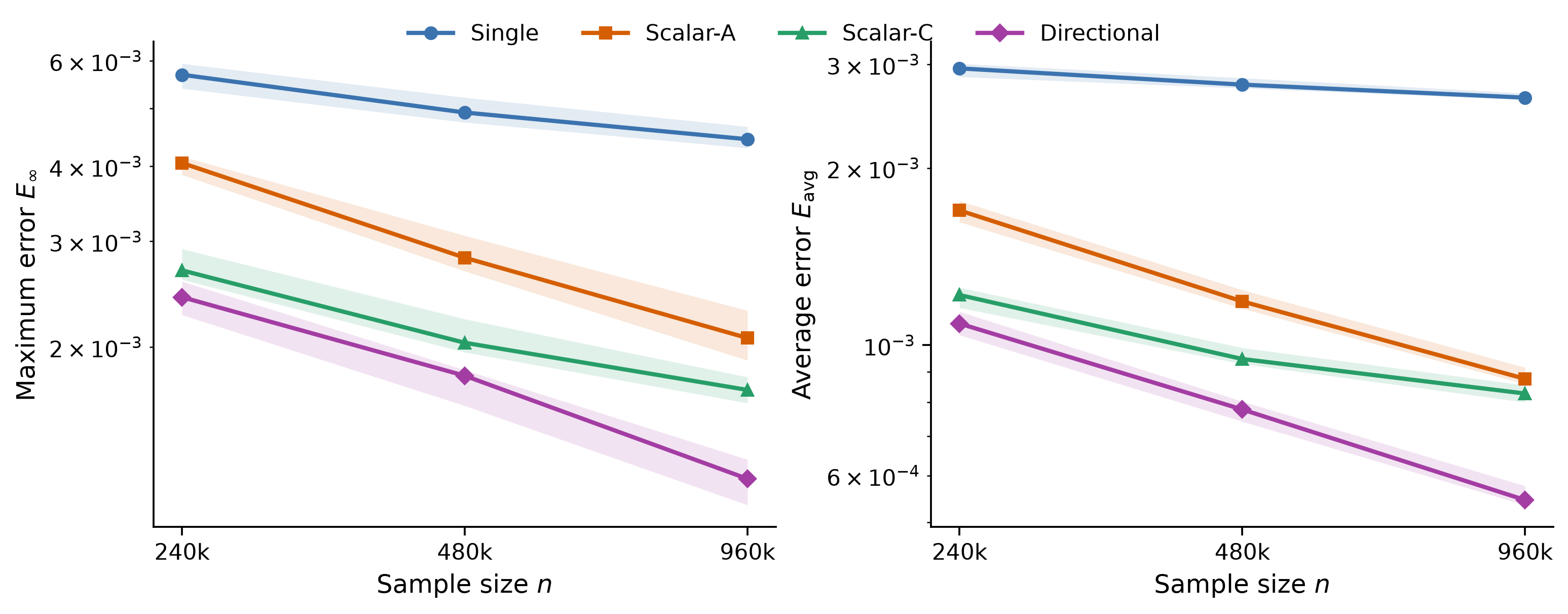}
    \caption{
Finite-sample directional ablation for YY at \(\sigma=0.05\) and
\((r_{1,1}, r_{1,2}, r_{1,3})/\sigma=(5,6,12)\). Left: maximum error \(E_\infty\); right:
average error \(E_{\mathrm{avg}}\). Single uses \(r_1=5\sigma\).
Curves show medians and bands show interquartile ranges over 20 paired
held-out repetitions with 256 queries. The radius configuration was chosen based on preliminary
experiments and kept fixed for the 20 independent
evaluation repetitions.
}
    \label{fig:e5-finite-sample}
\end{figure}

Figure~\ref{fig:e5-finite-sample} reports the median errors over the 20
repetitions, together with the middle $50\%$ range. Under both metrics,
\textsc{Directional} achieves the smallest median error at every sample size,
while the two scalar controls remain less accurate. This ordering shows that
the gain does not come merely from combining several scales, but from applying
different extrapolation weights to the axial and complementary normal
components. At $n=960{,}000$, for example, the median $E_{\infty}$ decreases
from $4.44\times10^{-3}$ for \textsc{Single} to $1.21\times10^{-3}$ for
\textsc{Directional}, while the median $E_{\mathrm{avg}}$ decreases from
$2.63\times10^{-3}$ to $5.46\times10^{-4}$. These correspond to reductions
of approximately $72.8\%$ and $79.3\%$, respectively. The finite-sample
results therefore support the main motivation for the directional YY refinement:
in higher codimension, separating the axial and complementary components
before debiasing can provide a substantial improvement over both the
single-scale estimator and scalar multiscale cancellation.

\endgroup

\section{Conclusion}
\label{sec:conclusion}
We develop geometric bias cancellation (GBC),
a geometric and statistical framework for higher-order
manifold fitting under Gaussian noise.
The central results show that the leading normal bias
can be canceled through estimated-normal correspondence,
and that the resulting improvement persists under
finite-sample fluctuations.
For STSP and Gaussian density ridges, GBC achieves
fourth-order Hausdorff recovery guarantees with explicit
stochastic error bounds.
For STSP, we further establish an asymptotic separation
from the oracle bandwidth-tuned single-scale estimator,
showing that GBC provides a statistical improvement
beyond bandwidth selection.

We further extend geometric bias cancellation to a directional construction
for a smooth-kernel variant of the YY estimator,
where distinct normal-bias components require different
cancellation weights.
The resulting construction achieves fourth-order
population accuracy with a uniform positive-reach
guarantee, together with finite-sample refinement
bounds under the stated preliminary-manifold conditions.

Our framework relies on estimator-specific verification
of the leading normal-bias expansion, suitable geometric
regularity, and uniformly bounded extrapolation weights.
While we establish fourth-order recovery and embeddedness,
a uniform positive-reach guarantee is not established
for the general framework or its empirical estimators.
Such a guarantee is obtained for the population
GBC-YY construction.
The YY refinement requires a preliminary manifold
satisfying the geometric and nondegeneracy conditions
of Assumption~\ref{ass:ysl_full_population}.
Proposition~\ref{prop:yy-sample-split-pilot}
provides a conditional sample-split guarantee,
while establishing these pilot conditions with
competitive total sample complexity remains open.
Developing automatic reference-branch identification and coverage,
as well as stronger geometric guarantees for empirical GBC estimators, remains an interesting direction for future work.

\section*{Acknowledgments}
Zhigang Yao was supported by the Singapore Ministry of Education Tier 2 grant 
A-8001562-00-00 and Tier 1 grants (A8004146-00-00 and A-8002931-00-00) at the National
University of Singapore.  Part of this work originated as a pilot project at the IDSG Summer School held at the Shanghai Institutefor Mathematics and Interdisciplinary Sciences (SIMIS). The authors gratefully acknowledge the
resources and facilities provided by SIMIS, which supported the completion of this work.

\FloatBarrier
\appendix
\numberwithin{equation}{section}
\numberwithin{figure}{section}
\numberwithin{table}{section}
\numberwithin{algorithm}{section}
\renewcommand{\theHsection}{appendix.\Alph{section}}
\section{General Theory and Proofs}
\label{app:general_ms}

This section provides the proofs of the population and finite-sample
results in Section~\ref{Multiscale Bias Cancellation}.

Throughout this section, write
\(P_x=P_{T_x\mathcal M}, Q_x=I-P_x=P_{N_x\mathcal M}, m=D-d\),
and let \(\iota:\mathcal M\hookrightarrow\mathbb R^D\) denote the
inclusion. We use \(\|f\|_\infty=\sup_{x\in\mathcal M}\|f(x)\|\)
and \(h_{\max}=\max_{1\leq k\leq K}h_k\).

We will use the following image-distance comparison. For a map
\(F:\mathcal M\to\mathbb R^D\), put
\(\varepsilon=\|F-\iota\|_\infty\). For every \(x\in\mathcal M\),
\(F(x)\in F(\mathcal M)\), so
\(\operatorname{dist}(x,F(\mathcal M))\leq\|x-F(x)\|\leq\varepsilon\).
Conversely, every \(y\in F(\mathcal M)\) has the form \(y=F(x)\), and
\(\operatorname{dist}(y,\mathcal M)\leq\|y-x\|\leq\varepsilon\).
Taking both suprema gives
\begin{equation}
\label{eq:app-general-image-distance}
d_H(F(\mathcal M),\mathcal M)\leq\|F-\iota\|_\infty.
\end{equation}

\subsection{Proof of content in subsection~\ref{known_direction}}
\subsubsection{Proof of Proposition~\ref{thm:manifold_recovery_1}}
\begin{proof}
For each scale \(h_i\), define the remainder
\(R_{h_i}(x)=u_{h_i}(x)-c(h_i,\sigma)b(x)\),
where \(X_{h_i}=\iota+u_{h_i}\) and
\(u_{h_i}\in\Gamma^1(N\mathcal M)\) by
Assumption~\ref{ass:normal-section}.
By Assumption~\ref{assump:normal_error_t},
\(\|R_{h_i}\|_\infty\leq C_Rh_i^4\).
Using the weight constraints
\(\sum_{i=1}^K w_i=1\) and
\(\sum_{i=1}^K w_ic(h_i,\sigma)=0\), we obtain
\[
\widetilde u(x)
:=\widetilde X(x)-x
=\sum_{i=1}^K w_i u_{h_i}(x)
=\sum_{i=1}^K w_i R_{h_i}(x).
\]
Consequently,
\[
\sup_{x\in\mathcal M}\|\widetilde X(x)-x\|
=\|\widetilde u\|_\infty
\leq\sum_{i=1}^K|w_i|\,\|R_{h_i}\|_\infty
\leq C_R\sum_{i=1}^K|w_i|h_i^4.
\]
This proves the asserted uniform approximation bound.
We next prove that \(\widetilde X\) is a \(C^1\) embedding. Being a
finite linear combination of the \(C^1\) normal sections \(u_{h_i}\),
\(\widetilde u\) belongs to \(\Gamma^1(N\mathcal M)\). Thus
\(\widetilde u(x)\in N_x\mathcal M\) for every \(x\), and
\(\widetilde X=\iota+\widetilde u\) is \(C^1\).

Consider the open normal disk bundle
\(
N_\tau\mathcal M
=
\left\{
(x,v):
x\in\mathcal M,\;
v\in N_x\mathcal M,\;
\|v\|<\tau
\right\}.
\)
Define the normal-bundle map
\(\Phi_\tau: N_\tau\mathcal M \longrightarrow \operatorname{Tub}_\tau(\mathcal M), \Phi_\tau(x,v)=x+v\).
Because
\(
0<\tau<\operatorname{reach}(\mathcal M),
\)
the tubular-neighborhood theorem implies that \(\Phi_\tau\) is a
\(C^1\) diffeomorphism onto
\(\operatorname{Tub}_\tau(\mathcal M)\). Its inverse is
\(\Phi_\tau^{-1}(y) = \left( \pi_{\mathcal M}(y), y-\pi_{\mathcal M}(y) \right)\),
where \(\pi_{\mathcal M}\) is the nearest-point projection. In
particular,
\(
\pi_{\mathcal M}\bigl(\Phi_\tau(x,v)\bigr)=x.
\)

By the assumed smallness condition,
\(
C_R\sum_{i=1}^{K}|w_i|h_i^4<\tau,
\)
and the uniform bound already proved, we have
\(
\|\widetilde u(x)\|<\tau,
\text{for every }x\in\mathcal M.
\)
Therefore,
\(
\bigl(x,\widetilde u(x)\bigr)\in N_\tau\mathcal M
\)
and
\(
\widetilde X(x)
=
\Phi_\tau\bigl(x,\widetilde u(x)\bigr).
\)
It follows that
\(
\pi_{\mathcal M}(\widetilde X(x))
=
x,
\)
or equivalently,
\(
\pi_{\mathcal M}\circ\widetilde X
=
\operatorname{Id}_{\mathcal M}.
\)

This identity immediately implies injectivity. Indeed, if
\(
\widetilde X(x_1)=\widetilde X(x_2),
\)
then applying \(\pi_{\mathcal M}\) to both sides yields
\(
x_1
=
\pi_{\mathcal M}(\widetilde X(x_1))
=
\pi_{\mathcal M}(\widetilde X(x_2))
=
x_2.
\)

Differentiating \(\pi_{\mathcal M}\circ\widetilde X
=\operatorname{Id}_{\mathcal M}\) gives
\(D\pi_{\mathcal M}|_{\widetilde X(x)}\circ D\widetilde X|_x =I_{T_x\mathcal M}\).
Set \(A_x=D\widetilde X|_x\in\mathcal L(T_x\mathcal M,\mathbb R^D)\)
and \(B_x=D\pi_{\mathcal M}|_{\widetilde X(x)}
\in\mathcal L(\mathbb R^D,T_x\mathcal M)\). Since
\(B_xA_x=I_{T_x\mathcal M}\), \(A_xv=0\) implies
\(v=B_xA_xv=0\). Hence \(\ker A_x=\{0\}\),
\(\operatorname{rank}D\widetilde X|_x=d\), and \(\widetilde X\) is an
immersion. Quantitatively,
\[
\|v\|=\|B_xA_xv\|\leq\|B_x\|_{\mathrm{op}}\|A_xv\|,\qquad
\|D\widetilde X|_xv\|
\geq\frac{\|v\|}{\|D\pi_{\mathcal M}|_{\widetilde X(x)}\|_{\mathrm{op}}}.
\]
Compactness of \(\mathcal M\) and continuity of \(x\mapsto B_x\) give
\(C_\pi:=\sup_{x\in\mathcal M}\|B_x\|_{\mathrm{op}}<\infty\).
Thus \(\|D\widetilde X|_xv\|\geq C_\pi^{-1}\|v\|\) uniformly in \(x\),
so the differential is uniformly left-invertible.

Because \(\mathcal M\) is compact and \(\mathbb R^D\) is Hausdorff,
the continuous injective map \(\widetilde X\) is a homeomorphism onto
its image. Together with the immersion property, this proves that
\(\widetilde X\) is a \(C^1\) embedding and
\(\widetilde{\mathcal M}=\widetilde X(\mathcal M)\) is a
\(d\)-dimensional \(C^1\) embedded submanifold of \(\mathbb R^D\).
Its inverse is \(\pi_{\mathcal M}|_{\widetilde{\mathcal M}}\), which is
\(C^1\) because \(\pi_{\mathcal M}\) is \(C^1\) on the tube. Hence
\(\widetilde X:\mathcal M\to\widetilde{\mathcal M}\) is a \(C^1\)
diffeomorphism.

Finally, applying \eqref{eq:app-general-image-distance} with
\(F=\widetilde X\) and the uniform approximation bound above gives the
asserted Hausdorff bound and completes the proof.
\end{proof}

\subsection{Proof of content in subsection~\ref{sec:normal_guided_extrapolation}}
\label{app:general_enf_proofs}

For the correspondence proofs, use the reference points and projector fields
\[
y_1(x)=X_{h_1}(x),\qquad
\widehat y_1(x)=\widehat X_{h_1}(x),\qquad
\widetilde Q(x)=P_{\widetilde N(x)},\qquad
\widehat Q(x)=P_{\widehat N(x)}.
\]
Whenever the corresponding intersections are defined, write
\[
x_k=\pi_{\mathcal M}(Z_{h_k}(x)),\qquad
v_k(x)=Z_{h_k}(x)-y_1(x).
\]
\subsubsection{Proof of Lemma~\ref{lem:pca_normal_accuracy}}
\begin{proof}
Fix \(x\in\mathcal M\).
We work uniformly with centers
\[
y=x+h^2\beta,
\qquad
\beta\in N_x\mathcal M,
\qquad
\|\beta\|\le B,
\]
where \(B<\infty\) is fixed.  The case \(y=x\) corresponds to
\(\beta=0\).

Set \(\lambda=\sigma/h\in\Lambda\Subset(0,\infty)\), where
\(\Lambda\) is fixed. All \(C^2\) bounds in this proof are taken in
local normal-bundle coordinates \((x,\beta)\), with \(h,\lambda\)
held fixed, and are uniform over \(x\in\mathcal M\),
\(\lambda\in\Lambda\), and \(\|\beta\|\leq B\).

Choose a finite collection of local orthonormal frames adapted to
\(\mathbb R^D = T_x\mathcal M\oplus N_x\mathcal M\),
and write
\(v=v_T+v_N, \qquad v_T\in T_x\mathcal M, \quad v_N\in N_x\mathcal M\).

Let \(p_\sigma\) denote the density of the noisy observation
\(Y=X+\xi\).  After the rescaling
\(Y=y+hv, \qquad \|v\|\le1\),
the kernel-weighted law of \(v\) has density proportional to
\(K(\|v\|^2)p_\sigma(y+hv)\) on the unit ball, where
\(K(t)=(1-t)_+^6\). We write \(\mathbb E_h\) for expectation under
this weighted law.

We use the second fundamental form convention
\[
\Psi_x(t)=x+t+\frac12\mathrm{II}_x(t,t)+O(\|t\|^3),
\qquad
H(x)=\frac1d\operatorname{tr}\mathrm{II}_x,
\qquad t\in T_x\mathcal M.
\]
Using this \(C^4\) local graph representation, together with
\(\sigma=\lambda h, y-x=h^2\beta\),
the Gaussian convolution admits the expansion
\begin{equation}
\label{eq:pca_density_expansion}
p_\sigma(y+hv)
=
C_{h,x,\lambda}\,
\rho_\lambda(v_N)
\left[
1+hL_{x,\beta,\lambda}(v)
+O_{C^2}(h^2)
\right],
\end{equation}
where
\(\rho_\lambda(v_N) = \exp\!\left( -\frac{\|v_N\|^2}{2\lambda^2} \right)\).
Compactness of \(\mathcal M\), the assumed \(C^6\) regularity,
and Gaussian smoothing allow the local density expansion to
be differentiated twice with respect to the base point and the bounded
normal coordinate \(\beta\), without changing the displayed powers of
\(h\).

The first-order correction satisfies
\[
L_{x,\beta,\lambda}(-v_T,v_N)
=
L_{x,\beta,\lambda}(v_T,v_N),
\]
and
\[
L_{x,\beta,\lambda}(v_T,-v_N)
=
-
L_{x,\beta,\lambda}(v_T,v_N).
\]
More precisely,
\[
L_{x,\beta,\lambda}(v)
=
\frac{1}{\lambda^2}
\left\langle
v_N,
A_{x,\beta,\lambda}(v_T)
\right\rangle,
\]
where
\[
A_{x,\beta,\lambda}(v_T)
=
-\beta+\frac12\mathrm{II}_x(v_T,v_T)
+\frac{d\lambda^2}{2}H(x)
\in N_x\mathcal M
\]
is even in \(v_T\).  These parity relations are preserved after taking
up to two derivatives in \((x,\beta)\).

Introduce the probability measure
\[
d\nu_\lambda(v)
=
\frac{
K(\|v\|^2)\rho_\lambda(v_N)\,dv
}{
Z_\lambda
},
\qquad
Z_\lambda
=
\int_{\|v\|\le1}
K(\|v\|^2)\rho_\lambda(v_N)\,dv .
\]
Write \(\mathbb E_\lambda\) for expectation under \(\nu_\lambda\).
The radial kernel preserves the reflection symmetries. Since
\(L_{x,\beta,\lambda}\) is odd in \(v_N\),
\(\mathbb E_\lambda \left[ L_{x,\beta,\lambda}(v) \right] =0\).
Therefore normalization of
\eqref{eq:pca_density_expansion} introduces no \(O(h)\) correction.
For every polynomially bounded \(f\) on the unit ball,
\begin{equation}
\label{eq:pca_expectation_expansion}
\mathbb E_h[f(v)]
=
\mathbb E_\lambda[f(v)]
+
h\,
\mathbb E_\lambda
\left[
f(v)L_{x,\beta,\lambda}(v)
\right]
+
O_{C^2}(h^2),
\end{equation}

By rotational invariance of \(\nu_\lambda\) under
\(O(d)\times O(m)\), there exist positive quantities
\(\tau_\lambda\) and \(a_\lambda\) such that
\[
\mathbb E_\lambda[v_Tv_T^\top]
=
\tau_\lambda P_x,
\qquad
\mathbb E_\lambda[v_Nv_N^\top]
=
a_\lambda Q_x,
\]
and
\(\mathbb E_\lambda[v_Nv_T^\top]=0\).
Moreover, the parity relations imply
\[
\begin{aligned}
\mathbb E_\lambda[v_Tv_T^\top L_{x,\beta,\lambda}(v)]&=0,\\
\mathbb E_\lambda[v_Nv_N^\top L_{x,\beta,\lambda}(v)]&=0,\\
\mathbb E_\lambda[v_Nv_T^\top L_{x,\beta,\lambda}(v)]&=0.
\end{aligned}
\]
Applying \eqref{eq:pca_expectation_expansion} therefore gives
\[
\begin{aligned}
\mathbb E_h[v_Tv_T^\top]&=\tau_\lambda P_x+O_{C^2}(h^2),\\
\mathbb E_h[v_Nv_N^\top]&=a_\lambda Q_x+O_{C^2}(h^2),\\
\mathbb E_h[v_Nv_T^\top]&=O_{C^2}(h^2).
\end{aligned}
\]

The corresponding first moments satisfy
\(\mathbb E_h[v_T]=O_{C^2}(h^2), \mathbb E_h[v_N]=O_{C^2}(h)\).
Hence subtracting the weighted mean changes the second moments only
by \(O_{C^2}(h^2)\).  Therefore
\[
\operatorname{Cov}_h(v)
=
\tau_\lambda P_x
+
a_\lambda Q_x
+
O_{C^2}(h^2).
\]
Returning to the original variable \(Y=y+hv\), we obtain
\begin{equation}
\label{eq:pca_covariance_expansion}
\Sigma_h(x+h^2\beta)
=
h^2
\left(
\tau_\lambda P_x
+
a_\lambda Q_x
\right)
+
O_{C^2}(h^4).
\end{equation}
In particular, its \(C^0\) part gives
\(P_x\Sigma_h(y)P_x = h^2\tau_\lambda P_x+O(h^4)\),
\(Q_x\Sigma_h(y)Q_x = h^2a_\lambda Q_x+O(h^4)\),
and
\begin{equation}
\label{eq:pca_cross_block}
Q_x\Sigma_h(y)P_x
=
O(h^4).
\end{equation}

It remains to verify that the leading tangent and normal eigenvalues
are uniformly separated.  Let
\(
S=\|v_N\|^2.
\)
By rotational invariance in the normal coordinates,
\(
a_\lambda
=
\frac{1}{m}\mathbb E_\lambda[S].
\)
Conditionally on \(v_N\), the tangent component \(v_T\) has density
proportional to \((1-S-\|v_T\|^2)_+^6\) on the \(d\)-dimensional ball
of radius \(\sqrt{1-S}\). Radial integration gives
\(\mathbb E_\lambda[v_{T,1}^2\mid v_N]=(1-S)/(d+14)\). Hence
\(
\tau_\lambda
=
\frac{1}{d+14}
\mathbb E_\lambda[1-S].
\)
Thus
\(
\tau_\lambda-a_\lambda>0
\)
is equivalent to
\(
\mathbb E_\lambda[S]
<
\frac{m}{D+14}.
\)

For the radial density proportional to \(K(\|v\|^2)\) on the unit
\(D\)-ball, corresponding formally to \(\lambda=\infty\), beta
integration gives
\(
\mathbb E_\infty[S]
=
\frac{m}{D+14}.
\)
For every finite \(\lambda>0\), the measure
\(\nu_\lambda\) is obtained from this kernel-weighted measure by the
strictly decreasing weight
\(
S\longmapsto
\exp\!\left(-\frac{S}{2\lambda^2}\right).
\)
Since \(S\) is nonconstant, this reweighting strictly decreases its
expectation.  Therefore
\(
\mathbb E_\lambda[S]
<
\frac{m}{D+14},
\)
and consequently
\(
\tau_\lambda>a_\lambda
\text{for every }\lambda\in(0,\infty).
\)
Continuity in \(\lambda\), together with compactness of \(\Lambda\),
gives
\(
g_*
:=
\inf_{\lambda\in\Lambda}
(\tau_\lambda-a_\lambda)
>0.
\)

It follows from \eqref{eq:pca_covariance_expansion} and Weyl's
inequality that, for all sufficiently small \(h\),
\[
\lambda_d(\Sigma_h(y))
-
\lambda_{d+1}(\Sigma_h(y))
\ge
\frac{g_*}{2}h^2.
\]
Thus the tangent--normal population eigengap is uniformly of order
\(h^2\).

For the projector estimate, define
\(\overline\Sigma_h(x,\beta):=h^{-2}\Sigma_h(x+h^2\beta),\qquad \Sigma_0(x,\lambda):=\tau_\lambda P_x+a_\lambda Q_x\).
By \eqref{eq:pca_covariance_expansion},
\(\overline\Sigma_h=\Sigma_0+O_{C^2}(h^2)\), and the limiting matrix
has the uniformly positive gap \(\tau_\lambda-a_\lambda\geq g_*>0\).
Hence the spectral projector onto its leading \(d\)-dimensional
eigenspace is a smooth function of the matrix in a fixed neighborhood
of this family.  Smooth spectral-projector calculus therefore yields
\(P_h(x+h^2\beta) = P_x+O_{C^2}(h^2)\),
uniformly in
\(x\in\mathcal M\),
\(\lambda\in\Lambda\), and bounded \(\beta\).
Since \(Q_h=I-P_h\),
\(Q_h(x+h^2\beta) = Q_x+O_{C^2}(h^2)\).

Taking only the \(C^0\) part recovers, in particular,
\[
\left\|
P_h(y)-P_{T_x\mathcal M}
\right\|_{\mathrm{op}}
+
\left\|
Q_h(y)-P_{N_x\mathcal M}
\right\|_{\mathrm{op}}
\le Ch^2.
\]
This proves the result.
\end{proof}
\subsubsection{Proof of Lemma~\ref{lem:wellposed_correspondence}}

\begin{proof}
We first work in codimension one, using compatibly oriented local
unit normals \(n\) and \(\widetilde n\).
Since the codimension is one, there exists a scalar
\(\alpha_k(x)\) such that
\(
u_{h_k}(x)=\alpha_k(x)n(x).
\)
Similarly, \(y_1(x)=x+\alpha_1(x)n(x)\).
Therefore,
\(
X_{h_k}(x)-y_1(x)
=
\bigl(\alpha_k(x)-\alpha_1(x)\bigr)n(x).
\)
By Assumption~\ref{assump:normal_error_t},
\(|\alpha_k(x)|=O(h_k^2), |\alpha_1(x)|=O(h_1^2)\),
uniformly in \(x\), and hence
\(\|X_{h_k}(x)-y_1(x)\| = O(h_k^2+h_1^2) = O(h_{\max}^2)\).

We first show that the estimated normal direction
\(\widetilde n(x)\) remains non-tangent to
\(\mathcal M_{h_k}\). Denote
\(P_{k,x}=P_{T_{X_{h_k}(x)}\mathcal M_{h_k}}\).
Since \(P_xn(x)=0\), we have
\[
\|P_{k,x}\widetilde n(x)\|
\leq
\|(P_{k,x}-P_x)\widetilde n(x)\|
+
\|P_x(\widetilde n(x)-n(x))\|
\leq
C_T h_k^2+C_N h_1^2.
\]
Hence, for sufficiently small \(h_{\max}\),
\(
\|P_{k,x}\widetilde n(x)\|<1.
\)
Since \(\widetilde n(x)\) is a unit vector, it follows that
\(
\widetilde n(x)
\notin
T_{X_{h_k}(x)}\mathcal M_{h_k}.
\)
Therefore,
\[
T_{X_{h_k}(x)}\mathcal M_{h_k}
\oplus
\operatorname{span}\{\widetilde n(x)\}
=
\mathbb R^D.
\]

Choose a local \(C^2\) parametrization
\(\psi_x:B_{r_0}(0)\subset\mathbb R^d\to\mathcal M\) with
\(\psi_x(0)=x\) and \(D\psi_x(0)\) an isometry onto
\(T_x\mathcal M\), and set
\(\varphi_k(q):=X_{h_k}(\psi_x(q))\).
The radius \(r_0>0\) and the first- and second-derivative bounds
for \(\psi_x\) can be chosen uniformly in \(x\), independently of
the scales. In particular, \(\varphi_k(0)=X_{h_k}(x)\).
The intersection between the estimated normal line
\(\widetilde L_{h_1}(x;s) = y_1(x)+s\widetilde n(x)\)
and \(\mathcal M_{h_k}\) is characterized by
\[
F_k(q,s;x)
:=
\varphi_k(q)
-
y_1(x)
-
s\widetilde n(x)
=
0.
\]
The derivative with respect to \((q,s)\) is
\(
D_{(q,s)}F_k
=
\begin{bmatrix}
D\varphi_k(q)
&
-\widetilde n(x)
\end{bmatrix}.
\)
At \(q=0\),
\[
\operatorname{Ran}(D\varphi_k(0))
=
T_{X_{h_k}(x)}\mathcal M_{h_k}.
\]
The direct-sum relation above therefore implies that
\(
\begin{bmatrix}
D\varphi_k(0)
&
-\widetilde n(x)
\end{bmatrix}
\)
is invertible.

For the quantitative argument, the normal-graph identity,
\(\|u_{h_k}\|_{C^0}=O(h_k^2)\), and
Assumption~\ref{ass:tangent-stability} give
\[
\sup_{z\in\mathcal M}
\|DX_{h_k}(z)-D\iota(z)\|_{\mathrm{op}}\leq Ch_k^2.
\]
Indeed, differentiating the normality of \(u_{h_k}\) bounds the
tangential part of \(Du_{h_k}\) by \(C\|u_{h_k}\|\), and the
tangent-projector gap then bounds its normal part by \(Ch_k^2\).
Consequently, the inverse of the displayed matrix has norm at most
\(C_I\), uniformly in \(x\) and the sufficiently small scales.
Moreover,
\(\|D\varphi_k(q)-D\varphi_k(0)\|_{\mathrm{op}} \leq C\bigl(\|q\|+h_k^2\bigr)\).

To locate a nearby solution, define
\(
s_k^0(x)
=
\alpha_k(x)-\alpha_1(x).
\)
Since
\(
X_{h_k}(x)
=
y_1(x)+s_k^0(x)n(x),
\)
we obtain
\[
F_k(0,s_k^0(x);x)
=
X_{h_k}(x)
-y_1(x)
-s_k^0(x)\widetilde n(x)
=
s_k^0(x)
\bigl(n(x)-\widetilde n(x)\bigr).
\]
Using Assumption~\ref{ass:estimated-normal-bundle},
\[
\|F_k(0,s_k^0(x);x)\|
\leq
|s_k^0(x)|
\,
\|n(x)-\widetilde n(x)\|
=
O(h_{\max}^2h_1^2).
\]
Thus, for sufficiently small \(h_{\max}\), the estimated-normal
intersection equation is a sufficiently small perturbation of the
true-normal intersection equation.

Write \(\zeta=(q,s)\), \(\zeta_0=(0,s_k^0(x))\), and
\(A_{k,x}=D_{(q,s)}F_k(\zeta_0;x)\).
The preceding derivative estimate, together with the affine dependence
on \(s\), allows us to choose a fixed \(r>0\) such that
\[
\sup_{\|\zeta-\zeta_0\|\leq r}
\left\|I-A_{k,x}^{-1}D_{(q,s)}F_k(\zeta;x)\right\|_{\mathrm{op}}
\leq \frac12
\]
for all sufficiently small scales, uniformly in \(x\).
Define
\(\mathcal T_{k,x}(\zeta) :=\zeta-A_{k,x}^{-1}F_k(\zeta;x)\).
This is a contraction on the closed ball
\(\overline B_r(\zeta_0)\). Since
\[
\|\mathcal T_{k,x}(\zeta_0)-\zeta_0\|
\leq C_I\|F_k(\zeta_0;x)\|
=O(h_1^2h_{\max}^2),
\]
its displacement at the center is at most \(r/2\) for sufficiently
small \(h_{\max}\), so it maps this ball into itself.
The contraction mapping theorem therefore gives a unique zero
\(\zeta_k=(q_k(x),s_k(x))\) of \(F_k\) in this neighborhood, with
\[
\|\zeta_k-\zeta_0\|
\leq 2C_I\|F_k(\zeta_0;x)\|
=O(h_1^2h_{\max}^2).
\]
Hence the estimated normal line
\(\widetilde L_{h_1}(x;\cdot)\) has a unique local intersection with
\(\mathcal M_{h_k}\), given by
\(
Z_{h_k}(x)
=
\varphi_k(q_k(x))
=
y_1(x)+s_k(x)\widetilde n(x).
\)

The bounds in Assumptions~\ref{ass:estimated-normal-bundle} and
\ref{ass:tangent-stability} are uniform in \(x\), and
\(\mathcal M\) is compact. Therefore the above construction holds
uniformly over \(x\in\mathcal M\) for sufficiently small
\(h_{\max}\).

Moreover, if \(x\mapsto P_{\widetilde N(x)}\) is \(C^1\), then
\(\widetilde n\) can be chosen \(C^1\) locally. In these local
coordinates, \(F_k(q,s;x)\) is \(C^1\) in \((q,s,x)\), since
\(y_1=X_{h_1}\) and the local parametrizations of
\(\mathcal M_{h_k}\) are \(C^1\). Since
\(D_{(q,s)}F_k\) is invertible, the parameterized implicit function
theorem implies that \(q_k(x)\) and \(s_k(x)\) depend \(C^1\)-smoothly
on \(x\). Consequently,
\(
Z_{h_k}(x)
=
y_1(x)+s_k(x)\widetilde n(x)
\)
defines a \(C^1\) map
\(
Z_{h_k}:\mathcal M\rightarrow\mathcal M_{h_k}.
\)

For \(m>1\), use the affine normal fiber
\(y_1(x)+\widetilde N(x)\). At each fixed \(x\), choose an
orthonormal frame \(\widetilde E_\alpha(x)\in\mathbb R^{D\times m}\)
of \(\widetilde N(x)\), choosing \(C^1\) local frames when the
projector field is \(C^1\). The intersection equation in this
trivialization is
\[
F_{k,\alpha}(q,s;x)
:=\varphi_k(q)-y_1(x)-\widetilde E_\alpha(x)s=0,
\qquad s\in\mathbb R^m.
\]
The projector and tangent-space bounds imply the uniform transversality
\[
T_{X_{h_k}(x)}\mathcal M_{h_k}\oplus\widetilde N(x)
=\mathbb R^D,
\]
so \([D\varphi_k(0),-\widetilde E_\alpha(x)]\) has a uniformly
bounded inverse. Set
\[
\Delta_k(x):=X_{h_k}(x)-y_1(x)\in N_x\mathcal M,
\qquad
s_{k,\alpha}^0(x):=\widetilde E_\alpha(x)^\top\Delta_k(x).
\]
Then
\[
F_{k,\alpha}(0,s_{k,\alpha}^0(x);x)
=\bigl(P_{N_x\mathcal M}-P_{\widetilde N(x)}\bigr)\Delta_k(x),
\]
and hence
\[
\|F_{k,\alpha}(0,s_{k,\alpha}^0(x);x)\|
\leq Ch_1^2(h_k^2+h_1^2).
\]
The same quantitative contraction argument gives a unique nearby
intersection. When the projector field is \(C^1\), the parameterized
implicit-function theorem gives its \(C^1\) dependence on the base
point in each local trivialization.

On an overlap, two local orthonormal frames satisfy
\(\widetilde E_\beta(x)=\widetilde E_\alpha(x)R_{\alpha\beta}(x)\)
with \(R_{\alpha\beta}(x)\in O(m)\). Their coordinates are related by
\(s_\beta=R_{\alpha\beta}(x)^\top s_\alpha\), so they describe the
same ambient intersection point. Local uniqueness therefore glues the
local constructions into the globally defined correspondence
\(Z_{h_k}\), without assuming a global normal frame.

\end{proof}

\subsubsection{Proof of Lemma~\ref{lem:normal_guided_correspondence}}
\begin{lemma}[Error induced by normal-guided correspondence]
\label{lem:normal_guided_correspondence}
Under Assumptions~\ref{ass:geometry-observation},
\ref{assump:normal_error_t}, and
\ref{ass:estimated-normal-bundle},
\(\|x_k-x\| = O(h_1^2h_k^2)\)
uniformly in \(x\in\mathcal M\). In particular, if the scales are of
the same order, then
\(\|x_k-x\| = O(h_{\max}^4)\).
\end{lemma}

\begin{proof}
By construction, \(v_k(x)\in\widetilde N(x)\).
By Assumption~\ref{ass:estimated-normal-bundle},
\(
\|P_xv_k(x)\|
\leq
C_Nh_1^2\|v_k(x)\|.
\)

Now define
\(
z_k^0
=
y_1(x)+Q_x v_k(x).
\)
Since \(y_1(x)-x=u_{h_1}(x)\in N_x\mathcal M\),
\(z_k^0-x\in N_x\mathcal M\).
For sufficiently small scales, \(z_k^0\) lies in the tubular
neighborhood of \(\mathcal M\), and hence
\(
\pi_{\mathcal M}(z_k^0)=x.
\)
On the other hand,
\(
\pi_{\mathcal M}(Z_{h_k}(x))=x_k.
\)
Since the nearest-point projection is smooth, and therefore uniformly
Lipschitz on a sufficiently small tubular neighborhood of the compact
manifold \(\mathcal M\), there exists \(C_\pi>0\) such that
\begin{align*}
\|x_k-x\|
&=
\left\|
\pi_{\mathcal M}(Z_{h_k}(x))
-
\pi_{\mathcal M}(z_k^0)
\right\|\\
&\leq
C_\pi
\|Z_{h_k}(x)-z_k^0\|\\
&=
C_\pi\|P_xv_k(x)\|\\
&\leq
C_\pi C_Nh_1^2\|v_k(x)\|.
\end{align*}

It remains to control \(\|v_k(x)\|\). Since
\(Z_{h_k}(x) = x_k+u_{h_k}(x_k), \qquad y_1(x) = x+u_{h_1}(x)\),
we have
\[
\|v_k(x)\|
\leq
\|x_k-x\|
+
\|u_{h_k}(x_k)\|
+
\|u_{h_1}(x)\|.
\]
By Assumption~\ref{assump:normal_error_t},
\(\sup_{z\in\mathcal M}\|u_{h_k}(z)\|=O(h_k^2), \sup_{z\in\mathcal M}\|u_{h_1}(z)\|=O(h_1^2)\).
Therefore,
\(\|v_k(x)\| \leq \|x_k-x\| + O(h_k^2+h_1^2)\).
Combining the preceding inequalities gives
\[
\|x_k-x\|
\leq
C h_1^2
\left(
\|x_k-x\|+h_k^2+h_1^2
\right).
\]
For sufficiently small \(h_1\), the first term on the right-hand side
can be absorbed into the left-hand side, yielding
\[
\|x_k-x\|
\leq
C h_1^2(h_k^2+h_1^2).
\]
Since \(h_1\leq h_k\),
\(
\|x_k-x\|
=
O(h_1^2h_k^2).
\)
If all scales are of the same order, this further gives
\(
\|x_k-x\|
=
O(h_{\max}^4).
\)
\end{proof}
\subsubsection{Proof of Lemma~\ref{lem:normal_field_stability}}
\begin{lemma}[Stability of normal displacement under normal-guided correspondence]
\label{lem:normal_field_stability}
Suppose Assumptions~\ref{ass:geometry-observation},
\ref{ass:normal-section}, \ref{assump:normal_error_t},
\ref{ass:estimated-normal-bundle}, and
\ref{ass:tangent-stability} hold. Let
\(x_k=\pi_{\mathcal M}(Z_{h_k}(x))\) be as in
Lemma~\ref{lem:normal_guided_correspondence}. Then
\[
\|u_{h_k}\|_{C^1(\mathcal M)}\leq Ch_k^2,
\qquad
\|u_{h_k}(x_k)-u_{h_k}(x)\|
\leq Ch_1^2h_k^4.
\]
The bounds are uniform in \(x\) and the admissible scales.
\end{lemma}

\begin{proof}
Assumption~\ref{assump:normal_error_t} and the boundedness of \(b\)
give \(\|u_h\|_\infty\leq Ch^2\).
Write \(P_x=P_{T_x\mathcal M}\), \(Q_x=I-P_x\), and
\(P_{h,x}=P_{T_{X_h(x)}\mathcal M_h}\).
For \(v\in T_x\mathcal M\), differentiating the normality of
\(u_h\) and using the bounded second fundamental form gives
\[
\|P_xD_vu_h\|\leq C\|u_h\|_\infty\|v\|
\leq Ch^2\|v\|.
\]
Let \(w=D_vX_h=v+D_vu_h\in T_{X_h(x)}\mathcal M_h\).
Since \(P_{h,x}w=w\) and \(Q_xv=0\),
Assumption~\ref{ass:tangent-stability} yields
\[
\begin{aligned}
\|Q_xD_vu_h\|=\|Q_xw\|
&=\|Q_x(P_{h,x}-P_x)w\|\\
&\leq Ch^2\bigl(\|v\|+\|D_vu_h\|\bigr).
\end{aligned}
\]
Combining the tangent and normal bounds and absorbing the
\(Ch^2\|D_vu_h\|\) term for small \(h\), we obtain
\(\|D_vu_h\|\leq Ch^2\|v\|\), uniformly in \(x\).
Together with the position bound, this proves
\(\|u_h\|_{C^1(\mathcal M)}\leq Ch^2\).

For sufficiently close points on the compact manifold,
geodesic distance is bounded by a constant times Euclidean distance.
Lemma~\ref{lem:normal_guided_correspondence} therefore gives
\[
\begin{aligned}
\|u_{h_k}(x_k)-u_{h_k}(x)\|
&\leq Ch_k^2d_{\mathcal M}(x_k,x)\\
&\leq Ch_k^2\|x_k-x\|\\
&\leq Ch_1^2h_k^4.
\end{aligned}
\]
\end{proof}
\subsubsection{Proof of Proposition~\ref{thm:normal_guided_manifold_recovery}}

\begin{proof}[Proof of the moreover part]
Since \(\widetilde Q\) is \(C^1\), the total space
\(\widetilde{\mathcal N} :=\{(x,v):x\in\mathcal M,\ v\in\widetilde N(x)\}\)
is a rank-\((D-d)\) \(C^1\) vector bundle over \(\mathcal M\).
For \(\rho>0\), define its open disk bundle and the affine-fiber map by
\[
\widetilde{\mathcal N}_\rho
:=\{(x,v)\in\widetilde{\mathcal N}:\|v\|<\rho\},
\qquad
\Phi_{h_1}(x,v):=y_1(x)+v.
\]
These definitions do not require a global orthonormal frame.

The normal-graph and tangent-space estimates used in the proof of
Lemma~\ref{lem:wellposed_correspondence} give
\(\|X_{h_1}-\iota\|_{C^1}=O(h_1^2)\). Thus, for all
sufficiently small \(h_1\), there is a uniform \(c_0>0\) such that
\[
\|(I-\widetilde Q(x))DX_{h_1}(x)a\|
\geq c_0\|a\|,
\qquad a\in T_x\mathcal M.
\]
If \((a,w)\in T_{(x,v)}\widetilde{\mathcal N}\), differentiating
\(\widetilde Q(x)v=v\) gives
\((I-\widetilde Q(x))w=D\widetilde Q(x)[a]v\).
Consequently, \(D\Phi_{h_1}(x,v)(a,w)=0\) implies
\[
c_0\|a\|
\leq\|D\widetilde Q(x)[a]v\|
\leq C\rho\|a\|.
\]
Choosing \(\rho\) sufficiently small, independently of \(h_1\),
forces \(a=0\) and then \(w=0\). Since the domain and codomain have
the same dimension, \(\Phi_{h_1}\) is a local \(C^1\)
diffeomorphism on \(\widetilde{\mathcal N}_\rho\).

We also verify global injectivity on a uniformly small disk bundle.
Suppose
\(y_1(x)+v=y_1(x')+v'\), with \(\|v\|,\|v'\|<\rho\).
The \(C^0\) normal-graph bound gives
\(\|x'-x\|\leq2\rho+Ch_1^2\), so the base points lie in a common
fixed-size local chart when \(\rho\) and \(h_1\) are small.
The geometry of \(\mathcal M\), the \(C^1\) bound on
\(X_{h_1}-\iota\), and the projector approximation imply
\[
\|(I-\widetilde Q(x))(y_1(x')-y_1(x))\|
\geq c_1\|x'-x\|
\]
for a uniform \(c_1>0\). On the other hand, the assumed equality and
\(v\in\widetilde N(x)\), \(v'\in\widetilde N(x')\) give
\[
\begin{aligned}
\|(I-\widetilde Q(x))(y_1(x')-y_1(x))\|
&=\|(I-\widetilde Q(x))v'\|\\
&\leq\|\widetilde Q(x')-\widetilde Q(x)\|_{\mathrm{op}}\|v'\|\\
&\leq C\rho\|x'-x\|.
\end{aligned}
\]
Here the last bound follows from the uniform \(C^1\) projector bound
and local equivalence of intrinsic and Euclidean distances.
Reducing \(\rho\) so that \(C\rho<c_1\) forces \(x'=x\), and then
\(v'=v\). Hence \(\Phi_{h_1}\) is an injective local diffeomorphism,
and therefore a \(C^1\) embedding onto its image.

By Lemma~\ref{lem:wellposed_correspondence}, the sections \(v_k\)
are \(C^1\) and satisfy \(v_k(x)\in\widetilde N(x)\).
The normal-direction bias expansion and
Lemma~\ref{lem:normal_guided_correspondence} give
\(\sup_x\|v_k(x)\|=O(h_{\max}^2)\).
Since \(\sum_{k=1}^K w_k=1\), we have
\[
\widetilde v(x):=\sum_{k=1}^K w_kv_k(x)\in\widetilde N(x),
\qquad
\widetilde X(x)
=y_1(x)+\widetilde v(x)
=\Phi_{h_1}(x,\widetilde v(x)).
\]
For fixed or uniformly bounded extrapolation weights,
\(\sup_x\|\widetilde v(x)\|=O(h_{\max}^2)<\rho\) for sufficiently
small \(h_{\max}\).
The graph of the section,
\(x\mapsto(x,\widetilde v(x))\), is a \(C^1\) embedding into
\(\widetilde{\mathcal N}_\rho\), since the bundle projection is its
left inverse. Composing it with \(\Phi_{h_1}\) shows that
\(\widetilde X\) is a \(C^1\) embedding. Consequently,
\(\widetilde{\mathcal M}=\widetilde X(\mathcal M)\) is a
\(d\)-dimensional \(C^1\) embedded manifold that is
\(C^1\)-diffeomorphic to \(\mathcal M\).
\end{proof}

For the following two finite-sample proofs, let \(E_X,E_N\) denote
the events in Assumptions~\ref{assump:normal_guided_empirical} and
\ref{assump:finite_normal_estimation}, respectively, and abbreviate
\(e_{n,k}^{\mathrm{corr}}:=r_n(h_1)+r_n(h_k)+s_n(h_1)h_k^2\).
The assumed probability bounds and the union bound give
\[
\mathbb P(E_X)\geq1-\delta_n^X,\qquad
\mathbb P(E_N)\geq1-\delta_n^N,\qquad
\mathbb P(E_X\cap E_N)\geq1-\delta_n^X-\delta_n^N.
\]

\begin{proof}[Proof of Proposition~\ref{prop:enf_local_conditioning}]
Write $P_x=P_{T_x\mathcal M}$ and $Q_x=I_D-P_x$.
By Assumption~\ref{assump:normal_error_t} and the tangent-space stability in
Assumption~\ref{ass:tangent-stability}, the argument of Lemma~\ref{lem:normal_field_stability} gives
\begin{equation}\label{eq:enf_sup_C1_graph}
 \|DX_{h_k}-D\iota\|_{\infty}\le Ch_k^2,
\end{equation}
where $\iota:\mathcal M\hookrightarrow\mathbb R^D$
is the inclusion. Choose a uniformly controlled local chart
$\psi_x$ with $D\psi_x(0)$ an isometry.
For $\|u\|\le C_{\mathrm{loc}}h_k^2$,
\[
 \|D(X_{h_k}\circ\psi_x)(u)-D\psi_x(0)\|_{\mathrm{op}}
 \le C(h_k^2+\|u\|)\le Ch_k^2.
\]
By Assumption~\ref{ass:estimated-normal-bundle}, local orthonormal frames may be
aligned so that an orthonormal true-normal basis
$E_x^0$ satisfies
\[
 \|E_x-E_x^0\|_{\mathrm{op}}\le Ch_1^2.
\]
This alignment is local and does not require the normal
bundle of $\mathcal M$ to be trivial. The $D\times D$
matrix
\[
 O_x=[\,D\psi_x(0),-E_x^0\,]
\]
is orthogonal. Since
\[
 D_{(u,t)}\mathcal F_{k,x}(u,t)
 =[\,D(X_{h_k}\circ\psi_x)(u),-E_x\,],
\]
the two preceding estimates prove the claimed Jacobian approximation.
The singular-value perturbation inequality applied to
$O_x$ yields the claimed condition-number bound.

Let $\zeta=(u,t)$ and $\zeta^*=(u_k^*,t_k^*)$ be
points of the same convex local search region, with
$\mathcal F_{k,x}(\zeta^*)=0$. Integration of the
Jacobian along their connecting segment gives
\[
 \begin{split}
 \|\mathcal F_{k,x}(\zeta)-\mathcal F_{k,x}(\zeta^*)\|
 &=\left\|O_x(\zeta-\zeta^*)
    +\int_0^1\{D\mathcal F_{k,x}(\zeta^*+s(\zeta-\zeta^*))
              -O_x\}(\zeta-\zeta^*)\,ds\right\|\\
 &\ge(1-\delta_k)\|\zeta-\zeta^*\|.
 \end{split}
\]
This proves the residual-to-solution estimate.

For the empirical statement, replace $X_{h_k}$
by $\widehat X_{h_k}$ and $E_x$ by an aligned empirical
frame $\widehat E_x$. Uniform empirical $C^1$ branch
error $e_{n,k}$ and frame error $s_{n,1}$ add at most
$C(e_{n,k}+s_{n,1})$ to the Jacobian perturbation
in the displayed Jacobian estimate. The same
singular-value and line-segment arguments apply on
any event where the local empirical intersection is
well defined. For STSP and Gaussian density ridges,
the differentiated branch and reference-normal estimates
in Appendices~\ref{app:stsp_finite_sample_proof} and~\ref{app:dr_finite_sample_proof} yield
$e_{n,k}=O(\eta_n(h_k))$ and
$s_{n,1}=O(\eta_n(h_1))$ on their respective
high-probability events.
\end{proof}

\begin{proposition}[Fourth-order closest-point localization
and well-posedness]\label{prop:cp_supplement}
Suppose Assumptions~\ref{ass:geometry-observation}, \ref{ass:normal-section}, \ref{assump:normal_error_t}, and~\ref{ass:tangent-stability} hold.
For $h_1<\cdots<h_K$ and $x\in\mathcal M$, set
\[
 \mathcal Z^{\mathrm{CP}}_{h_k}(x)
 :=\arg\min_{z\in\mathcal M_{h_k}}
      \|z-X_{h_1}(x)\|.
\]
For all sufficiently small $h_{\max}$ this set is
nonempty, and every $Z\in\mathcal Z^{\mathrm{CP}}_{h_k}(x)$
satisfies
\[
 \|\pi_{\mathcal M}(Z)-x\|\le Ch_k^4,
 \qquad \|Z-X_{h_k}(x)\|\le Ch_k^4,
 \qquad
 \|P_{N_x\mathcal M}(Z-X_{h_k}(x))\|\le Ch_k^6.
\]
Any selections $Z^{\mathrm{CP}}_{h_k}(x)$ and weights
satisfying $\sum_k w_k=1$ and
$\sum_k w_kc(h_k,\sigma)=0$ obey
\begin{equation}\label{eq:cp_supp_recovery}
 \sup_{x\in\mathcal M}
 \left\|\sum_{k=1}^K w_k Z^{\mathrm{CP}}_{h_k}(x)-x\right\|
 \le C\sum_{k=1}^K|w_k|h_k^4.
\end{equation}
If, additionally,
\begin{equation}\label{eq:cp_supp_curvature}
 \sup_{h\,\mathrm{admissible}}\|u_h\|_{C^2(\mathcal M)}\le B<\infty,
\end{equation}
then $\inf_h\operatorname{reach}(\mathcal M_h)>0$,
the closest-point selections are unique and $C^1$,
and $\|Z^{\mathrm{CP}}_{h_k}-\iota\|_{C^1}\le Ch_k^2$.
For fixed $K$ and uniformly bounded $\sum_k|w_k|$,
the multiscale image is a $C^1$ embedded manifold
$C^1$-diffeomorphic to $\mathcal M$.
\end{proposition}

\subsubsection*{Proof of Proposition~\ref{prop:cp_supplement}}
\begin{proof}
All constants below are uniform in the admissible scales.
Assumptions \ref{assump:normal_error_t} and \ref{ass:tangent-stability}, together with the argument of
Lemma \ref{lem:normal_field_stability}, imply
\begin{equation}\label{eq:cp-c1-rate}
 \|u_h\|_{C^1(\mathcal M)}\le Ch^2.
\end{equation}
Let $y_1=X_{h_1}(x)$. Since $\mathcal M_{h_k}$ is compact,
its closest-point set $\mathcal Z^{\mathrm{CP}}_{h_k}(x)$
is nonempty. Take any $Z$ in this set and write
$Z=X_{h_k}(z)$ for its unique normal-graph base point $z$.
By minimality,
\[
 \|X_{h_k}(z)-y_1\|
 \le\|X_{h_k}(x)-y_1\|\le Ch_k^2.
\]
The $C^0$ bounds on the normal sections then give
\begin{equation}\label{eq:cp-localization}
 \|z-x\|\le Ch_k^2.
\end{equation}
Set
\[
 d=z-x,\qquad
 \Delta=u_{h_k}(x)-u_{h_1}(x)\in N_x\mathcal M,
 \qquad e=u_{h_k}(z)-u_{h_k}(x).
\]
By \eqref{eq:cp-c1-rate} and $h_1\le h_k$,
\[
 \|\Delta\|\le Ch_k^2,\qquad
 \|e\|\le Ch_k^2\|d\|.
\]
Because $\mathcal M$ has reach at least $\tau$,
its chord--tangent estimate implies
\[
 \|P_{N_x\mathcal M}d\|\le C\|d\|^2.
\]
Since $\Delta$ is normal, it follows that
$|\langle d,\Delta\rangle|\le Ch_k^2\|d\|^2$.
Minimality relative to the competitor $X_{h_k}(x)$ gives
\begin{align*}
 0&\ge\|X_{h_k}(z)-y_1\|^2
     -\|X_{h_k}(x)-y_1\|^2\\
  &=\|d+\Delta+e\|^2-\|\Delta\|^2\\
  &=\|d\|^2+2\langle d,\Delta\rangle
    +2\langle d,e\rangle+2\langle\Delta,e\rangle
    +\|e\|^2\\
  &\ge (1-Ch_k^2)\|d\|^2-Ch_k^4\|d\|.
\end{align*}
For sufficiently small $h_k$, this forces
$\|d\|\le Ch_k^4$. Therefore,
\[
 \|\pi_{\mathcal M}(Z)-x\|=\|z-x\|\le Ch_k^4,
\]
and
\[
 \|Z-X_{h_k}(x)\|
 \le\|z-x\|+\|u_{h_k}(z)-u_{h_k}(x)\|
 \le Ch_k^4.
\]
Moreover, the chord--tangent estimate and the $C^1$ bound give
\[
 \|P_{N_x\mathcal M}\{Z-X_{h_k}(x)\}\|
 \le C\|z-x\|^2+Ch_k^2\|z-x\|
 \le Ch_k^6.
\]
This argument applies to every closest point, whether or not
the minimizer is unique. Writing
$Z^{\mathrm{CP}}_{h_k}(x)=X_{h_k}(x)+e_k(x)$,
$\|e_k\|_\infty\le Ch_k^4$, we obtain
\begin{align*}
 \sum_k w_k Z^{\mathrm{CP}}_{h_k}(x)-x
 &=\sum_k w_k\bigl(u_{h_k}(x)+e_k(x)\bigr)\\
 &=\sum_k w_k\bigl(
     u_{h_k}(x)-c(h_k,\sigma)b(x)+e_k(x)\bigr),
\end{align*}
and the cancellation constraints prove
\eqref{eq:cp_supp_recovery}.

Now assume \eqref{eq:cp_supp_curvature}. From
\eqref{eq:cp-c1-rate}, $X_h$ is uniformly bi-Lipschitz
and $C^1$-close to the inclusion. The uniform $C^2$ bound,
Taylor's theorem in fixed local charts, and the lower
bi-Lipschitz bound give
\[
 \operatorname{dist}\bigl(X_h(z),
  X_h(x)+T_{X_h(x)}\mathcal M_h\bigr)
 \le C\|X_h(z)-X_h(x)\|^2
\]
for nearby $x,z\in\mathcal M$. For separated base points,
compactness gives the same inequality after increasing $C$.
By Federer's chord--tangent characterization,
\begin{equation}\label{eq:cp-reach}
 \inf_{h\text{ small}}\operatorname{reach}(\mathcal M_h)
 \ge\tau_*>0.
\end{equation}
Since $\operatorname{dist}(y_1,\mathcal M_{h_k})\le Ch_k^2$,
all reference points lie in the $\tau_*$-tube for small
$h_{\max}$, and their closest points are therefore unique.
The nearest-point projections onto these $C^2$ manifolds
are $C^1$ on their fixed tubes.

Let $Z^{\mathrm{CP}}_{h_k}(x)$ be this unique projection.
Its normal projection residual has norm $O(h_k^2)$.
The differential of the tubular projection thus differs
from the tangent projector at the projected point by
$O(h_k^2)$, using the uniform bound on the second
fundamental form. Assumption \ref{ass:tangent-stability} and
$\|\pi_{\mathcal M}(Z^{\mathrm{CP}}_{h_k}(x))-x\|
=O(h_k^4)$ imply that the projected tangent projector
differs from $P_{T_x\mathcal M}$ by $O(h_k^2)$.
Composing $D\pi_{\mathcal M_{h_k}}$ with
$DX_{h_1}=D\iota+O(h_1^2)$ therefore gives
\[
 \|D Z^{\mathrm{CP}}_{h_k}-D\iota\|_\infty
 \le Ch_k^2.
\]
Together with the $C^0$ bound, this proves the $C^1$
claim. For fixed $K$ and a uniformly bounded absolute
weight sum, the corrected map is a $C^1$-small perturbation
of $\iota$. Local injectivity and immersion follow from
fixed coordinate charts, while compactness and positional
closeness preserve separation of distant base points.
Thus the map is a $C^1$ embedding and a $C^1$
diffeomorphism onto its image.
\end{proof}

\subsubsection*{Counterexample: closest-point correspondence
need not be unique under the general $C^1$ assumptions}
\label{ex:cp-nonunique}
We give an explicit family of smooth normal graphs satisfying
Assumptions \ref{ass:normal-section}, \ref{assump:normal_error_t} and \ref{ass:tangent-stability}, for which the finest-scale
reference manifold has an exact, uniformly regular normal
field, but closest-point matching fails to be single-valued.
Let $\mathcal M=\mathbb S^1\subset\mathbb R^2$, with
$x(\theta)=(\cos\theta,\sin\theta)$ and outward unit normal
$n(\theta)=x(\theta)$. For $0<s<1/2$, let
\[
 N_s=\lceil s^{-6}\rceil,\qquad
 \beta(s)=\sin^2\!\left(\frac{\pi}{2}\log_2 s\right),\qquad
 a_s=\frac{s^2\beta(s)}{N_s},
\]
and set
\begin{equation}\label{eq:cp-oscillatory-graph}
 X_s(x(\theta))
 =\bigl(1+s^2+a_s\cos(N_s\theta)\bigr)x(\theta).
\end{equation}
Each $\mathcal M_s=X_s(\mathcal M)$ is an embedded smooth
normal graph. Uniformly in $s$,
\[
 \|u_s-s^2 n\|_\infty\le s^2/N_s=O(s^8),
 \qquad \|u_s\|_{C^1}\le Cs^2.
\]
Thus the factorized expansion of Assumption \ref{assump:normal_error_t} holds with
$b=n$ and $c(s,\sigma)=s^2$, and Assumption \ref{ass:tangent-stability} follows
from the $C^1$ bound. However,
$\|u_s\|_{C^2}$ is not uniformly bounded.

Choose $h_n=2^{-2n}$ and the scales $h_1=h_n$,
$h_2=2h_n$. Then $\beta(h_1)=0$,
$\beta(h_2)=1$; in particular, $\mathcal M_{h_1}$ is a
concentric circle, whose geometric normal field equals the
true normal field $n$, and so the estimated-normal-fiber
construction is locally unique and uniformly transverse.
Take the reference point $y_1=(1+h_n^2,0)$, and put
$R=1+h_n^2$, $s=2h_n$, $N=N_s$, $a=s^2/N$, and
$r(\theta)=1+s^2+a\cos(N\theta)$. The squared distance from
$y_1$ to the target point $r(\theta)x(\theta)$ is
\[
 g(\theta)=r(\theta)^2+R^2-2Rr(\theta)\cos\theta.
\]
This is even in $\theta$ and
\[
 g''(0)
 =2\left[R(1+s^2+a)-(s^2-h_n^2+a)s^2N\right]
 =-\Theta(h_n^{-2})<0
\]
for all sufficiently large $n$. Thus $\theta=0$ is a strict
local maximum of the squared-distance function. Since
$g(0)=O(h_n^4)$ while $g(\pi)$ is bounded away from zero,
no global minimizer lies at $0$ or $\pi$. The symmetry
$g(\theta)=g(-\theta)$ then gives at least two distinct
global closest points. Consequently,
$\pi_{\mathcal M_{h_2}}(y_1)$ is not single-valued,
although the reference normal fiber through $y_1$ has a
unique local intersection with $\mathcal M_{h_2}$.

This counterexample concerns the abstract fitted-manifold
class of Section \ref{Multiscale Bias Cancellation}, not population STSP or density-ridge
fitting. It shows that the extra uniform-curvature condition
in Proposition~\ref{prop:cp_supplement} cannot simply
be discarded when asserting global closest-point uniqueness.

\subsubsection{Proof of Lemma~\ref{lem:finite_normal_correspondence}}

\begin{proof}
We first prove the codimension-one case using compatible local unit
normals. The higher-codimension extension uses aligned local frames,
as specified below.

Fix \(x\in\mathcal M\) and \(k\in\{1,\ldots,K\}\). Use the local
chart \(\psi_x\), the map \(\varphi_k=X_{h_k}\circ\psi_x\), and the
population equation \(F_k(q,s;x)=0\) from the proof of
Lemma~\ref{lem:wellposed_correspondence}, and set
\(\widehat\varphi_k=\widehat X_{h_k}\circ\psi_x\).
The population solution \((q_k,s_k)\) satisfies
\(Z_{h_k}(x)=\varphi_k(q_k)=y_1(x)+s_k\widetilde n(x)\).
By Lemma~\ref{lem:wellposed_correspondence}, the derivative
\(D_{(q,s)}F_k\) is uniformly invertible in a neighborhood of
\((q_k,s_k)\). Hence, by the local inverse-function theorem, there
exists a constant \(C_I>0\), uniform in \(x\), such that
\(\|q-q_k\|+|s-s_k| \leq C_I\|F_k(q,s;x)\|\)
for \((q,s)\) in this neighborhood.

Similarly, the empirical correspondence satisfies
\(\widehat\varphi_k(\widehat q_k) - \widehat y_1(x) - \widehat s_k\widehat n(x) = 0\),
with
\(
\widehat Z_{h_k}(x)
=
\widehat\varphi_k(\widehat q_k).
\)
Evaluating the population equation at
\((\widehat q_k,\widehat s_k)\) gives
\begin{align*}
F_k(\widehat q_k,\widehat s_k;x)
&=
\varphi_k(\widehat q_k)
-y_1(x)
-\widehat s_k\widetilde n(x)\\
&=
\bigl[
\varphi_k(\widehat q_k)
-\widehat\varphi_k(\widehat q_k)
\bigr]
+
\bigl[
\widehat y_1(x)-y_1(x)
\bigr]\\
&\quad+
\widehat s_k
\bigl[
\widehat n(x)-\widetilde n(x)
\bigr].
\end{align*}

On \(E_X\),
\[
\sup_q
\|\widehat\varphi_k(q)-\varphi_k(q)\|
\leq
r_n(h_k),
\qquad
\|\widehat y_1(x)-y_1(x)\|
\leq
r_n(h_1).
\]
On \(E_N\),
\(\|\widehat n(x)-\widetilde n(x)\| \leq s_n(h_1)\).
Therefore,
\[
\|F_k(\widehat q_k,\widehat s_k;x)\|
\leq
r_n(h_1)
+
r_n(h_k)
+
|\widehat s_k|s_n(h_1).
\]

Let
\(\Delta_k = \|\widehat q_k-q_k\| + |\widehat s_k-s_k|\).
The local inverse estimate yields
\[
\Delta_k
\leq
C_I
\left[
r_n(h_1)
+
r_n(h_k)
+
|\widehat s_k|s_n(h_1)
\right].
\]
Moreover,
\(
|\widehat s_k|
\leq
|s_k|+\Delta_k.
\)
By the population normal-direction expansion and
Lemma~\ref{lem:normal_guided_correspondence},
\(
|s_k|=O(h_k^2),
\)
where we use \(h_1\leq h_k\). Hence
\[
\Delta_k
\leq
C
\left[
e_{n,k}^{\mathrm{corr}}
+
s_n(h_1)\Delta_k
\right].
\]
Since \(s_n(h_1)\to0\), for sufficiently large \(n\) the last term can
be absorbed into the left-hand side, giving
\(\Delta_k \leq C e_{n,k}^{\mathrm{corr}}\).

Finally,
\begin{align*}
\|\widehat Z_{h_k}(x)-Z_{h_k}(x)\|
&=
\|
\widehat\varphi_k(\widehat q_k)
-
\varphi_k(q_k)
\|\\
&\leq
\|
\widehat\varphi_k(\widehat q_k)
-
\varphi_k(\widehat q_k)
\|
+
\|
\varphi_k(\widehat q_k)
-
\varphi_k(q_k)
\|.
\end{align*}
Since \(\varphi_k\) is \(C^1\) with uniformly bounded derivative in the
local chart,
\[
\|\widehat Z_{h_k}(x)-Z_{h_k}(x)\|
\leq
r_n(h_k)
+
C\|\widehat q_k-q_k\|.
\]
The constants can be chosen uniformly in \(x\) because
\(\mathcal M\) is compact and the local regularity bounds are uniform.
Taking the supremum over \(x\in\mathcal M\) proves
\[
\sup_{x\in\mathcal M}
\left\|
\widehat Z_{h_k}(x)-Z_{h_k}(x)
\right\|
\leq
C
e_{n,k}^{\mathrm{corr}}.
\]

For \(m>1\), use local normal coordinates rather
than choosing a global frame. Choose an orthonormal frame
\(\widetilde E_\alpha(x)\in\mathbb R^{D\times m}\) of
\(\widetilde N(x)\) on each coordinate neighborhood. For this position estimate,
these frames are needed only pointwise; \(C^1\) local frames are
used only when \(C^1\) regularity is invoked.
On the same high-probability
event, for sufficiently large \(n\), choose the empirical frame by
\[
\widehat E_\alpha(x)
:=\widehat Q(x)\widetilde E_\alpha(x)
\left[
\widetilde E_\alpha(x)^\top\widehat Q(x)\widetilde E_\alpha(x)
\right]^{-1/2}.
\]
The matrix in brackets is uniformly positive definite because the
projector gap is at most \(s_n(h_1)=o(1)\). This construction gives an
orthonormal frame of \(\widehat N(x)\) satisfying
\[
\|\widehat E_\alpha(x)-\widetilde E_\alpha(x)\|_{\mathrm{op}}
\leq C\|P_{\widehat N(x)}-P_{\widetilde N(x)}\|_{\mathrm{op}}
\leq Cs_n(h_1).
\]
Replacing scalar normal coordinates by vectors in \(\mathbb R^m\)
and the two unit normals by these aligned local frames yields exactly
the same residual and absorption estimates. Changes of local
orthonormal frame change only the fiber coordinates, not the ambient
intersection points or their error bounds. Thus the local estimates
apply globally without assuming triviality of either normal bundle.
\end{proof}

\subsubsection{Proof of Theorem~\ref{thm:finite_normal_guided_decomposition}}

\begin{proof}
We work on the event \(E_X\cap E_N\). For every \(x\in\mathcal M\),
\begin{align*}
\widehat X(x)-x
&=
\widehat X(x)
-
\widetilde X(x)
+
\widetilde X(x)-x\\
&=
\sum_{k=1}^{K}
w_k
\left[
\widehat Z_{h_k}(x)-Z_{h_k}(x)
\right]
+
\widetilde X(x)-x.
\end{align*}
Therefore, by the triangle inequality,
\begin{align*}
\|\widehat X(x)-x\|
&\leq
\sum_{k=1}^{K}|w_k|
\left\|
\widehat Z_{h_k}(x)-Z_{h_k}(x)
\right\|\\
&\quad+
\|\widetilde X(x)-x\|.
\end{align*}
Taking the supremum over \(x\in\mathcal M\) gives
\begin{equation}
\label{eq:finite-ng-decomposition}
\sup_{x\in\mathcal M}\|\widehat X(x)-x\|
\leq\sum_{k=1}^K|w_k|\sup_{x\in\mathcal M}
\|\widehat Z_{h_k}(x)-Z_{h_k}(x)\|
+\sup_{x\in\mathcal M}\|\widetilde X(x)-x\|.
\end{equation}

By Lemma~\ref{lem:finite_normal_correspondence},
\(\sup_{x\in\mathcal M} \left\| \widehat Z_{h_k}(x)-Z_{h_k}(x) \right\| \leq C e_{n,k}^{\mathrm{corr}}\).
Therefore,
\begin{align*}
&\sum_{k=1}^{K}|w_k|
\sup_{x\in\mathcal M}
\left\|
\widehat Z_{h_k}(x)-Z_{h_k}(x)
\right\|\\
&\qquad\leq
C\sum_{k=1}^{K}|w_k|
e_{n,k}^{\mathrm{corr}}.
\end{align*}

On the other hand, by
Proposition~\ref{thm:normal_guided_manifold_recovery},
the population GBC estimator satisfies
\[\sup_{x\in\mathcal M} \|\widetilde X(x)-x\| \leq C \sum_{k=1}^{K}|w_k|h_k^4\].
Substituting the preceding two bounds into
\eqref{eq:finite-ng-decomposition} yields
\[
\sup_{x\in\mathcal M}
\left\|
\widehat X(x)-x
\right\|
\leq
C\sum_{k=1}^{K}|w_k|h_k^4
+
C\sum_{k=1}^{K}|w_k|
e_{n,k}^{\mathrm{corr}}.
\]
Since this bound holds on \(E_X\cap E_N\), it holds with probability
at least
\(
1-\delta_n^X-\delta_n^N.
\)

Finally, since \(\widehat{\mathcal M}=\widehat X(\mathcal M)\),
\eqref{eq:app-general-image-distance} gives
\(d_H(\widehat{\mathcal M},\mathcal M)
\leq\varepsilon_n:=\sup_{x\in\mathcal M}\|\widehat X(x)-x\|\).
Combining this inequality with the uniform bound above gives
\[
d_H\left(
\widehat{\mathcal M},
\mathcal M
\right)
\leq
C\sum_{k=1}^{K}|w_k|h_k^4
+
C\sum_{k=1}^{K}|w_k|
\left[
r_n(h_1)+r_n(h_k)+s_n(h_1)h_k^2
\right].
\]
This completes the proof.
\end{proof}

\subsubsection*{Proof of Proposition~\ref{prop:noise_level_perturbation}}
\begin{proof}
Work on \(E_X\cap E_N\), whenever the local empirical correspondences
exist. Put \(c_k=c(h_k,\sigma)\) and
\(e_{n,k}^{\mathrm{corr}}=r_n(h_1)+r_n(h_k)+s_n(h_1)h_k^2\).
Writing \(x_k=\pi_{\mathcal M}(Z_{h_k}(x))\), define
\(R_{h_k}(x)=u_{h_k}(x)-c_kb(x)\). Then
\[
\begin{aligned}
Z_{h_k}(x)-x
={}&c_kb(x)+(x_k-x)\\
&+\bigl(u_{h_k}(x_k)-u_{h_k}(x)\bigr)+R_{h_k}(x).
\end{aligned}
\]
Here \(\|R_{h_k}\|_\infty\leq Ch_k^4\).
By Lemmas~\ref{lem:normal_guided_correspondence} and
\ref{lem:normal_field_stability}, the last three terms are bounded by
\(C(h_1^2h_k^2+h_1^2h_k^4+h_k^4)\leq Ch_k^4\).
Consequently,
\[
Z_{h_k}(x)-x=c_kb(x)+\mathcal E_k(x),\qquad
\|\mathcal E_k\|_\infty\leq Ch_k^4.
\]
This expansion does not require a cancellation constraint on the weights.
The empirical correspondence bound gives
\(\|\widehat Z_{h_k}-Z_{h_k}\|_\infty\leq Ce_{n,k}^{\mathrm{corr}}\).
Since the plug-in weights sum to one and cancel \(c(h_k,\widehat\sigma)\),
\[
\begin{aligned}
\sum_k\widehat w_kc_k&=\Delta_{n,c},\\
\widehat X_{\widehat\sigma}(x)-x
&=\Delta_{n,c}b(x)
+\sum_k\widehat w_k\{
\mathcal E_k(x)+\widehat Z_{h_k}(x)-Z_{h_k}(x)\}.
\end{aligned}
\]
Taking the uniform norm proves the positional bound, and
\eqref{eq:app-general-image-distance} gives the Hausdorff bound.
The mean value theorem on \(I_n\) gives, on \(E_\sigma\),
\[
|\Delta_{n,c}|
\leq\sum_k|\widehat w_k|L_{n,k}|\widehat\sigma-\sigma|
\leq\rho_n\sum_k|\widehat w_k|L_{n,k}.
\]
All inequalities hold for each realization of the weights, so data dependence
requires no independence assumption. The stated probability follows by the
union bound for \(E_X^c,E_N^c,E_\sigma^c\). The same calculation for
the population maps omits the \(e_{n,k}^{\mathrm{corr}}\) terms.
\end{proof}

\FloatBarrier
\section{Successive Tangent-Space Projection}
\label{sec:app_stsp}
The STSP calculations below use the fixed kernel
\(K(t)=(1-t)_+^6\) used in Section \ref{sec:stsp_validation}. For the STSP and density-ridge
proofs, write
\(P_x=P_{T_x\mathcal M}, Q_x=P_{N_x\mathcal M}\).
In the STSP proofs, let
\(\lambda_h=\sigma/h\in\Lambda_{\mathrm{STSP}}\Subset(0,\infty)\).
The noise level \(\sigma\) is common to the bandwidths in each multiscale
comparison; at a fixed bandwidth we abbreviate \(\lambda=\lambda_h\).
For the population analysis, setting \(\sigma=\lambda h\) defines the
two-parameter family
\[
\mathcal F_h^{\mathrm{STSP}}(x,\beta;\lambda)
=
h^{-2}Q_xR_{h}^{\mathrm{STSP}}(x+h^2\beta),
\qquad \beta\in N_x\mathcal M,
\]
where the residual is evaluated at noise level \(\lambda h\).
We write \(\mathcal F_h^{\mathrm{STSP}}(x,\beta)\) for this family at
\(\lambda=\lambda_h\); uniformity in \(\lambda\) permits this evaluation
without assuming a common noise-to-bandwidth ratio across scales.
All population estimates are uniform in \(x\in\mathcal M\),
\(\lambda\in\Lambda_{\mathrm{STSP}}\), and bounded \(\beta\).
The notation \(O_{C^2}\) refers to derivatives in local normal-bundle
coordinates \((x,\beta)\), with the scale parameters fixed.

\begin{proof}[Local density expansion]
Use a finite collection of reach-controlled graph charts,
\[
 \Psi_x(t)
 =x+t+\tfrac12\mathrm{II}_x(t,t)
      +\tfrac16\mathcal C_x[t,t,t]+\mathcal H_x(t).
\]
Taylor's formula and $C^6$ regularity give, for $|\alpha|\le2$,
\(\|D_x^\alpha\mathcal H_x(t)\|\le C\|t\|^4\),
and
\(J_x(t)=1+J_{2,x}(t)+\mathcal J_x(t), |D_x^\alpha\mathcal J_x(t)|\le C\|t\|^3\),
where $J_{2,x}$ is quadratic. Derivatives are taken in fixed
local base coordinates, with local frame coordinates for $t$
held fixed. The constants are uniform after reducing the
chart radii.

Substitute $t=hw$ into the Gaussian convolution.
Contributions outside the fixed local chart are exponentially
small, including the required derivatives. In the local
integral, write $w=v_T+\lambda g$ and integrate against a
standard Gaussian in $g$. Taylor expansion of the normal
Gaussian factor to second order in $h$, together with the
graph and volume bounds, leaves an $O_{C^2}(h^3)$ integrated
remainder.

To justify this uniformly, first restrict to
$\|g\|\le h^{-1/8}$. On this region the remainder and its first
two parameter derivatives are bounded by $Ch^3$ times an
integrable polynomial Gaussian envelope. The complementary
Gaussian tail is smaller than every fixed power of $h$.
Compactness of $\Lambda$ makes all constants uniform.

Gaussian averaging gives
\[
 \mathbb E\{
 \mathrm{II}_x(v_T+\lambda g,v_T+\lambda g)
 \}
 =
 \mathrm{II}_x(v_T,v_T)+d\lambda^2H(x),
\]
which yields $L$.

The second-order terms consist of the quadratic likelihood
correction, the quadratic volume term, and the cubic graph
term contracted with $v_N$. Each is even under
\[
 (v_T,v_N,g)\mapsto(-v_T,-v_N,-g).
\]
This proves the parity of $M$. Differentiation in $(x,\beta)$
preserves these algebraic symmetries and the remainder bound.
\end{proof}
\subsection{Proof of Lemma~\ref{lem:stsp_scaled_normal_equation}}
\begin{proof}
Fix \(x\in\mathcal M\), and let
\(y_h=x+h^2\beta, \beta\in N_x\mathcal M, \|\beta\|\le R\),
for a fixed \(R<\infty\). The assumed geometric regularity and Gaussian
convolution allow all local graph, Jacobian, and density expansions below
to be differentiated twice under the integral sign, uniformly over the
parameter range specified above.

Choose orthonormal coordinates adapted to
\(
\mathbb R^D=T_x\mathcal M\oplus N_x\mathcal M,
\)
and write \(v=v_T+v_N\) accordingly. Under the sign convention used in
this section, a local parametrization of \(\mathcal M\) around \(x\) is
\(\Psi_x(u) = x+u+\frac12 \mathrm{II}_x(u,u)+O(\|u\|^3)\).
Let \(p_\sigma\) denote the density of the noisy observation
\(Y=X+\xi\). After introducing
\(Y=y_h+hv, u=ht,  \sigma=\lambda h\),
the Gaussian convolution admits the expansion
\begin{equation}
\label{eq:stsp-local-density-expansion}
p_\sigma(y_h+hv)
=
C_{h,x,\lambda}\,
\rho_\lambda(v_N)
\left[
1+hL_{x,\beta,\lambda}(v)
+h^2M_{x,\beta,\lambda}(v)
+O_{C^2}(h^3)
\right],
\end{equation}
where
\(\rho_\lambda(v_N) = \exp\!\left( -\frac{\|v_N\|^2}{2\lambda^2} \right)\),
and the first-order term is
\begin{equation}
\label{eq:stsp-L1}
L_{x,\beta,\lambda}(v)
=
-\frac1{\lambda^2}
\left\langle
v_N,\,
\beta
-\frac12 \mathrm{II}_x(v_T,v_T)
-\frac{d\lambda^2}{2}H(x)
\right\rangle .
\end{equation}
The last term in \eqref{eq:stsp-L1} is produced by Gaussian averaging
of the quadratic graph term, rather than by the Jacobian correction.
Indeed, if \(Z\sim N(0,I_d)\), then
\[
\mathbb E
\left[
\mathrm{II}_x(v_T+\lambda Z,v_T+\lambda Z)
\right]
=
\mathrm{II}_x(v_T,v_T)
+
\lambda^2\sum_{j=1}^d \mathrm{II}_x(e_j,e_j)
=
\mathrm{II}_x(v_T,v_T)+d\lambda^2H(x).
\]
In particular,
\[
L_{x,\beta,\lambda}(-v_T,v_N)
=
L_{x,\beta,\lambda}(v_T,v_N),
\qquad
L_{x,\beta,\lambda}(v_T,-v_N)
=
-
L_{x,\beta,\lambda}(v_T,v_N).
\]
After collecting the second-order terms in the Gaussian exponent, the
local graph, and the Jacobian, one also obtains the simultaneous-reflection
parity
\begin{equation}
\label{eq:stsp-M-parity}
M_{x,\beta,\lambda}(-v_T,-v_N)
=
M_{x,\beta,\lambda}(v_T,v_N).
\end{equation}
The terms quadratic in the first-order normal displacement and the
second-order Jacobian corrections are even under this reflection. The only
second-order terms that are odd in \(v_N\) arise from the cubic normal graph;
after Gaussian averaging they are also odd in \(v_T\), and hence are even
under the simultaneous reflection. These parity identities are preserved
after taking derivatives with respect to \((x,\beta)\) up to order two.

Introduce the reference probability measure
\[
d\nu_{\lambda}(v)
=
\frac{
K(\|v\|^2)
\rho_\lambda(v_N)\,dv
}{
Z_{\lambda}
},
\qquad
Z_{\lambda}
=
\int_{\|v\|\le1}
K(\|v\|^2)\rho_\lambda(v_N)\,dv,
\]
and denote expectation with respect to \(\nu_{\lambda}\) by
\(\mathbb E_{\lambda}\). Here \(\mathbb E_{h}\) denotes expectation under
the rescaled law with density proportional to
\(K(\|v\|^2)p_\sigma(y_h+hv)\). The fixed radial weight preserves all
reflection symmetries of the density expansion.
The measure \(\nu_{\lambda}\) is invariant under
\(v\mapsto-v\). Therefore
\[
\mathbb E_{\lambda}[v]=0,
\qquad
\mathbb E_{\lambda}[L_{x,\beta,\lambda}(v)]=0,
\qquad
\mathbb E_{\lambda}
\left[
vM_{x,\beta,\lambda}(v)
\right]=0.
\]
More explicitly, for any integrable test function \(g\), normalization of
\eqref{eq:stsp-local-density-expansion} gives
\[
\mathbb E_{h}[g(v)]
=
\frac{
\mathbb E_{\lambda}
\left[
g(v)
\left\{
1+hL_{x,\beta,\lambda}(v)
+h^2M_{x,\beta,\lambda}(v)
+O_{C^2}(h^3)
\right\}
\right]
}{
\mathbb E_{\lambda}
\left[
1+hL_{x,\beta,\lambda}(v)
+h^2M_{x,\beta,\lambda}(v)
+O_{C^2}(h^3)
\right]
}.
\]
Taking \(g(v)=v\) and using the preceding cancellations yields the refined
scaled-mean expansion
\begin{equation}
\label{eq:stsp-scaled-mean}
\mathbb E_{h}[v]
=
h\,\mathbb E_{\lambda}
\left[
vL_{x,\beta,\lambda}(v)
\right]
+
O_{C^2}(h^3).
\end{equation}

We first compute its normal component. By rotational invariance of
\(\nu_{\lambda}\) under \(O(d)\times O(D-d)\),
\(\mathbb E_{\lambda}[v_Nv_N^\top] = a_{\lambda} I_{N_x\mathcal M}\),
where \(a_{\lambda}>0\). Moreover, by the definition of \(c_{\lambda}\),
\[
\mathbb E_{\lambda}
\left[
v_N
\left\langle
v_N,\mathrm{II}_x(v_T,v_T)
\right\rangle
\right]
=
c_{\lambda}
\sum_{j=1}^d \mathrm{II}_x(e_j,e_j)
=
d\,c_{\lambda}H(x).
\]
Substituting \eqref{eq:stsp-L1}, we therefore obtain
\begin{align}
Q_x\mathbb E_{\lambda}
\left[
vL_{x,\beta,\lambda}(v)
\right]
&=
-\frac1{\lambda^2}
\mathbb E_{\lambda}
\left[
v_N
\left\langle
v_N,\beta
-\frac12 \mathrm{II}_x(v_T,v_T)
-\frac{d\lambda^2}{2}H(x)
\right\rangle
\right]
\nonumber\\
&=
-\frac1{\lambda^2}
\left[
a_{\lambda}\beta
-
\frac d2 c_{\lambda}H(x)
-
\frac{d\lambda^2}{2}a_{\lambda}H(x)
\right]
\nonumber\\
&=
-\frac1{\lambda^2}
\left[
a_{\lambda}\beta
-
\frac d2
\left(
c_{\lambda}+\lambda^2a_{\lambda}
\right)
H(x)
\right].
\label{eq:stsp-normal-moment-computation}
\end{align}

Since \(\mu_{h}(y_h)-y_h=h\,\mathbb E_{h}[v]\), combining \eqref{eq:stsp-scaled-mean} and
\eqref{eq:stsp-normal-moment-computation} gives
\begin{equation}
\label{eq:stsp-local-normal-mean}
Q_x\{\mu_{h}(y_h)-y_h\}
=
-\frac{h^2}{\lambda^2}
\left[
a_{\lambda}\beta
-
\frac d2
\left(
c_{\lambda}+\lambda^2a_{\lambda}
\right)
H(x)
\right]
+
O_{C^2}(h^4).
\end{equation}

For the tangential component, the first-order contribution vanishes
because
\[
\mathbb E_{\lambda}
\left[
v_TL_{x,\beta,\lambda}(v)
\right]
=0,
\]
the integrand being odd under \(v_T\mapsto-v_T\). The second-order
contribution vanishes by \eqref{eq:stsp-M-parity}, since
\(\mathbb E_{\lambda}[vM_{x,\beta,\lambda}(v)]=0\). Consequently,
\begin{equation}
\label{eq:stsp-local-tangent-mean}
P_x\{\mu_{h}(y_h)-y_h\}
=
O_{C^2}(h^4).
\end{equation}

Set \(v_h:=y_h-\mu_{h}(y_h)\). Equations \eqref{eq:stsp-local-normal-mean} and
\eqref{eq:stsp-local-tangent-mean} imply
\begin{equation}
\label{eq:stsp-normal-v}
Q_xv_h
=
\frac{h^2}{\lambda^2}
\left[
a_{\lambda}\beta
-
\frac d2
\left(
c_{\lambda}+\lambda^2a_{\lambda}
\right)
H(x)
\right]
+
O_{C^2}(h^4),
\end{equation}
and \(P_xv_h=O_{C^2}(h^4)\). In particular,
\(\|v_h\|_{C^2}=O(h^2)\).

Let \(Q_{h}(y_h)=I-P_{h}(y_h)\) be the population PCA normal projector
entering the STSP operator. The same
\(C^2\)-uniform local expansion applied to the centered second moments gives
\[
\Sigma_{h}(y_h)
=
h^2\bigl(\tau_{\lambda} P_x+a_{\lambda} Q_x\bigr)
+
O_{C^2}(h^4).
\]
For completeness, the eigengap must be checked for the new kernel.
Write \(q=D-d\) and \(S=\|v_N\|^2\). Conditional on \(v_N\), the tangent
density is proportional to \((1-S-\|v_T\|^2)_+^6\), so radial integration
gives
\[
\tau_{\lambda}=\frac{\mathbb E_{\lambda}[1-S]}{d+14},\qquad
 a_{\lambda}=\frac{\mathbb E_{\lambda}[S]}q.
\]
For the radial density proportional to \((1-\|v\|^2)_+^6\), corresponding
to \(\lambda=\infty\), beta integration gives
\(\mathbb E_{\infty}[S]=q/(D+14)\). The finite-\(\lambda\) law is its
tilting by the strictly decreasing function
\(S\mapsto\exp(-S/(2\lambda^2))\). Since \(S\) is nonconstant, this strictly
decreases its expectation. Therefore
\(\mathbb E_{\lambda}[S]<\frac q{D+14}, \tau_{\lambda}>a_{\lambda}>0\).
Continuity and compactness imply
\(g:=\inf_{\lambda\in\Lambda_{\mathrm{STSP}}} (\tau_{\lambda}-a_{\lambda})>0\).
Thus Weyl's inequality gives a population eigengap at least
\(gh^2/2\) for sufficiently small \(h\). Smooth spectral-projector
calculus applied to the normalized covariance then yields
\(Q_{h}(y_h)=Q_x+E_h(x,\beta), \|E_h\|_{C^2}=O(h^2)\).
Since \(R_{h}^{\mathrm{STSP}}(y_h)=Q_{h}(y_h)v_h\), we obtain
\(Q_xR_{h}^{\mathrm{STSP}}(y_h) = Q_xv_h+Q_xE_hv_h\).
The preceding bounds give \(Q_xE_hv_h=O_{C^2}(h^4)\).
Together with \eqref{eq:stsp-normal-v}, this yields
\begin{equation}
\label{eq:stsp-normal-residual-expansion}
Q_xR_{h}^{\mathrm{STSP}}(x+h^2\beta)
=
\frac{h^2}{\lambda^2}
\left[
a_{\lambda}\beta
-
\frac d2
\left(
c_{\lambda}+\lambda^2a_{\lambda}
\right)
H(x)
\right]
+
O_{C^2}(h^4).
\end{equation}

Dividing \eqref{eq:stsp-normal-residual-expansion} by \(h^2\) gives
\[
\mathcal F_h^{\mathrm{STSP}}(x,\beta)
=
\mathcal F_0^{\mathrm{STSP}}(x,\beta;\lambda_h)
+
O_{C^2}(h^2),
\]
where, for \(\lambda\in\Lambda_{\mathrm{STSP}}\), we define
\[
\mathcal F_0^{\mathrm{STSP}}(x,\beta;\lambda)
=
\frac1{\lambda^2}
\left[
a_{\lambda}\beta
-
\frac d2
\left(
c_{\lambda}+\lambda^2a_{\lambda}
\right)
H(x)
\right].
\]
Since \(a_{\lambda}>0\), the limiting equation is affine and
nondegenerate in \(\beta\), with the unique zero
\[
\beta_0(x,\lambda)
=
\frac d2
\left(
\lambda^2+\frac{c_{\lambda}}{a_{\lambda}}
\right)
H(x).
\]
Furthermore,
\[
D_\beta
\mathcal F_0^{\mathrm{STSP}}(x,\beta;\lambda)
=
\frac{a_{\lambda}}{\lambda^2}
I_{N_x\mathcal M}.
\]
Since \(a_{\lambda}\) is continuous and strictly positive and
\(\Lambda_{\mathrm{STSP}}\Subset(0,\infty)\) is compact,
\(\gamma_0 := \inf_{\lambda\in\Lambda_{\mathrm{STSP}}} \frac{a_{\lambda}}{\lambda^2} >0\).
Hence
\[
\inf_{\substack{x\in\mathcal M\\
\lambda\in\Lambda_{\mathrm{STSP}}}}
\sigma_{\min}
\left[
D_\beta
\mathcal F_0^{\mathrm{STSP}}
\bigl(x,\beta_0(x,\lambda);\lambda\bigr)
\right]
\ge\gamma_0>0.
\]
This proves the claimed uniform normal nondegeneracy.
\end{proof}

\subsection{Proof of Lemma~\ref{lem:stsp_cubic_cancellation}}
\begin{proof}
Let
\(y_h(x)=x+h^2\beta_h(x)\)
be a \(C^2\) population STSP fixed-point section with uniformly bounded
\(C^2\) norm. Substituting its fixed-point identity
\(R_{h}^{\mathrm{STSP}}(y_h)=0\) into
\eqref{eq:stsp-normal-residual-expansion} gives, uniformly in the
admissible ratio \(\lambda_h\),
\[
0
=
\mathcal F_0^{\mathrm{STSP}}
\bigl(x,\beta_h(x);\lambda_h\bigr)
+
O_{C^2}(h^2).
\]
By definition,
\[
\mathcal F_0^{\mathrm{STSP}}
\bigl(x,\beta;\lambda_h\bigr)
=
\frac{a_{\lambda_h}}{\lambda_h^2}
\left[
\beta-\beta_0(x,\lambda_h)
\right],
\]
with \(\beta_0\) as in the preceding proof. The same lower bound
\(\gamma_0>0\) therefore gives
\[
\left\|
\beta_h-\beta_0(\cdot,\lambda_h)
\right\|_{C^2}
\leq
\gamma_0^{-1}O(h^2)
=
O(h^2),
\]
uniformly over all admissible pairs \((h,\sigma)\). Multiplication by
\(h^2\) yields
\[
y_h(x)
=
x+h^2\beta_0(x,\lambda_h)+O_{C^2}(h^4),
\]
which proves the claim.
\end{proof}

\subsection{Proof of Proposition~\ref{thm:stsp_bias_expansion}}
\begin{proof}
For the two-parameter family defined at the start of this section,
Lemma~\ref{lem:stsp_scaled_normal_equation} gives
\[
\mathcal F_h^{\mathrm{STSP}}(x,\beta;\lambda)
=
\mathcal F_0^{\mathrm{STSP}}(x,\beta;\lambda)
+
O_{C^2}(h^2)
\]
uniformly in \(x\in\mathcal M\), \(\lambda\in\Lambda_{\mathrm{STSP}}\),
and bounded \(\beta\). The reference map, its unique zero
\(\beta_0(x,\lambda)\), and the lower bound \(\gamma_0>0\) are those
established in the proof of that lemma. In particular,
\(D_\beta\mathcal F_0^{\mathrm{STSP}}\) is uniformly invertible.

By compactness of \(\mathcal M\) and
\(\Lambda_{\mathrm{STSP}}\), we may choose \(R<\infty\) such that
\[
\sup_{\substack{x\in\mathcal M\\
\lambda\in\Lambda_{\mathrm{STSP}}}}
\|\beta_0(x,\lambda)\|
<
\frac R2.
\]
The uniform \(C^2\) convergence and invertibility allow the use of the
parameterized implicit-function theorem. Thus there exists \(h_0>0\) such that, for every
\(0<h<h_0\) and \(\lambda\in\Lambda_{\mathrm{STSP}}\), there exists a
\(C^2\) section
\(\beta_h(\,\cdot\,;\lambda) \in\Gamma(N\mathcal M)\)
in a fixed neighborhood of \(\beta_0(\cdot,\lambda)\) satisfying
\(\mathcal F_h^{\mathrm{STSP}} \bigl(x,\beta_h(x;\lambda);\lambda\bigr) =0\).
Moreover,
\(\|\beta_h(\cdot;\lambda)-\beta_0(\cdot,\lambda)\|_{C^2} = O(h^2)\)
uniformly in \(\lambda\in\Lambda_{\mathrm{STSP}}\).

To exclude additional zeros in \(\|\beta\|\leq R\), use the same
\(\gamma_0\). Uniform \(C^1\) convergence allows us to decrease \(h_0\)
once more so that
\[
\sup_{\substack{x\in\mathcal M,\ \lambda\in\Lambda_{\mathrm{STSP}}\\
\|\beta\|\leq R}}
\left\|
D_\beta\mathcal F_h^{\mathrm{STSP}}(x,\beta;\lambda)
-
\frac{a_{\lambda}}{\lambda^2}I_{N_x\mathcal M}
\right\|_{\mathrm{op}}
\leq
\frac{\gamma_0}{2}.
\]
For any \(\beta_1,\beta_2\in N_x\mathcal M\) with
\(\|\beta_i\|\leq R\), integration along the segment joining them gives
\begin{align*}
&\left\langle
\mathcal F_h^{\mathrm{STSP}}(x,\beta_1;\lambda)
-
\mathcal F_h^{\mathrm{STSP}}(x,\beta_2;\lambda),
\beta_1-\beta_2
\right\rangle
\\
&\qquad\geq
\frac{\gamma_0}{2}
\|\beta_1-\beta_2\|^2.
\end{align*}
Hence \(\mathcal F_h^{\mathrm{STSP}}(x,\cdot;\lambda)\) is injective on
that ball. The implicit-function zero is therefore the unique zero in
\(\|\beta\|\leq R\).

Define
\[
u_h^{\mathrm{STSP}}(x;\lambda)
=
h^2\beta_h(x;\lambda),
\qquad
X_h^{\mathrm{STSP}}(x;\lambda)
=
x+u_h^{\mathrm{STSP}}(x;\lambda).
\]
Then
\(\|u_h^{\mathrm{STSP}}\|_{C^2}=O(h^2)\).

We next verify that solving the projected normal equation is equivalent
to solving the full STSP fixed-point equation. Let
\(y=x+h^2\beta, w=R_{h}^{\mathrm{STSP}}(y)\).
By definition of the STSP residual, \(w=Q_{h}(y)w\), where \(Q_{h}(y)=I-P_{h}(y)\) is the population PCA normal projector.
The population PCA consistency gives
\(\|Q_{h}(y)-Q_x\|_{\mathrm{op}}=O(h^2)\)
uniformly in the region under consideration.

If \(\mathcal F_h^{\mathrm{STSP}}(x,\beta;\lambda)=0\), then
\(Q_xw=0\). Since \(w=Q_{h}(y)w\),
\(w = \{Q_{h}(y)-Q_x\}w\).
For sufficiently small \(h\),
\(\|Q_{h}(y)-Q_x\|_{\mathrm{op}}<\frac12\),
so the preceding identity implies \(w=0\). Hence
\[
\mathcal F_h^{\mathrm{STSP}}(x,\beta;\lambda)=0
\quad\Longleftrightarrow\quad
R_{h}^{\mathrm{STSP}}(x+h^2\beta)=0
\]
throughout the \(Rh^2\) tubular neighborhood.

Consequently, all population STSP fixed points in
\[
\left\{
x+v:
x\in\mathcal M,\;
v\in N_x\mathcal M,\;
\|v\|\leq Rh^2
\right\}
\]
are exactly
\[
\mathcal M_h^{\mathrm{STSP}}
=
\left\{
x+h^2\beta_h(x;\lambda):
x\in\mathcal M
\right\}.
\]
Thus the local fixed-point set is a unique \(C^2\) normal graph.

Uniformity in \(\lambda\) now permits evaluation at
\(\lambda=\lambda_h\) in the common-noise multiscale setting. Defining
\(u_h^{\mathrm{STSP}}(x) := u_h^{\mathrm{STSP}}(x;\lambda_h)\),
Lemma~\ref{lem:stsp_cubic_cancellation}, equivalently the uniform
\(O_{C^2}(h^2)\) estimate for \(\beta_h-\beta_0\), gives
\[
u_h^{\mathrm{STSP}}(x)
=
h^2\beta_0(x,\lambda_h)
+
O_{C^2}(h^4).
\]
Substituting the expression for \(\beta_0\) yields
\[
u_h^{\mathrm{STSP}}(x)
=
\frac d2h^2
\left(
\lambda_h^2+
\frac{c_{\lambda_h}}{a_{\lambda_h}}
\right)
H(x)
+
O_{C^2}(h^4)
=
c_{\mathrm{STSP}}(h,\sigma)H(x)
+
O_{C^2}(h^4),
\]
where
\[
c_{\mathrm{STSP}}(h,\sigma)
=
\frac d2h^2
\left[
\left(\frac{\sigma}{h}\right)^2
+
\frac{c_{\sigma/h}}{a_{\sigma/h}}
\right].
\]

Equivalently, defining
\(\chi_{\mathrm{STSP}}(\lambda) = 1+ \frac{c_{\lambda}}{\lambda^2a_{\lambda}}\)
gives the equivalent representation
\[
c_{\mathrm{STSP}}(h,\sigma)
=
\frac d2
\sigma^2
\chi_{\mathrm{STSP}}
\left(
\frac{\sigma}{h}
\right).
\]

Because \(a_{\lambda}\), \(c_{\lambda}\), and hence
\(\chi_{\mathrm{STSP}}(\lambda)\), are uniformly bounded on the compact
set \(\Lambda_{\mathrm{STSP}}\), while
\(\sigma/h\in\Lambda_{\mathrm{STSP}}\), we have \(\sigma\asymp h\) and therefore
\(c_{\mathrm{STSP}}(h,\sigma) = O(h^2) = O(h^2+\sigma^2)\).

Finally, consider the map
\(X_h^{\mathrm{STSP}}(x)=x+u_h^{\mathrm{STSP}}(x)\). Since \(u_h^{\mathrm{STSP}}(x)\in N_x\mathcal M\) and
\(\|u_h^{\mathrm{STSP}}\|_{C^0}=O(h^2)<\tau\)
for sufficiently small \(h\), the tubular-neighborhood theorem gives
\(\pi_{\mathcal M}\circ X_h^{\mathrm{STSP}} = \operatorname{id}_{\mathcal M}\).
Hence \(X_h^{\mathrm{STSP}}\) is injective. Differentiating the above identity gives
\[
D\pi_{\mathcal M}(X_h^{\mathrm{STSP}}(x))
\circ
DX_h^{\mathrm{STSP}}(x)
=
I_{T_x\mathcal M},
\]
so \(DX_h^{\mathrm{STSP}}(x)\) has rank \(d\); thus \(X_h^{\mathrm{STSP}}\) is an immersion.
Since \(\mathcal M\) is compact, a continuous injective immersion into
the Hausdorff space \(\mathbb R^D\) is a \(C^2\) embedding onto its
image. Therefore
\(\mathcal M_h^{\mathrm{STSP}} = X_h^{\mathrm{STSP}}(\mathcal M)\)
is a \(d\)-dimensional \(C^2\) embedded manifold, and
\(\pi_{\mathcal M} \big|_{\mathcal M_h^{\mathrm{STSP}}} = \bigl(X_h^{\mathrm{STSP}}\bigr)^{-1}\)
is a \(C^2\) diffeomorphism.

This completes the proof.
\end{proof}

\begin{proof}[Proof of Proposition \ref{prop:stsp-naive}]
Write
\[
r(\lambda)=\frac{c_\lambda}{a_\lambda}.
\]
By Proposition \ref{thm:stsp_bias_expansion}, uniformly over $x\in\mathcal M$,
\[
X_h^{\mathrm{STSP}}(x)-x
=
\frac d2\left\{
\sigma^2+h^2r(\sigma/h)
\right\}H(x)+O_{C^2}(h^4).
\]

We first consider exact-normal correspondence.
Since $h_1=3\sigma$ and $h_2=6\sigma$,
\begin{align*}
\widetilde X_{\mathrm{sq}}(x)-x
&=\frac43 u_{h_1}(x)-\frac13 u_{h_2}(x)\\
&=\frac d2\sigma^2
\left\{
1+12r(1/3)-12r(1/6)
\right\}H(x)
+O(\sigma^4).
\end{align*}

We next establish strict positivity of the leading
coefficient. Recall that
\[
a_\lambda
=\mathbb E_\lambda[v_{N,1}^2],
\qquad
c_\lambda
=\mathbb E_\lambda[v_{N,1}^2v_{T,1}^2],
\]
where $\mathbb E_\lambda$ denotes expectation under
the kernel-weighted measure $\nu_\lambda$.
Put $S=\|v_N\|^2$. Conditional on $v_N$, the
tangent coordinates have density proportional to
\[
(1-S-\|v_T\|^2)_+^6.
\]
Radial integration over the tangent ball gives
\[
\mathbb E_\lambda[v_{T,1}^2\mid v_N]
=\frac{1-S}{d+14}.
\]
Therefore,
\[
r(\lambda)
=
\frac{
\mathbb E_\lambda[v_{N,1}^2(1-S)]
}{
(d+14)\mathbb E_\lambda[v_{N,1}^2]
}.
\]
Since $0<S<1$ almost surely under $\nu_\lambda$,
\[
0<r(\lambda)<\frac{1}{d+14}.
\]
Consequently,
\begin{align*}
\Gamma_{d,D}
&=\frac d2\{1+12r(1/3)-12r(1/6)\}\\
&>\frac d2\left(1-\frac{12}{d+14}\right)\\
&=\frac{d(d+2)}{2(d+14)}>0.
\end{align*}

For estimated-normal correspondence, put
$x_k=\pi_{\mathcal M}(Z_{h_k}(x))$.
The population PCA projector estimates in
Section \ref{sec:stsp_validation} verify the conditions of Lemmas \ref{lem:normal_guided_correspondence}
and \ref{lem:normal_field_stability}. Those lemmas give, uniformly in $x$,
\[
\|x_k-x\|=O(h_1^2h_k^2)=O(\sigma^4)
\]
and
\[
\|u_{h_k}(x_k)-u_{h_k}(x)\|
=O(h_1^2h_k^4)=O(\sigma^6).
\]
Since
\[
Z_{h_k}(x)=x_k+u_{h_k}(x_k),
\]
it follows that
\[
\sup_{x\in\mathcal M}
\|Z_{h_k}(x)-X_{h_k}^{\mathrm{STSP}}(x)\|
=O(\sigma^4).
\]
The extrapolation weights are fixed, so both
correspondence constructions satisfy
\[
\sup_{x\in\mathcal M}
\left\|
\widetilde X_{\mathrm{sq}}(x)-x
-\Gamma_{d,D}\sigma^2H(x)
\right\|
\le C\sigma^4.
\]

It remains to establish the Hausdorff asymptotics.
Define the normal graph
\[
F_\sigma(x)
=x+\Gamma_{d,D}\sigma^2H(x).
\]
For sufficiently small $\sigma$,
$\Gamma_{d,D}\sigma^2\|H\|_\infty
<\operatorname{reach}(\mathcal M)$.
Thus
\[
\pi_{\mathcal M}(F_\sigma(x))=x,
\]
and hence
\[
d_H(F_\sigma(\mathcal M),\mathcal M)
=\Gamma_{d,D}\sigma^2\|H\|_\infty.
\]
The uniform expansion implies
\[
d_H\!\left(
\widetilde{\mathcal M}_{\mathrm{sq}},
F_\sigma(\mathcal M)
\right)
\le C\sigma^4.
\]
By the triangle inequality,
\[
\left|
d_H(\widetilde{\mathcal M}_{\mathrm{sq}},\mathcal M)
-\Gamma_{d,D}\sigma^2\|H\|_\infty
\right|
\le C\sigma^4.
\]

Finally, $\|H\|_\infty>0$. Otherwise,
$H\equiv0$ would imply that every coordinate
function of the embedding
$\mathcal M\hookrightarrow\mathbb R^D$
is harmonic. On each compact connected component,
the maximum principle would force all coordinate
functions to be constant, contradicting $d\ge1$.

Therefore,
\[
d_H(\widetilde{\mathcal M}_{\mathrm{sq}},\mathcal M)
=\Theta(\sigma^2),
\]
as claimed.
\end{proof}

\subsection{Proof of Proposition~\ref{lem:stsp_empirical_rate}}
\begin{proof}
Write \(\eta=\eta_n(h)\), and work on the deterministic tube
\(\mathcal T_h\). All estimates below also hold with \(R\) replaced by
\(2R\), so boundary points of the tube cause no differentiation issue.
Derivatives \(D_y^j\) in this proof are ambient derivatives, not derivatives
in the scaled normal coordinates used for the population expansion.

\paragraph*{Uniform differentiated moment bounds}
For \(\ell=0,1,2\), put
\[
M_{\ell,h}(y)
=\mathbb E\!\left[(Y-y)^{\otimes\ell}K_h(y,Y)\right],\qquad
\widehat M_{\ell,h}(y)
=\frac1n\sum_{i=1}^n(Y_i-y)^{\otimes\ell}K_h(y,Y_i),
\]
where \(z^{\otimes0}=1\). In particular,
\(M_{0,h}=m_{0,h}\). The local density expansion and compactness of the
admissible \(\lambda\)-range give, uniformly on the enlarged tube,
\(c h^d\leq M_{0,h}(y)\leq C h^d, \mathbb P\{\|Y-y\|\leq h\}\leq Ch^d\).
Since \(K(t)=(1-t)_+^6\), the radial kernel is \(C^5\), including across
its support boundary. Each component of
\[
f_{\ell,j,y}(Y)
=D_y^j\!\left[(Y-y)^{\otimes\ell}K_h(y,Y)\right],
\qquad j=0,1,2,
\]
is supported on \(\{\|Y-y\|\leq h\}\) and satisfies
\[
\|f_{\ell,j,y}\|_\infty\leq Ch^{\ell-j},\qquad
\mathbb E|f_{\ell,j,y}(Y)|^2\leq Ch^{d+2\ell-2j},\qquad
\|D_y f_{\ell,j,y}\|_\infty\leq Ch^{\ell-j-1}.
\]
The last bound follows by differentiating the kernel once more; it is not
a consequence of a position-error estimate.

To make the concentration uniform, cover the deterministic tube by a
Euclidean net with mesh \(n^{-3}h^{d+1}\). The condition
\(nh^{d+2}/\log n\to\infty\) implies
\(\log(1/h)=O(\log n)\), so this net has at most polynomially many points
in \(n\). Bernstein's inequality at the net points, followed by a union
bound over the finitely many tensor components and derivative orders,
gives, with probability at least \(1-n^{-2}\) after increasing constants,
\[
\sup_{y\in\mathcal T_h}
\bigl\|D_y^j(\widehat M_{\ell,h}-M_{\ell,h})(y)\bigr\|
\leq C\!\left\{
\sqrt{\frac{h^{d+2\ell-2j}\log n}{n}}
+\frac{h^{\ell-j}\log n}{n}
\right\}.
\]
Indeed, the Lipschitz bound above makes the off-net error at most
\(Cn^{-3}h^{d+\ell-j}\), which is absorbed in this display.
Since \(\eta=o(h)\), the Bernstein envelope contribution is no larger
than the square-root contribution for sufficiently large \(n\). Thus
\begin{equation}
\label{eq:stsp-kernel-moment-concentration}
\sup_{y\in\mathcal T_h}
\bigl\|D_y^j(\widehat M_{\ell,h}-M_{\ell,h})(y)\bigr\|
\leq Ch^{d+\ell-j}\eta,
\qquad \ell,j=0,1,2.
\end{equation}
The deterministic moment derivatives obey the corresponding bounds
\(\sup_y\|D_y^jM_{\ell,h}(y)\|\leq Ch^{d+\ell-j}\).
We work on this single high-probability event from now on.
In particular, \(\widehat M_{0,h}(y)\geq ch^d\) throughout the tube.

The identities
\[
\mu_{h}(y)=y+\frac{M_{1,h}(y)}{M_{0,h}(y)},\qquad
\Sigma_{h}(y)
=\frac{M_{2,h}(y)}{M_{0,h}(y)}
-\frac{M_{1,h}(y)M_{1,h}(y)^\top}{M_{0,h}(y)^2}
\]
and their empirical analogues allow differentiation of the normalized
moments. Applying the quotient and product rules to
\eqref{eq:stsp-kernel-moment-concentration} yields
\begin{equation}
\label{eq:stsp-kernel-normalized-concentration}
\begin{aligned}
\sup_{y\in\mathcal T_h}
\bigl\|D_y^j(\widehat\mu_{h}-\mu_{h})(y)\bigr\|
&\leq Ch^{1-j}\eta,\\
\sup_{y\in\mathcal T_h}
\bigl\|D_y^j(\widehat\Sigma_{h}-\Sigma_{h})(y)\bigr\|_{\mathrm{op}}
&\leq Ch^{2-j}\eta,
\qquad j=0,1,2.
\end{aligned}
\end{equation}
For example, the reciprocal-mass derivatives of order \(j\) have size
\(O(h^{-d-j})\); their empirical perturbations have an additional factor
\(\eta\). Every product in the quotient rules therefore has the displayed
power of \(h\). Products of two empirical errors are absorbed because
\(\eta=o(1)\).

\paragraph*{Empirical eigengap and residual derivatives}
The population covariance has eigengap at least \(gh^2/2\), as proved
above. By \eqref{eq:stsp-kernel-normalized-concentration} and Weyl's
inequality, the empirical eigengap is at least \(gh^2/4\) for sufficiently
large \(n\). Hence the empirical PCA projector and residual are well
defined and at least \(C^2\) on a neighborhood of the tube.

Put \(S_h=h^{-2}\Sigma_{h}\) and
\(\widehat S_h=h^{-2}\widehat\Sigma_{h}\). These matrices range in a
bounded set with a gap bounded away from zero. Let \(\mathcal Q\) denote
the complementary spectral-projector map. Its matrix derivatives through
order three are uniformly bounded on a neighborhood of this set. For one
ambient derivative, the chain rule gives the exact identity
\[
\begin{aligned}
D_y(\widehat Q_{h}-Q_{h})
={}&D\mathcal Q(\widehat S_h)
   [D_y(\widehat S_h-S_h)]\\
&+\{D\mathcal Q(\widehat S_h)-D\mathcal Q(S_h)\}[D_yS_h].
\end{aligned}
\]
Its two terms are bounded by \(C\eta/h\). For two derivatives, apply
\[
D_y^2\mathcal Q(S_h)
=D^2\mathcal Q(S_h)[D_yS_h,D_yS_h]
+D\mathcal Q(S_h)[D_y^2S_h]
\]
and the empirical analogue. The bounds
\(D_y^j(\widehat S_h-S_h)=O(h^{-j}\eta)\) and
\(D_y^jS_h=O(h^{-j})\) make their difference \(O(\eta/h^2)\).
Together with the zeroth-order projector perturbation estimate, this gives
\begin{equation}
\label{eq:stsp-kernel-projector-concentration}
\sup_{y\in\mathcal T_h}
\bigl\|D_y^j(\widehat Q_{h}-Q_{h})(y)\bigr\|_{\mathrm{op}}
\leq Ch^{-j}\eta,
\qquad j=0,1,2.
\end{equation}
The same chain rules give the deterministic bounds
\(\sup_{y\in\mathcal T_h}\|D_y^jQ_{h}(y)\|\leq Ch^{-j}, j=0,1,2\),
and the empirical counterparts satisfy these bounds on the same event.

Abbreviate
\(R_h=R_{h}^{\mathrm{STSP}}\),
\(\widehat R_h=\widehat R_{h}^{\mathrm{STSP}}\), and
\(v_h(y)=y-\mu_{h}(y)\).
The population expansion gives \(\sup_{\mathcal T_h}\|v_h\|\leq Ch^2\),
while the differentiated moment bounds give
\(\|D_yv_h\|\leq C\) and \(\|D_y^2v_h\|\leq C/h\).
Use the exact identity
\[
\widehat R_h-R_h
=(\widehat Q_{h}-Q_{h})v_h
-\widehat Q_{h}(\widehat\mu_{h}-\mu_{h}).
\]
Its derivative is
\[
\begin{aligned}
D_y(\widehat R_h-R_h)
={}&D_y(\widehat Q_{h}-Q_{h})v_h
+(\widehat Q_{h}-Q_{h})D_yv_h\\
&-D_y\widehat Q_{h}(\widehat\mu_{h}-\mu_{h})
-\widehat Q_{h}D_y(\widehat\mu_{h}-\mu_{h}).
\end{aligned}
\]
The preceding bounds therefore imply
\begin{equation}
\label{eq:stsp-kernel-residual-C1}
\sup_{\mathcal T_h}\|\widehat R_h-R_h\|\leq Ch\eta,\qquad
\sup_{\mathcal T_h}\|D_y\widehat R_h-D_yR_h\|_{\mathrm{op}}\leq C\eta.
\end{equation}
They also give the deterministic estimates
\begin{equation}
\label{eq:stsp-population-ambient-derivatives}
\sup_{\mathcal T_h}\|D_yR_h\|\leq C,\qquad
\sup_{\mathcal T_h}\|D_y^2R_h\|\leq C/h.
\end{equation}
Thus the residual derivative control has been derived from differentiated
kernel moments and spectral projectors, independently of any assumed
empirical branch.

\paragraph*{Existence, uniqueness, and gluing of the empirical branch}
Work in a fixed local chart on \(\mathcal M\), with a smooth orthonormal
normal frame \(E(x)\in\mathbb R^{D\times q}\), \(q=D-d\).
For \(z\in\mathbb R^q\), \(\|z\|\leq Rh^2\), define the unscaled normal
equations
\[
f_h(x,z)=E(x)^\top R_h(x+E(x)z),\qquad
\widehat f_h(x,z)=E(x)^\top\widehat R_h(x+E(x)z).
\]
Let \(z_h(x)=E(x)^\top u_h^{\mathrm{STSP}}(x)\) be the population zero.
Increasing the fixed \(R\), if necessary, ensures
\(\sup_x\|z_h(x)\|\leq Rh^2/2\). Put
\[
\gamma_{\lambda}=\frac{a_{\lambda}}{\lambda^2},\qquad
\inf_{\lambda\in\Lambda_{\mathrm{STSP}}}\gamma_{\lambda}>0.
\]
The population scaled normal equation, expressed in this frame, gives
\[
\sup_{x,\,\|z\|\leq Rh^2}
\|D_zf_h(x,z)-\gamma_{\lambda}I_q\|\leq Ch^2.
\]
Equation~\eqref{eq:stsp-kernel-residual-C1} implies
\[
\sup_{x,z}\|\widehat f_h-f_h\|\leq Ch\eta,\qquad
\sup_{x,z}\|D_{(x,z)}\widehat f_h-D_{(x,z)}f_h\|\leq C\eta.
\]
Consequently, for sufficiently small \(h\) and large \(n\),
\(\left\|I_q-\gamma_{\lambda}^{-1}D_z\widehat f_h(x,z)\right\| \leq\frac12\)
uniformly on these normal balls.

Consider the contraction
\(\mathcal T_{h,x}(z) =z-\gamma_{\lambda}^{-1}\widehat f_h(x,z)\).
At \(z=z_h(x)\), its displacement is at most \(C_0h\eta\).
It therefore maps the closed ball
\(\overline B_{2C_0h\eta}(z_h(x))\) into itself and has a unique fixed
point \(\widehat z_h(x)\). This ball lies inside \(B_{Rh^2}(0)\), because
\(h\eta=o(h^2)\). Moreover,
\begin{equation}
\label{eq:stsp-empirical-root-position}
\sup_x\|\widehat z_h(x)-z_h(x)\|\leq Ch\eta.
\end{equation}
The derivative estimate on the whole convex ball \(\overline B_{Rh^2}(0)\) also
implies strict monotonicity there, so no additional zero of
\(\widehat f_h(x,\cdot)\) occurs in that ball. The parameterized
implicit-function theorem makes \(\widehat z_h\) at least \(C^2\).

The projected equation is equivalent to the full residual equation.
Indeed, for \(y=x+E(x)z\), the population projector expansion and
\eqref{eq:stsp-kernel-projector-concentration} give
\(\|\widehat Q_{h}(y)-Q_x\|\leq C(h^2+\eta)<\frac12\).
If \(\widehat f_h(x,z)=0\), then
\(Q_x\widehat R_h(y)=0\), whereas
\(\widehat R_h(y)=\widehat Q_{h}(y)\widehat R_h(y)\). Thus
\(\widehat R_h(y)=(\widehat Q_{h}(y)-Q_x)\widehat R_h(y)\), forcing
\(\widehat R_h(y)=0\). The reverse implication is immediate.
On overlapping normal frames, uniqueness therefore identifies the same
ambient point \(x+E(x)\widehat z_h(x)\). The local sections glue into a
global \(C^2\) normal section \(\widehat u_h^{\mathrm{STSP}}\), and all
empirical zeros in \(\mathcal T_h\) belong to this graph.

\paragraph*{Derivative stability and embedding}
The smooth fixed charts and
\eqref{eq:stsp-population-ambient-derivatives} imply
\[\sup\|D_{(x,z)}^2f_h\|\leq C/h\].
Differentiate the two normal equations to obtain
\[
D_xz_h=-B_h^{-1}C_h,\qquad
D_x\widehat z_h=-\widehat B_h^{-1}\widehat C_h,
\]
where
\[
\begin{aligned}
B_h(x)&=D_zf_h(x,z_h(x)),&
C_h(x)&=D_xf_h(x,z_h(x)),\\
\widehat B_h(x)&=D_z\widehat f_h(x,\widehat z_h(x)),&
\widehat C_h(x)&=D_x\widehat f_h(x,\widehat z_h(x)).
\end{aligned}
\]
Both inverse matrices are uniformly bounded. Combining the residual
\(C^1\) bound with \eqref{eq:stsp-empirical-root-position} gives
\[
\|\widehat B_h-B_h\|_\infty+
\|\widehat C_h-C_h\|_\infty
\leq C\eta+\frac Ch\|\widehat z_h-z_h\|_\infty
\leq C\eta.
\]
The population section has uniformly bounded \(C^1\) norm (in fact,
\(\|z_h\|_{C^1}=O(h^2)\)). The inverse-matrix identity now yields
\(\|D_x\widehat z_h-D_xz_h\|_\infty\leq C\eta\).
Returning to the ambient normal sections, and using the bounded frame
derivatives, proves
\[
\|\widehat u_h^{\mathrm{STSP}}-u_h^{\mathrm{STSP}}\|_\infty
\leq Ch\eta,\qquad
\|\widehat u_h^{\mathrm{STSP}}-u_h^{\mathrm{STSP}}\|_{C^1}
\leq C\eta.
\]
Since \(\widehat u_h^{\mathrm{STSP}}(x)\in N_x\mathcal M\) and
\(\|\widehat u_h^{\mathrm{STSP}}\|_\infty\leq Rh^2<\tau\),
\(\pi_{\mathcal M}\circ\widehat X_h^{\mathrm{STSP}} =\operatorname{id}_{\mathcal M}\).
As in the population proof, this identity gives injectivity and immersion;
compactness then makes \(\widehat X_h^{\mathrm{STSP}}\) a \(C^2\)
embedding. The derivative estimate
\(\|D\widehat X_h^{\mathrm{STSP}}-D\iota\|_\infty
\leq C(h^2+\eta)\) also gives the asserted tangent-projector bound.
This proves the proposition.

\paragraph*{PCA at the empirical reference point}
We record the derivative estimates used subsequently for correspondence.
Set
\[
\widetilde Q_h^{\mathrm{ref}}(x)
=Q_{h}\!\left(X_h^{\mathrm{STSP}}(x)\right),\qquad
\widehat Q_h^{\mathrm{ref}}(x)
=\widehat Q_{h}\!\left(\widehat X_h^{\mathrm{STSP}}(x)\right).
\]
The population normal-coordinate expansion, composed with the population
section, gives
\(\|\widetilde Q_h^{\mathrm{ref}}-Q\|_{C^1}\leq Ch^2\), where
\(Q(x)=Q_x\). To compare the two reference fields, use
\[
\begin{aligned}
\widehat Q_h^{\mathrm{ref}}-\widetilde Q_h^{\mathrm{ref}}
={}&(\widehat Q_{h}-Q_{h})\circ\widehat X_h^{\mathrm{STSP}}\\
&+Q_{h}\circ\widehat X_h^{\mathrm{STSP}}
-Q_{h}\circ X_h^{\mathrm{STSP}}.
\end{aligned}
\]
Equation~\eqref{eq:stsp-kernel-projector-concentration}, the position and
first-derivative estimates just proved, and the ambient bounds
\(\|D_yQ_{h}\|\leq C/h\),
\(\|D_y^2Q_{h}\|\leq C/h^2\), give
\begin{equation}
\label{eq:stsp-reference-PCA-stability}
\|\widehat Q_h^{\mathrm{ref}}-\widetilde Q_h^{\mathrm{ref}}\|_\infty
\leq C\eta,\qquad
\|\widehat Q_h^{\mathrm{ref}}-\widetilde Q_h^{\mathrm{ref}}\|_{C^1}
\leq C\eta/h.
\end{equation}
For the second estimate, differentiating the last display introduces
\(D_y(\widehat Q_{h}-Q_{h})\), a change of evaluation point in
\(D_yQ_{h}\), and the difference
\(D\widehat X_h^{\mathrm{STSP}}-DX_h^{\mathrm{STSP}}\).
Their respective bounds are
\(C\eta/h\), \(Ch^{-2}(h\eta)\), and \(Ch^{-1}\eta\).
This also verifies the reference-field claims in Section~\ref{sec:stsp_validation}, including the
uniform \(C^1\) bound under \(\eta=o(h)\).
\end{proof}

\subsection{Proof of Theorem~\ref{thm:finite_sample_stsp}}
\label{app:stsp_finite_sample_proof}
\begin{proof}
Proposition~\ref{thm:stsp_bias_expansion} gives the population expansion
\[
u_h^{\mathrm{STSP}}(x)
=
c_{\mathrm{STSP}}(h,\sigma)H(x)
+
O_{C^2}(h^4)
\]
and verifies the required tangent-space stability.

Apply Proposition~\ref{lem:stsp_empirical_rate} simultaneously at the fixed
number of selected scales. On the resulting event of probability tending
to one, the local empirical branches exist uniquely and
\[
\|\widehat X_{h_k}^{\mathrm{STSP}}-X_{h_k}^{\mathrm{STSP}}\|_\infty
\leq Ch_k\eta_n(h_k),\qquad
\|\widehat X_{h_k}^{\mathrm{STSP}}-X_{h_k}^{\mathrm{STSP}}\|_{C^1}
\leq C\eta_n(h_k).
\]
The population and empirical reference projector fields satisfy
\eqref{eq:stsp-reference-PCA-stability} at \(h=h_1\).
In particular, \(\widehat Q_1\) has a uniformly bounded \(C^1\) norm and
is uniformly close to \(Q_x\).

We first verify the local empirical correspondence. In a local chart \(\psi\) of \(\mathcal M\), choose a
\(C^1\) orthonormal frame \(E_1(x)\) of the population reference normal
space. Use the aligned empirical frame
\[
\widehat E_1(x)
=\widehat Q_1(x)E_1(x)
\bigl[E_1(x)^\top\widehat Q_1(x)E_1(x)\bigr]^{-1/2}.
\]
The reference-field bounds imply
\[
\|\widehat E_1-E_1\|_\infty\leq C\eta_n(h_1),\qquad
\|\widehat E_1-E_1\|_{C^1}\leq C\eta_n(h_1)/h_1.
\]
Define the local intersection equations
\[
\begin{aligned}
F_k(q,s;x)
&=X_{h_k}^{\mathrm{STSP}}(\psi(q))
-X_{h_1}^{\mathrm{STSP}}(x)-E_1(x)s,\\
\widehat F_k(q,s;x)
&=\widehat X_{h_k}^{\mathrm{STSP}}(\psi(q))
-\widehat X_{h_1}^{\mathrm{STSP}}(x)-\widehat E_1(x)s.
\end{aligned}
\]
By Lemma~\ref{lem:wellposed_correspondence}, the population solution
\((q_k(x),s_k(x))\) exists, \(\|s_k(x)\|\leq Ch_k^2\), and
\(D_{(q,s)}F_k\) has a uniformly bounded inverse near that solution.
The empirical \(C^1\) branch and frame bounds give
\[
\|D_{(q,s)}\widehat F_k-D_{(q,s)}F_k\|
\leq C\{\eta_n(h_k)+\eta_n(h_1)\}=o(1)
\]
on fixed local neighborhoods. At the population solution,
\[
\|\widehat F_k(q_k,s_k;x)\|
\leq C e_{n,k}^{\mathrm{corr}},\qquad
e_{n,k}^{\mathrm{corr}}=h_1\eta_n(h_1)+h_k\eta_n(h_k)
+h_k^2\eta_n(h_1).
\]
The uniform quantitative implicit-function argument therefore gives a
unique nearby empirical intersection, with displacement at most
\(Ce_{n,k}^{\mathrm{corr}}\). Its differential in \((q,s)\) remains invertible, proving
transversality. Since the empirical equation is \(C^1\) in \(x\), these
intersections depend \(C^1\)-smoothly on the reference point. Local
uniqueness glues them across charts, independently of the chosen frames.
For \(k=1\), the intersection is simply
\(\widehat X_{h_1}^{\mathrm{STSP}}(x)\).
Thus every empirical correspondence required in the construction is now
well defined, and
\[
\sup_{x\in\mathcal M}
\|\widehat Z_{h_k}^{\mathrm{STSP}}(x)-Z_{h_k}^{\mathrm{STSP}}(x)\|
\leq Ce_{n,k}^{\mathrm{corr}}.
\]

The same event verifies
Assumption~\ref{assump:normal_guided_empirical} with the deterministic bound
\(r_n^{\mathrm{STSP}}(h_k)=Ch_k\eta_n(h_k)\),
and Assumption~\ref{assump:finite_normal_estimation} with
\(s_n^{\mathrm{STSP}}(h_1)=C\eta_n(h_1)\).
Substituting these rates into
Theorem~\ref{thm:finite_normal_guided_decomposition} gives
\[
d_H\!\left(\widehat{\mathcal M}^{\mathrm{STSP}},\mathcal M\right)
=O_p\!\left(
\sum_{k=1}^K|w_k|
\left[h_k^4+h_1\eta_n(h_1)+h_k\eta_n(h_k)
+h_k^2\eta_n(h_1)\right]\right).
\]

Finally, we verify the empirical embedding. The preceding bounds give
\[
\|\widehat X_{h_1}^{\mathrm{STSP}}-\iota\|_{C^1}
\leq C\{h_1^2+\eta_n(h_1)\},\qquad
\|\widehat Q_1-Q\|_\infty\leq C\{h_1^2+\eta_n(h_1)\},\qquad
\|\widehat Q_1\|_{C^1}\leq C.
\]
Consider the empirical affine-fiber map
\[
\widehat\Phi_{\mathrm{fib}}(x,v)=\widehat X_{h_1}^{\mathrm{STSP}}(x)+v,
\qquad v\in\widehat N_{h_1}^{\mathrm{STSP}}(x),\quad \|v\|<\rho.
\]
For sufficiently small fixed \(\rho>0\), its local differential is
invertible by the same affine-bundle argument as in the proof of
Proposition~\ref{thm:normal_guided_manifold_recovery}: differentiating
\(\widehat Q_1v=v\) gives
\((I-\widehat Q_1)\dot v=D\widehat Q_1[\dot x]v\), whereas
\((I-\widehat Q_1)D\widehat X_{h_1}^{\mathrm{STSP}}\) has a uniformly
positive lower singular value. The \(C^1\) projector bound makes these
two identities incompatible with a nonzero kernel vector when \(\rho\)
is small.
For global injectivity, equality of two images first implies that their
base points are in a common local chart, using the \(C^0\) closeness to
\(\iota\). In that chart the lower bound for the projected base chord,
together with
\(\|\widehat Q_1(x)-\widehat Q_1(x')\|\leq C\|x-x'\|\), gives
\(c\|x-x'\|\leq C\rho\|x-x'\|\). Reducing \(\rho\) forces
\(x=x'\), and then the fiber vectors agree. Thus \(\widehat\Phi_{\mathrm{fib}}\) is a
\(C^1\) embedding on a uniform disk bundle.

Every corresponded displacement
\(\widehat v_k(x) =\widehat Z_{h_k}^{\mathrm{STSP}}(x) -\widehat X_{h_1}^{\mathrm{STSP}}(x)\)
lies in the empirical reference fiber and has uniform norm \(O(h_k^2)\)
on this event; here \(\eta_n(h_k)=o(h_k)\) and \(h_1\leq h_k\).
Since the weights sum to one and their absolute sum is bounded,
\[
\widehat X^{\mathrm{STSP}}(x)
=\widehat\Phi_{\mathrm{fib}}\!\left(x,\sum_{k=1}^K w_k\widehat v_k(x)\right),\qquad
\sup_x\left\|\sum_{k=1}^K w_k\widehat v_k(x)\right\|
\leq CW_0h_{\max}^2<\rho
\]
for sufficiently small \(h_{\max}\). The graph of this \(C^1\) section
is an embedding into the disk bundle. Composition with
\(\widehat\Phi_{\mathrm{fib}}\) proves that \(\widehat X^{\mathrm{STSP}}\) is a \(C^1\)
embedding and its image is \(C^1\)-diffeomorphic to \(\mathcal M\).

\paragraph*{Numerical correspondence error}
Work on the same event. Write
\(z=\widehat X_{h_1}^{\mathrm{STSP}}(x)\),
\(\widehat E(z)=\widehat E_1(x)\), and let \(t_k^*(z)\) be the coordinates
of the exact empirical intersection, so that
\[
\widehat Z_{h_k}^{\mathrm{STSP}}(x)
=z+\widehat E(z)t_k^*(z),\qquad g_{k,z}(t_k^*(z))=0.
\]
For \(C_{\mathrm{loc}}\) sufficiently large, these coordinates lie in
\(\|t\|\leq C_{\mathrm{loc}}h_k^2\). All moment and derivative estimates
remain valid on a fixed enlargement of \(\mathcal T_{h_k}\) containing
these search balls. Put
\[
\gamma_k=\frac{a_{\sigma/h_k}}{(\sigma/h_k)^2},\qquad
\gamma_{\mathrm{ENF},*}:=\inf_{\lambda\in\Lambda_{\mathrm{STSP}}}
\frac{a_{\lambda}}{\lambda^2}>0.
\]
At the population fixed point \(y_k^0=X_{h_k}^{\mathrm{STSP}}(x)\),
differentiating
\(R_{h_k}^{\mathrm{STSP}}\circ X_{h_k}^{\mathrm{STSP}}=0\) and
\(R_{h_k}^{\mathrm{STSP}}=Q_{h_k}R_{h_k}^{\mathrm{STSP}}\),
and using the scaled normal equation and
\(DX_{h_k}^{\mathrm{STSP}}=D\iota+O(h_k^2)\), gives
\[
D_yR_{h_k}^{\mathrm{STSP}}(y_k^0)
=\gamma_k Q_x+O(h_k^2).
\]
For \(\|t\|\leq C_{\mathrm{loc}}h_k^2\), the point
\(z+\widehat E(z)t\) is \(O(h_k^2)\) from \(y_k^0\).
Using \eqref{eq:stsp-population-ambient-derivatives},
\eqref{eq:stsp-kernel-residual-C1}, and the reference projector bound
therefore yields, uniformly on these balls,
\[
\begin{aligned}
D_tg_{k,z}(t)
&=\widehat E(z)^\top
D_y\widehat R_{h_k}^{\mathrm{STSP}}(z+\widehat E(z)t)
\widehat E(z),\\
\|D_tg_{k,z}(t)-\gamma_k I_q\|_{\mathrm{op}}
&\leq C\{h_k+h_1^2+\eta_n(h_k)+\eta_n(h_1)\}=o(1).
\end{aligned}
\]
Consequently, for sufficiently small scales and large \(n\), integration
along the segment between any two points \(s,t\) in the same search ball
gives
\[
\langle g_{k,z}(t)-g_{k,z}(s),t-s\rangle
\geq\frac{\gamma_{\mathrm{ENF},*}}{2}\|t-s\|^2.
\]
At each retained query \(z\in\mathcal Z_N\), taking \(s=t_k^*(z)\),
the residual tolerance and Cauchy--Schwarz give
\[
\|\widehat t_k(z)-t_k^*(z)\|
\leq\frac{2}{\gamma_{\mathrm{ENF},*}}\|g_{k,z}(\widehat t_k(z))\|
\leq C\varepsilon_{n,k}.
\]
Set \(t_1^*(z)=\widehat t_1(z)=0\). Since \(\widehat E(z)\) is an
isometry, the numerical and exact extrapolated outputs satisfy
\[
\max_{1\leq j\leq N}
\left\|
z_j+\widehat E(z_j)\sum_{k=1}^K w_k\widehat t_k(z_j)
-\widehat F^{\mathrm{STSP}}(z_j)
\right\|
\leq C\sum_{k=1}^K|w_k|\varepsilon_{n,k}.
\]

\paragraph*{Finite query coverage}
The exact reference-branch map \(\widehat F^{\mathrm{STSP}}\) has a
Lipschitz constant bounded uniformly on the same event. Indeed, the
intersection equations above have uniformly bounded
\(D_x\widehat F_k\) and uniformly bounded inverses of
\(D_{(q,s)}\widehat F_k\): the branch and reference-frame derivatives
are bounded, and \(\|\widehat s_k\|=O(h_k^2)\).
Implicit differentiation therefore bounds
\(D_x\widehat Z_{h_k}^{\mathrm{STSP}}\) uniformly. Bounded weights give
a uniform derivative bound for \(\widehat X^{\mathrm{STSP}}\).
Moreover, \(\widehat X_{h_1}^{\mathrm{STSP}}\) is a uniformly
bi-Lipschitz normal graph, since its \(C^1\) distance from the inclusion
tends to zero. In fixed local charts these derivative bounds give
Lipschitz estimates in the ambient distance; separated base points are
handled by boundedness of the maps and the uniform separation of their
reference images. Consequently,
\[
\|\widehat F^{\mathrm{STSP}}(z)
  -\widehat F^{\mathrm{STSP}}(z')\|
\leq C\|z-z'\|,
\qquad z,z'\in\widehat{\mathcal M}_{h_1}^{\mathrm{STSP}}.
\]
Define the exact sampled image
\(\mathcal P_N^*=\{\widehat F^{\mathrm{STSP}}(z_j):1\leq j\leq N\}\).
It is a subset of the continuous estimator, and every point of that
estimator is within \(C\delta_N\) of \(\mathcal P_N^*\). Thus
\[
\begin{aligned}
d_H\!\left(\widehat{\mathcal P}_N^{\mathrm{STSP,num}},
           \mathcal P_N^*\right)
&\leq C\sum_{k=1}^K|w_k|\varepsilon_{n,k},\\
d_H\!\left(\mathcal P_N^*,\widehat{\mathcal M}^{\mathrm{STSP}}\right)
&\leq C\delta_N.
\end{aligned}
\]
The triangle inequality proves the additional numerical and covering
error bound. For comparable scales and bounded weights,
\(\max_k\varepsilon_{n,k}=O(h_{\max}^4)\) and
\(\delta_N=O(h_{\max}^4)\) make these contributions \(O(h_{\max}^4)\).
The proof only uses the retained nonempty query set and does not require
a continuous numerical solution map. No derivative bound on numerical
outputs, or global branch-selection guarantee, is used or inferred here.
\end{proof}

\subsection{Proof of noise-level sensitivity of the STSP coefficient}
\label{app:stsp_noise_sensitivity}
\begin{proof}
Write \(\psi(\lambda)=\lambda^2+c_\lambda/a_\lambda\), so that
\(c_{\mathrm{STSP}}(h,s)=(d/2)h^2\psi(s/h)\).
Let \(\Lambda_+\Subset(0,\infty)\) be the fixed compact interval containing
both noise-to-bandwidth ratios. For a bounded function \(f(v)\) independent
of \(\lambda\), differentiating the normalized kernel integral gives
\[
\frac{d}{d\lambda}\mathbb E_{\lambda}[f]
=\lambda^{-3}\operatorname{Cov}_{\lambda}
\bigl(f,\|v_N\|^2\bigr).
\]
Differentiation under the integral is justified by the fixed compact support
and the positive lower bound for \(\lambda\). The mass denominator is
uniformly positive on \(\Lambda_+\). Apply the identity to
\(f(v)=v_{N,1}^2\) and \(f(v)=v_{N,1}^2v_{T,1}^2\).
The resulting moments and their first derivatives are uniformly bounded;
positivity and continuity also give
\(\inf_{\lambda\in\Lambda_+}a_{\lambda}>0\). Hence
\[
\psi'(\lambda)
=2\lambda+
\frac{c_{\lambda}'a_{\lambda}
-c_{\lambda}a_{\lambda}'}{a_{\lambda}^2}
\]
is uniformly bounded. It follows that
\[
\partial_s c_{\mathrm{STSP}}(h,s)
=\frac d2h\,\psi'(s/h),\qquad
|\partial_s c_{\mathrm{STSP}}(h,s)|\leq Ch
\quad\text{when }s/h\in\Lambda_+.
\]
The interval contains the segment between \(\sigma/h\) and
\(\widehat\sigma/h\), so the mean value theorem yields the asserted
\(Ch|\widehat\sigma-\sigma|\) coefficient bound.
Proposition~\ref{prop:noise_level_perturbation}, together with the existing
STSP empirical correspondence estimates, then gives the additional error
\(C|\widehat\sigma-\sigma|\sum_k|\widehat w_k|h_k\).
The numerical root and query-coverage bounds are also pathwise and use
only the absolute weight sum; their proofs therefore apply unchanged to
these spatially constant plug-in weights. Their randomness introduces no
spatial derivatives of the weights.
\end{proof}

\subsection{Noise-level estimation from the STSP covariance spectrum}
\label{app:stsp_noise_estimation}

\begin{proposition}[Covariance-based STSP noise estimation]
\label{prop:stsp_noise_estimation}
Suppose the conditions of Proposition~\ref{lem:stsp_empirical_rate} hold at the
reference bandwidth \(h_1\), and define \(\widehat\sigma\) as described in Section \ref{sec:stsp_validation} using a fixed compact interval \(\Lambda_+\) containing
\(\Lambda_{\mathrm{STSP}}\). Then \(a:\lambda\mapsto a_\lambda\) is
strictly increasing, with \(\inf_{\lambda\in\Lambda_+}a'(\lambda)>0\).
For sufficiently small \(h_1\) and large \(n\), on an event of probability
at least \(1-n^{-2}\) where the reference branch exists,
\[
|\widehat\sigma-\sigma|
\leq C\{h_1^3+h_1\eta_n(h_1)\}.
\]
The bound holds uniformly over measurable selections of
\(\widehat y_1\) from that branch, including selections based on the same
observations. If the remaining assumptions of
Theorem~\ref{thm:finite_sample_stsp} hold at all fitted scales, replacing
its weights by spatially constant, uniformly bounded plug-in weights at a
fixed number of comparable deterministic scales preserves its recovery
rate and embedding conclusion, as well as the sufficient sample-size
condition of Corollary~\ref{cor:stsp_fourth_order_sample_size}.
\end{proposition}

\begin{proof}
Write \(q=D-d\), \(S=\|v_N\|^2\), and \(h=h_1\).
Rotational invariance gives \(a_\lambda=q^{-1}\mathbb E_\lambda S\).
The differentiation identity in
Appendix~\ref{app:stsp_noise_sensitivity} therefore yields
\[
a'(\lambda)
=\frac{\operatorname{Var}_\lambda(S)}{q\lambda^3}>0.
\]
The variance is strictly positive because the kernel-weighted law has
positive density in the interior of the unit ball and \(S\) is not
constant there. Continuity and compactness give
\(m_a:=\inf_{\lambda\in\Lambda_+}a'(\lambda)>0\), so \(a^{-1}\) is
\(m_a^{-1}\)-Lipschitz on \(a(\Lambda_+)\).
In particular, the calibration depends only on \(d,D\) and the fixed
kernel. Integrating out the tangent coordinates gives the
one-dimensional formula
\[
a_\lambda
=\frac{1}{q}
\frac{\displaystyle\int_0^1
s^{q/2}(1-s)^{d/2+6}e^{-s/(2\lambda^2)}\,ds}
{\displaystyle\int_0^1
s^{q/2-1}(1-s)^{d/2+6}e^{-s/(2\lambda^2)}\,ds}.
\]

Work on the uniform concentration event used in the proof of
Proposition~\ref{lem:stsp_empirical_rate}; its probability is at least
\(1-n^{-2}\). That proof also establishes existence of the reference
branch on this event. By Lemma~\ref{lem:pca_normal_accuracy} and
\eqref{eq:stsp-kernel-normalized-concentration}, uniformly on
\(\mathcal T_h\),
\[
\left\|h^{-2}\widehat\Sigma_h(y)
-\{\tau_{\sigma/h}P_x+a_{\sigma/h}Q_x\}\right\|_{\mathrm{op}}
\leq C\{h^2+\eta_n(h)\},\qquad x=\pi_{\mathcal M}(y).
\]
The population tangent eigenvalues exceed the normal eigenvalues
uniformly in the admissible regime. Weyl's inequality consequently gives
\[
\left|\frac{1}{qh^2}\sum_{j=d+1}^D
\lambda_j(\widehat\Sigma_h(y))-a_{\sigma/h}\right|
\leq C\{h^2+\eta_n(h)\}.
\]
Since \(\widehat y_1\in\widehat{\mathcal M}_h^{\mathrm{STSP}}
\subset\mathcal T_h\) on this event, the inequality applies to the
selected, possibly data-dependent center. Projection onto
\(a(\Lambda_+)\) cannot increase its distance from \(a_{\sigma/h}\).
The inverse Lipschitz bound proves
\[
|\widehat\sigma-\sigma|
\leq h m_a^{-1}|\widehat a_1-a_{\sigma/h}|
\leq C\{h^3+h\eta_n(h)\}.
\]
Outside the branch-existence event, the estimate may be assigned any
fixed value in \(h\Lambda_+\). This does not affect the asserted rate.

For the multiscale conclusion, take \(h_k\asymp h_1\).
All true and estimated ratios \(\sigma/h_k\) and
\(\widehat\sigma/h_k\) then lie in a common fixed compact subinterval
of \((0,\infty)\). The coefficient-sensitivity bound in
Appendix~\ref{app:stsp_noise_sensitivity} and
Proposition~\ref{prop:noise_level_perturbation} give the additional error
\[
C|\widehat\sigma-\sigma|\sum_k|\widehat w_k|h_k
\leq CW_0h_{\max}\{h_1^3+h_1\eta_n(h_1)\}
=O\!\left(h_{\max}^4+h_{\max}^2\eta_n(h_{\max})\right)
\]
on the joint concentration event. This is absorbed into the original
\(O(h_{\max}^4+h_{\max}\eta_n(h_{\max}))\) bound.
The embedding argument in the proof of
Theorem~\ref{thm:finite_sample_stsp} is pathwise: the empirical affine
fibers and correspondence maps are unchanged, and the weights are
spatially constant, sum to one, and have bounded absolute sum.
It therefore also applies to the plug-in weights. The numerical and
covering estimates are likewise unchanged apart from the weight
replacement.
Under the sample-size condition of
Corollary~\ref{cor:stsp_fourth_order_sample_size},
\(\eta_n(h_{\max})=O(h_{\max}^3)\); hence both the original and added
errors are \(O(h_{\max}^4)\) on the same joint event of probability
at least \(1-Kn^{-2}\), after adjusting constants.
\end{proof}

\subsection{Uniform-in-bandwidth STSP stability and oracle comparison}
\label{subsec:stsp_uniform_oracle}
Put $\mathcal H_n=[a\sigma,b\sigma]$ and
$\varepsilon_n=(\log n/(n\sigma^d))^{1/2}$, as in Corollary~\ref{cor:stsp_oracle_separation}. We first make the concentration argument of the proof of
Proposition~\ref{lem:stsp_empirical_rate} uniform in $h\in\mathcal H_n$.

For $\ell=0,1,2$ and $j=0,1,2$, let
\[
  f_{\ell,j,h,y}(Y)
  :=D_y^j\!\left[(Y-y)^{\otimes\ell}
          K\!\left(\frac{\|Y-y\|^2}{h^2}\right)\right],
  \qquad K(t)=(1-t)_+^6.
\]
Fix $C_T<\infty$ large enough that the enlarged local tubes for all
$h\in\mathcal H_n$ lie in
$T_\sigma^*:=\{y:\operatorname{dist}(y,\mathcal M)\le C_T\sigma^2\}$.
Since $h/\sigma\in[a,b]$, the kernel bounds used for \eqref{eq:stsp-kernel-moment-concentration} hold
uniformly over $h\in\mathcal H_n$ and $y\in T_\sigma^*$:
\[
\begin{aligned}
 \|f_{\ell,j,h,y}\|_\infty
   &\le C\sigma^{\ell-j},\\
 \mathbb E\|f_{\ell,j,h,y}(Y)\|^2
   &\le C\sigma^{d+2\ell-2j},\\
 \|D_y f_{\ell,j,h,y}\|_\infty
  +\|\partial_h f_{\ell,j,h,y}\|_\infty
   &\le C\sigma^{\ell-j-1}.
\end{aligned}
\]
The last bound is available because the selected polynomial kernel
is $C^5$, and all derivative orders used here are at most three.
Cover $T_\sigma^*\times\mathcal H_n$ by a Euclidean product net of
mesh $n^{-3}\sigma^{d+1}$ in both $y$ and $h$. The condition
$n\sigma^{d+2}/\log n\to\infty$ implies
$\log(1/\sigma)=O(\log n)$, so the cardinality of this
net is polynomial in $n$. Bernstein's inequality, a union bound over
the net, tensor components and derivative orders, and the preceding
Lipschitz bound give, with probability at least $1-Cn^{-2}$,
\begin{equation}\label{eq:uniform_stsp_moments}
 \sup_{h\in\mathcal H_n}\sup_{y\in T_\sigma^*}
 \left\|D_y^j\bigl(\widehat M_{\ell,h}(y)-M_{\ell,h}(y)\bigr)\right\|
 \le C\sigma^{d+\ell-j}\varepsilon_n,
 \qquad 0\le\ell,j\le2.
\end{equation}
Here $M_{\ell,h}$ and $\widehat M_{\ell,h}$ are the raw local
kernel moments defined in the proof of Proposition~\ref{lem:stsp_empirical_rate}. In
particular, the off-net error is at most
$Cn^{-3}\sigma^{d+\ell-j}$ and is absorbed in the right-hand
side of \eqref{eq:uniform_stsp_moments}.

The population kernel-mass lower bound, covariance eigengap and
normal-equation nondegeneracy are uniform in
$\sigma/h\in[1/b,1/a]$ by Lemmas~\ref{lem:pca_normal_accuracy} and~\ref{lem:stsp_scaled_normal_equation}. Repeating the
quotient, spectral-projector and residual estimates \eqref{eq:stsp-kernel-normalized-concentration}--\eqref{eq:stsp-population-ambient-derivatives}
with \eqref{eq:uniform_stsp_moments} yields, simultaneously for all
$h\in\mathcal H_n$, the same bounds as in \eqref{eq:stsp-kernel-residual-C1}, with
$h\eta_n(h)$ replaced by $C\sigma\varepsilon_n$ and
$\eta_n(h)$ replaced by $C\varepsilon_n$. Since
$\varepsilon_n=o(\sigma)$, the contraction and implicit-function
argument leading to \eqref{eq:stsp-empirical-root-position} is uniform in $h$. Thus all local empirical
STSP branches exist simultaneously on this event and satisfy
\begin{equation}\label{eq:uniform_stsp_branches}
\sup_{h\in\mathcal H_n}
\|\widehat u_h^{\mathrm{STSP}}-u_h^{\mathrm{STSP}}\|_\infty
  \le C\sigma\varepsilon_n,
\qquad
\sup_{h\in\mathcal H_n}
\|\widehat u_h^{\mathrm{STSP}}-u_h^{\mathrm{STSP}}\|_{C^1(\mathcal M)}
  \le C\varepsilon_n.
\end{equation}

We next identify the oracle single-scale error. By Proposition~\ref{thm:stsp_bias_expansion},
uniformly in $h=t\sigma$, $t\in[a,b]$,
\[
  \left\|u_h^{\mathrm{STSP}}-
     \sigma^2\gamma(t)H\right\|_\infty\le C\sigma^4,
  \qquad
  \gamma(t)=\frac d2\left(1+t^2\frac{c_{1/t}}{a_{1/t}}\right).
\]
For completeness, the function $\gamma$ is strictly increasing.
Indeed, write $S=\|v_N\|^2$ under the measure $\nu_\lambda$
of Section~\ref{sec:stsp_validation} and set $r(\lambda)=c_\lambda/a_\lambda$.
The conditional second-moment calculation in the proof of
Proposition~\ref{prop:stsp-naive} gives
\[
 r(\lambda)=\frac1{d+14}
 \left(1-\frac{\mathbb E_\lambda[S^2]}
                {\mathbb E_\lambda[S]}\right).
\]
Let $\widetilde\nu_\lambda$ be the size-biased law with density
proportional to $S\,d\nu_\lambda$. Differentiating its Gaussian
tilt shows
\[
 r'(\lambda)=-\frac{\operatorname{Var}_{\widetilde\nu_\lambda}(S)}
 {(d+14)\lambda^3}<0,
 \qquad
 \gamma'(t)=\frac d2\{2t\,r(1/t)-r'(1/t)\}>0.
\]
In particular, $\gamma_*=\gamma(a)>0$, and the two distinct
scales $a_1<a_2$ yield distinct coefficients and fixed, bounded
cancellation weights.

For a normal graph contained in the reach tube,
\[
 d_H(\widehat{\mathcal M}_h^{\mathrm{STSP}},\mathcal M)
   =\|\widehat u_h^{\mathrm{STSP}}\|_\infty.
\]
Combining this identity with
\eqref{eq:uniform_stsp_branches} and the population expansion
gives
\[
 \sup_{h\in\mathcal H_n}
 \left|d_H(\widehat{\mathcal M}_h^{\mathrm{STSP}},\mathcal M)
        -\sigma^2\gamma(h/\sigma)\|H\|_\infty\right|
 \le C\{\sigma^4+\sigma\varepsilon_n\}.
\]
Taking the infimum over $h\in\mathcal H_n$ proves the first
assertion of the corollary. Moreover, $\|H\|_\infty>0$ on a
compact boundaryless positive-dimensional Euclidean submanifold:
otherwise all ambient coordinate functions would be harmonic on
each compact connected component, which is impossible for an
immersion. Finally, for fixed $h_k=a_k\sigma$, Theorem~\ref{thm:finite_sample_stsp} gives
\[
 d_H(\widehat{\mathcal M}^{\mathrm{STSP}},\mathcal M)
   =O_p(\sigma^4+\sigma\varepsilon_n).
\]
Dividing by the oracle error, which is asymptotic to the positive
constant $\gamma_*\|H\|_\infty\sigma^2$, yields
\[
 \frac{d_H(\widehat{\mathcal M}^{\mathrm{STSP}},\mathcal M)}
 {\inf_{h\in\mathcal H_n}
       d_H(\widehat{\mathcal M}_h^{\mathrm{STSP}},\mathcal M)}
 =O_p\left(\sigma^2+\frac{\varepsilon_n}{\sigma}\right)
 =o_p(1).
\]
This completes the proof.

\FloatBarrier
\section{Gaussian Density-Ridge Fitting}
\label{sec:gaussian_density_ridge}
\subsection{Population and Finite-Sample Guarantees}

We next consider Gaussian density-ridge fitting as a second, structurally
different verification of the factorized bias expansion in
Section~\ref{Multiscale Bias Cancellation}.
Unlike STSP, the fitted manifold is characterized through the gradient and
Hessian of a smoothed density
\citep{genovese2014nonparametric,chen2015asymptotic}.
Under the present Gaussian observation model, the population smoothing scale
can be combined explicitly with the noise level, leading to a particularly
simple normal-bias expansion. Let
\begin{equation}
\label{eq:density_ridge_population_density}
 p_{h,\sigma}(y)
 :=
 \mathbb E\{\phi_h(y-Y)\}
 =
 \int_{\mathcal M}\phi_{\rho_h}(y-z)\,d\omega(z),
 \qquad
 \rho_h^2=h^2+\sigma^2,
\end{equation}
where the second equality follows from Gaussian convolution.

Let
\(
 \lambda_1^{h,\sigma}(y)
 \geq\cdots\geq
 \lambda_D^{h,\sigma}(y)
\)
be the eigenvalues of \(\nabla^2p_{h,\sigma}(y)\), and let
\(Q_{h,\sigma}^{\mathrm{DR}}(y)\) be the spectral projector associated
with its \(D-d\) smallest eigenvalues.  The \(d\)-dimensional population
density ridge is
\begin{equation}
\label{eq:density_ridge_definition}
 \mathcal R_{h,\sigma}
 :=
 \left\{
 y\in\mathbb R^D:
 Q_{h,\sigma}^{\mathrm{DR}}(y)
 \nabla p_{h,\sigma}(y)=0,\quad
 \lambda_{d+1}^{h,\sigma}(y)<0
 \right\}.
\end{equation}

\begin{proposition}[Population expansion for Gaussian density ridges]
\label{prop:density_ridge_population_expansion}
Suppose the geometric conditions of Assumption~\ref{ass:geometry-observation} hold and
\(\mathcal M\) is \(C^6\).  Under the scale regime
\(\sigma\lesssim h\), there exist constants \(R,C,h_0>0\) such that,
for \(0<h<h_0\),
\(\mathcal M_{h,\sigma}^{\mathrm{DR}} := \mathcal R_{h,\sigma} \cap \operatorname{Tub}_{R\rho_h^2}(\mathcal M)\)
is a unique \(C^2\) normal graph over \(\mathcal M\):
\[
 \mathcal M_{h,\sigma}^{\mathrm{DR}}
 =
 \left\{
 X_{h,\sigma}^{\mathrm{DR}}(x):
 x\in\mathcal M
 \right\},
 \qquad
 X_{h,\sigma}^{\mathrm{DR}}(x)
 =
 x+u_{h,\sigma}^{\mathrm{DR}}(x),
 \quad
 u_{h,\sigma}^{\mathrm{DR}}(x)\in N_x\mathcal M.
\]
Uniformly over \(x\in\mathcal M\),
\begin{equation}
\label{eq:density_ridge_population_expansion}
 \left\|
 u_{h,\sigma}^{\mathrm{DR}}
 -
 \frac d2\rho_h^2 H
 \right\|_{C^2(\mathcal M)}
 \leq
 C\rho_h^4.
\end{equation}
Moreover,
\(
 \left\|
 u_{h,\sigma}^{\mathrm{DR}}
 \right\|_{C^2(\mathcal M)}
 \leq
 C\rho_h^2,
\)
and \(X_{h,\sigma}^{\mathrm{DR}}\) is a \(C^2\) embedding whose image is
\(C^2\)-diffeomorphic to \(\mathcal M\).
\end{proposition}

The tubular restriction above is used only to identify the ridge branch
associated with \(\mathcal M\); the density-ridge construction itself
depends only on \(p_{h,\sigma}\).  Since
\(\rho_h\asymp h\) under \(\sigma\lesssim h\),
Proposition~\ref{prop:density_ridge_population_expansion} verifies
Assumptions~\ref{ass:normal-section},
\ref{assump:normal_error_t}, and
\ref{ass:tangent-stability} with
\(b(x)=H(x), \quad c_{\mathrm{DR}}(h,\sigma) = \frac d2(h^2+\sigma^2)\).

For admissible scales \(h_1<\cdots<h_K\), choose bounded weights satisfying
\begin{equation}
\label{eq:density_ridge_multiscale_weights}
 \sum_{k=1}^K w_k=1,
 \qquad
 \sum_{k=1}^K
 w_k(h_k^2+\sigma^2)=0.
\end{equation}
Use \(\mathcal M_{h_1,\sigma}^{\mathrm{DR}}\) as the reference manifold,
and let
\(Z_{h_k}^{\mathrm{DR}}(x)\in
\mathcal M_{h_k,\sigma}^{\mathrm{DR}}\)
be the estimated-normal correspondence points constructed in
Section~\ref{sec:normal_guided_extrapolation}.
Define
\begin{equation}
\label{eq:density_ridge_multiscale_map}
 \widetilde X^{\mathrm{DR}}(x)
 =
 \sum_{k=1}^K
 w_k Z_{h_k}^{\mathrm{DR}}(x),
 \qquad
 \widetilde{\mathcal M}^{\mathrm{DR}}
 =
 \widetilde X^{\mathrm{DR}}(\mathcal M).
\end{equation}
Proposition~\ref{prop:density_ridge_population_expansion} and
Remark~\ref{rem:reference_normal_space} verify the normal-accuracy and
regularity conditions required by
Proposition~\ref{thm:normal_guided_manifold_recovery}.
\begin{equation}
\label{eq:density_ridge_multiscale_rate}
 \sup_{x\in\mathcal M}
 \left\|
 \widetilde X^{\mathrm{DR}}(x)-x
 \right\|
 \leq
 C\sum_{k=1}^K |w_k|h_k^4,
 \qquad
 d_H\!\left(
 \widetilde{\mathcal M}^{\mathrm{DR}},\mathcal M
 \right)
 \leq
 C\sum_{k=1}^K |w_k|h_k^4.
\end{equation}
For fixed \(K\), uniformly bounded weights, and sufficiently small
\(h_{\max}\), the corrected map is a \(C^1\) embedding and
\(\widetilde{\mathcal M}^{\mathrm{DR}}\) is \(C^1\)-diffeomorphic to
\(\mathcal M\).


\paragraph*{Finite-sample consistency of GBC-DR}
We next instantiate the finite-sample estimated-normal theory for Gaussian
 density-ridge fitting.

For \(y\in\mathbb R^D\), define the Gaussian kernel density estimator by
\(\widehat p_h(y)=\frac1n\sum_{i=1}^n\phi_h(y-Y_i)\).
Let \(\widehat\lambda_1^h(y)\geq\cdots\geq
\widehat\lambda_D^h(y)\) be the eigenvalues of
\(\nabla^2\widehat p_h(y)\). Wherever the empirical \(d\)-th eigengap is
positive, let \(\widehat Q_h^{\mathrm{DR}}(y)\) be the orthogonal projector
onto the eigenspace corresponding to the \(D-d\) smallest eigenvalues.
The population and empirical normalized density-ridge residuals are
\begin{equation}\label{eq:dr_normalized_residuals}
R_{h,\sigma}^{\mathrm{DR}}(y)
=-h^2Q_{h,\sigma}^{\mathrm{DR}}(y)
\frac{\nabla p_{h,\sigma}(y)}{p_{h,\sigma}(y)}, \quad
\widehat R_h^{\mathrm{DR}}(y)
=-h^2\widehat Q_h^{\mathrm{DR}}(y)
\frac{\nabla\widehat p_h(y)}{\widehat p_h(y)}.
\end{equation}
This normalization leaves the ridge equation unchanged and makes its normal
 derivative uniformly nondegenerate in the regime \(\sigma\lesssim h\).
The projectors in these residuals are those of the density Hessians, not
of the Hessians of the log densities.

On every open region where the empirical density and its \(d\)-th Hessian
eigengap are positive, the spectral projector and normalized residual are
smooth. For a fixed sufficiently large \(R\), put
\[
\mathcal T_h^{\mathrm{DR}}
=
\{x+v:x\in\mathcal M,\ v\in N_x\mathcal M,
\ \|v\|\leq R\rho_h^2\},\qquad
\rho_h^2=h^2+\sigma^2.
\]
The empirical fitted set is the local ridge zero set of the computable residual,
\[
\widehat{\mathcal M}_h^{\mathrm{DR}}
=
\{y\in\mathcal T_h^{\mathrm{DR}}:
\widehat R_h^{\mathrm{DR}}(y)=0,\ 
\widehat\lambda_{d+1}^h(y)<0\}.
\]
The tube identifies the branch associated with \(\mathcal M\) for the
analysis; neither the true normal coordinates nor \(\mathcal M\) are inputs
to the residual calculation. Recall
\(\eta_n(h)=\sqrt{\frac{\log n}{nh^d}}, q=D-d\).

\begin{proposition}[Empirical density-ridge branch existence and stability]
\label{lem:dr_empirical_rate}
Suppose Assumption~\ref{ass:geometry-observation} holds, \(\mathcal M\)
is \(C^6\), and the observations are i.i.d. under the model in~\eqref{eq:model}.
Assume \(0\leq\sigma\leq C_\sigma h\), where \(C_\sigma<\infty\) is fixed,
and
\(
\eta_n(h)=o(h), \text{equivalently},
\frac{nh^{d+2}}{\log n}\longrightarrow\infty.
\)
Choose \(R\) large enough to contain the population ridge branch in the
interior of \(\mathcal T_h^{\mathrm{DR}}\). For sufficiently small \(h\)
and large \(n\), with probability at least \(1-n^{-2}\), uniformly on
\(\mathcal T_h^{\mathrm{DR}}\),
\begin{equation}\label{eq:dr_empirical_eigengap}
\begin{aligned}
\widehat p_h(y)&\geq ch^{-q},\\
\widehat\lambda_{d+1}^h(y)&\leq-ch^{-q-2},\\
\widehat\lambda_d^h(y)-\widehat\lambda_{d+1}^h(y)
&\geq ch^{-q-2}.
\end{aligned}
\end{equation}
On the same event, the entire local empirical ridge is a unique \(C^2\)
normal graph
\[
\widehat{\mathcal M}_h^{\mathrm{DR}}
=\{\widehat X_h^{\mathrm{DR}}(x):x\in\mathcal M\},\qquad
\widehat X_h^{\mathrm{DR}}(x)=x+\widehat u_h^{\mathrm{DR}}(x),\quad
\widehat u_h^{\mathrm{DR}}(x)\in N_x\mathcal M,
\]
and \(\widehat X_h^{\mathrm{DR}}\) is a \(C^2\) embedding. Moreover,
\begin{equation}\label{eq:dr_empirical_c2}
\bigl\|\widehat u_h^{\mathrm{DR}}-u_{h,\sigma}^{\mathrm{DR}}\bigr\|_{C^j(\mathcal M)}
\leq Ch^{1-j}\eta_n(h),\qquad j=0,1,2.
\end{equation}
Consequently, the uniform empirical approximation condition holds with
\(
r_n^{\mathrm{DR}}(h)=Ch\eta_n(h),
\)
and the empirical tangent spaces satisfy
\[
\sup_{x\in\mathcal M}
\left\|
P_{T_{\widehat X_h^{\mathrm{DR}}(x)}
\widehat{\mathcal M}_h^{\mathrm{DR}}}
-P_{T_x\mathcal M}
\right\|_{\mathrm{op}}
\leq C\bigl\{h^2+\eta_n(h)\bigr\}.
\]
The constants are uniform over the stated noise--bandwidth regime.
\end{proposition}

For the estimated-normal correspondence, use the reference ridge's
geometric normal spaces at scale \(h_1\). Define the population and
empirical reference projector fields by
\[
\widetilde Q_1^{\mathrm{DR}}(x)
=P_{N_{X_{h_1,\sigma}^{\mathrm{DR}}(x)}
\mathcal M_{h_1,\sigma}^{\mathrm{DR}}},\qquad
\widehat Q_1^{\mathrm{DR}}(x)
=P_{N_{\widehat X_{h_1}^{\mathrm{DR}}(x)}
\widehat{\mathcal M}_{h_1}^{\mathrm{DR}}},
\]
and put
\(\widetilde N_{h_1}^{\mathrm{DR}}(x)
=\operatorname{Ran}\widetilde Q_1^{\mathrm{DR}}(x)\)
and
\(\widehat N_{h_1}^{\mathrm{DR}}(x)
=\operatorname{Ran}\widehat Q_1^{\mathrm{DR}}(x)\).
These geometric normal projectors need not coincide with the
 density-Hessian spectral projectors defining the ridge. They can be
computed from the empirical ridge residual: at any empirical ridge point
\(y\),
\begin{equation}\label{eq:dr_geometric_normals}
T_y\widehat{\mathcal M}_h^{\mathrm{DR}}
=\ker D_y\widehat R_h^{\mathrm{DR}}(y),\qquad
N_y\widehat{\mathcal M}_h^{\mathrm{DR}}
=\operatorname{Ran}\!\left\{
\bigl[D_y\widehat R_h^{\mathrm{DR}}(y)\bigr]^\top\right\},
\end{equation}
and this Jacobian has rank \(q\). Thus the geometric normal calculation
requires no access to \(\mathcal M\).
The population expansion and empirical branch stability give
\(\bigl\|\widetilde Q_1^{\mathrm{DR}} -P_{N\mathcal M}\bigr\|_{C^1}\leq Ch_1^2\),
where \(P_{N\mathcal M}\) denotes the field
\(x\mapsto P_{N_x\mathcal M}\). The differentiated residual and
empirical-branch estimates in the proof of
Proposition~\ref{lem:dr_empirical_rate} also give
\begin{equation}\label{eq:dr_reference_projector_bounds}
\begin{aligned}
\bigl\|\widehat Q_1^{\mathrm{DR}}-
\widetilde Q_1^{\mathrm{DR}}\bigr\|_\infty
&\leq C\eta_n(h_1),\\
\bigl\|\widehat Q_1^{\mathrm{DR}}-
\widetilde Q_1^{\mathrm{DR}}\bigr\|_{C^1}
&\leq C\frac{\eta_n(h_1)}{h_1}.
\end{aligned}
\end{equation}
Thus the empirical reference field is \(C^1\), its \(C^1\) norm is bounded
on the same high-probability event under \(\eta_n(h_1)=o(h_1)\), and the
finite-sample normal-estimation condition holds with
\(s_n^{\mathrm{DR}}(h_1)=C\eta_n(h_1)\).

\paragraph*{PCA-based correspondence in the experiments}
The numerical experiments instead use kernel-weighted PCA with
$K(t)=(1-t)_+^6$ for cross-scale correspondence. This changes only the
correspondence rule, not the underlying density-ridge estimator.
Let $Q_h^{\mathrm{PCA}}$ and $\widehat Q_h^{\mathrm{PCA}}$ be the
population and empirical normal projectors from the centered
kernel-weighted covariance in Lemma~\ref{lem:pca_normal_accuracy}, and set
\[
\begin{aligned}
\widetilde Q_1^{\mathrm{PCA}}(x)
&=Q_{h_1}^{\mathrm{PCA}}
  \bigl(X_{h_1,\sigma}^{\mathrm{DR}}(x)\bigr),\\
\widehat Q_1^{\mathrm{PCA}}(x)
&=\widehat Q_{h_1}^{\mathrm{PCA}}
  \bigl(\widehat X_{h_1}^{\mathrm{DR}}(x)\bigr).
\end{aligned}
\]
Suppose additionally that $\sigma/h_1$ lies in a fixed compact subset
of $(0,\infty)$. The population reference field satisfies
$\|\widetilde Q_1^{\mathrm{PCA}}-P_{N\mathcal M}\|_{C^1}\leq Ch_1^2$,
and the differentiated moment and ridge-branch bounds give
\[
\begin{aligned}
\|\widehat Q_1^{\mathrm{PCA}}
  -\widetilde Q_1^{\mathrm{PCA}}\|_\infty
&\leq C\eta_n(h_1),\\
\|\widehat Q_1^{\mathrm{PCA}}
  -\widetilde Q_1^{\mathrm{PCA}}\|_{C^1}
&\leq C\eta_n(h_1)/h_1.
\end{aligned}
\]
These bounds hold on the joint ridge and PCA concentration event,
whose probability is at least $1-(K+1)n^{-2}$; their proof is given
at the end of Appendix~\ref{app:dr_empirical_rate_proof}.
Under $\max_k\eta_n(h_k)/h_k\to0$, the empirical PCA reference field
has uniformly bounded $C^1$ norm, and the corresponding local
intersections are unique, transverse, and $C^1$ in the reference point.
Appendix~\ref{app:dr_finite_sample_proof} shows that the recovery,
embedding, and numerical-error conclusions of
Theorem~\ref{thm:finite_sample_dr}, and the fourth-order guarantee of
Corollary~\ref{cor:dr_fourth_order_sample_size}, extend to this PCA
construction on the same joint event, under the additional
noise-to-bandwidth restriction above.

\medskip
\noindent
For the geometric-normal construction of
Theorem~\ref{thm:finite_sample_dr}, use
$\widehat{\mathcal M}_{h_1}^{\mathrm{DR}}$
as the reference estimator and
$\widehat N_{h_1}^{\mathrm{DR}}$
as its geometric normal space. Let
\(
\widehat Z_{h_k}^{\mathrm{DR}}(x)
\in\widehat{\mathcal M}_{h_k}^{\mathrm{DR}}
\)
denote the empirical correspondence point, and let
\(
Z_{h_k}^{\mathrm{DR}}(x)
\in\mathcal M_{h_k,\sigma}^{\mathrm{DR}}
\)
denote its population counterpart.

For fixed $K$, if
$\max_k\eta_n(h_k)/h_k\to0$, the empirical branch
and geometric-normal projector bounds imply that
these local intersections exist uniquely, are
transverse, and depend $C^1$-smoothly on the
reference point, with probability at least
$1-Kn^{-2}$. The argument is given in the proof
of Theorem~\ref{thm:finite_sample_dr}.
The finite-sample correspondence bound gives,
on the same event,
\[
\sup_{x\in\mathcal M}
\left\|
\widehat Z_{h_k}^{\mathrm{DR}}(x)
-Z_{h_k}^{\mathrm{DR}}(x)
\right\|
\le C\left[
h_1\eta_n(h_1)
+h_k\eta_n(h_k)
+h_k^2\eta_n(h_1)
\right].
\]
The same bound holds for PCA-based correspondence on the joint
ridge and PCA concentration event described above, with both
intersection points defined using the PCA reference fields.

Define the finite-sample GBC-DR map by
\[
\widehat X^{\mathrm{DR}}(x)
=\sum_{k=1}^K w_k\widehat Z_{h_k}^{\mathrm{DR}}(x),
\qquad
\widehat{\mathcal M}^{\mathrm{DR}}
=\widehat X^{\mathrm{DR}}(\mathcal M).
\]
This formal estimator is the continuous image of the smooth empirical
ridge branches, not a finite point cloud. On the event where these
branches and correspondences are defined, its reference-branch
parameterization is
\[
\widehat F^{\mathrm{DR}}
=\widehat X^{\mathrm{DR}}\circ
\bigl(\widehat X_{h_1}^{\mathrm{DR}}\bigr)^{-1},
\qquad
\widehat{\mathcal M}^{\mathrm{DR}}
=\widehat F^{\mathrm{DR}}
\bigl(\widehat{\mathcal M}_{h_1}^{\mathrm{DR}}\bigr).
\]

\begin{theorem}[Finite-sample recovery of GBC-DR]
\label{thm:finite_sample_dr}
Suppose the conditions of Proposition~\ref{lem:dr_empirical_rate} hold at
 a fixed number \(K\) of fitting scales \(h_1<\cdots<h_K\). Assume
\(\max_k\eta_n(h_k)/h_k\to0\), and let the deterministic extrapolation
weights satisfy \(\sum_{k=1}^K|w_k|\leq W_0<\infty\) and
\[
\sum_{k=1}^K w_k=1,\qquad
\sum_{k=1}^K w_k(h_k^2+\sigma^2)=0.
\]
Then, for sufficiently small \(h_{\max}=\max_k h_k\) and sufficiently
large \(n\), with probability at least \(1-Kn^{-2}\),
\begin{equation}\label{eq:dr_multiscale_finite_sample}
\begin{aligned}
d_H\!\left(\widehat{\mathcal M}^{\mathrm{DR}},\mathcal M\right)
&\leq\sup_{x\in\mathcal M}
\bigl\|\widehat X^{\mathrm{DR}}(x)-x\bigr\|\\
&\leq C\sum_{k=1}^K|w_k|
\left[
h_k^4+h_1\eta_n(h_1)+h_k\eta_n(h_k)+h_k^2\eta_n(h_1)
\right].
\end{aligned}
\end{equation}
For sufficiently small \(h_{\max}\), the empirical local ridge branches
and estimated-normal correspondence points are well defined on the same
event, and \(\widehat{\mathcal M}^{\mathrm{DR}}\) is a \(C^1\) embedded
manifold \(C^1\)-diffeomorphic to \(\mathcal M\).

For the numerical implementation, suppose that the supplied reference
branch is the branch \(\widehat{\mathcal M}_{h_1}^{\mathrm{DR}}\)
identified above. Let
\(\mathcal Z_N=\{z_1,\ldots,z_N\}\subset
\widehat{\mathcal M}_{h_1}^{\mathrm{DR}}\)
be a nonempty set of retained query points, with fill distance
\(\delta_N =\sup_{z\in\widehat{\mathcal M}_{h_1}^{\mathrm{DR}}} \min_{1\leq j\leq N}\|z-z_j\|\).
Let \(q=D-d\), and choose an orthonormal basis
\(\widehat E(z)\in\mathbb R^{D\times q}\) of the geometric reference
normal space at each query point. For \(k=2,\ldots,K\), define
\(g_{k,z}(t) =\widehat E(z)^\top\widehat R_{h_k}^{\mathrm{DR}} \bigl(z+\widehat E(z)t\bigr)\).
Fix a sufficiently large constant \(C_{\mathrm{loc}}\), independent of
\(n\) and the scales. If numerical solutions \(\widehat t_k(z)\) satisfy
\[
\|\widehat t_k(z)\|\leq C_{\mathrm{loc}}h_k^2,\qquad
\|g_{k,z}(\widehat t_k(z))\|\leq\varepsilon_{n,k}
\]
for all \(z\in\mathcal Z_N\), for deterministic tolerances
\(\varepsilon_{n,k}\geq0\), put \(\widehat t_1(z)=0\),
\(\varepsilon_{n,1}=0\), and define the finite numerical output
\[
\widehat{\mathcal P}_N^{\mathrm{DR,num}}
=\left\{
z_j+\widehat E(z_j)\sum_{k=1}^K w_k\widehat t_k(z_j):
1\leq j\leq N
\right\}.
\]
Then, on the same high-probability event,
\begin{equation}\label{eq:dr_numerical_error}
\begin{aligned}
d_H\!\left(\widehat{\mathcal P}_N^{\mathrm{DR,num}},\mathcal M\right)
&\leq
d_H\!\left(\widehat{\mathcal M}^{\mathrm{DR}},\mathcal M\right)
+C\sum_{k=1}^K|w_k|\varepsilon_{n,k}+C\delta_N.
\end{aligned}
\end{equation}
For comparable scales \(h_k\asymp h\), fixed \(K\), and bounded weights,
\eqref{eq:dr_multiscale_finite_sample} simplifies to
\[
d_H\!\left(\widehat{\mathcal M}^{\mathrm{DR}},\mathcal M\right)
=O_p\!\left(h^4+h\sqrt{\frac{\log n}{nh^d}}\right).
\]
Thus \(\eta_n(h)=O(h^3)\), equivalently
\(nh^{d+6}\gtrsim\log n\), is sufficient for fourth-order stochastic
accuracy. For the finite numerical output, it also suffices that
\(\max_k\varepsilon_{n,k}=O(h^4)\) and \(\delta_N=O(h^4)\).
The embedding conclusion above concerns the continuous exact-intersection
estimator, not the finite point cloud; the residual tolerance alone does
not imply regularity of a numerical interpolation of that cloud.
\end{theorem}

\begin{corollary}[A sufficient sample size for fourth-order density-ridge accuracy]
\label{cor:dr_fourth_order_sample_size}
Suppose the assumptions of Theorem~\ref{thm:finite_sample_dr}, other
than its sample-size condition, hold. Let
\(h_k=a_kh\), where \(0<a_1<\cdots<a_K<\infty\) are fixed and
\(h\to0\), with \(0\leq\sigma\leq C_\sigma h_1\).
There exist constants \(C_0,C_1>0\) such that, if
\[
 n h^{d+6}\geq C_0\log n,
\]
then, for sufficiently small \(h\) and large \(n\), with probability
at least \(1-Kn^{-2}\),
\[
 d_H\!\left(\widehat{\mathcal M}^{\mathrm{DR}},\mathcal M\right)
 \leq C_1h^4,
\]
and \(\widehat{\mathcal M}^{\mathrm{DR}}\) is a \(C^1\) embedded
manifold \(C^1\)-diffeomorphic to \(\mathcal M\).
In particular, when \(h\asymp\sigma\to0\), there exist constants
\(C_0',C_1'>0\) such that the sufficient sample-size condition
\[
 n\geq C_0'\sigma^{-(d+6)}\log(1/\sigma)
\]
yields
\(d_H(\widehat{\mathcal M}^{\mathrm{DR}},\mathcal M)
\leq C_1'\sigma^4\) on the same event.
Under the numerical assumptions of
Theorem~\ref{thm:finite_sample_dr}, if also
\(\max_k\varepsilon_{n,k}=O(h^4)\) and \(\delta_N=O(h^4)\), then
\(d_H\!\left(\widehat{\mathcal P}_N^{\mathrm{DR,num}},\mathcal M\right) \leq Ch^4\)
on that event.
\end{corollary}

\paragraph*{Estimated noise level}
For fixed deterministic bandwidths, select one point
\(\widehat y_1\in\widehat{\mathcal M}_{h_1}^{\mathrm{DR}}\) on the
reference ridge. Using the density-Hessian eigenvalues defined above, put
\(q=D-d\), \(\gamma_*=(1+C_\sigma^2)^{-1}\), and define
\[
\widehat\gamma_1
=\Pi_{[\gamma_*,1]}\!\left(
-\frac{h_1^2}{q\widehat p_{h_1}(\widehat y_1)}
\sum_{j=d+1}^D\widehat\lambda_j^{h_1}(\widehat y_1)
\right),
\qquad
\widehat\sigma^2=h_1^2(\widehat\gamma_1^{-1}-1).
\]
Here \(C_\sigma\) is the fixed upper bound in
Proposition~\ref{lem:dr_empirical_rate}, and \(\Pi_I\) denotes interval
projection. Under that proposition's conditions,
Proposition~\ref{prop:dr_noise_estimation} in
Appendix~\ref{app:dr_noise_estimation} establishes
\[
|\widehat\sigma^2-\sigma^2|
=O_p\!\left(h_1^4+h_1^2\eta_n(h_1)\right),
\]
uniformly for \(0\leq\sigma\leq C_\sigma h_1\).
Use this single variance estimate in spatially constant weights satisfying
\(\sum_k\widehat w_k=1\) and
\(\sum_k\widehat w_k(h_k^2+\widehat\sigma^2)=0\).
For fixed comparable scales and uniformly bounded weights,
Theorem~\ref{thm:finite_sample_dr} retains its error rate and embedding
conclusion, and Corollary~\ref{cor:dr_fourth_order_sample_size}
retains its fourth-order guarantee.

\begin{remark}[Extrapolation to zero total smoothing variance]
\label{rem:dr_deconvolution}
For Gaussian density ridges, the effective smoothing variance is
\(t_k=h_k^2+\sigma^2\), and
Proposition~\ref{prop:density_ridge_population_expansion} gives
\(X_{h_k,\sigma}^{\mathrm{DR}}(x) = x+\frac d2 t_k H(x)+O(t_k^2)\).
The constraints \(\sum_k w_k=1\) and \(\sum_k w_k t_k=0\)
therefore extrapolate the matched ridge locations to \(t=0\),
canceling the leading geometric effect of both observation
noise and kernel smoothing. Equivalently,
\(\sum_k w_k h_k^2=-\sigma^2\).
By contrast, extrapolation based only on
\(\sum_k w_k h_k^2=0\) leaves the leading term
\((d/2)\sigma^2H(x)\).

This provides a finite-order geometric interpretation of
partial deconvolution: the correction acts on matched ridge
locations, rather than inverting the convolution operator
to recover the entire latent distribution.
Under the conditions of Theorem~\ref{thm:finite_sample_dr},
with \(h_k\asymp\sigma\), fixed \(K\), and uniformly bounded
weights,
\[
d_H(\widehat{\mathcal M}^{\mathrm{DR}},\mathcal M)
=
O_p\!\left(
\sigma^4+\sigma\sqrt{\frac{\log n}{n\sigma^d}}
\right).
\]
Thus
\(n\gtrsim\sigma^{-(d+6)}\log(1/\sigma)\)
suffices for \(O_p(\sigma^4)\) accuracy in the small-noise
regime, improving on the second-order approximation of
the single-scale ridge.

This conclusion does not contradict the logarithmic minimax
lower bound for singular deconvolution established by
\citet{genovese2012manifold}. That bound concerns a fixed
Gaussian noise distribution, whereas our result concerns
the joint limit \(\sigma\to0\) and \(h_k\asymp\sigma\).
For fixed positive \(\sigma\), the present bound retains
an \(O(\sigma^4)\) approximation term and does not establish
consistency merely by increasing \(n\).
\end{remark}

\subsection{Proofs and Noise-Variance Estimation}

\label{sec:app_density_ridge}
\subsubsection{Proof of
Proposition~\ref{prop:density_ridge_population_expansion}}
\begin{proof}
Put
\(\rho=\rho_h=(h^2+\sigma^2)^{1/2}, p_\rho=p_{h,\sigma}, s_\rho=\nabla\log p_\rho\).
For \(y\) in an \(O(\rho^2)\) tubular neighborhood of \(\mathcal M\),
write
\(y=x+\rho^2\beta, x=\pi_{\mathcal M}(y), \beta\in N_x\mathcal M, \|\beta\|\leq R\).

Let \(\mathbb E_{\rho,y}\) denote expectation with respect to the
probability measure on \(\mathcal M\) whose density relative to
\(\omega\) is proportional to \(\phi_\rho(y-z)\).  Gaussian
differentiation gives the exact identities
\begin{equation}
\label{eq:density_ridge_score_covariance_identity}
 \begin{split}
 s_\rho(y)
 &=
 \rho^{-2}
 \left\{
 \mathbb E_{\rho,y}[Z]-y
 \right\},\\
 \nabla^2\log p_\rho(y)
 &=
 -\rho^{-2}I_D
 +
 \rho^{-4}
 \operatorname{Cov}_{\rho,y}(Z).
 \end{split}
\end{equation}

Use the same reach-controlled normal graph and Gaussian localization as
in the proof of Lemma~\ref{lem:pca_normal_accuracy}.  In coordinates adapted to
\(T_x\mathcal M\oplus N_x\mathcal M\), the local graph has the expansion
\begin{equation}
\label{eq:density_ridge_local_graph}
 \Psi_x(\rho t)
 =
 x+\rho t
 +\frac{\rho^2}{2}\mathrm{II}_x(t,t)
 +\frac{\rho^3}{6}C_x[t,t,t]
 +O_{C^2}\!\left(
 \rho^4(1+\|t\|^4)
 \right),
\end{equation}
and its volume density equals
\[
 1+\rho^2J_{2,x}(t)
 +O_{C^2}\!\left(
 \rho^3(1+\|t\|^3)
 \right).
\]
Contributions outside the fixed local chart are exponentially small.
After the change of variables \(z=\Psi_x(\rho t)\), all remainders are
integrable against a standard Gaussian envelope.

Reflection symmetry in \(t\) eliminates the odd terms.  Using
\[
 \mathbb E[t_it_j]=\delta_{ij},
 \qquad
 \sum_{i=1}^d
 \mathrm{II}_x(e_i,e_i)
 =
 dH(x),
\]
we obtain, uniformly for \(x\in\mathcal M\) and bounded \(\beta\),
\begin{equation}
\label{eq:density_ridge_local_moments}
 \begin{split}
 P_x\left\{
 \mathbb E_{\rho,y}[Z]-y
 \right\}
 &=
 O_{C^2}(\rho^4),\\
 Q_x\left\{
 \mathbb E_{\rho,y}[Z]-y
 \right\}
 &=
 -\rho^2
 \left\{
 \beta-\frac d2H(x)
 \right\}
 +O_{C^2}(\rho^4),\\
 \operatorname{Cov}_{\rho,y}(Z)
 &=
 \rho^2P_x+O_{C^2}(\rho^4).
 \end{split}
\end{equation}
The cubic graph term is odd under \(t\mapsto-t\), while the volume
correction begins at second order; consequently no third-order normal term
appears.  The same cancellations remain valid after taking up to two
derivatives in the local normal-bundle variables \((x,\beta)\).

Combining
\eqref{eq:density_ridge_score_covariance_identity} and
\eqref{eq:density_ridge_local_moments} yields
\begin{equation}
\label{eq:density_ridge_score_expansion}
 \begin{split}
 P_xs_\rho(x+\rho^2\beta)
 &=
 O_{C^2}(\rho^2),\\
 Q_xs_\rho(x+\rho^2\beta)
 &=
 -\beta+\frac d2H(x)
 +O_{C^2}(\rho^2),\\
 \rho^2\nabla^2\log p_\rho(x+\rho^2\beta)
 &=
 -Q_x+O_{C^2}(\rho^2).
 \end{split}
\end{equation}
Since
\(\frac{\nabla^2p_\rho}{p_\rho} = \nabla^2\log p_\rho+s_\rho s_\rho^\top\)
and \(s_\rho=O_{C^2}(1)\) in the region under consideration,
\begin{equation}
\label{eq:density_ridge_hessian_expansion}
 \frac{\rho^2}{p_\rho(y)}
 \nabla^2p_\rho(y)
 =
 -Q_x+O_{C^2}(\rho^2).
\end{equation}

Since \(\rho^2/p_\rho(y)>0\), this rescaling preserves eigenspaces and
eigenvalue signs. Smooth spectral-projector calculus applied to
\eqref{eq:density_ridge_hessian_expansion} therefore gives
\begin{equation}
\label{eq:density_ridge_projector_expansion}
 Q_{h,\sigma}^{\mathrm{DR}}
 (x+\rho^2\beta)
 =
 Q_x+O_{C^2}(\rho^2).
\end{equation}
Moreover, for sufficiently small \(\rho\),
\[
 \lambda_{d+1}^{h,\sigma}(y)
 \leq
 -c\frac{p_\rho(y)}{\rho^2}<0,
 \qquad
 \lambda_d^{h,\sigma}(y)
 -
 \lambda_{d+1}^{h,\sigma}(y)
 \geq
 c\frac{p_\rho(y)}{\rho^2}.
\]
Thus the ridge projector is well defined throughout the stated tubular
neighborhood.

Define the scaled normal ridge equation
\begin{equation}
\label{eq:density_ridge_scaled_equation}
 F_\rho(x,\beta)
 :=
 Q_x
 Q_{h,\sigma}^{\mathrm{DR}}
 (x+\rho^2\beta)
 s_\rho(x+\rho^2\beta).
\end{equation}
Equations
\eqref{eq:density_ridge_score_expansion} and
\eqref{eq:density_ridge_projector_expansion} imply
\begin{equation}
\label{eq:density_ridge_scaled_equation_expansion}
 F_\rho(x,\beta)
 =
 -\beta+\frac d2H(x)
 +O_{C^2}(\rho^2),
 \qquad
 D_\beta F_\rho(x,\beta)
 =
 -I_{N_x\mathcal M}+O_{C^1}(\rho^2).
\end{equation}
The limiting equation has the unique solution
\(
 \beta_0(x)
 =
 \frac d2H(x),
\)
and its derivative in \(\beta\) is
\(-I_{N_x\mathcal M}\).  The parameterized implicit-function theorem
therefore gives a unique \(C^2\) section
\(\beta_\rho\) satisfying
\[
 F_\rho(x,\beta_\rho(x))=0,
 \qquad
 \left\|
 \beta_\rho-\frac d2H
 \right\|_{C^2(\mathcal M)}
 \leq
 C\rho^2.
\]
Choosing \(R\) larger than
\(\sup_x\|\beta_0(x)\|+1\), the second estimate in
\eqref{eq:density_ridge_scaled_equation_expansion} also excludes any
additional zero with \(\|\beta\|\leq R\).

To identify the zeros with the density ridge, use
\eqref{eq:density_ridge_projector_expansion}: the restriction of
\(Q_x\) to
\(\operatorname{im}Q_{h,\sigma}^{\mathrm{DR}}(y)\)
is invertible for sufficiently small \(\rho\).  Since
\(Q_{h,\sigma}^{\mathrm{DR}}(y)s_\rho(y) \in \operatorname{im}Q_{h,\sigma}^{\mathrm{DR}}(y)\),
we have
\[
 F_\rho(x,\beta)=0
 \quad\Longleftrightarrow\quad
 Q_{h,\sigma}^{\mathrm{DR}}(y)s_\rho(y)=0.
\]
Because \(p_\rho(y)>0\), this is equivalent to
\(Q_{h,\sigma}^{\mathrm{DR}}(y) \nabla p_\rho(y)=0\).
The eigenvalue condition has already been verified, so the zeros of
\(F_\rho\) are exactly the points of
\(\mathcal R_{h,\sigma}\) in
\(\operatorname{Tub}_{R\rho^2}(\mathcal M)\).

Set
\[
 u_{h,\sigma}^{\mathrm{DR}}(x)
 =
 \rho^2\beta_\rho(x),
 \qquad
 X_{h,\sigma}^{\mathrm{DR}}(x)
 =
 x+u_{h,\sigma}^{\mathrm{DR}}(x).
\]
Then
\[
 \left\|
 u_{h,\sigma}^{\mathrm{DR}}
 -
 \frac d2\rho^2H
 \right\|_{C^2}
 \leq
 C\rho^4,
 \qquad
 \left\|
 u_{h,\sigma}^{\mathrm{DR}}
 \right\|_{C^2}
 \leq
 C\rho^2.
\]
The tubular-neighborhood argument used in the proof of
Proposition~\ref{thm:manifold_recovery_1}, therefore shows that
\(X_{h,\sigma}^{\mathrm{DR}}\) is a \(C^2\) embedding and that its image
is \(C^2\)-diffeomorphic to \(\mathcal M\).  Differentiating the normal
graph also gives
\[
 \sup_{x\in\mathcal M}
 \left\|
 P_{T_{X_{h,\sigma}^{\mathrm{DR}}(x)}
       \mathcal M_{h,\sigma}^{\mathrm{DR}}}
 -
 P_x
 \right\|_{\mathrm{op}}
 \leq
 C\rho^2.
\]
Finally, \(\rho^2=h^2+\sigma^2\asymp h^2\) under
\(\sigma\lesssim h\), so the remainder \(O_{C^2}(\rho^4)\) is
\(O_{C^2}(h^4)\).  This completes the proof.
\end{proof}


\subsubsection{Proof of Proposition~\ref{lem:dr_empirical_rate}}
\label{app:dr_empirical_rate_proof}
\begin{proof}
Abbreviate
\[
 \rho=\rho_h,\quad p=p_{h,\sigma},\quad
 Q_h=Q_{h,\sigma}^{\mathrm{DR}},\quad
 \widehat Q_h=\widehat Q_h^{\mathrm{DR}},\quad
 R_h=R_{h,\sigma}^{\mathrm{DR}},\quad
 \widehat R_h=\widehat R_h^{\mathrm{DR}},\quad
 \eta=\eta_n(h).
\]
All estimates are uniform on \(\mathcal T_h^{\mathrm{DR}}\) and on any
fixed enlargement of that tube. We choose the enlargement large enough
to include the local search regions used below. Derivatives \(D_y^j\)
are ambient derivatives, as distinguished from derivatives in the scaled
normal-bundle coordinates used for the population expansion.

\paragraph*{Uniform differentiated density bounds}
For every fixed \(a>0\), Gaussian integration conditional on \(X=z\)
gives
\[
 \mathbb E\!\left[
 e^{-a\|y-Y\|^2/h^2}\mid X=z\right]
 =\left(\frac{h^2}{h^2+2a\sigma^2}\right)^{D/2}
   \exp\!\left(-\frac{a\|y-z\|^2}{h^2+2a\sigma^2}\right).
\]
Since \(\sigma\leq C_\sigma h\), integration over the compact manifold
and its uniformly controlled graph charts yields
\begin{equation}
\label{eq:dr_gaussian_local_mass}
 \sup_{y\in\mathbb R^D}
 \mathbb E e^{-a\|y-Y\|^2/h^2}\leq C_a h^d.
\end{equation}
Indeed, the manifold volume in any ambient ball of sufficiently small
radius \(r\) is at most \(Cr^d\), uniformly in its center; integration
of the Gaussian radial envelope gives the displayed bound. This
argument does not require a positive lower bound on \(\sigma/h\).
Gaussian convolution with the latent manifold measure also gives
\begin{equation}
\label{eq:dr_population_density_size}
 ch^{-q}\leq p(y)\leq Ch^{-q}
 \quad\text{on the enlarged tube},\qquad q=D-d.
\end{equation}
For the lower bound, integrate over a tangent-coordinate ball of radius
proportional to \(\rho\) centered at \(\pi_{\mathcal M}(y)\), using
\(d(y,\mathcal M)=O(\rho^2)\) and \(\rho\asymp h\).

For a component of the order-\(m\) derivative tensor, put
\(f_{m,y}(u)=D_y^m\phi_h(y-u)\). Gaussian derivative formulas imply,
for \(0\leq m\leq5\),
\(|f_{m,y}(u)|\leq C_mh^{-D-m} e^{-c_m\|y-u\|^2/h^2}\).
Consequently, for \(0\leq m\leq4\),
\[
 \|f_{m,y}\|_\infty\leq Ch^{-D-m},\qquad
 \mathbb E|f_{m,y}(Y)|^2\leq Ch^{d-2D-2m},\qquad
 \|D_yf_{m,y}\|_\infty\leq Ch^{-D-m-1}.
\]
Here and below, bounds on tensor components imply bounds in any fixed
tensor norm, since \(D\) is fixed.

Cover the enlarged deterministic tube by a Euclidean net of mesh
\(n^{-3}h^{d+1}\). The sample-size condition implies
\(\log(1/h)=O(\log n)\), so the number of net points is polynomial in
\(n\). Bernstein's inequality, a union bound over the net and derivative
components, and the preceding Lipschitz bound give an event
\(\mathcal E_h\), of probability at least \(1-n^{-2}\), on which
\[
 \begin{split}
 \sup_y\|D_y^m(\widehat p_h-p)(y)\|
 &\leq C\left\{
 \sqrt{\frac{h^{d-2D-2m}\log n}{n}}
          +\frac{h^{-D-m}\log n}{n}\right\}\\
 &\leq Ch^{-q-m}\eta,\qquad 0\leq m\leq4.
 \end{split}
\]
The off-net error is at most \(Cn^{-3}h^{-q-m}\), and is absorbed
in this bound. The Bernstein envelope term equals
\(h^{-q-m}\eta^2\), which is also absorbed because \(\eta=o(1)\).
Thus, on this one event,
\begin{equation}
\label{eq:dr_density_derivative_concentration}
 \sup_y\|D_y^m(\widehat p_h-p)(y)\|
       \leq Ch^{-q-m}\eta,\qquad 0\leq m\leq4.
\end{equation}
The deterministic derivative bounds, including the additional order
needed later, are
\begin{equation}
\label{eq:dr_population_density_derivatives}
 \sup_y\|D_y^m p(y)\|\leq Ch^{-q-m},\qquad 0\leq m\leq5.
\end{equation}
They follow directly from \eqref{eq:dr_gaussian_local_mass} and the
Gaussian derivative envelope. Equations
\eqref{eq:dr_population_density_size} and
\eqref{eq:dr_density_derivative_concentration} imply
\(\widehat p_h\geq ch^{-q}\) on \(\mathcal E_h\). We work on this
event for the rest of the proof.

\paragraph*{Normalized Hessians, projectors, and residual derivatives}
Define
\[
 A_h=h^2\frac{\nabla^2p}{p},\qquad
 \widehat A_h=h^2\frac{\nabla^2\widehat p_h}{\widehat p_h},\qquad
 v_h=-h^2\frac{\nabla p}{p},\qquad
 \widehat v_h=-h^2\frac{\nabla\widehat p_h}{\widehat p_h}.
\]
Applying the quotient and product rules to
\eqref{eq:dr_density_derivative_concentration} gives
\begin{equation}
\label{eq:dr_normalized_derivative_concentration}
 \begin{split}
 \sup_y\|D_y^j(\widehat A_h-A_h)(y)\|
       &\leq Ch^{-j}\eta,\\
 \sup_y\|D_y^j(\widehat v_h-v_h)(y)\|
       &\leq Ch^{1-j}\eta,\qquad j=0,1,2.
 \end{split}
\end{equation}
For example, \(D_y^j(p^{-1})=O(h^{q-j})\); its empirical perturbation
has the additional factor \(\eta\). Each term in the differentiated
quotients therefore has the indicated scale. Terms with two or more
empirical perturbations are absorbed using \(\eta=o(1)\).

The population Hessian expansion
\eqref{eq:density_ridge_hessian_expansion} implies
\begin{equation}
\label{eq:dr_normalized_hessian_gap}
 A_h(x+\rho^2\beta)
 =-\gamma_hQ_x+O_{C^2}(h^2),\qquad
 \gamma_h=\frac{h^2}{\rho^2},\qquad
 \gamma_*:=\frac1{1+C_\sigma^2}\leq\gamma_h\leq1.
\end{equation}
Thus the lowest \(q\) eigenvalues of \(A_h\) are uniformly negative
and separated from the other \(d\) eigenvalues. Weyl's inequality and
\eqref{eq:dr_normalized_derivative_concentration} give the same properties
for \(\widehat A_h\). Multiplication by
\(\widehat p_h/h^2>0\) proves both empirical eigenvalue bounds in
\eqref{eq:dr_empirical_eigengap}.

Let \(\mathcal Q\) be the spectral-projector map onto these lowest
\(q\) eigenvalues. Its matrix derivatives are bounded on a neighborhood
of the bounded matrix family with this uniform gap. In particular,
\[
 \begin{split}
 D_y(\widehat Q_h-Q_h)
 ={}&D\mathcal Q(\widehat A_h)
                       [D_y(\widehat A_h-A_h)]\\
    &+\{D\mathcal Q(\widehat A_h)-D\mathcal Q(A_h)\}[D_yA_h].
 \end{split}
\]
Differentiating once more, and using
\(D_y^j A_h=O(h^{-j})\), yields
\begin{equation}
\label{eq:dr_projector_derivative_concentration}
 \sup_y\|D_y^j(\widehat Q_h-Q_h)(y)\|
       \leq Ch^{-j}\eta,\qquad j=0,1,2.
\end{equation}
The deterministic bounds from
\eqref{eq:dr_population_density_derivatives} similarly give
\(\sup_y\|D_y^jQ_h(y)\|\leq Ch^{-j},\qquad j=0,1,2,3\).
These bounds use only separation between the two eigenvalue clusters;
individual eigenvalues within either cluster may coincide.

By \eqref{eq:density_ridge_score_expansion},
\(\|v_h\|_\infty\leq Ch^2\) on the tube. The quotient rules also give
\(\|D_y^jv_h\|_\infty\leq Ch^{1-j},\qquad j=1,2,3\).
Use the exact identity
\[
 \widehat R_h-R_h
 =(\widehat Q_h-Q_h)v_h+\widehat Q_h(\widehat v_h-v_h).
\]
Leibniz's rule, \eqref{eq:dr_normalized_derivative_concentration}, and
\eqref{eq:dr_projector_derivative_concentration} now imply
\begin{equation}
\label{eq:dr_residual_derivative_concentration}
 \sup_y\|D_y^j(\widehat R_h-R_h)(y)\|
       \leq Ch^{1-j}\eta,\qquad j=0,1,2.
\end{equation}
For the deterministic residual, the sharper bound
\(v_h=O(h^2)\) and the preceding derivative estimates give
\begin{equation}
\label{eq:dr_population_residual_derivatives}
 \|D_yR_h\|_\infty\leq C,\qquad
 \|D_y^2R_h\|_\infty\leq C/h,\qquad
 \|D_y^3R_h\|_\infty\leq C/h^2.
\end{equation}
In particular, the empirical regularity used next has been obtained from
differentiated density bounds, without assuming an empirical ridge branch.

\paragraph*{Existence, uniqueness, and gluing of the empirical branch}
Work in a fixed local chart on \(\mathcal M\), with a smooth orthonormal
normal frame \(E(x)\in\mathbb R^{D\times q}\). For
\(z\in\mathbb R^q\), \(\|z\|\leq R\rho^2\), define
\[
 f_h(x,z)=E(x)^\top R_h(x+E(x)z),\qquad
 \widehat f_h(x,z)=E(x)^\top\widehat R_h(x+E(x)z).
\]
Let \(z_h(x)=E(x)^\top u_{h,\sigma}^{\mathrm{DR}}(x)\) be the population
zero. Increase \(R\), if necessary, so that
\(\sup_x\|z_h(x)\|\leq R\rho^2/2\). The scaled population normal
equation \eqref{eq:density_ridge_scaled_equation_expansion}, multiplied
by \(-h^2\) and written with \(z=\rho^2E(x)^\top\beta\), gives
\begin{equation}
\label{eq:dr_normal_jacobian}
 \sup_{x,z}\|D_zf_h(x,z)-\gamma_h I_q\|\leq Ch^2.
\end{equation}
The fixed chart and frame have bounded derivatives, so
\eqref{eq:dr_residual_derivative_concentration} implies
\begin{equation}
\label{eq:dr_chart_equation_concentration}
 \sup_{x,z}\|D_{(x,z)}^j(\widehat f_h-f_h)(x,z)\|
       \leq Ch^{1-j}\eta,\qquad j=0,1,2.
\end{equation}
For sufficiently small \(h\) and large \(n\),
\[
 \sup_{x,z}\|I_q-\gamma_h^{-1}D_z\widehat f_h(x,z)\|\leq\tfrac12.
\]
The map
\(\mathcal T_{h,x}(z)=z-\gamma_h^{-1}\widehat f_h(x,z)\)
is therefore a contraction. Its displacement at \(z_h(x)\) is at most
\(C_0h\eta\), so it maps
\(\overline B_{2C_0h\eta}(z_h(x))\) into itself. This ball lies inside
\(B_{R\rho^2}(0)\), because \(h\eta=o(h^2)\). It has a unique fixed
point \(\widehat z_h(x)\), and
\begin{equation}
\label{eq:dr_empirical_root_position}
 \sup_x\|\widehat z_h(x)-z_h(x)\|\leq Ch\eta.
\end{equation}
The derivative estimate on the entire convex normal ball also gives
\[
 \langle\widehat f_h(x,z)-\widehat f_h(x,z'),z-z'\rangle
 \geq\tfrac12\gamma_*\|z-z'\|^2.
\]
Hence there is no other zero in that ball. The empirical density is
smooth, its mass is positive, and the spectral gap persists on an open
neighborhood of the tube. The parameterized implicit-function theorem
therefore makes \(\widehat z_h\) at least \(C^2\).

The projected equation is equivalent to the full residual equation.
Indeed, for \(y=x+E(x)z\),
\(\|\widehat Q_h(y)-Q_x\|\leq C(h^2+\eta)<\tfrac12\).
If \(\widehat f_h(x,z)=0\), then \(Q_x\widehat R_h(y)=0\), whereas
\(\widehat R_h(y)=\widehat Q_h(y)\widehat R_h(y)\). Thus
\(\widehat R_h(y)=(\widehat Q_h(y)-Q_x)\widehat R_h(y)\),
which forces \(\widehat R_h(y)=0\). The reverse implication is immediate,
and the required eigenvalue sign holds throughout the tube. Uniqueness
on overlapping normal frames identifies the same ambient zero. The
local sections consequently glue into a global \(C^2\) normal section
\(\widehat u_h^{\mathrm{DR}}\), containing all empirical ridge points
in \(\mathcal T_h^{\mathrm{DR}}\).

\paragraph*{First- and second-derivative stability}
The deterministic residual bounds and bounded chart derivatives imply
\[
 \|D_{(x,z)}^2f_h\|_\infty\leq C/h,\qquad
 \|D_{(x,z)}^3f_h\|_\infty\leq C/h^2.
\]
Put
\[
 \begin{array}{ll}
 B_h(x)=D_zf_h(x,z_h(x)),&C_h(x)=D_xf_h(x,z_h(x)),\\
 \widehat B_h(x)=D_z\widehat f_h(x,\widehat z_h(x)),&
 \widehat C_h(x)=D_x\widehat f_h(x,\widehat z_h(x)).
 \end{array}
\]
Both inverse matrices \(B_h^{-1}\) and \(\widehat B_h^{-1}\) are
uniformly bounded. Equations
\eqref{eq:dr_chart_equation_concentration} and
\eqref{eq:dr_empirical_root_position} give
\[
 \|\widehat B_h-B_h\|_\infty+
 \|\widehat C_h-C_h\|_\infty
 \leq C\eta+(C/h)\|\widehat z_h-z_h\|_\infty
 \leq C\eta.
\]
Implicit differentiation yields
\[
 D_xz_h=-B_h^{-1}C_h,\qquad
 D_x\widehat z_h=-\widehat B_h^{-1}\widehat C_h,
\]
and hence
\begin{equation}
\label{eq:dr_empirical_root_C1}
 \|D_x\widehat z_h-D_xz_h\|_\infty\leq C\eta.
\end{equation}
Here the population section satisfies \(\|z_h\|_{C^2}\leq Ch^2\),
so both first derivatives are uniformly bounded.

For the second derivatives, first note that
\[
 \begin{split}
 &\sup_x\bigl\|
 D_{(x,z)}^2\widehat f_h(x,\widehat z_h(x))
       -D_{(x,z)}^2 f_h(x,z_h(x))\bigr\|\\
 &\hspace{2cm}\leq C\eta/h+(C/h^2)
                          \|\widehat z_h-z_h\|_\infty
 \leq C\eta/h.
 \end{split}
\]
For coordinate directions \(a,b\), differentiating the population
implicit equation twice gives
\[
 \begin{split}
 D_x^2z_h[a,b]=-B_h^{-1}\bigl\{
 &f_{xx}[a,b]+f_{xz}[a,D_xz_h b]\\
 &+f_{zx}[D_xz_h a,b]
      +f_{zz}[D_xz_h a,D_xz_h b]\bigr\},
 \end{split}
\]
with all derivatives of \(f_h\) evaluated at \((x,z_h(x))\). The same
identity holds for the empirical equation. Subtract the two identities.
The second-derivative coefficients differ by at most \(C\eta/h\);
their norms are at most \(C/h\); the first derivatives of the sections
differ by at most \(C\eta\); and the inverse matrices differ by at most
\(C\eta\). These bounds imply
\(\|D_x^2\widehat z_h-D_x^2z_h\|_\infty\leq C\eta/h\).
Returning to ambient normal sections and using the bounded derivatives
of \(E\) proves all three estimates in
\eqref{eq:dr_empirical_c2}.

Since \(\widehat u_h^{\mathrm{DR}}(x)\in N_x\mathcal M\) and
\(\|\widehat u_h^{\mathrm{DR}}\|_\infty\leq R\rho^2<\tau\),
\(\pi_{\mathcal M}\circ\widehat X_h^{\mathrm{DR}} =\operatorname{id}_{\mathcal M}\).
This identity gives injectivity and, after differentiation, immersion.
Compactness of \(\mathcal M\) makes \(\widehat X_h^{\mathrm{DR}}\) a
\(C^2\) embedding. Its \(C^1\) distance from the inclusion is at most
\(C(h^2+\eta)\), which proves the tangent-projector bound.

\paragraph*{Geometric normals at the empirical reference point}
We verify explicitly the normal-field claims used in
Appendix~\ref{sec:gaussian_density_ridge}.
At a zero \(y\) of \(\widehat R_h\), differentiating
\(\widehat R_h=\widehat Q_h\widehat R_h\) gives
\(D_y\widehat R_h(y)=\widehat Q_h(y)D_y\widehat R_h(y)\).
Its rank is at most \(q\), and is at least \(q\) because the normal
Jacobian \(D_z\widehat f_h\) is invertible. Thus its rank is exactly
\(q\). Differentiating the zero equation along the empirical branch
identifies its tangent space with the kernel, proving
\eqref{eq:dr_geometric_normals}.

In a fixed local base chart \(\psi\), write
\[
 J_h=D(X_{h,\sigma}^{\mathrm{DR}}\circ\psi),\qquad
 \widehat J_h=D(\widehat X_h^{\mathrm{DR}}\circ\psi).
\]
The smallest singular values of these \(D\times d\) matrices are
uniformly positive. Their geometric normal projectors are
\[
 \mathscr N(J_h),\quad \mathscr N(\widehat J_h),\qquad
 \mathscr N(J)=I_D-J(J^\top J)^{-1}J^\top.
\]
The map \(\mathscr N\) and its first two derivatives are bounded on the
relevant full-rank matrix family. Since
\[
 \|\widehat J_h-J_h\|_\infty\leq C\eta,\qquad
 \|D\widehat J_h-DJ_h\|_\infty\leq C\eta/h,
\]
the chain rule gives
\[
 \|\mathscr N(\widehat J_h)-\mathscr N(J_h)\|_\infty\leq C\eta,\qquad
 \|\mathscr N(\widehat J_h)-\mathscr N(J_h)\|_{C^1}\leq C\eta/h.
\]
Applying the same calculation to the population embedding and the base
inclusion, whose \(C^2\) difference is \(O(h^2)\), gives the population
\(C^1\) normal-projector error \(O(h^2)\). At \(h=h_1\), these are
exactly \eqref{eq:dr_reference_projector_bounds}. Since \(\eta=o(h)\),
the empirical reference field has uniformly bounded \(C^1\) norm on
\(\mathcal E_h\). This completes the proof.

\paragraph*{PCA at the empirical reference point}
For the PCA alternative, suppose additionally that
$\lambda=\sigma/h\in\Lambda_{\mathrm{PCA}}\Subset(0,\infty)$.
Let $Q_h^{\mathrm{PCA}}$ and $\widehat Q_h^{\mathrm{PCA}}$ denote
the population and empirical normal projectors obtained from the
centered covariance with $K_h(y,Y)=(1-\|Y-y\|^2/h^2)_+^6$.
Since $\rho_h\asymp h$, the enlarged density-ridge tube is contained
in a tube of radius $Ch^2$. The local density calculation used for
\eqref{eq:stsp-kernel-moment-concentration} gives, uniformly there,
$ch^d\leq M_{0,h}(y)\leq Ch^d$ and
$\mathbb P\{\|Y-y\|\leq h\}\leq Ch^d$.
Lemma~\ref{lem:pca_normal_accuracy} also gives a population covariance
eigengap of at least $gh^2$ for some $g>0$.
These facts depend on the local moments and the tube, not on the STSP
fixed-point equation. The arguments for
\eqref{eq:stsp-kernel-normalized-concentration} and
\eqref{eq:stsp-kernel-projector-concentration} therefore give an event
$\mathcal E_h^{\mathrm{PCA}}$, of probability at least $1-n^{-2}$,
on which the empirical PCA eigengap remains at least $gh^2/2$ and
\[
\sup_y\|D_y^j(\widehat Q_h^{\mathrm{PCA}}-Q_h^{\mathrm{PCA}})(y)\|_{\mathrm{op}}
\leq Ch^{-j}\eta,\qquad j=0,1.
\]
The same moment and projector calculations give the deterministic
bounds \[\sup_y\|D_y^jQ_h^{\mathrm{PCA}}(y)\|\leq Ch^{-j}\] for
$j=1,2$. All these suprema are over the enlarged deterministic tube.

Set
\[
\begin{aligned}
\widetilde Q_h^{\mathrm{PCA,ref}}(x)
&=Q_h^{\mathrm{PCA}}\bigl(X_{h,\sigma}^{\mathrm{DR}}(x)\bigr),\\
\widehat Q_h^{\mathrm{PCA,ref}}(x)
&=\widehat Q_h^{\mathrm{PCA}}\bigl(\widehat X_h^{\mathrm{DR}}(x)\bigr).
\end{aligned}
\]
The population ridge expansion makes
$h^{-2}u_{h,\sigma}^{\mathrm{DR}}$ uniformly bounded in $C^2$.
Composing the normal-coordinate PCA expansion in
Lemma~\ref{lem:pca_normal_accuracy} with this section yields
\[\|\widetilde Q_h^{\mathrm{PCA,ref}}-P_{N\mathcal M}\|_{C^1}\leq Ch^2\].
On $\mathcal E_h\cap\mathcal E_h^{\mathrm{PCA}}$, abbreviate
$X=X_{h,\sigma}^{\mathrm{DR}}$, $\widehat X=\widehat X_h^{\mathrm{DR}}$,
$Q=Q_h^{\mathrm{PCA}}$, and $\widehat Q=\widehat Q_h^{\mathrm{PCA}}$.
The branch bounds in~\eqref{eq:dr_empirical_c2} give
$\|\widehat X-X\|_\infty\leq Ch\eta$,
$\|D\widehat X-DX\|_\infty\leq C\eta$, and
$\|DX\|_\infty+\|D\widehat X\|_\infty\leq C$.
Using
\[
\widehat Q\circ\widehat X-Q\circ X
=(\widehat Q-Q)\circ\widehat X+Q\circ\widehat X-Q\circ X
\]
and the chain rule, we obtain
\[
\begin{aligned}
\|\widehat Q_h^{\mathrm{PCA,ref}}
  -\widetilde Q_h^{\mathrm{PCA,ref}}\|_\infty
&\leq C\eta,\\
\|D\widehat Q_h^{\mathrm{PCA,ref}}
  -D\widetilde Q_h^{\mathrm{PCA,ref}}\|_\infty
&\leq C\left\{\frac{\eta}{h}
  +\frac{\|\widehat X-X\|_\infty}{h^2}
  +\frac{\|D\widehat X-DX\|_\infty}{h}\right\}
\leq C\eta/h.
\end{aligned}
\]
Thus the reference fields have the same $C^0$ and $C^1$ stability
orders as in~\eqref{eq:stsp-reference-PCA-stability}. At $h=h_1$
these are the PCA fields in Appendix~\ref{sec:gaussian_density_ridge}.
In particular, $\|\widehat Q_1^{\mathrm{PCA}}-P_{N\mathcal M}\|_\infty
\leq C\{h_1^2+\eta_n(h_1)\}$ and
$\|\widehat Q_1^{\mathrm{PCA}}\|_{C^1}\leq C$ under
$\eta_n(h_1)=o(h_1)$. 
\end{proof}

\subsubsection{Proof of Theorem~\ref{thm:finite_sample_dr}}
\label{app:dr_finite_sample_proof}
\begin{proof}
Apply Proposition~\ref{lem:dr_empirical_rate} at all \(K\) scales, on fixed
enlargements of their tubes. A union bound gives an event
\(\mathcal E=\bigcap_{k=1}^K\mathcal E_{h_k}\) of probability at least
\(1-Kn^{-2}\). No independence between the scale-specific estimators is
needed. Write \(\eta_k=\eta_n(h_k)\). On this event,
\[
 \|\widehat X_{h_k}^{\mathrm{DR}}-X_{h_k,\sigma}^{\mathrm{DR}}\|_\infty
       \leq Ch_k\eta_k,\qquad
 \|\widehat X_{h_k}^{\mathrm{DR}}-X_{h_k,\sigma}^{\mathrm{DR}}\|_{C^1}
       \leq C\eta_k,
\]
and \eqref{eq:dr_reference_projector_bounds} holds at the reference scale.
We henceforth work on \(\mathcal E\).

\paragraph*{Existence and stability of empirical correspondence}
Choose a local base chart \(\psi\) and a \(C^1\) orthonormal frame
\(E_1(x)\) for the geometric normal space of the population reference
ridge. This frame has uniformly bounded \(C^1\) norm. An aligned
empirical frame is
\[
 \widehat E_1(x)
 =\widehat Q_1^{\mathrm{DR}}(x)E_1(x)
 \bigl[E_1(x)^\top\widehat Q_1^{\mathrm{DR}}(x)E_1(x)\bigr]^{-1/2}.
\]
The matrix inside the inverse square root is uniformly positive definite.
The reference-field bounds give
\begin{equation}
\label{eq:dr_aligned_reference_frames}
 \|\widehat E_1-E_1\|_\infty\leq C\eta_1,\qquad
 \|\widehat E_1-E_1\|_{C^1}\leq C\eta_1/h_1,
 \qquad \|\widehat E_1\|_{C^1}\leq C.
\end{equation}
For \(a\in\mathbb R^d\) and \(t\in\mathbb R^q\), define
\[
 \begin{split}
 F_k(a,t;x)
 &=X_{h_k,\sigma}^{\mathrm{DR}}(\psi(a))
       -X_{h_1,\sigma}^{\mathrm{DR}}(x)-E_1(x)t,\\
 \widehat F_k(a,t;x)
 &=\widehat X_{h_k}^{\mathrm{DR}}(\psi(a))
       -\widehat X_{h_1}^{\mathrm{DR}}(x)-\widehat E_1(x)t.
 \end{split}
\]
Lemma~\ref{lem:wellposed_correspondence}, with the geometric reference
normals, gives the population solution \((a_k(x),t_k(x))\), where
\(\|t_k(x)\|\leq Ch_k^2\). Its differential in \((a,t)\) has a
uniformly bounded inverse near that solution: the target tangent space
is close to \(T_x\mathcal M\) and the reference normal space is close
to \(N_x\mathcal M\).

On these local neighborhoods,
\(\|D_{(a,t)}\widehat F_k-D_{(a,t)}F_k\| \leq C(\eta_k+\eta_1)=o(1)\),
while at the population solution,
\[
 \|\widehat F_k(a_k,t_k;x)\|\leq Ce_{n,k}^{\mathrm{corr}},\qquad
 e_{n,k}^{\mathrm{corr}}=h_1\eta_1+h_k\eta_k+h_k^2\eta_1.
\]
The population second derivatives in \((a,t)\) are uniformly bounded.
Using the inverse population differential, the usual quantitative
implicit-function argument is therefore a contraction on a ball of
radius \(Ce_{n,k}^{\mathrm{corr}}\) about \((a_k,t_k)\), for a sufficiently large fixed
constant. This proves a unique nearby empirical solution
\((\widehat a_k,\widehat t_k)\), with displacement at most \(Ce_{n,k}^{\mathrm{corr}}\).
Since \(e_{n,k}^{\mathrm{corr}}=o(h_k^2)\), it stays on the selected local branches.
Its differential in \((a,t)\) is invertible, proving transversality.
The equation is \(C^1\) in \(x\), so the intersections depend
\(C^1\)-smoothly on the reference point. Local uniqueness glues the
solutions across charts, independently of the selected frames. For
\(k=1\), the intersection is simply
\(\widehat X_{h_1}^{\mathrm{DR}}(x)\).

The resulting corresponded points satisfy
\begin{equation}
\label{eq:dr_empirical_correspondence_bound}
 \sup_{x\in\mathcal M}
 \|\widehat Z_{h_k}^{\mathrm{DR}}(x)-Z_{h_k}^{\mathrm{DR}}(x)\|
       \leq Ce_{n,k}^{\mathrm{corr}}.
\end{equation}

\paragraph*{Multiscale error and embedding}
The population result \eqref{eq:density_ridge_multiscale_rate}, with
this choice of reference normals, gives
\[
 \sup_{x\in\mathcal M}\|\widetilde X^{\mathrm{DR}}(x)-x\|
       \leq C\sum_{k=1}^K|w_k|h_k^4.
\]
Combining it with \eqref{eq:dr_empirical_correspondence_bound} proves
\[
 \sup_{x\in\mathcal M}\|\widehat X^{\mathrm{DR}}(x)-x\|
 \leq C\sum_{k=1}^K|w_k|\{h_k^4+e_{n,k}^{\mathrm{corr}}\},
\]
and hence the Hausdorff bound. Equivalently, on \(\mathcal E\), the
general finite-sample theorem applies with deterministic quantities
\(r_n(h_k)=Ch_k\eta_k\) and \(s_n(h_1)=C\eta_1\).

For completeness, we verify the regularity of the multiscale image.
The empirical sections and reference frame satisfy
\[
 \|\widehat u_{h_k}^{\mathrm{DR}}\|_\infty\leq Ch_k^2,\qquad
 \|\widehat u_{h_k}^{\mathrm{DR}}\|_{C^1}
       \leq C(h_k^2+\eta_k),\qquad
 \|\widehat E_1\|_{C^1}\leq C.
\]
In the intersection equations, compare the empirical solution with
\((a_0(x),0)\), where \(\psi(a_0(x))=x\). The value of
\(\widehat F_k\) there is
\(\widehat u_{h_k}^{\mathrm{DR}}(x)
       -\widehat u_{h_1}^{\mathrm{DR}}(x)=O(h_k^2)\).
Uniform transversality gives
\(\|\widehat a_k-a_0\|+\|\widehat t_k\|\leq Ch_k^2\).
Implicit differentiation then yields
\begin{equation}
\label{eq:dr_correspondence_C1_inclusion}
 \|\widehat Z_{h_k}^{\mathrm{DR}}-\iota\|_{C^1}
       \leq C\{h_k^2+\eta_k+h_1^2+\eta_1\},
\end{equation}
where \(\iota:\mathcal M\hookrightarrow\mathbb R^D\) is the inclusion.
To see the scale in this estimate explicitly, the differential in
\((a,t)\) is
\[
 \bigl[D(\widehat X_{h_k}^{\mathrm{DR}}\circ\psi)(\widehat a_k),
                              -\widehat E_1(x)\bigr],
\]
whereas the base derivative of the equation contains
\(-D\widehat X_{h_1}^{\mathrm{DR}}(x)
       -D\widehat E_1(x)\widehat t_k\).
Their differences from the corresponding inclusion equation are bounded
by \(C(h_k^2+\eta_k+h_1^2+\eta_1)\); the frame-derivative term is
\(O(h_k^2)\). The inverses of the differentials are uniformly bounded,
which proves \eqref{eq:dr_correspondence_C1_inclusion}.

Since \(\sum_k w_k=1\) and \(\sum_k|w_k|\leq W_0\),
\[
 \|\widehat X^{\mathrm{DR}}-\iota\|_{C^1}
 \leq CW_0\{h_{\max}^2+\max_k\eta_k\}.
\]
This is sufficiently small for small scales and large \(n\). A
sufficiently small \(C^1\) perturbation of the inclusion of a fixed
compact embedded manifold is an embedding. Indeed, in fixed local
charts, projection onto the chart tangent coordinates gives local
injectivity and full differential rank; compactness and \(C^0\)
closeness preserve separation of pairs outside these local
neighborhoods. Thus \(\widehat X^{\mathrm{DR}}\) is a \(C^1\)
embedding, and its image is \(C^1\)-diffeomorphic to \(\mathcal M\).

\paragraph*{Numerical correspondence error}
Write \(z=\widehat X_{h_1}^{\mathrm{DR}}(x)\), and let \(t_k^*(z)\)
be the coordinates of the exact empirical intersection in the chosen
orthonormal reference frame \(\widehat E(z)\). Then
\[
 \widehat Z_{h_k}^{\mathrm{DR}}(x)
 =z+\widehat E(z)t_k^*(z),\qquad g_{k,z}(t_k^*(z))=0,
 \qquad \|t_k^*(z)\|\leq Ch_k^2.
\]
For large enough fixed \(C_{\mathrm{loc}}\), the exact intersection
lies in the stated search ball. Put
\[
 y_k^0=X_{h_k,\sigma}^{\mathrm{DR}}(x),\qquad
 \gamma_k=\frac{h_k^2}{h_k^2+\sigma^2}\geq\gamma_{\mathrm{ENF},*}
:=\frac{1}{1+C_\sigma^2}>0.
\]
At the population ridge point, the normal Jacobian estimate
\eqref{eq:dr_normal_jacobian}, together with
\(R_{h_k,\sigma}^{\mathrm{DR}}\circ X_{h_k,\sigma}^{\mathrm{DR}}=0\)
and \(R_{h_k,\sigma}^{\mathrm{DR}}
       =Q_{h_k,\sigma}^{\mathrm{DR}}R_{h_k,\sigma}^{\mathrm{DR}}\), gives
\begin{equation}
\label{eq:dr_ambient_jacobian_at_ridge}
 D_yR_{h_k,\sigma}^{\mathrm{DR}}(y_k^0)
       =\gamma_k Q_x+O(h_k^2).
\end{equation}
In detail, the derivative annihilates the tangent space of the
population ridge, which differs from \(T_x\mathcal M\) by \(O(h_k^2)\).
Its range lies in the Hessian normal eigenspace, which differs from
\(N_x\mathcal M\) by the same order. Its normal-to-normal block is
\(\gamma_k I_q+O(h_k^2)\). These three facts prove
\eqref{eq:dr_ambient_jacobian_at_ridge}.

Every point \(z+\widehat E(z)t\) with
\(\|t\|\leq C_{\mathrm{loc}}h_k^2\) is within \(Ch_k^2\) of
\(y_k^0\), and the connecting segment stays inside a fixed enlargement
of the tube. By \eqref{eq:dr_population_residual_derivatives},
\eqref{eq:dr_residual_derivative_concentration}, and the reference normal
bound,
\[
 \begin{split}
 D_tg_{k,z}(t)
 &=\widehat E(z)^\top D_y\widehat R_{h_k}^{\mathrm{DR}}
                   (z+\widehat E(z)t)\widehat E(z),\\
 \|D_tg_{k,z}(t)-\gamma_k I_q\|
 &\leq C\{h_k+h_1^2+\eta_k+\eta_1\}.
 \end{split}
\]
For sufficiently small scales and large \(n\), integration along line
segments in the convex search ball therefore gives
\[
 \langle g_{k,z}(t)-g_{k,z}(s),t-s\rangle
       \geq\tfrac12\gamma_{\mathrm{ENF},*}\|t-s\|^2.
\]
Taking \(s=t_k^*(z)\) and using the residual tolerance proves
\[
 \|\widehat t_k(z)-t_k^*(z)\|
       \leq\frac{2}{\gamma_{\mathrm{ENF},*}}\varepsilon_{n,k}.
\]
Set \(t_1^*(z)=\widehat t_1(z)=0\). Since \(\widehat E(z)\) is an
isometry, the numerical output at each retained query differs from
\(\widehat F^{\mathrm{DR}}(z)\) by at most
\(C\sum_k|w_k|\varepsilon_{n,k}\).

\paragraph*{Finite query coverage}
The reference-branch map
\(\widehat F^{\mathrm{DR}}
 =\widehat X^{\mathrm{DR}}\circ(\widehat X_{h_1}^{\mathrm{DR}})^{-1}\)
has a uniformly bounded Lipschitz constant on \(\mathcal E\). This
follows from \eqref{eq:dr_correspondence_C1_inclusion}, bounded weights,
and the uniform bi-Lipschitz property of the empirical reference normal
graph. The latter follows from its \(C^1\) closeness to the inclusion,
using fixed local charts and compactness for separated point pairs.
For the exact sampled image
\(\mathcal P_N^* =\{\widehat F^{\mathrm{DR}}(z_j):1\leq j\leq N\}\),
we therefore have
\[
 \begin{split}
 d_H(\widehat{\mathcal P}_N^{\mathrm{DR,num}},\mathcal P_N^*)
       &\leq C\sum_k|w_k|\varepsilon_{n,k},\\
 d_H(\mathcal P_N^*,\widehat{\mathcal M}^{\mathrm{DR}})
       &\leq C\delta_N.
 \end{split}
\]
The triangle inequality proves
\eqref{eq:dr_numerical_error}. Only the retained nonempty query
set is used. No derivative control on the numerical outputs or
regularity of an interpolation of the point cloud is inferred.
This completes the proof.
\end{proof}

\paragraph*{PCA-based correspondence}
Suppose additionally that $\sigma/h_1$ lies in a fixed compact subset
of $(0,\infty)$. For the PCA construction, work on
$\mathcal E\cap\mathcal E_{h_1}^{\mathrm{PCA}}$, which has probability
at least $1-(K+1)n^{-2}$ by the union bound.
The PCA reference-field estimates at the end of
Appendix~\ref{app:dr_empirical_rate_proof} give the frame bounds
\eqref{eq:dr_aligned_reference_frames} with $E_1$ and $\widehat E_1$
replaced by local population and aligned empirical PCA frames.
The population PCA field has $O(h_1^2)$ normal-projector error,
so Lemma~\ref{lem:wellposed_correspondence} supplies the corresponding
population intersections. With these frames in the same intersection
equations, the differential perturbation is still
$O(\eta_k+\eta_1)$ and the residual at the population solution is
$O(e_{n,k}^{\mathrm{corr}})$. The preceding quantitative implicit-function argument
therefore gives unique transverse empirical PCA-fiber intersections,
$C^1$ dependence on the reference point, and
\eqref{eq:dr_empirical_correspondence_bound} for the PCA-defined
correspondence points.
The uniform $C^1$ frame bound also gives
\eqref{eq:dr_correspondence_C1_inclusion} by the same implicit
differentiation. Consequently, the multiscale error bound and
embedding argument apply to the PCA-defined continuous estimator.

For its numerical implementation, let $\widehat E^{\mathrm{PCA}}(z)$
be an orthonormal PCA normal basis at the empirical reference point
$z=\widehat X_{h_1}^{\mathrm{DR}}(x)$, and put
\[
g_{k,z}^{\mathrm{PCA}}(t)
=\widehat E^{\mathrm{PCA}}(z)^\top
\widehat R_{h_k}^{\mathrm{DR}}
\bigl(z+\widehat E^{\mathrm{PCA}}(z)t\bigr).
\]
The ridge branches and residuals are unchanged. On the prescribed
search ball, the population Jacobian identity
\eqref{eq:dr_ambient_jacobian_at_ridge}, the residual derivative
bounds, and the PCA normal-projector error give
\[
\sup_{\|t\|\leq C_{\mathrm{loc}}h_k^2}
\|D_tg_{k,z}^{\mathrm{PCA}}(t)-\gamma_k I_q\|_{\mathrm{op}}
\leq C\{h_k+h_1^2+\eta_k+\eta_1\}=o(1).
\]
Thus the same strong-monotonicity argument bounds the distance
from an accepted numerical root to the exact PCA-fiber intersection
by $2\varepsilon_{n,k}/\gamma_{\mathrm{ENF},*}$.
The reference-branch map is uniformly Lipschitz by the PCA version
of~\eqref{eq:dr_correspondence_C1_inclusion}, so the query-coverage
argument is unchanged. Hence~\eqref{eq:dr_numerical_error} also
holds for the corresponding PCA-defined exact estimator and
numerical output. The fourth-order conclusion of
Corollary~\ref{cor:dr_fourth_order_sample_size} follows under its
sample-size and numerical conditions, on this joint event and with
the additional noise-to-bandwidth restriction above.

\subsubsection{Noise-variance estimation from the density-Hessian spectrum}
\label{app:dr_noise_estimation}

\begin{proposition}[Spectral noise-variance estimation for density ridges]
\label{prop:dr_noise_estimation}
Suppose the conditions of Proposition~\ref{lem:dr_empirical_rate} hold at the
reference bandwidth \(h_1\), with \(0\leq\sigma\leq C_\sigma h_1\).
Define \(\widehat\gamma_1\) and \(\widehat\sigma^2\) as in
Appendix~\ref{sec:gaussian_density_ridge}.
For sufficiently small \(h_1\) and large \(n\), on an event of probability
at least \(1-n^{-2}\) where the reference ridge exists,
\[
|\widehat\sigma^2-\sigma^2|
\leq C\{h_1^4+h_1^2\eta_n(h_1)\}.
\]
The constants are uniform over the stated noise--bandwidth regime and
measurable selections of \(\widehat y_1\) from that ridge.
If the remaining assumptions of Theorem~\ref{thm:finite_sample_dr}
hold at all fitted scales, replacing its weights by spatially constant,
uniformly bounded plug-in weights at fixed comparable deterministic scales
preserves its recovery rate and embedding conclusion, as well as the
sufficient sample-size condition of
Corollary~\ref{cor:dr_fourth_order_sample_size}, using the same
observations for fitting and noise estimation.
\end{proposition}

\begin{proof}
Write \(h=h_1\), \(\eta=\eta_n(h)\), and \(q=D-d\).
Use the normalized density Hessians
\[
A_h(y)=h^2\frac{\nabla^2p_{h,\sigma}(y)}{p_{h,\sigma}(y)},
\qquad
\widehat A_h(y)=h^2\frac{\nabla^2\widehat p_h(y)}{\widehat p_h(y)}.
\]
By \eqref{eq:dr_normalized_hessian_gap}, uniformly on
\(\mathcal T_h^{\mathrm{DR}}\),
\[
A_h(y)=-\gamma_hQ_x+O(h^2),\qquad
\gamma_h=\frac{h^2}{h^2+\sigma^2}\in[\gamma_*,1],
\qquad x=\pi_{\mathcal M}(y).
\]
On the event \(\mathcal E_h\) used in the proof of
Proposition~\ref{lem:dr_empirical_rate}, which has probability at least
\(1-n^{-2}\), the empirical ridge exists and
\eqref{eq:dr_normalized_derivative_concentration} gives
\[
\sup_{y\in\mathcal T_h^{\mathrm{DR}}}
\|\widehat A_h(y)-A_h(y)\|_{\mathrm{op}}\leq C\eta.
\]
The lowest \(q\) eigenvalues of \(-\gamma_hQ_x\) are all
\(-\gamma_h\). Weyl's inequality yields
\[
\left|-\frac1q\sum_{j=d+1}^D
\lambda_j(\widehat A_h(y))-\gamma_h\right|
\leq C(h^2+\eta).
\]
Since \(\widehat p_h(y)>0\), the eigenvalues of
\(\widehat A_h(y)\) are
\(h^2\widehat\lambda_j^h(y)/\widehat p_h(y)\), in the same order.
The uniform bound applies at \(\widehat y_1\), even when the point is
selected using the observations. Projection onto \([\gamma_*,1]\)
cannot increase distance from \(\gamma_h\), so
\(|\widehat\gamma_1-\gamma_h|\leq C(h^2+\eta)\).
Both arguments of the reciprocal are at least \(\gamma_*>0\), and
\(\sigma^2=h^2(\gamma_h^{-1}-1)\). Therefore
\[
|\widehat\sigma^2-\sigma^2|
=h^2|\widehat\gamma_1^{-1}-\gamma_h^{-1}|
\leq h^2\gamma_*^{-2}|\widehat\gamma_1-\gamma_h|
\leq C(h^4+h^2\eta).
\]
Outside the reference-branch existence event, define the estimate to be
zero. The estimator is nonnegative and at most \(C_\sigma^2h^2\).
No division by \(\sigma\) is involved, so the bound remains uniform
as \(\sigma/h\to0\). When \(\sigma\asymp h\), taking the square root
also gives \(|\widehat\sigma-\sigma|=O_p(h^3+h\eta)\).

For the plug-in weights, their affine constraint gives the exact identity
\[
\sum_k\widehat w_k
\{c_{\mathrm{DR}}(h_k,\sigma)-c_{\mathrm{DR}}(h_k,\widehat\sigma)\}
=\frac d2(\sigma^2-\widehat\sigma^2).
\]
Here \(\widehat\sigma=\sqrt{\widehat\sigma^2}\); the identity and the
pathwise error decomposition also hold when \(\widehat\sigma=0\),
because \(c_{\mathrm{DR}}\) is a polynomial in the noise level.
Thus the additional Hausdorff-error term is at most
\[
\frac d2\|H\|_\infty|\widehat\sigma^2-\sigma^2|
\leq C\{h_1^4+h_1^2\eta_n(h_1)\}.
\]
For comparable scales this is absorbed into the original
\(O(h_{\max}^4+h_{\max}\eta_n(h_{\max}))\) bound.
The branch and correspondence maps are unchanged. The near-identity
\(C^1\) bound in the proof of Theorem~\ref{thm:finite_sample_dr} uses
only spatially constant weights with sum one and bounded absolute sum,
so its embedding conclusion applies to these random weights as well.
The numerical and covering estimates remain valid with the same weight
replacement.
Under Corollary~\ref{cor:dr_fourth_order_sample_size}, the additional
term is \(O(h_{\max}^4)\) on the same joint event of probability at least
\(1-Kn^{-2}\). 
\end{proof}

\FloatBarrier
\section{Directional GBC for the YY Estimator}
\label{app:ysl23_full_pipeline}
Throughout, write
\(Q(x)=P_{N_x\mathcal M}\), \(T(x)=I_D-Q(x)\),
\(v_0(x)=V_0(X_0(x))\), and
\(\overline U_0(x)=U_0(X_0(x))\).
In local moment calculations, suppress the argument \(x\) and write \(z=X_0(x)\).
For sufficiently small \(\sigma\), \(Qv_0\neq0\); define
\begin{equation}
\label{eq:ysl_common_normal_split}
 \nu_0=\frac{Qv_0}{\|Qv_0\|},\qquad
 P_0=\nu_0\nu_0^\top,\qquad
 S_0=Q-P_0,\qquad
 b_\perp=S_0u_0.
\end{equation}
Here \(P_0\) projects onto a line in the true normal space, \(S_0\)
onto its \(m\)-dimensional normal complement, and \(b_\perp\) is the
exact preliminary offset in that complement.

All derivatives are taken on the fixed manifold \(\mathcal M\), with
\(\sigma\) and the radii held fixed.  Norms are computed in a fixed finite
atlas with ambient vector or matrix components.  Constants may depend on
this atlas, the fixed \(C^6\) geometry, the dimensions, the kernel profiles,
and the uniform constants in Assumption~\ref{ass:ysl_full_population}, but
not on \(\sigma\), the base point, or the selected admissible scale.

The analysis below concerns the smooth-kernel YY map
$G_{r_1}$ defined by \eqref{eq:ysl_full_population_G}
and its common-start directional combination
\eqref{eq:ysl_common_algorithm}.
This is a modification of the weighted map in
\citet[equation~(4.14)]{yao2023manifold};
the projected population update in their
equation~(4.12) is not analyzed here.

For a vector, matrix, or tensor field \(f\), write
\begin{equation}
\label{eq:ysl_common_scaled_norm}
 \|f\|_{2,\sigma}
 =\max_{0\leq j\leq2}\sigma^{j/2}\|f\|_{C^j}.
\end{equation}
Leibniz's rule gives
\(\|fg\|_{2,\sigma}\leq C\|f\|_{2,\sigma}\|g\|_{2,\sigma}\).
Composition with a function having bounded first two derivatives on a fixed
compact parameter set preserves boundedness of this norm.  We use this
notation only to account for derivatives of the preliminary offset and the
first-stage axis; it does not assert a stronger unscaled regularity
assumption.

\subsection{Proof of preliminary normal coordinates from reach}

We use the chord--tangent characterization of reach, recalled in
\citet[Lemma~2.3]{yao2023manifold}: if a compact \(C^2\) embedded manifold
\(\mathcal N\) has reach at least \(\tau_{\mathcal N}\), then
\begin{equation}
\label{eq:ysl_common_reach_paraboloid}
 \operatorname{dist}(z,y+T_y\mathcal N)
 \leq\frac{\|z-y\|^2}{2\tau_{\mathcal N}},
 \qquad y,z\in\mathcal N.
\end{equation}
The supremum of twice the left-hand side divided by \(\|z-y\|^2\), over
unequal pairs, is the reciprocal of the reach.

\begin{lemma}[A thin positive-reach manifold is a normal graph]
\label{lem:ysl_common_tube_upgrade}
Assume that $\mathcal M$ is connected. Let $\mathcal N$ be a
nonempty, compact, connected, boundaryless, $C^2$ embedded
$d$-manifold such that
\[
 \mathcal N\subset\operatorname{Tub}_\varepsilon(\mathcal M),
 \qquad \operatorname{reach}(\mathcal N)\geq\tau_{\mathcal N}>0.
\]
For a sufficiently small numerical $c>0$, if
$\varepsilon\leq c\min\{\tau,\tau_{\mathcal N}\}$, then
$\pi_{\mathcal M}|_{\mathcal N}$ is a diffeomorphism onto
$\mathcal M$. Moreover, with
$x=\pi_{\mathcal M}(y)$ and
$\rho=\min\{\tau,\tau_{\mathcal N}\}$,
\begin{equation}
\label{eq:ysl_common_angle_bound}
 \|P_{T_y\mathcal N}-P_{T_x\mathcal M}\|
 \leq C\sqrt{\varepsilon/\rho}.
\end{equation}
\end{lemma}

\begin{proof}
The case \(\varepsilon=0\) follows from the equal dimensions and the
open-and-closed argument below.  Suppose \(\varepsilon>0\).  Fix a unit
vector \(v\in T_y\mathcal N\) and a unit-speed geodesic \(\gamma\) with
\(\gamma(0)=y\), \(\gamma'(0)=v\).  The curvature bound furnished by reach
gives
\(\gamma(s)=y+sv+e(s),\qquad \|e(s)\|\leq Cs^2/\tau_{\mathcal N}\).
Choose \(x_s=\pi_{\mathcal M}(\gamma(s))\).  Then
\(x_s-x=sv+(y-x)+e(s)+(x_s-\gamma(s))\).
Both \(\|y-x\|\) and \(\|x_s-\gamma(s)\|\) are at most \(\varepsilon\).
Set \(s=\sqrt{\varepsilon\rho}\).  For sufficiently small \(c\), this
choice satisfies \(\varepsilon\leq s\ll\rho\), and hence
\(\|x_s-x\|\leq Cs\).  Applying \(Q(x)\) and
\eqref{eq:ysl_common_reach_paraboloid} on \(\mathcal M\) yields
\[
 s\|Q(x)v\|
 \leq C\left(\varepsilon+s^2/\tau_{\mathcal N}
                    +\|x_s-x\|^2/\tau\right)
 \leq C(\varepsilon+s^2/\rho).
\]
Thus \(\|Q(x)v\|\leq C\sqrt{\varepsilon/\rho}\).  The equal dimensions
identify this one-sided tangent gap with the projector gap, proving
\eqref{eq:ysl_common_angle_bound}.

Make \(c\) smaller so that the gap is less than \(1/4\).
The differential of the nearest-point projection is the tangent projection
followed by the inverse of \(I-A_{y-x}\), where \(A\) is the shape
operator of \(\mathcal M\).  It is therefore invertible on
\(T_y\mathcal N\).  The restriction is a local diffeomorphism.  Its image
is open, and is also closed by compactness.  Connectedness of \(\mathcal M\)
makes the image all of \(\mathcal M\).

For injectivity, suppose \(y_1\neq y_2\) have the same projection \(x\).
Their difference lies in \(N_x\mathcal M\) and has norm at most
\(2\varepsilon\).  The tangent gap gives
\(\operatorname{dist}(y_2,y_1+T_{y_1}\mathcal N) \geq\tfrac12\|y_2-y_1\|\).
Combining this with \eqref{eq:ysl_common_reach_paraboloid} on
\(\mathcal N\) implies \(\|y_2-y_1\|\geq\tau_{\mathcal N}\), a
contradiction.  A bijective local diffeomorphism is a diffeomorphism.
\end{proof}

\begin{proof}[Proof of Lemma~\ref{lem:ysl_full_preliminary_bounds}]
Apply Lemma~\ref{lem:ysl_common_tube_upgrade} with
\(\mathcal N=\mathcal M^{\mathrm{pre}}\) and
\(\varepsilon=C_0\sigma\).  This gives \(X_0=x+u_0\), the \(C^0\)
bound, and an \(O(\sigma^{1/2})\) tangent gap.  For a normal section,
\[
 DX_0(x)v=(I-A_{u_0(x)})v+\nabla_v^\perp u_0.
\]
The tangent block is uniformly invertible because \(u_0=O(\sigma)\).
The tangent gap therefore implies
\(\|\nabla^\perp u_0\|_{C^0}\leq C\sigma^{1/2}\); the difference from
ambient differentiation consists of \(A_{u_0}\), which is
\(O(\sigma)\).

To bound the second derivative, differentiate the displayed graph identity
and project onto \(N_{X_0(x)}\mathcal M^{\mathrm{pre}}\).  The resulting
normal component of \(D^2X_0\) is
\(\mathrm{II}_{\mathcal M^{\mathrm{pre}}}(DX_0,DX_0)\).
The restriction of this normal projection to \(N_x\mathcal M\) is
uniformly invertible by the tangent gap.  It thus solves for
\(\nabla^{\perp,2}u_0\) in terms of the second fundamental form of the
preliminary manifold and bounded quantities involving the fixed geometry,
\(u_0\), and \(Du_0\).  Since reach bounds
\(\|\mathrm{II}_{\mathcal M^{\mathrm{pre}}}\|\leq\tau_*^{-1}\), we
obtain \(\|u_0\|_{C^2}\leq C\).
Finally, \(X_0-\operatorname{Id}_{\mathcal M}\) is \(C^1\)-small.
Within a fixed local chart, a line-segment estimate gives a uniform lower
Lipschitz bound; for separated pairs the \(C^0\) bound gives the same
conclusion.  An upper Lipschitz bound follows by the same two-region
argument.  This proves the last assertion.
\end{proof}

\subsection{Proof of cutoffs, localization, and differentiated integral bounds}
\label{sec:ysl_common_localization_proof}

For \(s\geq0\), let
\[
 K_{0,r_0}(s)=\left(1-\frac{s^2}{r_0^2}\right)_+^p,
 \qquad
 K_{\perp,r_1}(s)=\left(1-\frac{s^2}{r_1^2}\right)_+^p.
\]
Unlike the truncated Gaussian-type weights in
\citet{yao2023manifold}, which lack sufficient
smoothness at the truncation boundary, these
polynomial kernels and the $C^3$ longitudinal
plateau cutoff provide the regularity required
for the differentiated higher-order expansions
and finite-sample analysis while preserving
the original localization mechanism.

Let \(K_{\parallel,r_2}(s)=k(s/r_2)\), where the fixed \(C^3\) profile
\(k\) takes values in \([0,1]\), equals one on \([0,1/2]\), and vanishes
on \([1,\infty)\). The explicit polynomial profile is
\begin{equation}
\label{eq:ysl_common_plateau}
 k(t)=
 \begin{cases}
 1,&0\leq t\leq1/2,\\
 \eta(2t-1),&1/2<t<1,\\
 0,&t\geq1,
 \end{cases}
 \qquad
 \eta(u)=1-35u^4+84u^5-70u^6+20u^7.
\end{equation}
Indeed,
\[
 \eta(0)=1,\quad\eta(1)=0,\quad
 \eta^{(j)}(0)=\eta^{(j)}(1)=0\ (j=1,2,3),
 \qquad \eta'(u)=-140u^3(1-u)^3\leq0.
\]
Consequently, \(0\leq k\leq1\), \(k\in C^3\), and
\(\|K_{\parallel,r_2}^{(j)}\|_\infty\leq C_jr_2^{-j},\qquad j=0,1,2,3\).
The spherical and transverse profiles are also \(C^3\) for \(p\geq4\).
Their radial extensions are smooth at the origin: the first two are
functions of the squared norm, and the longitudinal profile is locally
constant there.  Thus the norms in the cylindrical weight create no
additional singularity at zero.  On the nondegenerate first-stage domain,
all maps used below are at least \(C^2\).

For all population expansions, use a fixed, uniformly controlled normal
graph chart around \(x\),
\begin{equation}
\label{eq:ysl_common_local_graph}
 X_x(t)=x+t+h_x(t),\qquad t\in T_x\mathcal M,
 \qquad
 h_x(t)=h_{2,x}(t)+h_{3,x}(t)+\mathcal H_{4,x}(t),
\end{equation}
where
\[
 h_{2,x}(t)=\tfrac12\mathrm{II}_x(t,t),\qquad
 h_{3,x}(t)=\tfrac16\mathcal C_x[t,t,t],\qquad
 \|\mathcal H_{4,x}(t)\|\leq C\|t\|^4.
\]
Its volume density satisfies
\begin{equation}
\label{eq:ysl_common_graph_volume}
 J_x(t)=\sqrt{\det(I+(Dh_x(t))^\top Dh_x(t))}
       =1+J_{2,x}(t)+O(\|t\|^3),
\end{equation}
with \(J_{2,x}\) quadratic.  The same graph remainder bounds hold after
up to two base derivatives, using the fixed \(C^6\) geometry.  Local
orthonormal frames can be chosen in finitely many bundle charts; the
estimates do not depend on which such frame is used.

The original spherical or cylindrical support is contained in an
\(O(\sigma\sqrt{L_\sigma})\) ball about \(x\).  A latent point outside the
fixed local chart can contribute only through a Gaussian displacement bounded
away from zero.  Its contribution, including the first two derivatives
needed here, is \(O(\sigma^N)\) relative to the local mass for every fixed
\(N\).  This follows directly from the Gaussian exponential bound, which
absorbs the fixed powers of \(\sigma^{-1}\) produced by differentiation.
One may insert a smooth chart cutoff that is one near the origin before
making any change to the longitudinal weight.

In the localized integrals, write \(t=\sigma(s+g)\), with \(s\) bounded
and \(g\) Gaussian.  Restricting first to
\(\|g\|\leq\sigma^{-1/8}\) loses less than any fixed power of \(\sigma\).
On this region, \(h_x(t)/\sigma=o(1)\).  The transition of the longitudinal
cutoff then requires a standardized longitudinal Gaussian displacement at
least
\((C_2/2)\sqrt{L_\sigma}-C-o(1)\).
Its mass, with any fixed polynomial factors required by two derivatives, is
bounded by \(C\sigma^{cC_2^2}\) after decreasing the positive numerical
constant \(c\).  Choose \(C_2\) large enough that all such tail errors in
normalized means and their needed differentiated bounds are, for example,
\(O(\sigma^6)\).  We may therefore replace the longitudinal weight by one
\emph{inside the localized integral}.  We never make this replacement on
the whole manifold, which could admit unrelated branches along the
resulting infinite cylinder.

All Taylor remainder estimates below are integrated against Gaussian
moments of polynomials in \(g\).  They therefore give pure powers of
\(\sigma\), rather than powers multiplied by truncation logarithms.
For parameter derivatives, one can translate the center and rotate the axis
to a fixed cylinder, transferring differentiation to the smooth
Gaussian-convolved density.  This also avoids imposing fourth or higher
classical derivatives on the compact cutoff itself.

\subsection{Proof of Lemma~\ref{lem:ysl_common_first_stage_bounds}}
\label{app:ysl_first_stage_bounds_proof}

\begin{proof}
Put \(a=u_0/\sigma\).  Lemma~\ref{lem:ysl_full_preliminary_bounds} gives
\begin{equation}
\label{eq:ysl_common_a_bounds}
 \|a\|_{C^0}\leq C,
 \qquad\|Da\|_{C^0}\leq C\sigma^{-1/2},
 \qquad\|D^2a\|_{C^0}\leq C\sigma^{-1}.
\end{equation}
At a fixed base point, use \eqref{eq:ysl_common_local_graph}.  The flat
spherical model in coordinates relative to \(z=x+\sigma a\) has the
probability density
\begin{equation}
\label{eq:ysl_common_first_flat}
 d\Pi_a^0(s,b)
 =(Z_a^0)^{-1}
 \left(1-\frac{\|s\|^2+\|b\|^2}{c_0^2}\right)_+^p
 \phi_q(b+a)\,ds\,db,
 \quad (s,b)\in\mathbb R^d\times\mathbb R^q.
\end{equation}
In this flat model, the Gaussian variable \(g\in\mathbb R^d\) is standard
and independent of \((s,b)\).
The flat tangential mean is zero.  Rotational symmetry in the normal block
gives
\(\mathbb E_{\Pi_a^0}b=-\gamma_0(a)a\).
For \(a=\rho e\), \(\rho>0\), pair a point with \(b\cdot e>0\) with its
reflection across \(e^\perp\).  The contribution to the numerator in the
\(e\)-direction is proportional to
\((b\cdot e)\{e^{-\rho(b\cdot e)}-e^{\rho(b\cdot e)}\}<0\).
Thus \(\gamma_0(a)>0\).  The coefficient is smooth on compact nonzero
normal annuli; it also has a positive continuous extension at zero given by
the corresponding flat variance.  All mass denominators are bounded below
on bounded sets of \(a\).

In the curved model the normal Gaussian becomes
\(\phi_q(b+a-h_x(t)/\sigma)\).  The quadratic graph produces a relative
first-order correction in \(\sigma\), whereas the volume correction starts
at order \(\sigma^2\).  Tangential reflection
\((s,g)\mapsto(-s,-g)\) cancels both the flat tangential numerator and its
quadratic-graph correction.  It follows that
\begin{equation}
\label{eq:ysl_common_first_expansion}
 \begin{split}
 Qv_0&=-\gamma_0(a)u_0+\sigma^2B_0(x,a,\sigma),\\
 Tv_0&=\sigma^3B_T(x,a,\sigma),
 \end{split}
\end{equation}
where \(B_0,B_T\), and their partial derivatives of order at most two in
\(x,a\), are uniformly bounded.  These coefficient functions may be defined
by the exact differences in \eqref{eq:ysl_common_first_expansion}; they retain
all higher-order terms.  The localization argument above and dominated
Gaussian differentiation justify the uniform derivative bounds.

The signal condition and the second identity imply
\(\inf_x\|Qv_0(x)\|\geq c\sigma\).  Since the flat mean is smooth and
vanishes at \(a=0\), the first identity also implies \(\|a(x)\|\geq c'>0\)
for sufficiently small \(\sigma\).  Thus \(\gamma_0\) has a uniform
positive lower bound on the actual compact parameter range.
Normalizing \(Qv_0/\sigma\) is smooth with bounded derivatives on this
range.  Equation~\eqref{eq:ysl_common_a_bounds} gives
\(\|\nu_0\|_{2,\sigma}+\|P_0\|_{2,\sigma} +\|S_0\|_{2,\sigma}\leq C\).
As \(S_0Qv_0=0\), the first identity in
\eqref{eq:ysl_common_first_expansion} yields the exact relation
\begin{equation}
\label{eq:ysl_common_exact_residual}
 b_\perp=S_0u_0
 =\sigma^2\gamma_0(a)^{-1}S_0B_0(x,a,\sigma).
\end{equation}
This proves \(\|b_\perp\|_{2,\sigma}\leq C\sigma^2\).

Finally set
\begin{equation}
\label{eq:ysl_common_tilt_parameter}
 \theta(x)=\frac{Tv_0(x)}{\|Qv_0(x)\|}\in T_x\mathcal M,
 \qquad
 \overline U_0
 =\frac{(\nu_0+\theta)(\nu_0+\theta)^\top}{1+\|\theta\|^2}.
\end{equation}
Then \(\|\theta\|_{2,\sigma}\leq C\sigma^2\).
Expanding the smooth projector expression proves
\(\|\overline U_0-P_0\|_{2,\sigma}\leq C\sigma^2\), and all asserted
bounds follow by the definition of the scaled norm.
\end{proof}

A second-order approximation to \(B_0\) or \(P_0\) could leave
third-order positional terms.  The exact field \(b_\perp\) in
\eqref{eq:ysl_common_exact_residual} includes these terms, is shared
across transverse scales, and is canceled by the same weight equation.

\subsection{Proof of the transverse and mixed moments}
\label{sec:ysl_common_moments}

For \(m\geq1\), define a probability measure on
\((s,b)\in\mathbb R^d\times\mathbb R^m\) by
\begin{equation}
\label{eq:ysl_common_flat_measure}
 d\Pi_\kappa(s,b)
 =Z_\kappa^{-1}
 \left(1-\frac{\|s\|^2+\|b\|^2}{\kappa^2}\right)_+^p
 \phi_m(b)\,ds\,db,
 \qquad \kappa>0,
\end{equation}
where \(\phi_m\) is the standard \(m\)-dimensional Gaussian density.  Set
\(
 \alpha_m(\kappa)=\mathbb E_{\Pi_\kappa}[s_1^2],
 \beta_m(\kappa)=\mathbb E_{\Pi_\kappa}[b_1^2].
\)
Blockwise rotational symmetry makes the coordinate choices immaterial.

Write
\[
 K_\kappa(s,b)
 =\left(1-\frac{\|s\|^2+\|b\|^2}{\kappa^2}\right)_+^p,
 \qquad J=d+2p+2.
\]
For \(m\geq1\), integrating over \(s\) first shows that the marginal
normal density is proportional to
\((1-\|b\|^2/\kappa^2)_+^{p+d/2}\phi_m(b)\).
Conditional on \(b\), radial polynomial integration in \(\mathbb R^d\)
gives
\[
 \mathbb E_{\Pi_\kappa}[s_i s_\ell\mid b]
 =\frac{\kappa^2-\|b\|^2}{J}\,\delta_{i\ell}.
\]
Passing to the squared radial coordinate in the normal block yields
\begin{equation}
\label{eq:ysl_common_one_dimensional_moments}
 \begin{split}
 \beta_m(\kappa)
 &=\frac1m
 \frac{\displaystyle\int_0^{\kappa^2}
 t^{m/2}e^{-t/2}(1-t/\kappa^2)^{p+d/2}\,dt}
 {\displaystyle\int_0^{\kappa^2}
 t^{m/2-1}e^{-t/2}(1-t/\kappa^2)^{p+d/2}\,dt},\\
 \alpha_m(\kappa)&=\frac{\kappa^2-m\beta_m(\kappa)}{J}.
 \end{split}
\end{equation}
The denominators are strictly positive.  Gaussian integration by parts in
\(b_j\) gives
\[
 \beta_m
 =1+Z_\kappa^{-1}\int b_j(\partial_{b_j}K_\kappa)\phi_m(b)\,ds\,db<1.
\]
The strict inequality follows from radial decrease on a set of positive
measure; positivity of the variance gives \(\beta_m>0\).
These moments and any fixed number of their scale derivatives are uniformly
bounded on compact subsets of \((0,\infty)\).  One can see this by rescaling
\((s,b)=\kappa(u,v)\), which fixes the integration domain.

\begin{lemma}[Mixed moment]
\label{lem:ysl_common_mixed_moment}
For \(m\geq1\), every tangent coordinate \(i\), and every complementary
normal coordinate \(j\),
\begin{equation}
\label{eq:ysl_common_mixed_identity}
 \mathbb E_{\Pi_\kappa}[s_i^2b_j^2]
 =\alpha_m(\kappa)-\beta_m(\kappa),
 \qquad
 \mathbb E_{\Pi_\kappa}[b_j^2(1+s_i^2)]
 =\alpha_m(\kappa).
\end{equation}
\end{lemma}

\begin{proof}
Integrating \(\partial_{b_j}(s_i^2b_jK_\kappa\phi_m)\) gives
\[
 \int s_i^2b_j^2K_\kappa\phi_m
 =\int s_i^2K_\kappa\phi_m
  +\int s_i^2b_j(\partial_{b_j}K_\kappa)\phi_m.
\]
Radiality implies
\(s_i^2b_j\partial_{b_j}K_\kappa
=s_ib_j^2\partial_{s_i}K_\kappa\).
An integration by parts in \(s_i\) changes the last integral into
\(-\int b_j^2K_\kappa\phi_m\).  The boundary terms vanish because the
cutoff vanishes at its support boundary.  Divide by \(Z_\kappa\) and use
the definitions of \(\alpha_m\) and \(\beta_m\) to obtain the first identity; the second
follows by addition of \(\beta_m\).
\end{proof}

The mixed expectation in \eqref{eq:ysl_common_mixed_identity} is not a
product of two marginal expectations.  The transverse cutoff couples the
tangent and complementary normal variables, even though the ambient noise
coordinates are independent before weighting.

\subsection{Proof of normal expansion at a common preliminary point}

We first calculate with the normal axis \(P_0(x)\), and subsequently
restore the actual axis \(\overline U_0(x)\).  This replacement is for the
proof only.  Split the normal space as
\(\operatorname{im}P_0\oplus\operatorname{im}S_0\), and write
\[
 u_0=\Delta\nu_0+c,\qquad c=b_\perp,\qquad
 h_L=P_0h_x,\quad h_E=S_0h_x.
\]
Then \(\Delta=O(\sigma)\) and \(c=O(\sigma^2)\).
In the localized integral, use observed tangent coordinate \(\sigma s\)
and observed complementary normal coordinate \(c+\sigma b\).
After integrating the longitudinal Gaussian, the flat integration law is
\(\Pi_\kappa(ds,db)\phi_d(g)\,dg\), with \(\kappa=r_1/\sigma\).
Denote expectation under this law by \(\mathbb E_\kappa\) and define
\begin{equation}
\label{eq:ysl_common_likelihood_ratio}
 \delta(t)=\frac{h_E(t)-c}{\sigma},
 \qquad
 L(t,b)=\exp\{b\cdot\delta(t)-\tfrac12\|\delta(t)\|^2\}.
\end{equation}
Indeed,
\(\phi_m(b+(c-h_E(t))/\sigma)=\phi_m(b)L(t,b)\).
If \(A^{P}_{r_1}\) denotes the mass for this localized, normal-axis model,
then, up to the negligible localization and longitudinal tails,
\begin{equation}
\label{eq:ysl_common_exact_local_moments}
 \begin{split}
 A^{P}_{r_1}
 &=\frac{\sigma^d Z_\kappa}{\operatorname{vol}(\mathcal M)}D_\kappa,
 \qquad D_\kappa=\mathbb E_\kappa[J_x(t)L(t,b)],\\
 M_L&=\frac{\mathbb E_\kappa[h_L(t)J_x(t)L(t,b)]}{D_\kappa},\\
 M_E&=c+\sigma\frac{\mathbb E_\kappa[bJ_x(t)L(t,b)]}{D_\kappa},\\
 M_T&=\sigma\frac{\mathbb E_\kappa[sJ_x(t)L(t,b)]}{D_\kappa}.
 \end{split}
\end{equation}
Here \(M_L,M_E,M_T\) are the axial, complementary normal, and tangent
components relative to \(x\); suppressing the local chart cutoff incurs
the stated negligible error.

For clarity, write \(w=s+g\), \(c=\sigma^2\bar c\), and
\[
 h_x(\sigma w)
 =\sigma^2h_{2,x}(w)+\sigma^3h_{3,x}(w)
    +O(\sigma^4\|w\|^4).
\]
Thus
\[
 \delta
 =\sigma\{S_0h_{2,x}(w)-\bar c\}
   +\sigma^2S_0h_{3,x}(w)+O(\sigma^3\|w\|^4).
\]
The field \(\bar c\) is bounded; it is not expanded further in \(\sigma\).
On the Gaussian truncation described above,
\begin{equation}
\label{eq:ysl_common_gaussian_taylor}
 L=1+b\cdot\delta
   +\tfrac12\{(b\cdot\delta)^2-\|\delta\|^2\}
   +\mathcal E_3,
\end{equation}
where the integrated remainder is bounded by the third moment of
\(\|\delta\|\) times a fixed polynomial in \(b\).  In particular,
\(\mathbb E_\kappa[\|\mathcal E_3\|]\leq C\sigma^3\), and analogous
polynomially weighted estimates hold.  Reflection in \(b\), followed by
\eqref{eq:ysl_common_graph_volume}, gives
\begin{equation}
\label{eq:ysl_common_mass_expansion}
 D_\kappa=1+O(\sigma^2),
 \qquad A^P_{r_1}\asymp\sigma^d.
\end{equation}
There is no first-order mass term: \(b\cdot\delta\) is odd in \(b\)
for fixed \(s,g\).

\paragraph*{The axial mean}
The cubic graph term has zero expectation under
\((s,g)\mapsto(-s,-g)\).  The linear likelihood term multiplied by
\(h_L(t)\) has zero expectation under \(b\mapsto-b\).
The fourth-order graph remainder, the volume correction, and the quadratic
likelihood correction therefore give
\[
 \mathbb E_\kappa[h_L(t)J_x(t)L(t,b)]
 =\sigma^2\mathbb E_\kappa[P_0h_{2,x}(w)]+O(\sigma^4).
\]
Since \(g\) is independent of \((s,b)\), centered, and has identity
covariance,
\(\mathbb E_\kappa[w_iw_\ell] =(1+\alpha_m(\kappa))\delta_{i\ell}\).
Using \(\operatorname{tr}\mathrm{II}_x=dH(x)\) and
\eqref{eq:ysl_common_mass_expansion} gives
\begin{equation}
\label{eq:ysl_common_axial_mean}
 M_L=\frac d2\sigma^2(1+\alpha_m(\kappa))P_0H+O(\sigma^4).
\end{equation}

\paragraph*{The complementary normal mean}
The constant term in \(\sigma\mathbb E_\kappa[bJ_xL]\) vanishes by
normal reflection.  Its quadratic likelihood term also vanishes: the
quadratic expression in \eqref{eq:ysl_common_gaussian_taylor} is even in
\(b\), and the numerator contains an additional \(b\).
The first nonzero contribution is the linear one,
\[
 \sigma\mathbb E_\kappa[b(b\cdot\delta)]
 =\mathbb E_\kappa[bb^\top\{h_E(t)-c\}].
\]
The cubic graph part vanishes by tangent reflection.  For normal coordinate
\(j\), all off-diagonal normal and tangent terms vanish by coordinate
reflections, so the retained term is
\[
 -\beta_m c_j
 +\frac12\sigma^2\sum_{i=1}^d
     (S_0\mathrm{II}_{ii})_j
     \mathbb E_{\Pi_\kappa}[b_j^2(1+s_i^2)].
\]
Lemma~\ref{lem:ysl_common_mixed_moment} identifies the mixed factor with
\(\alpha_m\).  The cubic likelihood remainder contributes
\(\sigma O(\sigma^3)=O(\sigma^4)\); the volume correction to the linear
term and the fourth-order graph contribution are also \(O(\sigma^4)\).
Dividing by \(D_\kappa\) changes a second-order numerator by only
\(O(\sigma^4)\).  Adding back the exact \(c\) in
\eqref{eq:ysl_common_exact_local_moments}, we conclude
\begin{equation}
\label{eq:ysl_common_complement_mean}
 M_E=(1-\beta_m(\kappa))c
      +\frac d2\sigma^2\alpha_m(\kappa)S_0H+O(\sigma^4).
\end{equation}
In particular, evaluating the mixed moment leaves no additional
\(\sigma^2\beta_m S_0H\) term.

\paragraph*{The tangent mean and differentiated bounds}
The flat tangent mean is zero.  The linear likelihood term again vanishes
by normal reflection.  The volume and quadratic likelihood corrections are
\(O(\sigma^2)\), and the numerator contains the factor \(\sigma s\).
Thus \(\|M_T\|\leq C\sigma^3\).  This bound does not require the
additional tangent cancellations that may occur in special geometries.

For derivatives, the bounded rescaled parameters are the base point, the
normal axis, \(\Delta/\sigma\), and \(c/\sigma^2\).  Their first two
base derivatives cost at most \(\sigma^{-1/2}\) and \(\sigma^{-1}\), by
Lemma~\ref{lem:ysl_common_first_stage_bounds} and
\eqref{eq:ysl_common_a_bounds}.  The preceding cancellations hold throughout
a fixed compact set of these parameters.  Gaussian differentiation and the
\(C^6\) graph bounds show that the exact tangent and normal means have the
forms \(\sigma^3\mathcal T\) and \(\sigma^2\mathcal V\), with
\(\mathcal T,\mathcal V\) uniformly \(C^2\) in these parameters.
Consequently,
\begin{equation}
\label{eq:ysl_common_normal_axis_derivatives}
 \|M_T\|_{C^j}\leq C\sigma^{3-j/2},
 \qquad
 \|M_L+M_E\|_{C^j}\leq C\sigma^{2-j/2},\quad j=0,1,2.
\end{equation}
Only the complete second-order normal mean is differentiated here.  No
fourth-order unscaled \(C^2\) estimate for its positional remainder is
being assumed.

\subsection{Proof of restoring the estimated axis}

At a fixed base point use the tilt parameter in
\eqref{eq:ysl_common_tilt_parameter}.  For
\(\theta\in T_x\mathcal M\), put
\[
 \nu_\theta=\frac{\nu_0+\theta}{\sqrt{1+\|\theta\|^2}},
 \qquad U_\theta=\nu_\theta\nu_\theta^\top,
\]
and let \(G_{r_1}^{\theta}\) be the local mean with the same center \(z\)
and radii but cylindrical axis \(U_\theta\).  The following estimate
controls the effect on the fourth-order normal remainder.

\begin{lemma}[Axis replacement for the full weighted mean]
\label{lem:ysl_common_axis_replacement}
Fix \(x\in\mathcal M\).  There are
constants \(\theta_*>0\), \(c>0\), and \(C<\infty\), uniform in \(x\),
\(r_1/\sigma\in[c_-,c_+]\), and sufficiently small \(\sigma\), such that,
for \(\|\theta\|\leq\theta_*\),
\begin{equation}
\label{eq:ysl_common_axis_response}
 \begin{split}
 \|T(G^{\theta}_{r_1}-G^{0}_{r_1})\|
 &\leq C\sigma\|\theta\|,\\
 \|Q(G^{\theta}_{r_1}-G^{0}_{r_1})\|
 &\leq C\{\sigma^2\|\theta\|+\sigma\|\theta\|^2\},
 \end{split}
\end{equation}
up to the longitudinal and localization tails already absorbed into the
fourth-order remainder.

For the numerator--denominator estimates underlying this conclusion, define
\begin{equation}
\label{eq:ysl_common_axis_mass_numerator}
 \begin{split}
 \Psi_\theta(y)
 &=K_{\parallel,r_2}(\|U_\theta(y-z)\|)
   K_{\perp,r_1}(\|(I_D-U_\theta)(y-z)\|),\\
 A_\theta&=\mathbb E[\Psi_\theta(Y)],\\
 N_\theta&=\mathbb E[(Y-x)\Psi_\theta(Y)],
 \qquad
 G_{r_1}^{\theta}=x+\frac{N_\theta}{A_\theta}.
 \end{split}
\end{equation}
Uniformly on the same parameter range,
\begin{equation}
\label{eq:ysl_common_axis_natural_bounds}
 c\sigma^d\leq A_\theta\leq C\sigma^d,
 \qquad
 \|N_\theta\|\leq C\sigma^{d+1}.
\end{equation}
For \(v,v_1,v_2\in T_x\mathcal M\), the first derivatives at the
normal axis satisfy
\begin{equation}
\label{eq:ysl_common_axis_first_derivatives}
 \begin{split}
 |D_\theta A_0[v]|
 &\leq C\sigma^{d+1}\|v\|,\\
 \|T D_\theta N_0[v]\|
 &\leq C\sigma^{d+1}\|v\|,\\
 \|Q D_\theta N_0[v]\|
 &\leq C\sigma^{d+2}\|v\|,
 \end{split}
\end{equation}
and
\begin{equation}
\label{eq:ysl_common_axis_base_numerator}
 \|TN_0\|\leq C\sigma^{d+3},
 \qquad
 \|QN_0\|\leq C\sigma^{d+2}.
\end{equation}
The second derivatives obey the natural bounds
\begin{equation}
\label{eq:ysl_common_axis_second_derivatives}
 \begin{split}
 |D_\theta^2A_\theta[v_1,v_2]|
 &\leq C\sigma^d\|v_1\|\|v_2\|,\\
 \|D_\theta^2N_\theta[v_1,v_2]\|
 &\leq C\sigma^{d+1}\|v_1\|\|v_2\|.
 \end{split}
\end{equation}
The corresponding first derivatives away from zero satisfy
\(|D_\theta A_\theta[v]|\leq C\sigma^d\|v\|\) and
\(\|D_\theta N_\theta[v]\|\leq C\sigma^{d+1}\|v\|\).
If \(A_\theta^{\mathrm{flat}}\) and
\(N_\theta^{\mathrm{flat}}\) denote the same integrals with the local
manifold replaced by its tangent plane and the graph volume density replaced
by one, then tangent reflection gives the exact identities
\begin{equation}
\label{eq:ysl_common_axis_flat_cancellation}
 D_\theta A_0^{\mathrm{flat}}[v]=0,
 \qquad
 QD_\theta N_0^{\mathrm{flat}}[v]=0.
\end{equation}
These identities and bounds imply \eqref{eq:ysl_common_axis_response}.
\end{lemma}

\begin{proof}
We work in the localized coordinates of
\eqref{eq:ysl_common_local_graph}.  Put \(w=y-z\),
\(w_T=Tw\), and \(w_L=\langle\nu_0,w\rangle\).  Since
\(D_\theta\nu_\theta|_{\theta=0}[v]=v\),
\begin{equation}
\label{eq:ysl_common_axis_distance_derivative}
 \begin{split}
 \left.D_\theta\|U_\theta w\|^2\right|_{\theta=0}[v]
 &=2w_L\langle v,w_T\rangle,\\
 \left.D_\theta\|(I_D-U_\theta)w\|^2\right|_{\theta=0}[v]
 &=-2w_L\langle v,w_T\rangle.
 \end{split}
\end{equation}
After the scaled changes of variables, the derivative of the cylindrical
weight is therefore an integrable even profile multiplied by
\(w_L\langle v,w_T\rangle/\sigma^2\).  It is odd under the tangent
reflection \(w_T\mapsto-w_T\).

The flat quantities in
\eqref{eq:ysl_common_axis_flat_cancellation} retain the possibly nonzero
axial center and the complementary center \(c\).  The tangential numerator
contains one additional factor \(Tw\), so the same
reflection does not force it to vanish.  Since \(Tw=O(\sigma)\) on the
moment scale and the effective local mass is \(O(\sigma^d)\),
\begin{equation}
\label{eq:ysl_common_axis_flat_tangent_scale}
 \|TD_\theta N_0^{\mathrm{flat}}[v]\|
 \leq C\sigma^{d+1}\|v\|.
\end{equation}
The complementary center does not alter these cancellations because it does
not break tangent reflection symmetry.

In the curved model, the first correction to the flat
Gaussian-convolved density and volume factor is
\(O(\sigma)\) with a fixed integrable polynomial envelope.  Combining this
factor with the differentiated weight gives
\(|D_\theta A_0[v]|\leq C\sigma^{d+1}\|v\|\).
For the normal numerator there is one further centered normal coordinate of
moment size \(O(\sigma)\), which yields
\(\|QD_\theta N_0[v]\| \leq C\sigma^{d+2}\|v\|\).
The tangential bound in
\eqref{eq:ysl_common_axis_first_derivatives} already holds at flat order.
The estimates in \eqref{eq:ysl_common_axis_base_numerator} are precisely the
unnormalized versions of the tangent and normal mean bounds in
\eqref{eq:ysl_common_normal_axis_derivatives}, together with
\(A_0\asymp\sigma^d\).

Two tilt derivatives of either squared distance in
\eqref{eq:ysl_common_axis_distance_derivative} have bounded rescaled
polynomial envelopes.  Thus a twice differentiated mass has size
\(O(\sigma^d)\), while a twice differentiated numerator has one additional
coordinate moment and has size \(O(\sigma^{d+1})\).  This proves
\eqref{eq:ysl_common_axis_second_derivatives}; the uniform first-derivative
bounds follow in the same way.  The mass lower bound follows by restricting
to a fixed inner part of the rescaled cylinder on which both cutoffs are
bounded below.  These arguments justify differentiation under the integral
sign; alternatively, the fixed-cylinder argument in
Appendix~\ref{sec:ysl_common_localization_proof} avoids derivatives of the
compact cutoff.

Differentiate the quotient in
\eqref{eq:ysl_common_axis_mass_numerator} at \(\theta=0\):
\begin{equation}
\label{eq:ysl_common_axis_quotient_first}
 D_\theta G^0_{r_1}[v]
 =\frac{D_\theta N_0[v]}{A_0}
  -\frac{N_0D_\theta A_0[v]}{A_0^2}.
\end{equation}
Using \eqref{eq:ysl_common_axis_natural_bounds}--
\eqref{eq:ysl_common_axis_base_numerator} gives
\begin{equation}
\label{eq:ysl_common_axis_mean_first_scale}
 \|TD_\theta G^0_{r_1}[v]\|
 \leq C\sigma\|v\|,
 \qquad
 \|QD_\theta G^0_{r_1}[v]\|
 \leq C\sigma^2\|v\|.
\end{equation}
The quotient second-derivative formula, together with the natural first- and
second-derivative bounds above, gives uniformly for
\(\|\theta\|\leq\theta_*\),
\begin{equation}
\label{eq:ysl_common_axis_mean_second_scale}
 \|D_\theta^2G^\theta_{r_1}[v_1,v_2]\|
 \leq C\sigma\|v_1\|\|v_2\|.
\end{equation}
Taylor's formula between \(0\) and \(\theta\) now gives
\[
 \|T(G^\theta_{r_1}-G^0_{r_1})\|
 \leq C\{\sigma\|\theta\|+\sigma\|\theta\|^2\},
\]
and
\[
 \|Q(G^\theta_{r_1}-G^0_{r_1})\|
 \leq C\{\sigma^2\|\theta\|+\sigma\|\theta\|^2\}.
\]
After reducing \(\theta_*\), the quadratic term in the first display is
absorbed into \(C\sigma\|\theta\|\), proving
\eqref{eq:ysl_common_axis_response}.
\end{proof}

For the actual first-stage axis,
\(\|\theta\|_{C^0}\leq C\sigma^2\), and hence
Lemma~\ref{lem:ysl_common_axis_replacement} gives an additional normal error
of order \(O(\sigma^4)\) and an additional tangent term of order
\(O(\sigma^3)\).  For base derivatives,
\(\theta/\sigma^2\) has bounded \(\|\cdot\|_{2,\sigma}\) norm.  Rotating to
a fixed cylinder gives bounded first two partial derivatives of all
normalized coefficient functions, so the chain rule preserves
\eqref{eq:ysl_common_normal_axis_derivatives} after the actual axis is
restored.

\begin{proof}[Proof of Proposition~\ref{thm:ysl_common_raw_expansion}]
Combine \eqref{eq:ysl_common_axial_mean},
\eqref{eq:ysl_common_complement_mean}, and
\eqref{eq:ysl_common_axis_response}, with \(c=b_\perp\) as defined above.  Define \(t_{r_1}\) by exact tangent
projection and define \(R_{r_1}^N\) by subtracting the three displayed
normal fields from the exact normal projection.  Thus the remainder is
normal by definition.  The preceding calculations give its \(C^0\)
fourth-order bound and all the differentiated bounds in
\eqref{eq:ysl_common_raw_bounds}.
When \(m=0\), the complementary Gaussian block and its likelihood ratio
are absent.  The same axial calculation uses the ordinary radial polynomial
second moment \(\alpha_0(\kappa)=\kappa^2/(d+2p+2)\), proving the
hypersurface statement as well.
\end{proof}

\subsection{Proof of near-identity images and normal reparameterization}

\begin{lemma}[A quantitative near-identity image criterion]
\label{lem:ysl_common_image_geometry}
Let \(\mathcal F:\mathcal M\to\mathbb R^D\) be \(C^2\), and write
\[
 \mathcal F(x)=x+t(x)+v(x),\qquad
 t(x)\in T_x\mathcal M,\quad v(x)\in N_x\mathcal M.
\]
If \(\|\mathcal F-\operatorname{Id}_{\mathcal M}\|_{C^1}\) is sufficiently
small, then \(\mathcal F\) is an embedding and
\(\varphi=\pi_{\mathcal M}\circ\mathcal F\) is a diffeomorphism of
\(\mathcal M\).  Its image is a unique normal graph \(x+u(x)\), and
\begin{equation}
\label{eq:ysl_common_image_distance}
 \|u\|_{C^0}
 \leq C\{\|v\|_{C^0}+\|t\|_{C^0}^2\}.
\end{equation}
If additionally \(\|\mathcal F-\operatorname{Id}_{\mathcal M}\|_{C^2}\)
is sufficiently small, the image has reach at least a fixed positive
constant, and
\begin{equation}
\label{eq:ysl_common_image_c2}
 \|u\|_{C^j}\leq
 C\|\mathcal F-\operatorname{Id}_{\mathcal M}\|_{C^j},\qquad j=1,2.
\end{equation}
Furthermore,
\begin{equation}
\label{eq:ysl_common_reparameterization_difference}
 \|u-v\|_{C^0}
 \leq C\{\|t\|_{C^0}\|v\|_{C^0}
           +\|t\|_{C^0}^2+\|Dv\|_{C^0}\|t\|_{C^0}\}.
\end{equation}
The difference in the last display compares the two ambient vector fields
at the same point of \(\mathcal M\).
\end{lemma}

\begin{proof}
Set \(R=\mathcal F-\operatorname{Id}_{\mathcal M}\).
For two points in a fixed small chart, their intrinsic distance is at most
a constant times their ambient distance.  Integrating \(DR\) along a short
chart path gives
\[
 \|\mathcal F(x)-\mathcal F(y)\|
 \geq(1-C\|DR\|_{C^0})\|x-y\|.
\]
For pairs outside these small charts, \(\|x-y\|\) has a fixed positive
lower bound, and
\(\|\mathcal F(x)-\mathcal F(y)\|\geq\|x-y\|-2\|R\|_{C^0}\)
gives another uniform lower bound.  This proves global injectivity and a
uniform lower Lipschitz constant; the differential is also injective.
The same argument applies to
\(\varphi-\operatorname{Id}_{\mathcal M}\), which is \(C^1\)-small
because the nearest-point projection is smooth in a fixed tube.
Its image is open by the inverse function theorem and closed by compactness,
hence equals \(\mathcal M\).  Thus \(\varphi\) is a diffeomorphism and
\(u=\mathcal F\circ\varphi^{-1}-\operatorname{Id}_{\mathcal M}\) is a
normal section.

For each \(x\), the point on \(\mathcal M\) with tangent graph coordinate
\(t(x)\) equals \(x+t(x)+O(\|t(x)\|^2)\).  It follows that
\[
 \operatorname{dist}(\mathcal F(x),\mathcal M)
 \leq\|v(x)\|+C\|t(x)\|^2.
\]
As \(\varphi\) is onto, this bound gives
\eqref{eq:ysl_common_image_distance} and both directed Hausdorff bounds.
For the differentiated estimates, the map
\((x,r)\mapsto\pi_{\mathcal M}(x+r)-x\) is smooth and vanishes at
\(r=0\).  The chain rule gives \(C^j\) bounds in terms of \(R\) for
\(j\leq2\).  Its near-identity inverse has the corresponding bounds.
Applying the same argument to
\(\mathcal F-\pi_{\mathcal M}\circ\mathcal F\) proves
\eqref{eq:ysl_common_image_c2}.

To verify a uniform reach bound, let \(a=\mathcal F(x)\),
\(b=\mathcal F(y)\).  For nearby base points, Taylor expansion in fixed
charts, the bounded \(C^2\) norm of \(\mathcal F\), and the lower
Lipschitz constant yield
\[
 \operatorname{dist}(b,a+T_a\mathcal F(\mathcal M))
 \leq C\|x-y\|^2\leq C'\|a-b\|^2.
\]
For separated base points, \(\|a-b\|\geq c'>0\), and the same bound
follows from \(\operatorname{dist}(b,a+T_a\mathcal F(\mathcal M))
\leq\|a-b\|\).  The chord--tangent characterization of reach therefore
gives a uniform positive lower bound.  This is the smooth-image mechanism
used in \citet[Theorem~5.4]{yao2023manifold}, with global injectivity and the
lower Lipschitz bound established explicitly here.

Finally the expansion of the nearest-point projection about a normal fiber
is
\[
 \varphi(x)=x+t(x)
 +O\bigl(\|t(x)\|\|v(x)\|+\|t(x)\|^2\bigr).
\]
Thus \(u(\varphi(x))-v(x)=x+t(x)-\varphi(x)\) is bounded by the first two
terms in \eqref{eq:ysl_common_reparameterization_difference}.
Comparing \(v(x)\) and \(v(\varphi(x))\) adds at most
\(C\|Dv\|_{C^0}\|t\|_{C^0}\), proving the final estimate.
\end{proof}

\begin{proof}[Proof of Corollary~\ref{thm:ysl_full_single_scale_section}]
Apply Lemma~\ref{lem:ysl_common_image_geometry} with
\(\mathcal F=Y_{r_1}\).  Proposition~\ref{thm:ysl_common_raw_expansion} gives
\[
 \|Y_{r_1}-\operatorname{Id}_{\mathcal M}\|_{C^1}
 \leq C\sigma^{3/2},\qquad
 \|Y_{r_1}-\operatorname{Id}_{\mathcal M}\|_{C^2}
 \leq C\sigma.
\]
Both tend to zero.  We obtain the single-scale embedding, unique normal
section, its differentiated bounds, and uniform reach.  Since \(X_0\) is
a diffeomorphism, the same embedding conclusion holds for the restriction
of \(G_{r_1}\) to \(\mathcal M^{\mathrm{pre}}\).
Let \(v_{r_1}=Q(Y_{r_1}-\operatorname{Id}_{\mathcal M})\).  In
\eqref{eq:ysl_common_reparameterization_difference}, the three contributions
have orders \(\sigma^5\), \(\sigma^6\), and \(\sigma^{9/2}\), respectively.
Therefore
\(\|u_{r_1}-v_{r_1}\|_{C^0}\leq C\sigma^{9/2}\leq C\sigma^4\).
The raw normal expansion consequently remains valid for the normal section,
with the claimed fourth-order positional remainder.
\end{proof}

\subsection{Proof of admissible scale triples and cancellation}

For the multiscale arguments below, write
\(G_k=G_{r_{1,k}}\), \(Y_k=G_k\circ X_0\), and \(d_k=a_k-b_k\).
Since both directional weight vectors sum to one, \(\sum_{k=1}^K d_k=0\).

The axial and complementary curvature coefficients are
\[
 c_L(r,\sigma)=\sigma^2(1+\alpha_m(r/\sigma)),
 \qquad c_E(r,\sigma)=\sigma^2\alpha_m(r/\sigma).
\]
If a single weight vector independently canceled both coefficients while
summing to one, subtraction of the cancellation equations would give
\(0=\sigma^2\).  This observation rules out that independent-field scalar
argument; it does not exclude special geometric relations among the fields.
The construction here does not assume any such relation.

\begin{lemma}[Existence of nondegenerate scale triples]
\label{lem:ysl_common_scale_triples}
For fixed \(d\geq1\), \(m\geq1\), and \(p\geq4\), the three functions
\(1\), \(\alpha_m(\kappa)\), and \(1-\beta_m(\kappa)\) are linearly
independent on every nonempty open interval in \((0,\infty)\).
Consequently, that interval contains three radii for which the weight matrix is nonsingular.
\end{lemma}

\begin{proof}
Rescale \((s,b)=\kappa(u,v)\) in
\eqref{eq:ysl_common_flat_measure}.  The integration domain is the fixed
unit ball, and the remaining Gaussian factor is
\(\exp(-\kappa^2\|v\|^2/2)\).  Its integrals, and therefore the moment
ratios, are real analytic for \(\kappa>0\).
The marginal formula preceding
\eqref{eq:ysl_common_one_dimensional_moments} and dominated convergence give
\(\beta_m(\kappa)\to1\) as \(\kappa\to\infty\).
The second identity in \eqref{eq:ysl_common_one_dimensional_moments} then
shows that \(\alpha_m(\kappa)\to\infty\).
If
\[
 A+B\alpha_m(\kappa)+C(1-\beta_m(\kappa))=0
\]
on a nonempty open interval, analyticity extends the identity to all
\(\kappa>0\).  Sending \(\kappa\) to infinity first gives \(B=0\), then
\(A=0\).  Since \(\beta_m(\kappa)<1\) for finite \(\kappa\), also
\(C=0\).  Linear independence implies the existence of three evaluation
vectors spanning \(\mathbb R^3\), which proves the final assertion.
The lemma is an existence statement, not a uniform numerical conditioning
bound for arbitrary triples; the latter is separately imposed in \eqref{eq:ysl_common_weights}.
\end{proof}

\begin{proof}[Proof of Proposition~\ref{thm:ysl_common_population}]
Put \(t_k=T(Y_k-x)\) and define the ideal directional matrices
\(B_k(x)=a_kP_0(x)+b_k(I_D-P_0(x))\).
They are used for comparison only.  As both weight vectors sum to one,
\(\sum_kB_k=I_D\).  On \(T_x\mathcal M\), \(B_k\) acts as the scalar
\(b_k\); on the two normal blocks it acts as \(a_k\) and \(b_k\).
Insert \eqref{eq:ysl_common_raw_expansion}.  The three cancellation equations
in \eqref{eq:ysl_common_weights} give the exact identity
\begin{equation}
\label{eq:ysl_common_ideal_cancellation}
 \sum_kB_k(x)Y_k(x)-x
 =\sum_k b_k t_k(x)
   +\sum_k\{a_kP_0(x)+b_kS_0(x)\}R_k^N(x).
\end{equation}
In particular, the ideal tangent term is \(O(\sigma^3)\) and its normal
term is \(O(\sigma^4)\).  The exact residual field \(b_\perp\), not a
truncated approximation to it, has canceled.

The algorithm uses \(\overline U_0\) rather than \(P_0\).  Using
\(\sum_k d_k=0\), the difference in the combined output is
\begin{equation}
\label{eq:ysl_common_axis_replacement_in_sum}
 \begin{split}
 \widetilde\Phi(x)-\sum_kB_k(x)Y_k(x)
 &=\sum_k d_k(\overline U_0(x)-P_0(x))(Y_k(x)-x).
 \end{split}
\end{equation}
The first-stage axis bound is \(O(\sigma^2)\), and
\(\max_k\|Y_k-\operatorname{Id}_{\mathcal M}\|_{C^0}\leq C\sigma^2\).
Hence \eqref{eq:ysl_common_axis_replacement_in_sum} is \(O(\sigma^4)\).
Together with \eqref{eq:ysl_common_ideal_cancellation}, this proves
\[
\left\|T\bigl(\widetilde\Phi-\operatorname{Id}_{\mathcal M}\bigr)\right\|_{C^0}
= O(\sigma^3),
\quad
\left\|Q\bigl(\widetilde\Phi-\operatorname{Id}_{\mathcal M}\bigr)\right\|_{C^0}
= O(\sigma^4).
\]
It is important to subtract the common base point in this calculation;
multiplying the axis error by an uncentered position would give an
unnecessarily poor bound.

For regularity, put
\(\overline A_k=A_k\circ X_0
=a_k\overline U_0+b_k(I_D-\overline U_0)\).
Lemma~\ref{lem:ysl_common_first_stage_bounds} implies
\(\|\overline U_0\|_{2,\sigma}\leq C\), while
Proposition~\ref{thm:ysl_common_raw_expansion} gives
\(\|Y_k-\operatorname{Id}_{\mathcal M}\|_{2,\sigma}\leq C\sigma^2\).
Use the exact centering identity
\[
 \widetilde\Phi-\operatorname{Id}_{\mathcal M}
 =\sum_k\overline A_k(Y_k-\operatorname{Id}_{\mathcal M}).
\]
The product estimate for \eqref{eq:ysl_common_scaled_norm} and bounded
weights yield
\(\|\widetilde\Phi-\operatorname{Id}_{\mathcal M}\|_{2,\sigma} \leq C\sigma^2\).
Apply Lemma~\ref{lem:ysl_common_image_geometry}.  The map and its
nearest-point projection are respectively an embedding and a diffeomorphism,
and the image has uniformly positive reach.  In the distance estimate,
\(\|t\|_{C^0}\leq C\sigma^3\) and \(\|v\|_{C^0}\leq C\sigma^4\), so
\[
 \|\widetilde u\|_{C^0}
 \leq C(\sigma^4+\sigma^6)\leq C'\sigma^4.
\]
Because the nearest-point projection is onto, this controls both directed
Hausdorff distances.  The \(C^2\) bound for the normal section follows
from \eqref{eq:ysl_common_image_c2}.  This proves
\eqref{eq:ysl_common_ms_image} and the proposition.
\end{proof}

For Remark~\ref{rem:ysl_common_hypersurface}, take \(P_0=Q\), delete the
complementary normal terms, and set \(b_k=a_k\).  The matrix weights reduce
to scalar weights and the axis-replacement term
\eqref{eq:ysl_common_axis_replacement_in_sum} vanishes identically.  The
same tangent--normal image argument gives the stated two-scale result.

\subsection{Proof of uniform finite-sample approximation of the common-start map}
\label{sec:ysl_common_finite_sample_proof}

Throughout this subsection, \(\mathcal M^{\mathrm{pre}}\) and \(X_0\)
are fixed.  Write \(\mathbb P_n f=n^{-1}\sum_{i=1}^n f(Y_i)\),
\(\mathbb P f=\mathbb E[f(Y)]\), \(\ell=\ell_\delta(\sigma)\), and
\(\zeta=\zeta_{n,\delta}(\sigma)\).

\paragraph*{Proof of
Theorem~\ref{thm:ysl_common_finite_sample}}

We first construct a single high-probability event controlling both stages.
Choose the fixed constant \(\theta_*>0\) from the axis-replacement estimate
above.  By the population first-stage bounds, there is a constant \(C_J\)
such that
\(\|DU_0\|_{C^0(\mathcal M^{\mathrm{pre}})} \leq \frac{C_J}{2}\sigma^{-1/2}\).
In each coordinate chart of the fixed finite atlas, let \(B\) denote one
coordinate derivative of a rank-one projector.  Consider the deterministic
projector-jet class
\begin{equation}
\label{eq:ysl_one_sample_projector_jet_class}
 \mathfrak J_\sigma
 =
 \left\{
 (z,k,P,B):
 \begin{array}{l}
 z\in\mathcal M^{\mathrm{pre}},\quad 1\leq k\leq K,\\
 P=P^\top=P^2,\quad \operatorname{rank}(P)=1,\\
 \|P-U_0(z)\|_{\mathrm{op}}\leq\theta_*,\\
 \|B\|_{\mathrm{op}}\leq C_J\sigma^{-1/2}
 \end{array}
 \right\}.
\end{equation}
For \((z,k,P,B)\in\mathfrak J_\sigma\), define
\[
 \Psi_{k,z,P}(y)
 =
 K_{\parallel,r_2}\bigl(\|P(y-z)\|\bigr)
 K_{\perp,r_{1,k}}\bigl(\|(I_D-P)(y-z)\|\bigr).
\]
Its derivative in the selected base coordinate is written as
\(D_{z,B}\Psi_{k,z,P} = \partial_z\Psi_{k,z,P} + D_P\Psi_{k,z,P}[B]\),
where \(\partial_z\) differentiates the center with the projector held
fixed.

Let \(w\) denote either a spherical first-stage weight or a cylindrical
weight indexed by \(\mathfrak J_\sigma\).  Set
\(A_w=\mathbb P w(Y), M_w=\frac{\mathbb P\{Yw(Y)\}}{A_w}\).
For \(j=0\), let \(D^0w=w\).  For \(j=1\), let \(D^1w\) denote any
coordinate derivative of the spherical weight or the corresponding
derivative \(D_{z,B}w\) of a cylindrical weight.

The localization and differentiated moment estimates established above give,
uniformly over both weight classes and for \(j=0,1\),
\begin{equation}
\label{eq:ysl_one_sample_moment_bounds}
 \begin{split}
 A_w
 &\asymp
 \sigma^d,\\
 \mathbb E\!\left[
 \bigl\|(Y-M_w)D^jw(Y)\bigr\|^2
 \right]
 &\leq
 C\sigma^{d+2-2j},\\
 \bigl\|(Y-M_w)D^jw(Y)\bigr\|_\infty
 &\leq
 C\sigma^{1-j}\sqrt{L_\sigma},\\
 \mathbb E|D^jw(Y)|^2
 &\leq
 C\sigma^{d-2j},
 \qquad
 \|D^jw\|_\infty
 \leq
 C\sigma^{-j}.
 \end{split}
\end{equation}
For the spherical class, the factor \(\sqrt{L_\sigma}\) in the envelope
bound may be omitted.

After rescaling the local supports, these classes are indexed by the base
point, finitely many radii, and the finite-dimensional projector jet
\((P,B)\).  Their normalized covering numbers satisfy
\begin{equation}
\label{eq:ysl_one_sample_entropy_bound}
 \sup_Q
 \log N\!\left(
 \varepsilon\|F_j\|_{Q,2},
 \mathscr W_{j,\sigma},
 L_2(Q)
 \right)
 \leq
 C\log\!\left(\frac{e}{\varepsilon\sigma}\right),
 \qquad
 0<\varepsilon<1,
\end{equation}
where \(\mathscr W_{j,\sigma}\) denotes the corresponding derivative class
and \(F_j\) is its envelope.  The constant \(C_{\mathrm{ent}}\) in
\(\ell_\delta(\sigma)\) is chosen large enough to dominate these
fixed-dimensional entropy bounds.

Uniform Bernstein inequalities applied to
\eqref{eq:ysl_one_sample_moment_bounds} and
\eqref{eq:ysl_one_sample_entropy_bound} imply that, with probability at
least \(1-\delta\), the following bounds hold simultaneously for both
stages, all admissible projector jets, and \(j=0,1\):
\begin{equation}
\label{eq:ysl_one_sample_empirical_process_bounds}
 \begin{split}
 \sup_w
 \left|
 (\mathbb P_n-\mathbb P)D^jw
 \right|
 &\leq
 C\sigma^{d-j}
 \left\{
 \sqrt{\frac{\ell}{n\sigma^d}}
 +
 \frac{\ell}{n\sigma^d}
 \right\},\\
 \sup_w
 \left\|
 (\mathbb P_n-\mathbb P)
 \bigl\{(Y-M_w)D^jw(Y)\bigr\}
 \right\|
 &\leq
 C\sigma^{d+1-j}\zeta .
 \end{split}
\end{equation}
Work throughout on this event, which is uniform over the deterministic
class \(\mathfrak J_\sigma\).

Under \(n\sigma^d\geq C\ell\), the first bound in
\eqref{eq:ysl_one_sample_empirical_process_bounds}, together with the
population mass lower bound, gives
\begin{equation}
\label{eq:ysl_one_sample_empirical_mass_bound}
 \inf_w\mathbb P_nw
 \geq
 c\sigma^d.
\end{equation}
Thus all empirical denominators are bounded away from zero.

For a fixed admissible weight field, define
\(\widehat M_w = \frac{\mathbb P_n\{Yw(Y)\}} {\mathbb P_nw(Y)}\).
The centered quotient identity is
\[
 \widehat M_w-M_w
 =
 \frac{\mathbb P_n\{(Y-M_w)w(Y)\}}
      {\mathbb P_nw(Y)}.
\]
For any coordinate derivative \(D\),
\[
 \begin{split}
 D\widehat M_w-DM_w
 &=
 \frac{\mathbb P_n\{(Y-M_w)Dw(Y)\}}
      {\mathbb P_nw(Y)}
 -
 \frac{\mathbb P\{(Y-M_w)Dw(Y)\}}
      {\mathbb Pw(Y)}\\
 &\quad
 -
 (\widehat M_w-M_w)
 \frac{\mathbb P_nDw(Y)}
      {\mathbb P_nw(Y)}.
 \end{split}
\]
Using \eqref{eq:ysl_one_sample_empirical_process_bounds} and
\eqref{eq:ysl_one_sample_empirical_mass_bound} in these two identities
gives
\begin{equation}
\label{eq:ysl_one_sample_uniform_mean_bound}
 \|\widehat M_w-M_w\|_{C^j}
 \leq
 C\sigma^{1-j}\zeta,
 \qquad
 j=0,1,
\end{equation}
uniformly over both stages and all projector jets in
\(\mathfrak J_\sigma\).  Terms involving products of empirical errors are
absorbed using \(\zeta\leq c\sigma^{1/2}\).

Applying \eqref{eq:ysl_one_sample_uniform_mean_bound} to the spherical
first-stage class yields
\begin{equation}
\label{eq:ysl_one_sample_first_stage_rate}
 \|\widehat F_{r_0}-F_{r_0}\|_
 {C^j(\mathcal M^{\mathrm{pre}})}
 \leq
 C\sigma^{1-j}\zeta,
 \qquad
 j=0,1.
\end{equation}
The population signal condition and the \(C^0\) bound imply
\[
 \begin{split}
 \inf_{z\in\mathcal M^{\mathrm{pre}}}
 \|\widehat F_{r_0}(z)-z\|
 &\geq
 \inf_{z\in\mathcal M^{\mathrm{pre}}}
 \|F_{r_0}(z)-z\|
 -
 \|\widehat F_{r_0}-F_{r_0}\|_{C^0}\\
 &\geq
 c_F\sigma-C\sigma\zeta
 \geq
 \frac{c_F}{2}\sigma,
 \end{split}
\]
after reducing the constant \(c\) in the condition
\(\zeta\leq c\sigma^{1/2}\).  In particular,
\(\widehat U_0\) is well defined throughout
\(\mathcal M^{\mathrm{pre}}\).

For \(\mathcal P(v)=vv^\top/\|v\|^2\) and \(r=1,2\),
\begin{equation}
\label{eq:ysl_one_sample_projector_derivatives}
 \|D^r\mathcal P(v)\|_{\mathrm{op}}
 \leq
 C_r\|v\|^{-r}.
\end{equation}
On the population and empirical signal regions, this gives
\(\|D^r\mathcal P(v)\|_{\mathrm{op}}\leq C\sigma^{-r}\), \(r=1,2\).
Moreover, the population first-stage bounds give
\(\|DV_0\|_{C^0(\mathcal M^{\mathrm{pre}})} \leq C\sigma^{1/2}\).

Since
\(U_0=\mathcal P(V_0), \qquad \widehat U_0=\mathcal P(\widehat V_0)\),
the \(C^0\) case of
\eqref{eq:ysl_one_sample_first_stage_rate} gives
\[
 \|\widehat U_0-U_0\|_{C^0}
 \leq
 C\sigma^{-1}
 \|\widehat V_0-V_0\|_{C^0}
 \leq
 C\zeta.
\]
For the first derivative, the chain rule and
\eqref{eq:ysl_one_sample_projector_derivatives} give
\[
 \begin{split}
 \|D\widehat U_0-DU_0\|_{C^0}
 &\leq
 C\sigma^{-1}
 \|D\widehat V_0-DV_0\|_{C^0}\\
 &\quad
 +
 C\sigma^{-2}
 \|\widehat V_0-V_0\|_{C^0}
 \left(
 \|DV_0\|_{C^0}
 +
 \|D\widehat V_0-DV_0\|_{C^0}
 \right)\\
 &\leq
 C\sigma^{-1}\zeta.
 \end{split}
\]
Consequently,
\begin{equation}
\label{eq:ysl_one_sample_axis_rate}
 \|\widehat U_0-U_0\|_
 {C^j(\mathcal M^{\mathrm{pre}})}
 \leq
 C\sigma^{-j}\zeta,
 \qquad
 j=0,1.
\end{equation}

After reducing \(c\) if necessary,
\eqref{eq:ysl_one_sample_axis_rate} and
\(\zeta\leq c\sigma^{1/2}\) imply
\[
 \|\widehat U_0-U_0\|_{C^0}
 \leq
 \theta_*,
 \qquad
 \|D\widehat U_0\|_{C^0}
 \leq
 C_J\sigma^{-1/2}.
\]
Hence \((z,k,\widehat U_0(z),D\widehat U_0(z))\in\mathfrak J_\sigma\)
simultaneously for every \(z\) and \(k\).

For a \(C^1\) rank-one projector field \(P\) whose pointwise jets belong to
\(\mathfrak J_\sigma\), denote the corresponding population and empirical
cylindrical means by \(G_{k,P}\) and \(\widehat G_{k,P}\).  Thus
\(G_{k,U_0}=G_k, \widehat G_{k,\widehat U_0}=\widehat G_k\).
Because \eqref{eq:ysl_one_sample_uniform_mean_bound} holds uniformly over
the entire projector-jet class, it may be evaluated at the data-dependent
axis \(\widehat U_0\).  Therefore,
\begin{equation}
\label{eq:ysl_one_sample_random_axis_empirical_rate}
 \max_{1\leq k\leq K}
 \|\widehat G_{k,\widehat U_0}
      -G_{k,\widehat U_0}\|_
 {C^j(\mathcal M^{\mathrm{pre}})}
 \leq
 C\sigma^{1-j}\zeta,
 \qquad
 j=0,1.
\end{equation}

For the full rank-one projector neighborhood, use uniformly regular local coordinates \(\vartheta\), reducing \(\theta_*\) if necessary, and write \(g_k(z,\vartheta)=G_{k,P(\vartheta)}(z)=z+N/A\),
where \(A=\mathbb E[\Psi_{k,z,P(\vartheta)}(Y)]\) and \(N=\mathbb E[(Y-z)\Psi_{k,z,P(\vartheta)}(Y)]\).
The rescaling \(Y-z=\sigma v\) and the localization argument in Appendix~\ref{sec:ysl_common_localization_proof} give \(A\geq c\sigma^d\) and
\(\|\partial_z^\ell\partial_\vartheta^j A\|+\sigma^{-1}\|\partial_z^\ell\partial_\vartheta^j N\|\leq C\sigma^{d-\ell}\) for nonnegative \(j,\ell\) with \(j+\ell\leq2\) and \(\ell\leq1\), by uniformly integrable polynomial Gaussian envelopes.
The quotient rule yields \(\|\partial_\vartheta g_k\|+\|\partial_\vartheta^2g_k\|\leq C\sigma\) and \(\|\partial_z\partial_\vartheta g_k\|\leq C\).
Integrate along local coordinate segments and differentiate once in \(z\).
Using \(\|DP\|_{C^0}+\|DU_0\|_{C^0}\leq C\sigma^{-1/2}\), the additional second-projector-derivative terms are at most \(C\sigma^{1/2}\|P-U_0\|_{C^0}\).
Absorbing these terms gives, uniformly over the same projector class,
\begin{equation}
\label{eq:ysl_one_sample_axis_stability}
 \begin{split}
 \|G_{k,P}-G_{k,U_0}\|_{C^0}
 &\leq
 C\sigma\|P-U_0\|_{C^0},\\
 \|G_{k,P}-G_{k,U_0}\|_{C^1}
 &\leq
 C\left\{
 \|P-U_0\|_{C^0}
 +
 \sigma\|P-U_0\|_{C^1}
 \right\}.
 \end{split}
\end{equation}
Taking \(P=\widehat U_0\) and applying
\eqref{eq:ysl_one_sample_axis_rate} gives
\begin{equation}
\label{eq:ysl_one_sample_population_axis_replacement}
 \max_{1\leq k\leq K}
 \|G_{k,\widehat U_0}-G_k\|_
 {C^j(\mathcal M^{\mathrm{pre}})}
 \leq
 C\sigma^{1-j}\zeta,
 \qquad
 j=0,1.
\end{equation}
Combining
\eqref{eq:ysl_one_sample_random_axis_empirical_rate} and
\eqref{eq:ysl_one_sample_population_axis_replacement} proves
\[
 \max_{1\leq k\leq K}
 \|\widehat G_k-G_k\|_
 {C^j(\mathcal M^{\mathrm{pre}})}
 \leq
 C\sigma^{1-j}\zeta,
 \qquad
 j=0,1.
\]
Together with \eqref{eq:ysl_one_sample_axis_rate}, this proves the stated
single-scale estimates.

For the combined common-start map, set \(\Delta U=\widehat U_0-U_0\).
Then \(A_k=b_kI_D+d_kU_0\) and
\(\widehat A_k=b_kI_D+d_k\widehat U_0\).  Using \(\sum_k d_k=0\),
for any fixed reference index \(k_*\) we obtain the exact centering identity
\begin{equation}
\label{eq:ysl_one_sample_combined_identity}
 \widehat G-\widetilde G
 =
 \sum_{k=1}^K
 \widehat A_k(\widehat G_k-G_k)
 +
 \Delta U
 \sum_{k=1}^K
 d_k(G_k-G_{k_*}).
\end{equation}
The centering by \(G_{k_*}\) prevents the axis error from multiplying an
order-one ambient position.

The population expansion established above gives
\begin{equation}
\label{eq:ysl_one_sample_cross_scale_difference}
 \max_{1\leq k\leq K}
 \|G_k-G_{k_*}\|_
 {C^j(\mathcal M^{\mathrm{pre}})}
 \leq
 C\sigma^{2-j/2},
 \qquad
 j=0,1.
\end{equation}
Moreover,
\[
 \sum_{k=1}^K\|\widehat A_k\|_{C^0}
 \leq
 CW,
 \qquad
 \sum_{k=1}^K\|D\widehat A_k\|_{C^0}
 \leq
 CW\sigma^{-1/2}.
\]

For the first term in
\eqref{eq:ysl_one_sample_combined_identity}, the single-scale bounds give
\[
 \begin{aligned}
 \left\|\sum_{k=1}^K\widehat A_k(\widehat G_k-G_k)\right\|_{C^0}
 &\leq CW\sigma\zeta,\\
 \left\|\sum_{k=1}^K\widehat A_k(\widehat G_k-G_k)\right\|_{C^1}
 &\leq CW\left(\zeta+\sigma^{1/2}\zeta\right)\leq CW\zeta.
 \end{aligned}
\]
For the second term, using
\eqref{eq:ysl_one_sample_axis_rate} and
\eqref{eq:ysl_one_sample_cross_scale_difference},
\[
 \begin{aligned}
 \left\|\Delta U\sum_{k=1}^K d_k(G_k-G_{k_*})\right\|_{C^0}
 &\leq CW\sigma^2\zeta,\\
 \left\|\Delta U\sum_{k=1}^K d_k(G_k-G_{k_*})\right\|_{C^1}
 &\leq CW\left(\sigma^{-1}\zeta\,\sigma^2
                     +\zeta\,\sigma^{3/2}\right)\\
 &\leq CW\sigma\zeta.
 \end{aligned}
\]
It follows from
\eqref{eq:ysl_one_sample_combined_identity} that
\begin{equation}
\label{eq:ysl_one_sample_combined_rate}
 \|\widehat G-\widetilde G\|_
 {C^j(\mathcal M^{\mathrm{pre}})}
 \leq
 CW\sigma^{1-j}\zeta,
 \qquad
 j=0,1.
\end{equation}
Since \(X_0\) has uniformly bounded \(C^1\) norm, composition with \(X_0\)
gives
\[
 \|\widehat\Phi
      -\widetilde\Phi\|_{C^j(\mathcal M)}
 \leq
 CW\sigma^{1-j}\zeta,
 \qquad
 j=0,1.
\]

By Proposition~\ref{thm:ysl_common_population},
\(d_H\!\left( \widetilde G(\mathcal M^{\mathrm{pre}}), \mathcal M \right) \leq C\sigma^4\).
Using the common parameter domain and
\eqref{eq:ysl_one_sample_combined_rate},
\[
 \begin{split}
 d_H\!\left(
 \widehat{\mathcal M}^{\mathrm{YY}},
 \mathcal M
 \right)
 &\leq
 d_H\!\left(
 \widehat G(\mathcal M^{\mathrm{pre}}),
 \widetilde G(\mathcal M^{\mathrm{pre}})
 \right)
 +
 d_H\!\left(
 \widetilde G(\mathcal M^{\mathrm{pre}}),
 \mathcal M
 \right)\\
 &\leq
 \|\widehat G-\widetilde G\|_
 {C^0(\mathcal M^{\mathrm{pre}})}
 +
 C\sigma^4\\
 &\leq
 CW\sigma\zeta+C\sigma^4.
 \end{split}
\]
This proves the claimed Hausdorff bound.

\paragraph*{Numerical evaluation and finite query coverage}
Work on the same event.  The population estimates and
\eqref{eq:ysl_one_sample_combined_rate} give uniform \(C^1\) bounds for
\(\widehat\Phi=\widehat G\circ X_0\), since \(W\leq W_0\).
In fixed local charts of \(\mathcal M\), integration along short chart
paths gives a Lipschitz bound in the ambient distance.  For separated
base points, the uniform \(C^0\) bound gives the same conclusion.
Since \(X_0\) is uniformly bi-Lipschitz by
Lemma~\ref{lem:ysl_full_preliminary_bounds},
\(\|\widehat G(z)-\widehat G(z')\|\leq C\|z-z'\|, z,z'\in\mathcal M^{\mathrm{pre}}\),
with \(C\) uniform in the stated regime.  Let
\(\mathcal P^*_{\mathcal S}
=\{\widehat G(z):z\in\mathcal S\}\).
This is a subset of \(\widehat{\mathcal M}^{\mathrm{YY}}\), and every
point of the continuous image is within \(C\delta_{\mathcal S}\) of it.
Pairing the numerical and exact outputs at each retained query also gives
\[
 \begin{aligned}
 d_H\!\left(\mathcal P^*_{\mathcal S},
            \widehat{\mathcal M}^{\mathrm{YY}}\right)
 &\leq C\delta_{\mathcal S},\\
 d_H\!\left(\widehat{\mathcal P}^{\mathrm{YY,num}}_N,
            \mathcal P^*_{\mathcal S}\right)
 &\leq\varepsilon_{\mathrm{num}}.
 \end{aligned}
\]
The triangle inequality and the preceding continuous-image bound yield
\[
 d_H\!\left(\widehat{\mathcal P}^{\mathrm{YY,num}}_N,\mathcal M\right)
 \leq C\sigma^4+CW\sigma\zeta
       +\varepsilon_{\mathrm{num}}+C\delta_{\mathcal S}.
\]
These estimates are pathwise on the original event and use only the
retained nonempty query set.  They require no independence of the queries
from the observations and no continuous numerical interpolation.

Finally, the empirical signal lower bound makes \(\widehat U_0\) well
defined on the entire preliminary manifold.  Since the spherical,
transverse, and longitudinal kernel profiles are \(C^3\), and all empirical
denominators are positive, the map
\(\widehat G|_{\mathcal M^{\mathrm{pre}}}\) is \(C^2\).  The population map
\(\widetilde G|_{\mathcal M^{\mathrm{pre}}}\) is a \(C^2\) embedding, and
\(\|\widehat G-\widetilde G\|_ {C^1(\mathcal M^{\mathrm{pre}})} \leq CW\zeta\).
For sufficiently small fixed \(\sigma\), the condition
\(\zeta\leq c\sigma^{1/2}\) makes this a sufficiently small \(C^1\)
perturbation of the population embedding.  Therefore
\(\widehat G|_{\mathcal M^{\mathrm{pre}}}\) is a \(C^2\) embedding. This argument establishes the empirical embedding through
\(C^1\) closeness; it does not invoke the additional second-order
condition used for the reach conclusion in
Lemma~\ref{lem:ysl_common_image_geometry}. This
completes the proof.

\paragraph*{Proof of
Corollary~\ref{cor:ysl_common_fourth_order_sample_size}}
Suppose that
\[
 n
 \geq
 CW^2\sigma^{-(d+6)}\ell.
\]
Then
\[
 \sqrt{
 \frac{\ell}
      {n\sigma^d}
 }
 \leq
 CW^{-1}\sigma^3
\]
and
\[
 \sqrt{L_\sigma}\,
 \frac{\ell}
      {n\sigma^d}
 \leq
 CW^{-2}\sigma^6\sqrt{L_\sigma}
 \leq
 CW^{-1}\sigma^3
\]
for sufficiently small \(\sigma\).  Consequently,
\[
 \zeta
 \leq
 CW^{-1}\sigma^3,
 \qquad
 W\sigma\zeta
 \leq
 C\sigma^4.
\]
The same sample-size condition also implies
\(n\sigma^d \geq C\ell, \qquad \zeta \leq c\sigma^{1/2}\)
for sufficiently small \(\sigma\).  Applying
Theorem~\ref{thm:ysl_common_finite_sample} gives
\[
 d_H\!\left(
 \widehat{\mathcal M}^{\mathrm{YY}},
 \mathcal M
 \right)
 \leq
 C\sigma^4,
\]
and the same event makes
\(\widehat{\mathcal M}^{\mathrm{YY}}\) a \(C^2\) embedded manifold.

For the numerical point cloud, the bound just proved adds only
\(\varepsilon_{\mathrm{num}}+C\delta_{\mathcal S}\).
If both quantities are \(O(\sigma^4)\), the same event therefore gives
\(d_H(\widehat{\mathcal P}^{\mathrm{YY,num}}_N,\mathcal M)
\leq C\sigma^4\).

\subsection{Proof of the sample-split normal-offset proposition}
\begin{proof}
Work on the pilot event $\mathcal E_{\rm pil}$.
Let $F_\sigma(p)=p+A\sigma\nu_\sigma(p)$ for
$p\in\mathcal P_\sigma$. Since $\nu_\sigma(p)$ is a unit normal
and $A\sigma<\operatorname{reach}(\mathcal P_\sigma)$ for small
$\sigma$, the tubular-neighborhood theorem shows that
$\pi_{\mathcal P_\sigma}(F_\sigma(p))=p$. Hence
$F_\sigma$ is a $C^6$ embedding and
\[
 d(F_\sigma(p),\mathcal P_\sigma)=A\sigma.
\]
Furthermore, $F_\sigma$ is uniformly $C^2$-close to the inclusion
of $\mathcal P_\sigma$, by the uniform $C^2$ bound on
$\nu_\sigma$. The reach bound for $\mathcal P_\sigma$ gives
uniform local charts with controlled second derivatives; after
composition with $F_\sigma$, the image has a uniform chord--tangent
quadratic bound and a uniform lower bi-Lipschitz constant.
The chord--tangent characterization of reach therefore gives
\[
 \operatorname{reach}(\mathcal M^{\mathrm{pre}})\ge\tau_*>0,
\]
where $\tau_*$ is independent of $\sigma$. The normal offset also
preserves compactness, connectedness, and the absence of boundary.

The Hausdorff distance bound on the pilot implies, for every
$z=F_\sigma(p)$,
\[
 |d(z,\mathcal M)-d(z,\mathcal P_\sigma)|
 \le d_H(\mathcal P_\sigma,\mathcal M)\le B\sigma.
\]
Consequently,
\[
 (A-B)\sigma\le d(z,\mathcal M)
 \le(A+B)\sigma,\qquad z\in\mathcal M^{\mathrm{pre}}.
\]
In particular, $\mathcal M^{\mathrm{pre}}\subset
\operatorname{Tub}_{(A+B)\sigma}(\mathcal M)$.
Lemma \ref{lem:ysl_common_tube_upgrade} now provides a normal-graph parameterization
$X_0(x)=x+u_0(x)$ over $\mathcal M$ with
\[
 A-B\le \|u_0(x)/\sigma\|\le A+B,
 \quad \|u_0\|_{C^1}\le C\sigma^{1/2},
 \quad \|u_0\|_{C^2}\le C.
\]

It remains to check the first-stage signal, which cannot be inferred
from the distance bound alone without the local Gaussian calculation.
For $z=X_0(x)$, the flat-model expansion established in Section \ref{app:ysl_first_stage_bounds_proof}
before normalization of the first-stage axis gives
\[
 P_{N_x\mathcal M}\{F_{r_0}(z)-z\}
 =-\gamma_0(u_0(x)/\sigma)u_0(x)+O(\sigma^2),
\]
uniformly over the compact normal annulus above. The flat-model
coefficient $\gamma_0(a)$ is strictly positive for $a\ne0$
and continuous in $a$; hence its infimum $\gamma_*>0$ on
$\{a:A-B\le\|a\|\le A+B\}$. Therefore,
\[
 \inf_{z\in\mathcal M^{\mathrm{pre}}}
 \|F_{r_0}(z)-z\|
 \ge \gamma_*(A-B)\sigma-C\sigma^2
 \ge \tfrac12\gamma_*(A-B)\sigma
\]
for sufficiently small $\sigma$. Thus the signal condition in
Assumption \ref{ass:ysl_full_population} holds with a uniform constant.
All other radii and regularity conditions in that assumption are
unchanged.

Conditional on the pilot sample and on $\mathcal E_{\rm pil}$,
$\mathcal M^{\mathrm{pre}}$ is fixed, satisfies Assumption \ref{ass:ysl_full_population} with
uniform constants, and is independent of the refinement sample.
Corollary \ref{cor:ysl_common_fourth_order_sample_size} therefore gives, under its sample-size condition,
a conditional probability at least $1-\delta_1$ for the stated
fourth-order Hausdorff error and $C^2$ embeddedness.
Integrating over the pilot event yields
\[
 \mathbb P(\text{refinement succeeds})
 \ge (1-\delta_1)\mathbb P(\mathcal E_{\rm pil})
 \ge 1-\delta_0-\delta_1.
\]
No finite-sample reach lower bound for the refined empirical image
is asserted.
\end{proof}
\FloatBarrier

\clearpage
\begingroup
\raggedbottom
\emergencystretch=2em
\setlength{\floatsep}{8pt plus 2pt minus 2pt}
\setlength{\textfloatsep}{10pt plus 2pt minus 2pt}
\setlength{\intextsep}{8pt plus 2pt minus 2pt}
\captionsetup{font=small,skip=4pt}
\captionsetup[subfigure]{font=small,skip=2pt}
\renewcommand{\topfraction}{0.95}
\renewcommand{\bottomfraction}{0.90}
\renewcommand{\textfraction}{0.05}
\renewcommand{\floatpagefraction}{0.70}
\makeatletter
\setlength{\@fptop}{0pt}
\setlength{\@fpsep}{10pt plus 2pt minus 2pt}
\makeatother

\section{Additional Experiments and Implementation Details}
\label{app:simulation-experiments}

\paragraph*{Implementation and code}
All numerical experiments are implemented in Python,
following the fitting operators and extrapolation
constructions in Sections~\ref{Multiscale Bias Cancellation}--\ref{sec:yao_contraction}
and Appendix~\ref{sec:gaussian_density_ridge}.
The computationally intensive simulations are conducted on a 64-bit
Linux server running Ubuntu 24.04.4 LTS, equipped with AMD EPYC 7763
processors, approximately 504\,GiB of system-visible RAM, and two
NVIDIA L40 GPUs. The GPU production environment uses Python~3.12.7,
NumPy~2.1.3, SciPy~1.18.1, and CuPy~14.2.0 with CUDA~12.9.
CuPy accelerates the GPU-based fitting computations, while independent
simulation tasks are executed in parallel. Numerical computations use
double-precision floating-point arithmetic (\texttt{float64}).
Result post-processing and figure generation are performed on
64-bit Windows~11 Professional using Python~3.12.9 and
Matplotlib~3.11.1. Experiment-specific settings are given below,
and detailed numerical procedures and additional results are provided
in Appendix~\ref{app:simulation-experiments}. The Python implementation,
experiment configurations, dependency specifications, and scripts
for reproducing the reported figures are hosted at
\url{https://github.com/zzjpp/Multiscale-Bias-Cancellation-for-Accurate-Manifold-Estimation}.

We first supplement the population bias-order and finite-sample results,
then examine how scale choice affects fitting error, numerical reliability,
and coverage. Directional ablation isolates the role of component-specific
weights in the refinement construction. Shared implementation details follow
the results; experiment-specific settings are recorded separately.

\Needspace{4\baselineskip}
\subsection{Population Verification}
\label{app:population-verification}

The population experiments supplement Section~\ref{sec:population-bias} using \(\sigma\in\{0.008,\allowbreak 0.01,\allowbreak 0.0125,\allowbreak 0.016,\allowbreak 0.02,\allowbreak 0.025,\allowbreak 0.03\}\).
STSP and DR use \(h_1=3\sigma\), \(h_2=6\sigma\), and Single-scale at
\(h_2\); their two-scale weights satisfy \(w_1+w_2=1\) and
\(w_1c(h_1,\sigma)+w_2c(h_2,\sigma)=0\), with the method-specific
coefficient \(c\). Equivalently,
\(w_1=\frac{c(h_2,\sigma)}{c(h_2,\sigma)-c(h_1,\sigma)}\) and
\(w_2=-\frac{c(h_1,\sigma)}{c(h_2,\sigma)-c(h_1,\sigma)}\).
YY uses transverse radii \((1.25,5,8)\sigma\), \(r_0=3\sigma\),
\(r_2=12\sigma\sqrt{\log(1/\sigma)}\), and Single-scale at \(r_1=5\sigma\).
Its GBC-YY output uses \(A_k=a_kU_0+b_k(I-U_0)\) with
\eqref{eq:ysl_common_weights}, setting \(b_k=a_k\) for hypersurfaces.
Figure annotations are least-squares slopes of log error against
\(\log\sigma\) over all seven levels; they summarize finite-range scaling,
not independently established asymptotic limits.

Figure~\ref{fig:e1-stsp-population-appendix} shows lower STSP average
errors after correction. On the five symmetric manifolds, the Single-scale and
GBC slopes are close to two and four, respectively. On the corrugated
torus, the GBC average-error slope is approximately \(3.95\), compared
with \(3.82\) for maximum error in
Figure~\ref{fig:intro_stsp_population}. Thus, the evidence for fourth-order
scaling is stronger for the average error than for the worst queried location
over this noise range.

\begin{figure}[!htbp]
    \centering
    \includegraphics[width=0.94\linewidth]{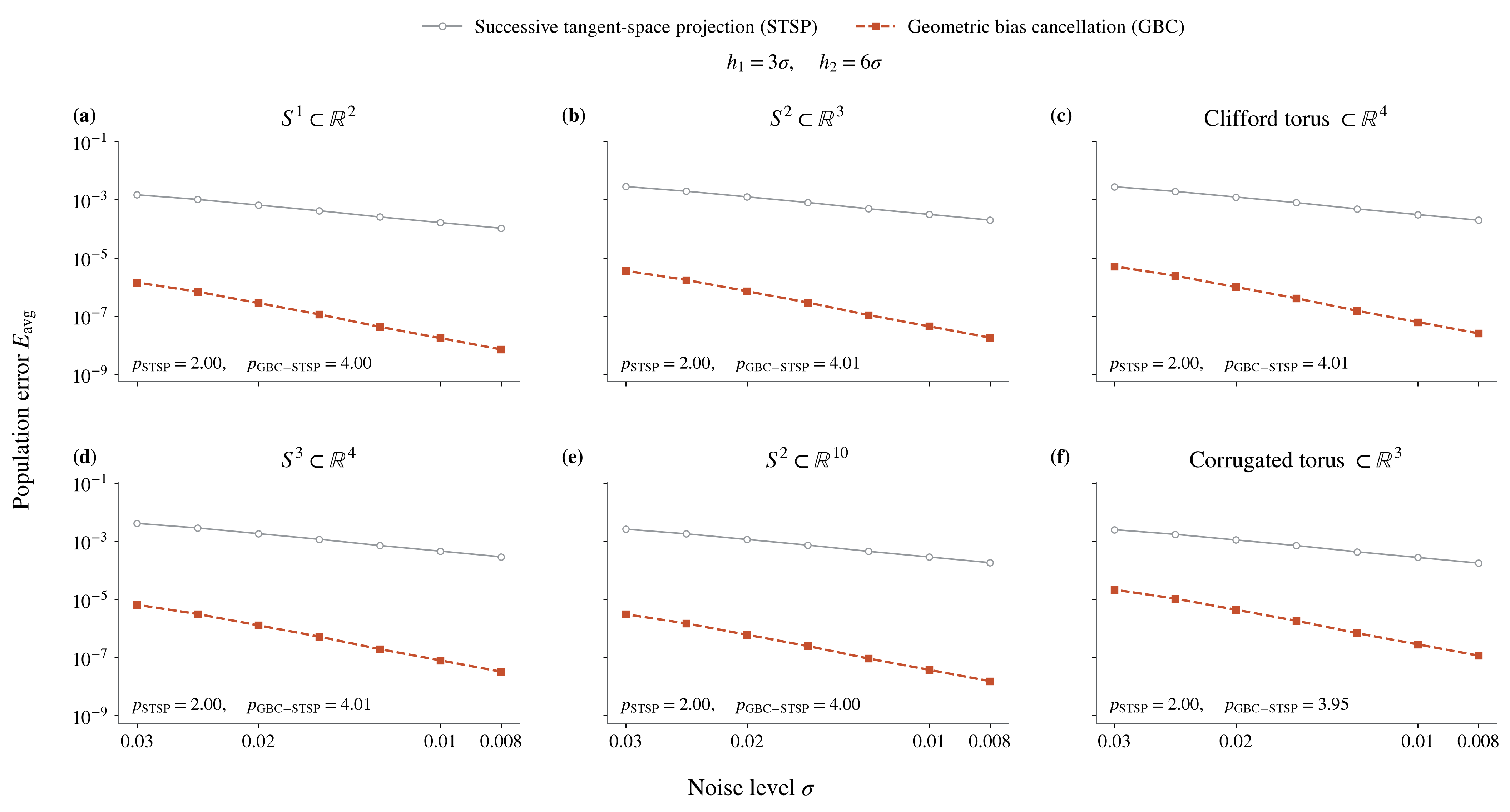}
    \caption{Population average error of STSP on six manifolds along
    \(h_1=3\sigma\), \(h_2=6\sigma\). Lines connect tested noise levels, and
    annotations give the Single-scale and GBC full-grid log--log slopes.}
    \label{fig:e1-stsp-population-appendix}
\end{figure}

Figures~\ref{fig:e1-dr-population-main} and
\ref{fig:e1-dr-population-average} show lower DR errors after correction,
with GBC slopes of approximately \(4.04\)--\(4.09\) on the five symmetric
manifolds. On the corrugated torus, the slopes are \(3.09\) for maximum
error and \(3.66\) for average error. This surface therefore demonstrates
bias reduction, but the tested range does not establish fourth-order
scaling there.

\begin{figure}[!htbp]
    \centering
    \begin{subfigure}{0.90\linewidth}
        \centering
        \includegraphics[width=\linewidth]{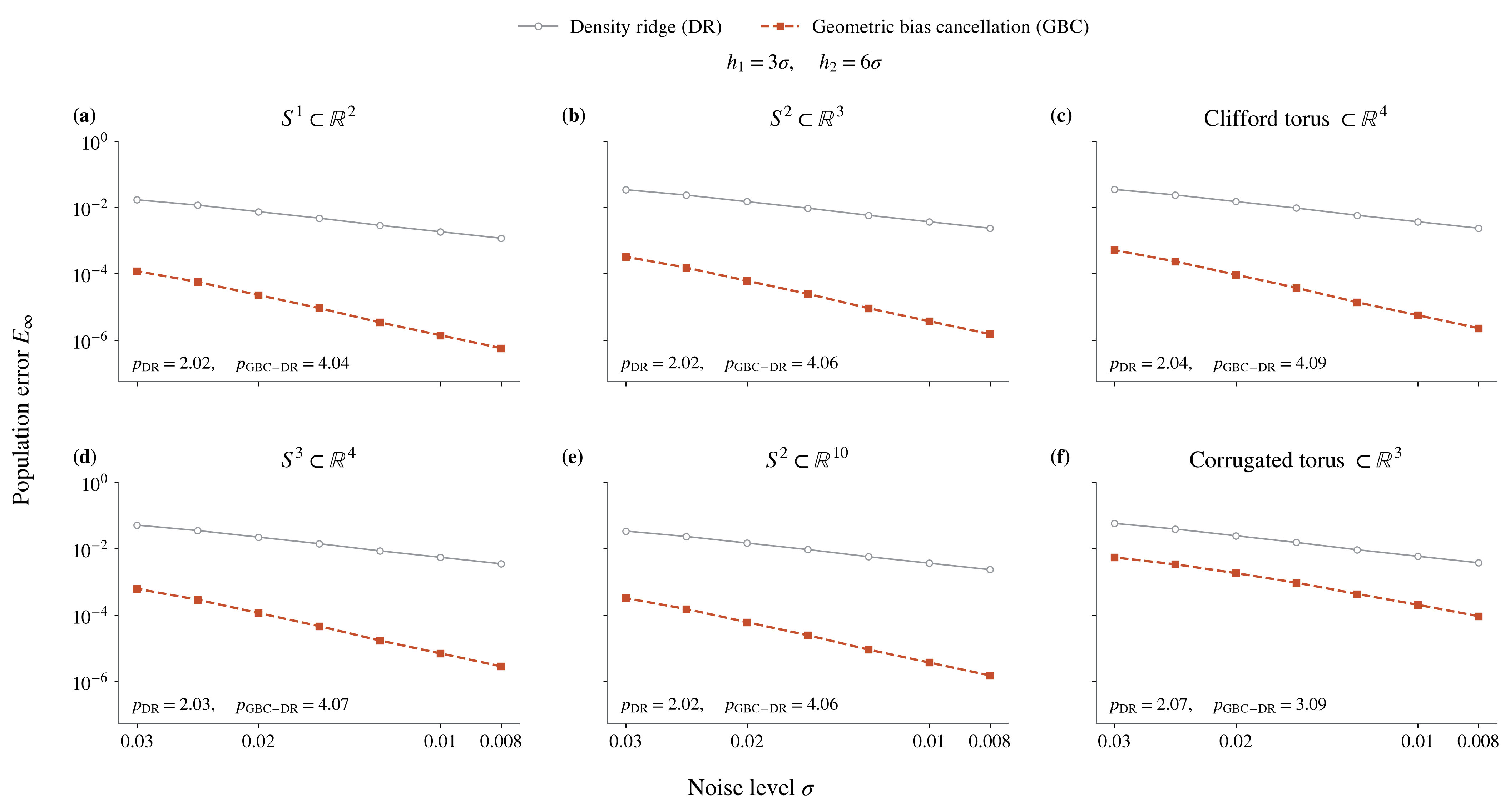}
        \caption{Maximum error \(E_\infty\).}
        \label{fig:e1-dr-population-main}
    \end{subfigure}
    \par\vspace{3pt}
    \begin{subfigure}{0.90\linewidth}
        \centering
        \includegraphics[width=\linewidth]{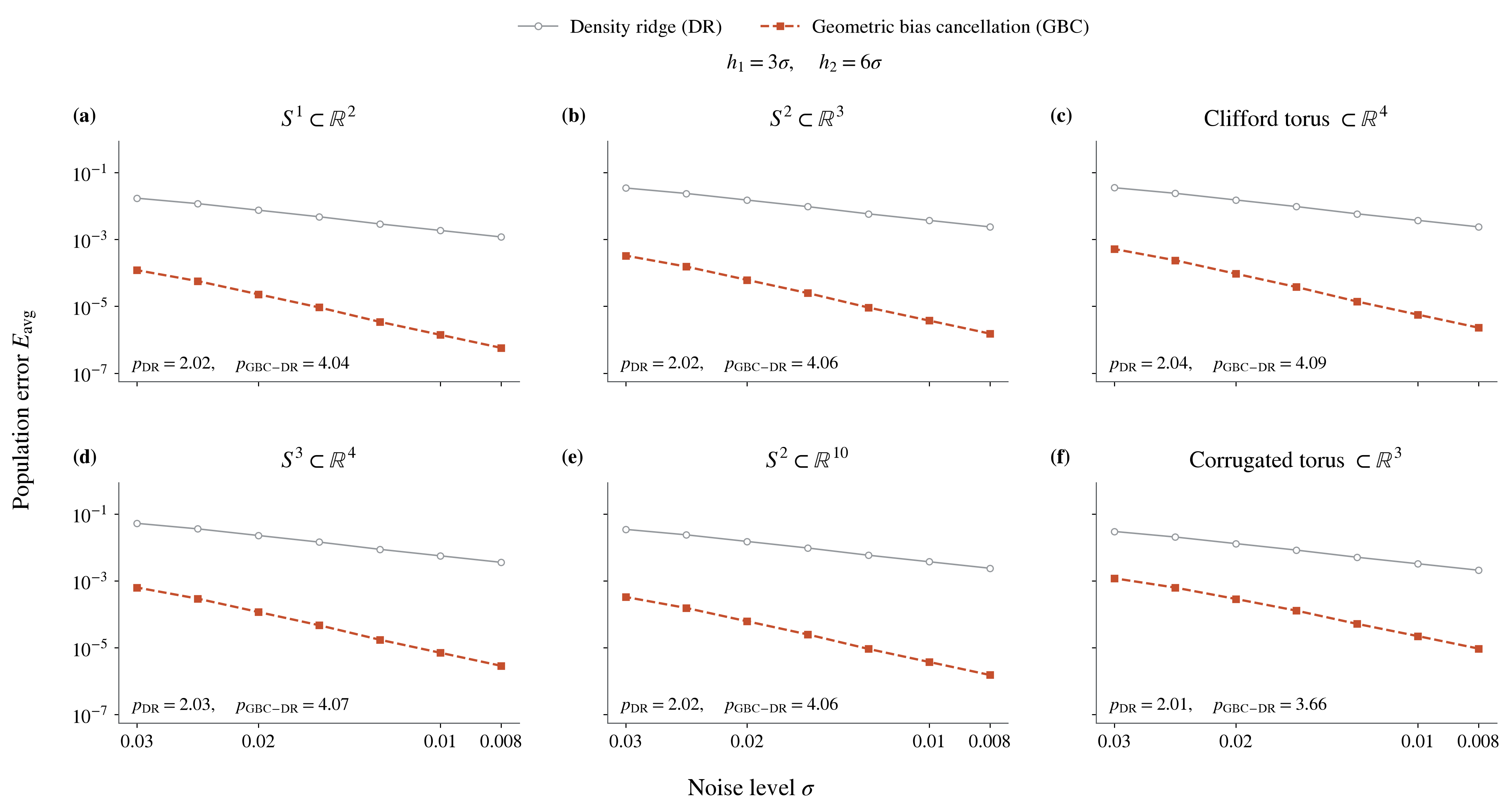}
        \caption{Average error \(E_{\mathrm{avg}}\).}
        \label{fig:e1-dr-population-average}
    \end{subfigure}
    \caption{Population DR errors along \(h_1=3\sigma\), \(h_2=6\sigma\):
    (a) maximum and (b) average point-to-manifold error.
    Lines connect noise levels, and annotations give the Single-scale and GBC
    full-grid log--log slopes with respect to \(\sigma\).}
    \label{fig:e1-dr-population-pair}
\end{figure}

Figures~\ref{fig:e1-yy-population-main} and
\ref{fig:e1-yy-population-average} show approximately second-order Single-scale
errors and smaller GBC-YY errors. The GBC-YY slopes are near four on the
symmetric manifolds and are \(3.89\) and \(4.03\) on the corrugated torus
for maximum and average error, respectively. Together, these experiments
support the predicted bias cancellation across the three constructions,
while showing that finite-range scaling depends on the geometry and error
metric. The YY results concern refinement of the supplied preliminary
manifold; the population calculations and resolution checks are described
in Appendix~\ref{app:e1-population}.

\begin{figure}[!htbp]
    \centering
    \begin{subfigure}{0.90\linewidth}
        \centering
        \includegraphics[width=\linewidth]{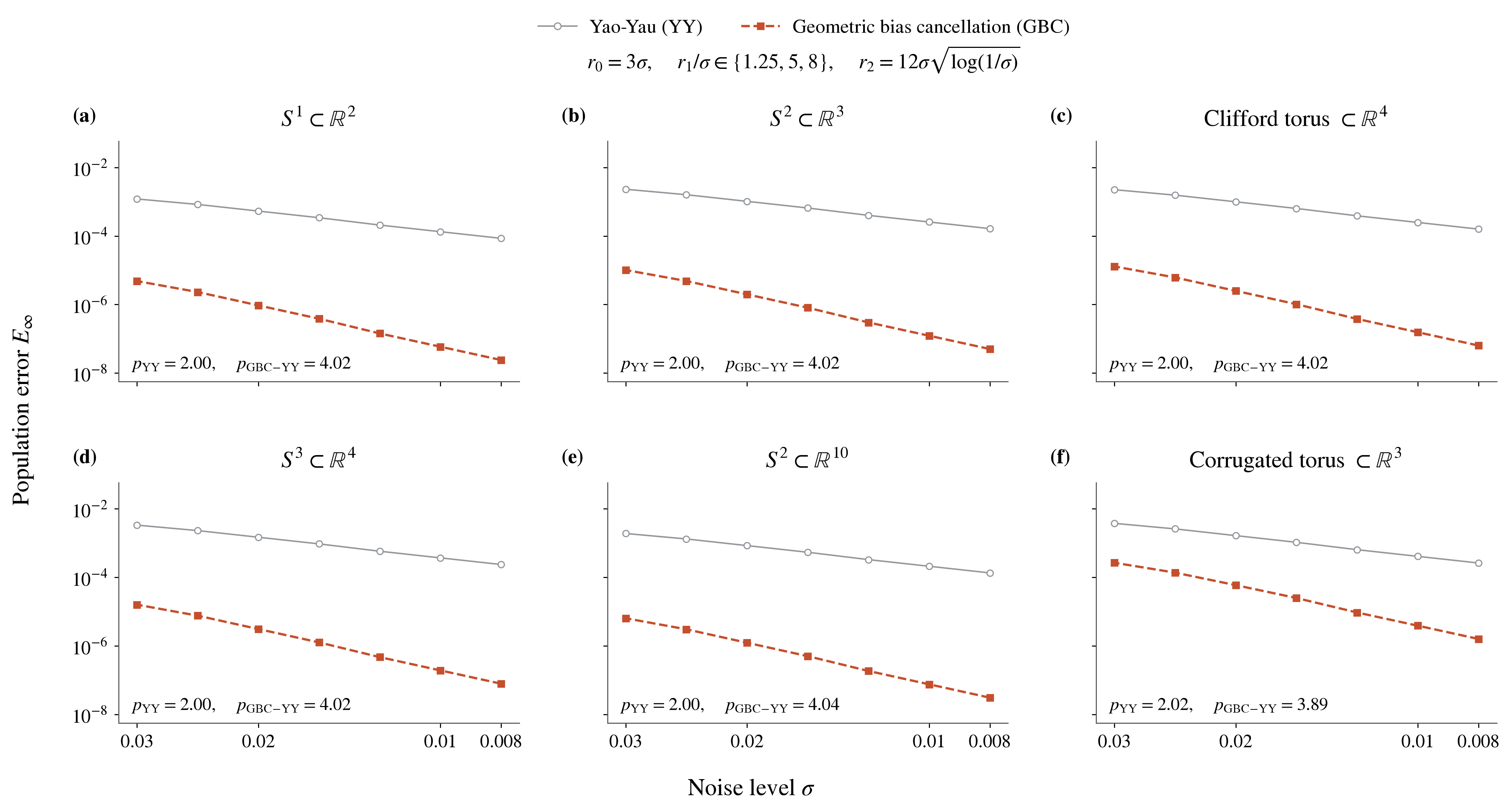}
        \caption{Maximum error \(E_\infty\).}
        \label{fig:e1-yy-population-main}
    \end{subfigure}
    \par\vspace{3pt}
    \begin{subfigure}{0.90\linewidth}
        \centering
        \includegraphics[width=\linewidth]{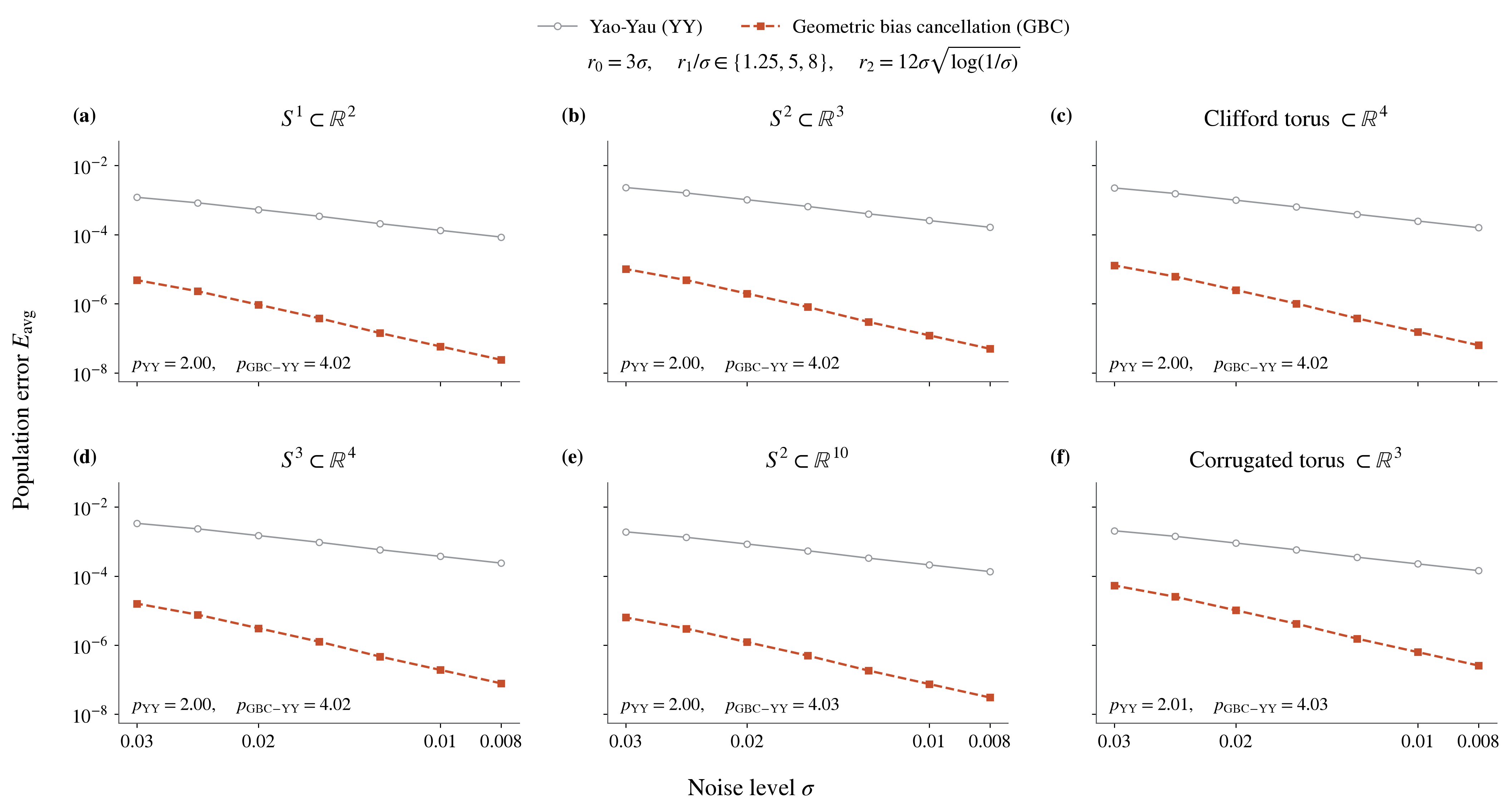}
        \caption{Average error \(E_{\mathrm{avg}}\).}
        \label{fig:e1-yy-population-average}
    \end{subfigure}
    \caption{Population YY errors: (a) maximum and (b) average point-to-manifold
    error. The transverse radii are \((1.25,5,8)\sigma\), with
    \(r_0=3\sigma\), \(r_2=12\sigma\sqrt{\log(1/\sigma)}\), and
    Single-scale at \(r_1=5\sigma\). Annotations give full-grid slopes
    with respect to \(\sigma\).}
    \label{fig:e1-yy-population-pair}
\end{figure}

\FloatBarrier
\Needspace{4\baselineskip}
\subsection{Finite-Sample Behavior at Fixed Scales}
\label{app:asymptotic-behavior}

We next examine whether the population improvement translates into
finite-sample accuracy at fixed fitting scales. The experiments use
\(\sigma=0.005\), sample sizes from \(2{,}000\) to \(500{,}000\), and
20 repetitions, with \(h_1=0.1\), \(h_2=0.2\) for STSP and DR and the
YY radii in Appendix~\ref{app:e2-finite-sample}. Curves show medians of the
finite repetition-level errors and bands show interquartile ranges;
hollow markers denote configurations that fail the reliability criterion
and are excluded from claims of an accuracy advantage. Because the scales
are fixed, these experiments assess sample-size effects at those scales,
not joint bandwidth--sample-size consistency.

Figure~\ref{fig:e2-stsp-finite-sample-appendix} complements the STSP
maximum-error results in Figure \ref{fig:e2-finite-sample-main}. At large sample sizes, the Single-scale
average errors level off while GBC continues to decrease. On several
manifolds, however, GBC initially has larger errors, and some small-sample
configurations are unreliable. The population bias reduction therefore
does not imply an improvement at every sample size.

\begin{figure}[!htbp]
    \centering
    \includegraphics[width=0.92\linewidth]{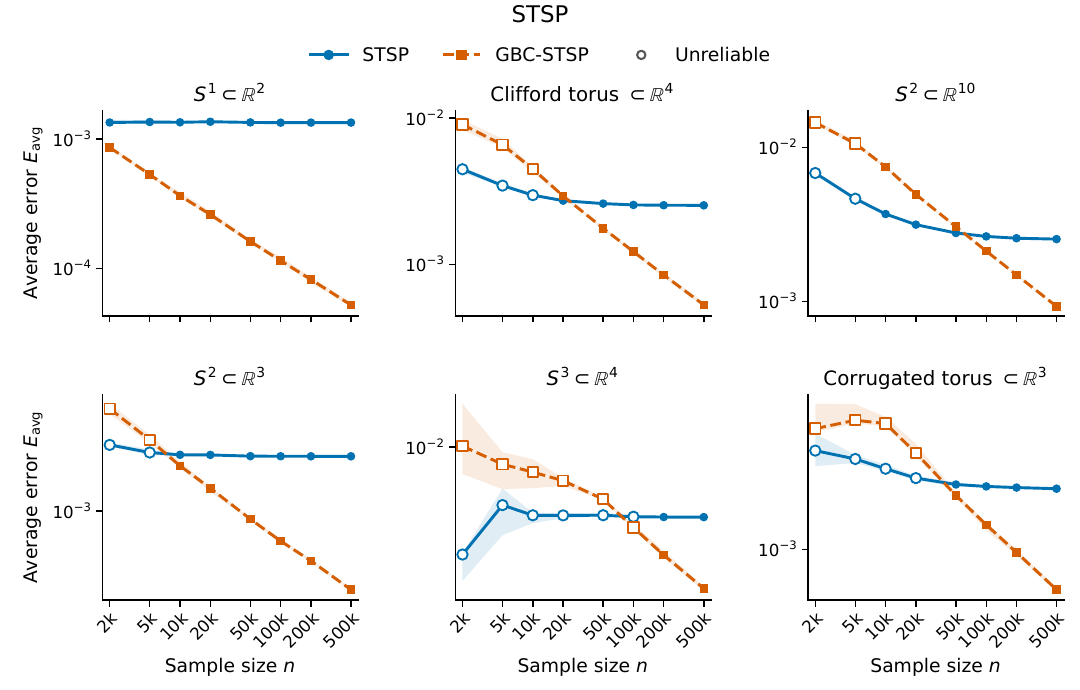}
    \caption{Finite-sample average STSP error at \(\sigma=0.005\),
    \(h_1=0.1\), and \(h_2=0.2\). Curves show medians and bands
    show interquartile ranges over 20 repetitions; hollow markers denote
    settings that fail the reliability criterion.}
    \label{fig:e2-stsp-finite-sample-appendix}
\end{figure}

Figures~\ref{fig:e2-dr-finite-sample-main} and
\ref{fig:e2-dr-finite-sample-average} show lower DR GBC errors in the
jointly reliable range. The GBC average errors begin to level off at
large \(n\), whereas the maximum errors retain more sample-size variation.
These trends are consistent with reduced deterministic bias, but do not
separately estimate the bias and stochastic contributions.

\begin{figure}[!htbp]
    \centering
    \begin{subfigure}{0.90\linewidth}
        \centering
        \includegraphics[width=\linewidth]{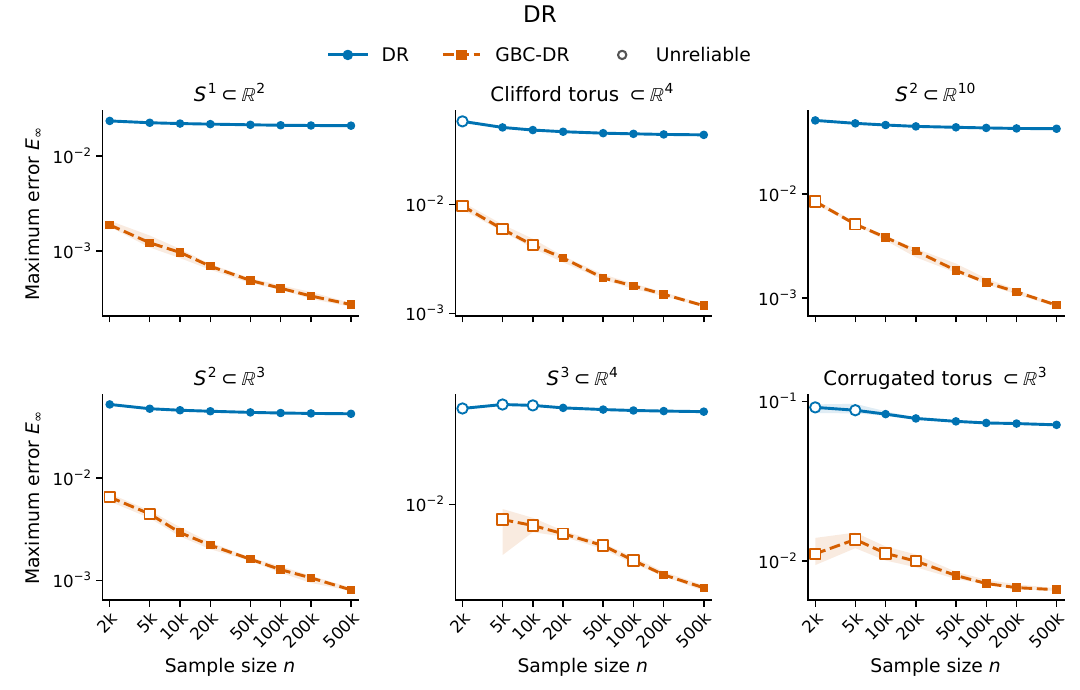}
        \caption{Maximum error \(E_\infty\).}
        \label{fig:e2-dr-finite-sample-main}
    \end{subfigure}
    \par\vspace{3pt}
    \begin{subfigure}{0.90\linewidth}
        \centering
        \includegraphics[width=\linewidth]{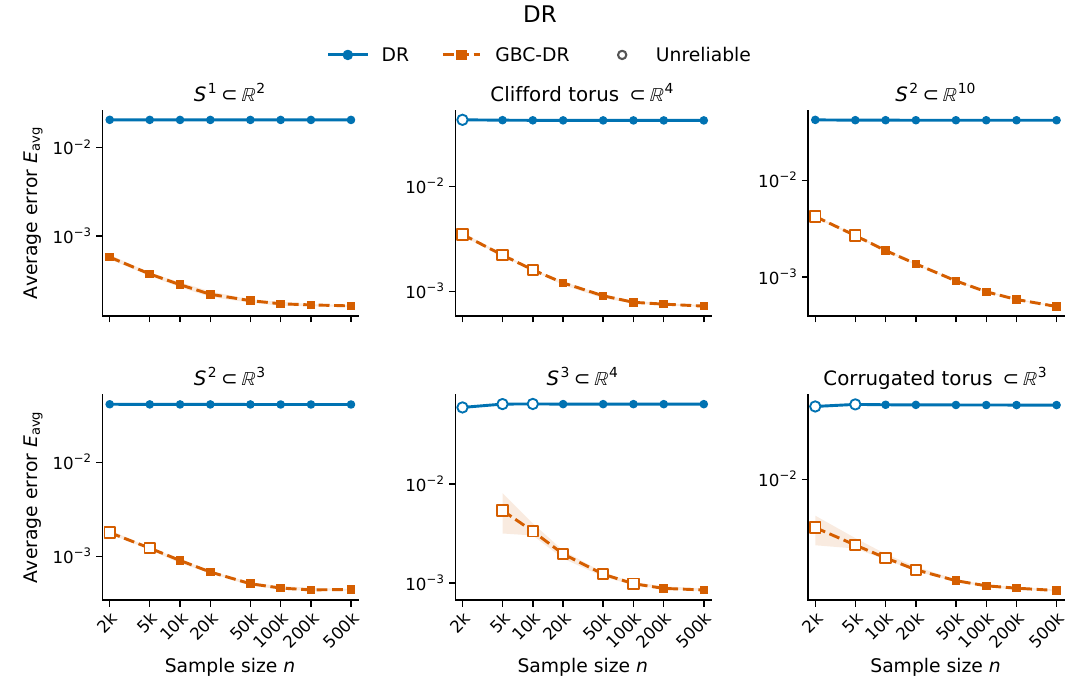}
        \caption{Average error \(E_{\mathrm{avg}}\).}
        \label{fig:e2-dr-finite-sample-average}
    \end{subfigure}
    \caption{Finite-sample DR errors at \(\sigma=0.005\), \(h_1=0.1\), and
    \(h_2=0.2\): (a) maximum and (b) average point-to-manifold error.
    Curves show medians and bands show interquartile ranges over
    20 repetitions; hollow markers denote unreliable settings.}
    \label{fig:e2-dr-finite-sample-pair}
\end{figure}

Figure~\ref{fig:e2-finite-sample-dr-yy} shows that GBC-YY improves average
error more consistently than maximum error. On the Clifford torus, even at
the largest sample size, GBC-YY has lower average error but higher maximum
error than Single-scale. This suggests that a lower overall error need not mean
better performance at the worst queried location. A possible explanation,
examined further in Appendix~\ref{app:bandwidth-robustness}, is that signed
extrapolation weights may amplify local fitting errors. 
\begin{figure}[!htbp]
    \centering
    \begin{subfigure}{0.90\linewidth}
        \centering
        \includegraphics[width=\linewidth]{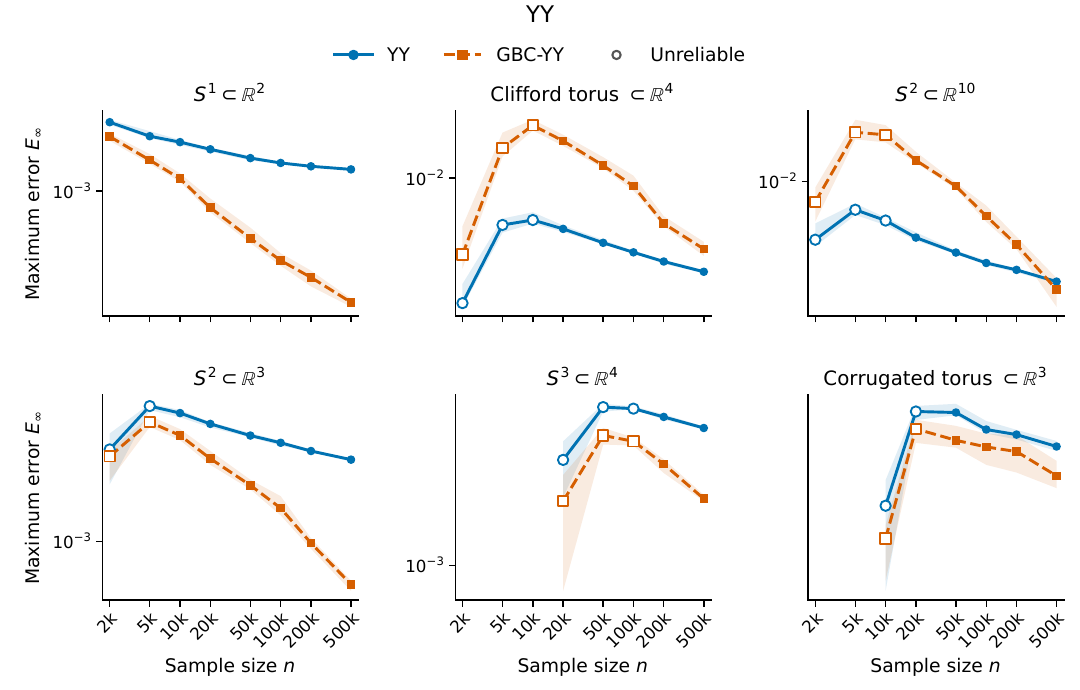}
        \caption{Maximum error \(E_\infty\).}
        \label{fig:e2-yy-maximum}
    \end{subfigure}
    \par\vspace{4pt}
    \begin{subfigure}{0.90\linewidth}
        \centering
        \includegraphics[width=\linewidth]{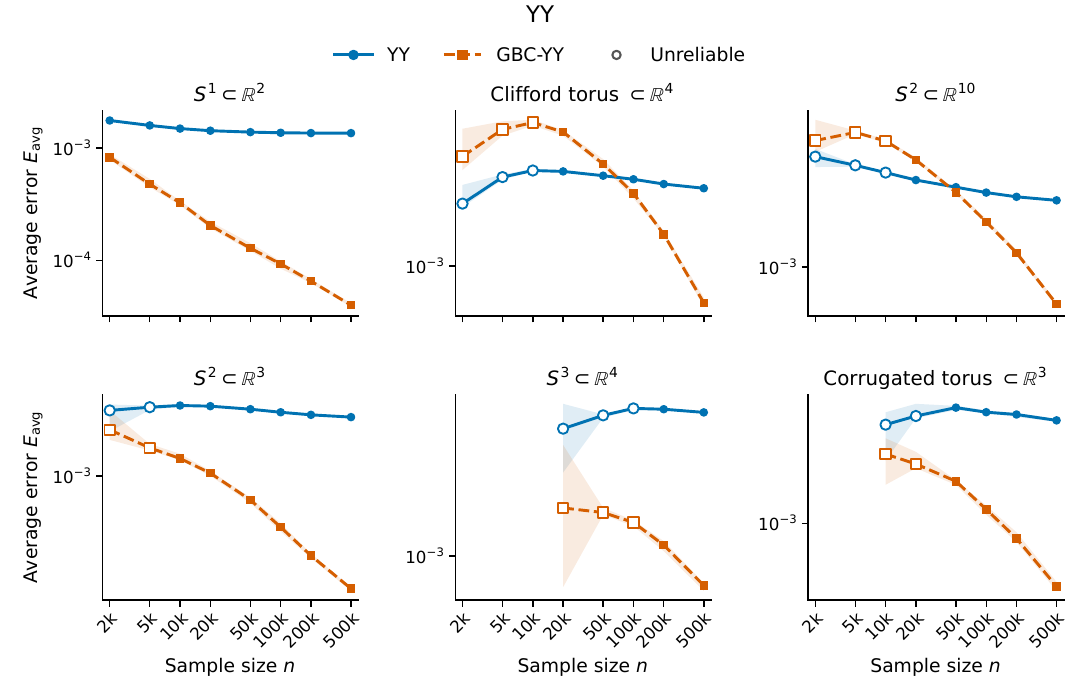}
        \caption{Average error \(E_{\mathrm{avg}}\).}
        \label{fig:e2-yy-average}
    \end{subfigure}
    \caption{Finite-sample YY errors at \(\sigma=0.005\) and the fixed radii
    in Appendix~\ref{app:e2-finite-sample}: (a) maximum and (b) average
    point-to-manifold error. Curves show medians and bands show
    interquartile ranges over 20 repetitions; hollow markers denote
    unreliable settings.}
    \label{fig:e2-finite-sample-dr-yy}
\end{figure}

\FloatBarrier
\Needspace{4\baselineskip}
\subsection{Bandwidth Robustness and Reverse Coverage}
\label{app:bandwidth-robustness}
This section complements Section~\ref{sec:e3-bandwidth-robustness} by examining
both the accuracy gain from multiscale cancellation and the range of
reliable scales over which low error is maintained. We evaluate STSP,
DR, and YY at fixed \(\sigma=0.015\); reverse coverage is assessed
separately.
The six-panel order is \(\mathbb S^1\subset\mathbb R^2\),
\(\mathbb S^2\subset\mathbb R^3\), the Clifford torus
\(\mathbb T_{\mathrm C}^2\subset\mathbb R^4\),
\(\mathbb S^3\subset\mathbb R^4\), \(\mathbb S^2\subset\mathbb R^{10}\), and the
corrugated torus \(\mathbb T_{\mathrm{corr}}^2\subset\mathbb R^3\),
with sample sizes \(n=(12000,24000,24000,60000,24000,36000)\).

The sample sizes were chosen based on preliminary
assessments of numerical stability and computational
feasibility. Within each manifold, all methods were
evaluated using the same sample-size budget and
reliability criteria. 
For STSP and DR, the parent bandwidth is
\(h_2\), with \(h_1=h_2/2\); for YY, it is \(r_1\), with transverse
radii \((0.4,1,1.6)r_1\), \(r_0=24\sigma\), and
\(r_2=12\sigma\sqrt{\log(1/\sigma)}\). All use the 28-point grid in
Appendix~\ref{app:experiment-specific}. We compare GBC with Single-scale at
each constituent scale, retaining the existing fitting and failure rules:
STSP/DR matching uses centered local \(h_1\)-PCA, and YY uses its
empirical first-stage axis and prescribed preliminary manifold.

\paragraph*{Forward fitting error}
Each setting uses 20 repetitions and 256 base locations for the first
five manifolds or 512 for the corrugated torus. For successful queries
\(\mathcal S\ne\varnothing\), the errors are
\(E_{\mathrm{avg}}=|\mathcal S|^{-1}\sum_{j\in\mathcal S}
 d(\widehat z_j,\mathcal M)\) and
\(E_\infty=\max_{j\in\mathcal S}d(\widehat z_j,\mathcal M)\).
Each estimator retains its own success set. Plotted medians include all
repetitions with finite errors, including partial successes; an empty
success set gives an undefined error. Filled markers require an aggregate
query-failure rate of at most \(1\%\) and at least 18 complete repetitions.
Other finite values use hollow markers, and undefined values leave gaps.
Comparisons use jointly reliable configurations at the same parent scale;
reported ranges refer to tested grid points. We assess scale robustness
from sustained low error across neighboring reliable scales, not just a
lower minimum or the total number of successful configurations. The
maximum error is over retained queries, not a continuous fitted manifold.

For STSP, Figure~\ref{fig:e3-stsp-average-appendix}
complements the maximum-error results in
Figure~\ref{fig:e3-bandwidth-robustness-main}:
GBC has lower intermediate-scale errors and retains an advantage over
both Single-scale baselines across several reliable bandwidths on each manifold.
This is consistent with bias cancellation allowing larger neighborhoods
without the same leading geometric bias as Single-scale. On
\(\mathbb S^3\subset\mathbb R^4\), the average-error advantage extends over the
tested jointly reliable bandwidths \(h_2=0.45\)--\(1.8\). At
\(h_2=0.82\), the median error is \(1.58\times10^{-3}\) for GBC,
compared with \(1.53\times10^{-2}\) and \(6.27\times10^{-2}\) for
Single-scale at \(h_1\) and \(h_2\). GBC is reliable at 10 of 28 scales on this
manifold; the substantial low-error range lies within this reliable region.

\begin{figure}[!htbp]
    \centering
    \includegraphics[width=0.94\linewidth]{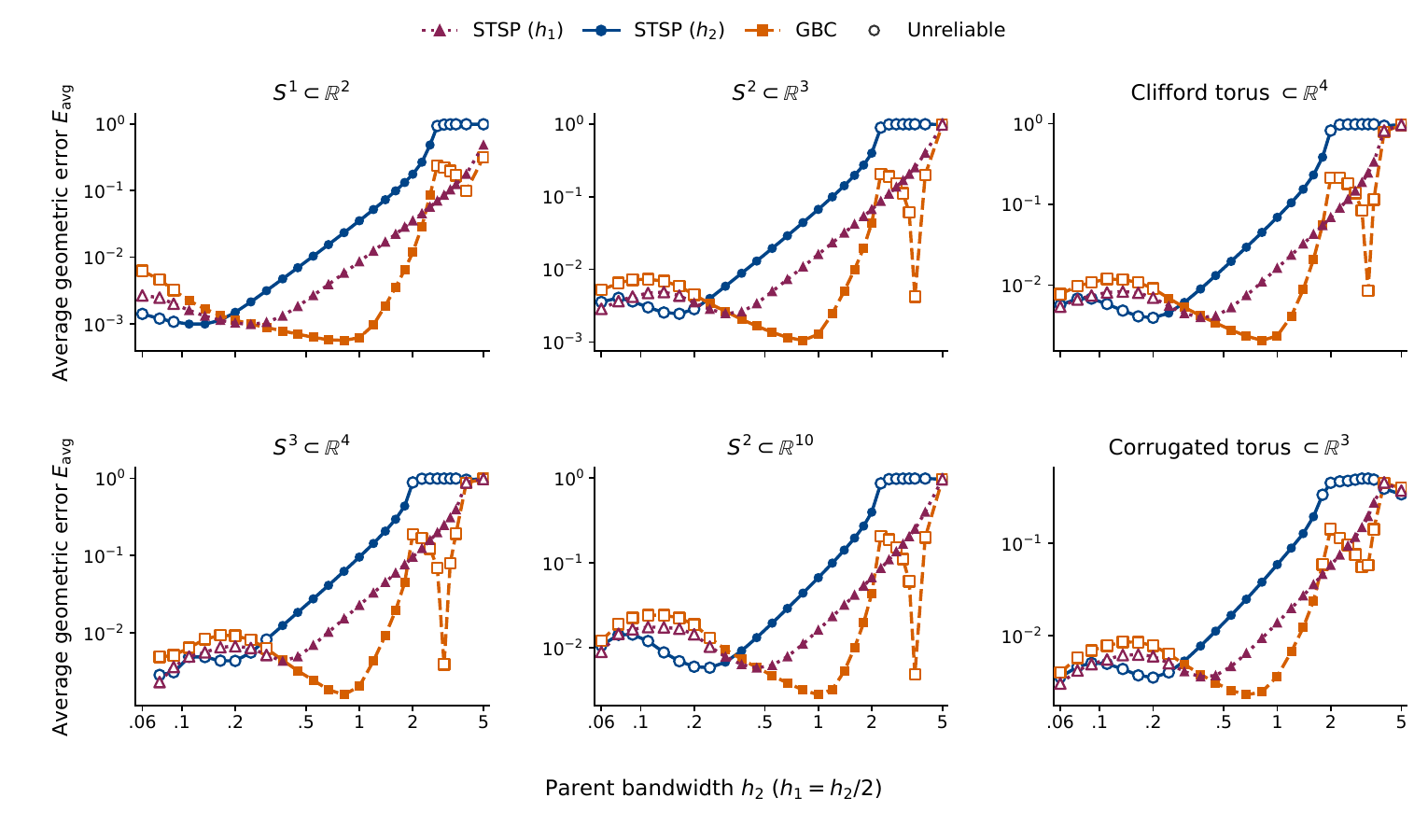}
    \caption{STSP average forward error as \(h_2\) varies, with \(h_1=h_2/2\).
    Curves compare Single-scale at \(h_1\), Single-scale at \(h_2\), and GBC
    under the sample sizes specified in Section \ref{app:bandwidth-robustness}. Values are medians over repetitions
    with finite errors; hollow markers denote unreliable settings.}
    \label{fig:e3-stsp-average-appendix}
\end{figure}

DR also improves accuracy across neighboring reliable scales
(Figure~\ref{fig:e3-dr-bandwidth-appendix}). GBC has lower average and
maximum errors than both Single-scale baselines at the jointly reliable
bandwidths \(h_2=0.2\)--\(0.45\) on \(\mathbb S^2\subset\mathbb R^3\) and
\(0.3\)--\(0.45\) on \(\mathbb S^2\subset\mathbb R^{10}\).
The benefit is more restricted on \(\mathbb S^3\subset\mathbb R^4\): despite
\(n=60000\), GBC is reliable at only 2 of 28 scales; eight further
values are finite but unreliable and 18 are undefined. Only
\(h_2=0.37\) is jointly reliable for all three curves, so this geometry
does not support a broad reliable low-error region for DR. Recorded
failures include insufficient local support, poorly separated PCA
eigenvalues, loss of normal transversality, and unsuccessful Newton
solves; their relative contributions are not isolated here.

\begin{figure}[!htbp]
    \centering
    \begin{subfigure}{0.90\linewidth}
        \centering
        \includegraphics[width=\linewidth]{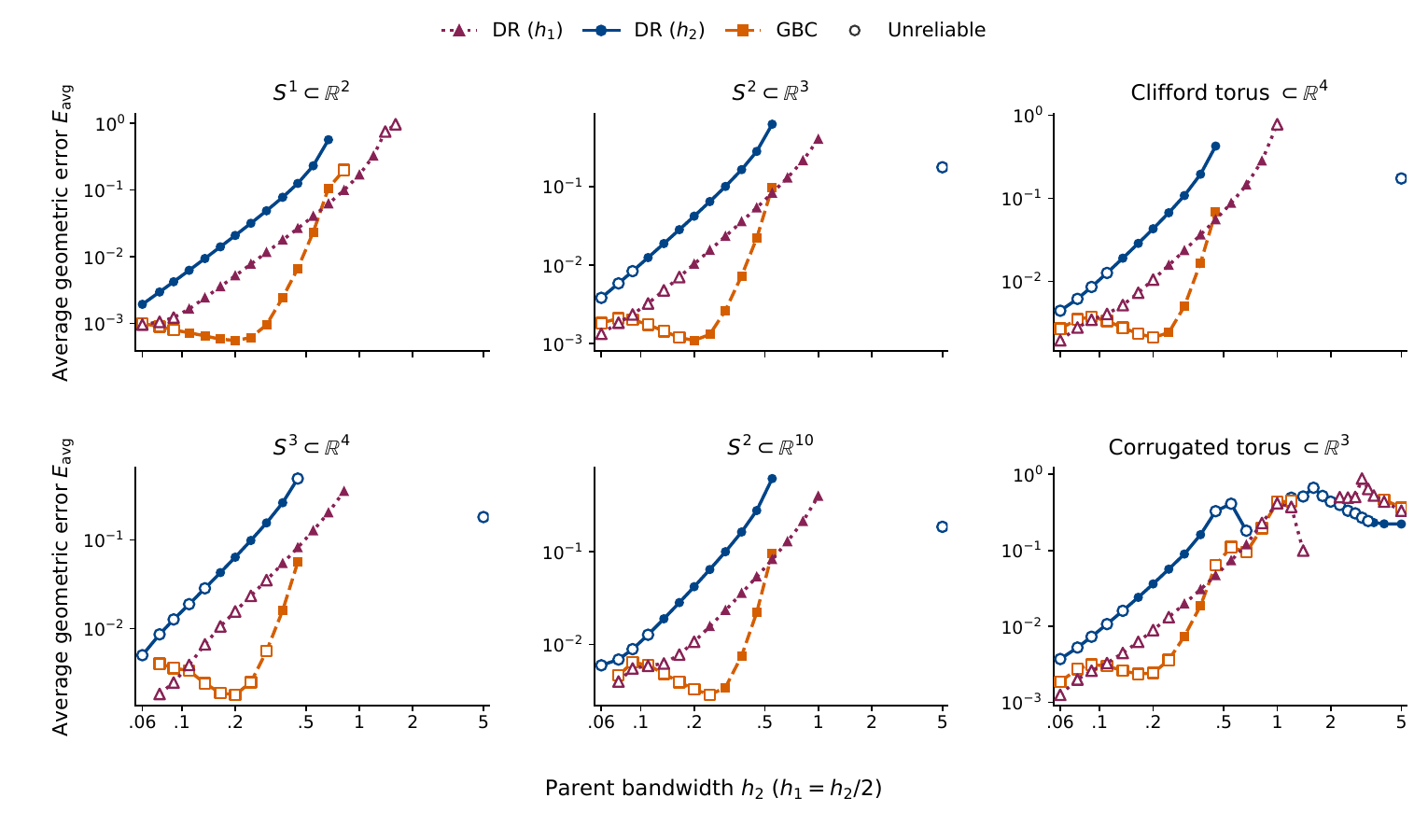}
        \caption{Average forward error \(E_{\mathrm{avg}}\).}
        \label{fig:e3-dr-average-appendix}
    \end{subfigure}
    \par\vspace{3pt}
    \begin{subfigure}{0.90\linewidth}
        \centering
        \includegraphics[width=\linewidth]{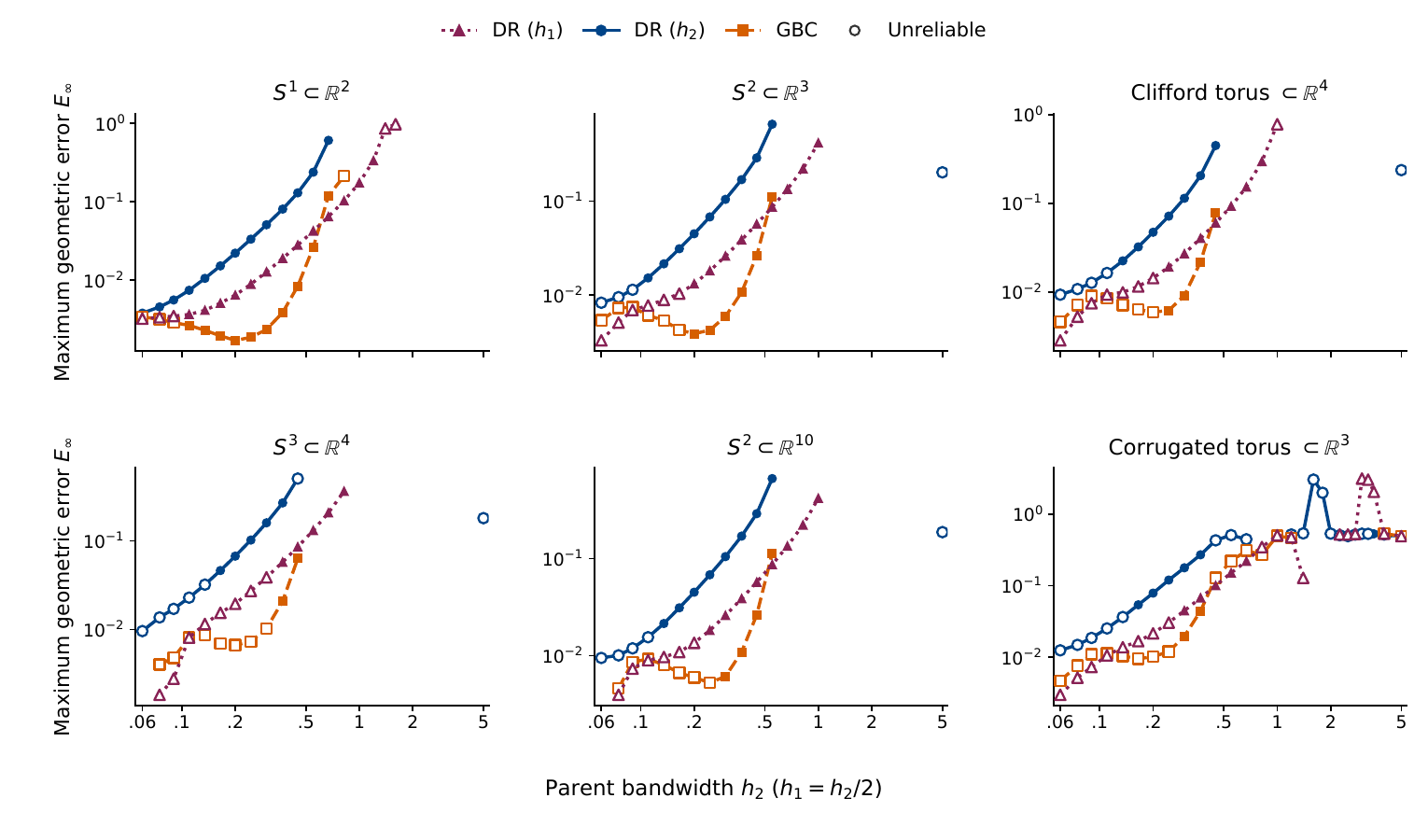}
        \caption{Maximum forward error \(E_\infty\).}
        \label{fig:e3-dr-max-appendix}
    \end{subfigure}
    \caption{DR forward errors under the sample sizes in Appendix~\ref{app:bandwidth-robustness}: (a) average and
    (b) maximum error as \(h_2\) varies, with \(h_1=h_2/2\).
    Curves compare both Single-scale baselines with GBC. Values are medians
    over repetitions with finite errors; hollow markers denote unreliable
    settings and gaps denote undefined values.}
    \label{fig:e3-dr-bandwidth-appendix}
\end{figure}

\paragraph*{YY accuracy and directional weights}
For YY, Figure~\ref{fig:e3-yy-bandwidth-appendix} shows reliable low-error
ranges, with more consistent gains in average than maximum error. On
\(\mathbb S^3\subset\mathbb R^4\), GBC-YY is reliable at 17 of 28 scales and
outperforms all three Single-scale fits at the jointly reliable radii
\(r_1=0.55\)--\(1\) for average error and \(0.55\)--\(0.82\) for maximum
error. Average error also improves over some scales on the Clifford torus
and \(\mathbb S^2\subset\mathbb R^{10}\), but neither has two neighboring
jointly reliable radii where GBC-YY beats every Single-scale fit in maximum
error. At \(r_1=0.55\), all six comparisons have zero query failures.
Relative to Single-scale at \(0.4r_1\), GBC-YY lowers the median maximum error
from \(0.00361\) to \(0.00227\) on \(\mathbb S^1\), from \(0.00807\) to
\(0.00496\) on \(\mathbb S^2\subset\mathbb R^3\), and from \(0.01350\) to
\(0.00873\) on \(\mathbb S^3\subset\mathbb R^4\). It instead raises this error
from \(0.00986\) to \(0.02056\) on the Clifford torus and from
\(0.01008\) to \(0.01682\) on \(\mathbb S^2\subset\mathbb R^{10}\), although
the latter's average error falls from \(5.69\times10^{-3}\) to
\(4.51\times10^{-3}\). These are losses of maximum-error accuracy, not
numerical failures or effects of excluding failed queries.

To quantify the sensitivity introduced by signed extrapolation, write
\(a_k=w_k^\parallel\) and \(b_k=w_k^\perp\), consistently with
\eqref{eq:ysl_common_weights}, and define
\begin{equation}
\begin{aligned}
W_\parallel(r_1)&=\sum_{k=1}^{3}|a_k|,\qquad
W_\perp(r_1)=\sum_{k=1}^{3}|b_k|,\\
W(r_1)&=W_\parallel(r_1)+W_\perp(r_1).
\end{aligned}
\label{eq:e3-yy-weight-sums}
\end{equation}
Table~\ref{tab:e3-yy-weight-sums} reports \(W\) for the bandwidth-robustness grid,
using \(\sigma=0.015\), the polynomial exponent \(p=6\), and transverse
radii \((0.4,1,1.6)r_1\). The weights follow the moment constraints and
the minimum-Euclidean-norm axial rule specified in
Appendix~\ref{app:e1-population}. These entries are deterministic
coefficient diagnostics, not fitting errors or matrix condition numbers;
they are defined even when an empirical fit fails the reliability test.

The role of \(W\) is a worst-case error bound. At a fixed preliminary
point and common axis \(U\), put \(A_k=a_kU+b_k(I-U)\). For perturbations
\(e_k\) of the three fitted outputs,
\begin{equation}
\left\|\sum_{k=1}^{3}A_ke_k\right\|
\leq \sum_{k=1}^{3}\max\{|a_k|,|b_k|\}\,\|e_k\|
\leq W(r_1)\max_k\|e_k\|.
\label{eq:e3-yy-weight-amplification}
\end{equation}
For hypersurfaces, \(b_k=a_k\) and \(A_k=a_kI\), so the sharper scalar
bound uses \(\sum_k|a_k|=W/2\). The doubled value in the table follows
the common two-block definition of \(W\); it does not mean that the
algorithm applies the scalar correction twice. This display holds for a fixed axis. The empirical axis affects both the
scale-specific fits and the matrices \(A_k\); these effects are treated
jointly in the finite-sample refinement analysis. Moreover,
all scales reuse the same observations, so their errors can be correlated:
\(W\) is not an estimated variance-inflation factor.

Across the tested grid, $W$ decreases and then levels off.
On the Clifford torus, it decreases from $12.2164$ at
$r_1=0.06$ to $4.5890$ at $r_1=5$; on
$\mathbb S^2\subset\mathbb R^{10}$, the corresponding values
are $13.1544$ and $4.5885$.
At $r_1=0.55$, the total weight sums are $4.6445$
and $4.6102$, respectively, with complementary weight
sums of approximately $3.20$ and $3.17$.
By comparison, the hypersurface values are approximately
$2.89$--$2.90$.
Thus, in higher codimension, the complementary weights
can amplify scale-specific fitting errors more strongly,
even when the leading population bias is canceled.

The effect of this amplification also depends on the
underlying geometry. The sphere $\mathbb S^2\subset\mathbb R^3$
and the corrugated torus have identical weights throughout
the grid, yet at $r_1=0.55$, GBC-YY reduces the maximum
error on the sphere while increasing it from $0.01525$
to $0.01637$ on the corrugated torus.
Together with the stronger improvements in average error,
these results suggest that signed extrapolation weights
interact with geometry-dependent fitting errors,
particularly at locations where local errors are large.

\begin{figure}[!htbp]
    \centering
    \begin{subfigure}{0.90\linewidth}
        \centering
        \includegraphics[width=\linewidth]{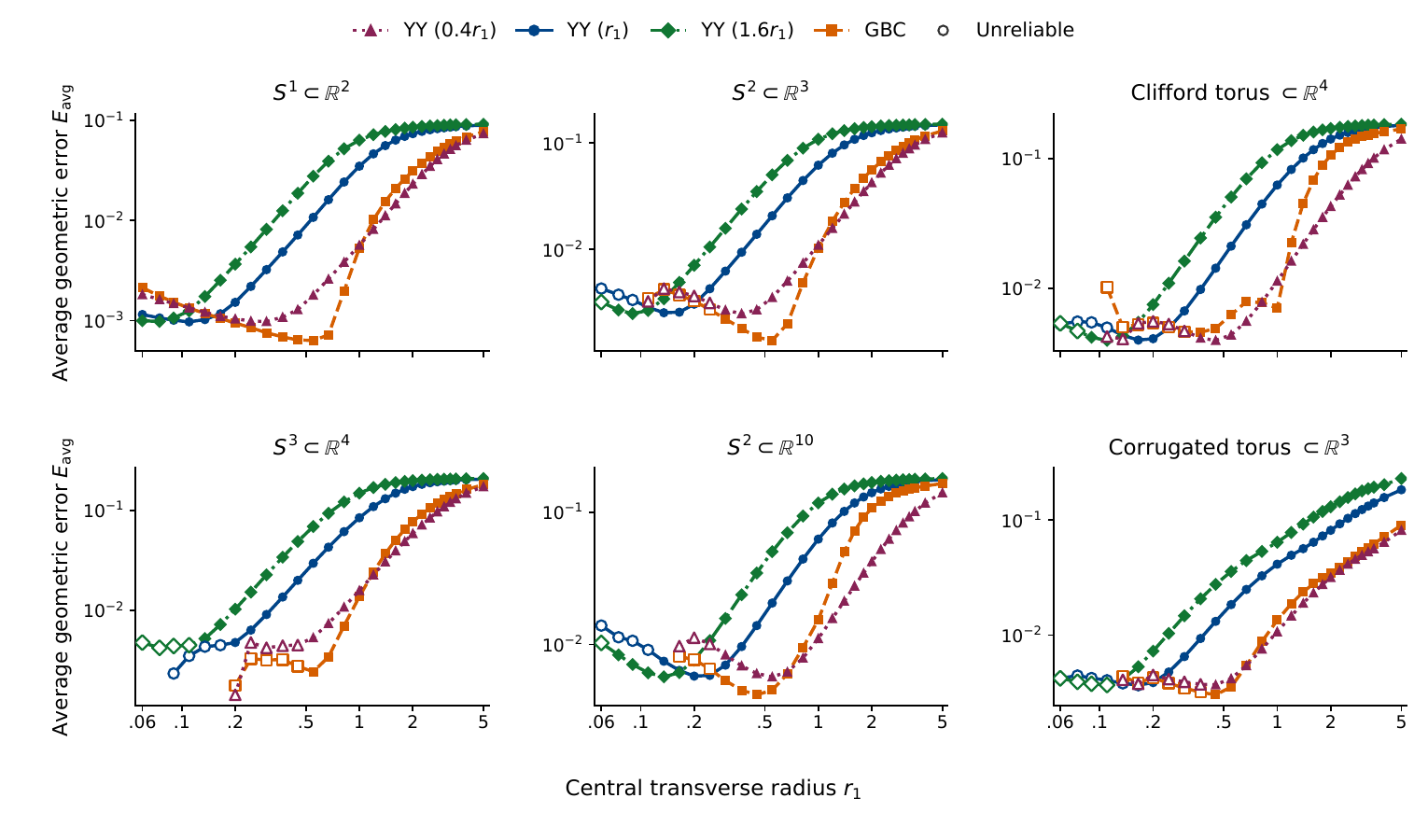}
        \caption{Average forward error \(E_{\mathrm{avg}}\).}
        \label{fig:e3-yy-average-appendix}
    \end{subfigure}
    \par\vspace{3pt}
    \begin{subfigure}{0.90\linewidth}
        \centering
        \includegraphics[width=\linewidth]{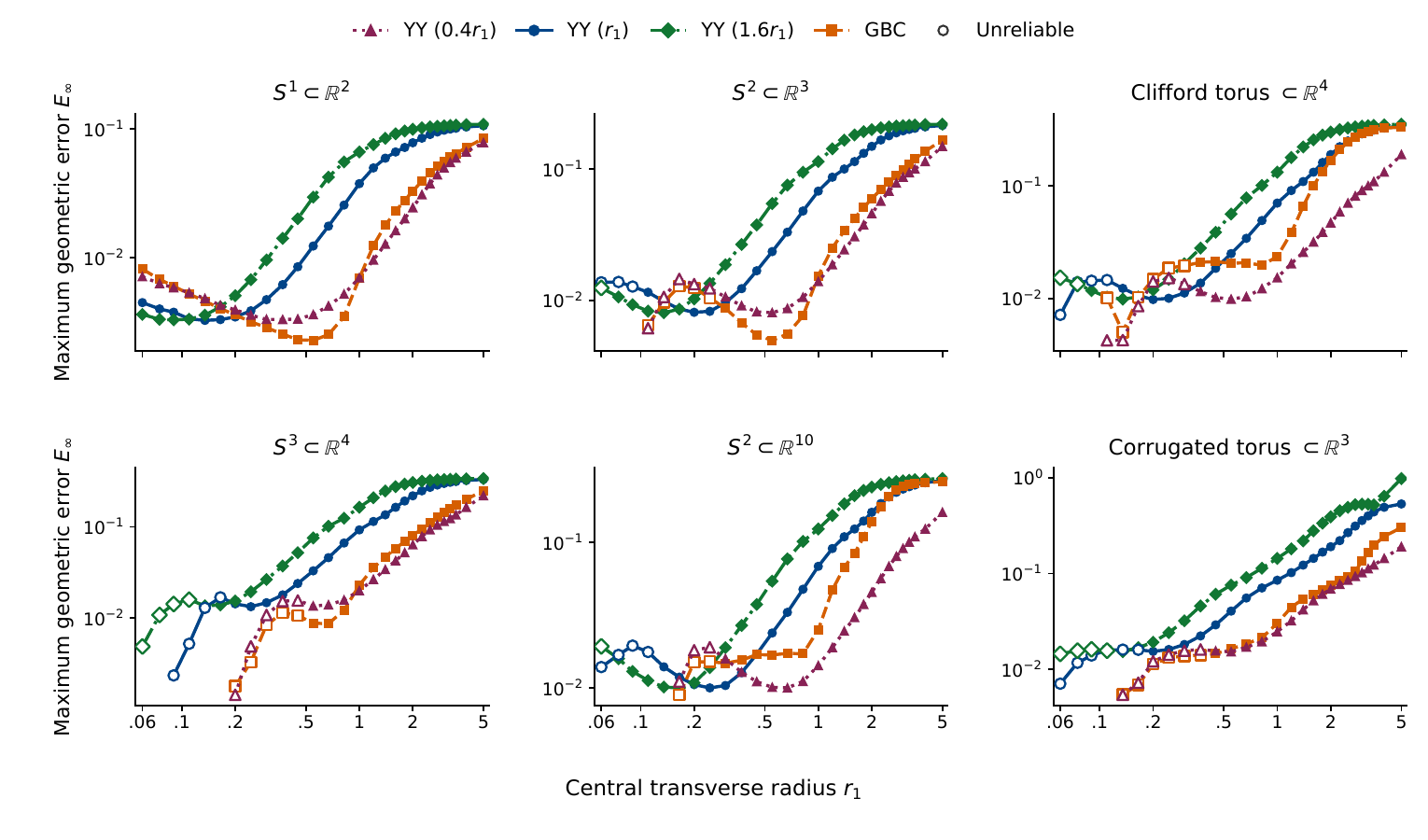}
        \caption{Maximum forward error \(E_\infty\).}
        \label{fig:e3-yy-max-appendix}
    \end{subfigure}
    \caption{YY forward errors under the sample sizes specified in Appendix~\ref{app:bandwidth-robustness}: (a) average and
    (b) maximum error as \(r_1\) varies. Curves compare Single-scale at
    \(0.4r_1\), \(r_1\), and \(1.6r_1\) with GBC-YY. Values are medians
    over repetitions with finite errors; hollow markers denote
    unreliable settings.}
    \label{fig:e3-yy-bandwidth-appendix}
\end{figure}

\begin{table}[!htbp]
    \centering
    \small
    \caption{Total absolute YY extrapolation weights
    on the bandwidth-robustness grid, with $\sigma=0.015$, $p=6$,
    and transverse radii $(0.4,1,1.6)r_1$.
    Values are rounded to four decimal places.
    For hypersurfaces, the sharper scalar
    perturbation bound uses $W/2$.}
    \label{tab:e3-yy-weight-sums}
    \setlength{\tabcolsep}{3.5pt}
    \renewcommand{\arraystretch}{1.05}
    \begin{tabular*}{\linewidth}{@{\extracolsep{\fill}}r*{6}{r}@{}}
    \toprule
    \(r_1\) & \(\mathbb S^1\subset\mathbb R^2\)
    & \(\mathbb S^2\subset\mathbb R^3\)
    & \shortstack{Clifford\\torus}
    & \(\mathbb S^3\subset\mathbb R^4\)
    & \(\mathbb S^2\subset\mathbb R^{10}\)
    & \shortstack{Corrugated\\torus}\\
    \midrule
0.060 & 4.5425 & 4.6537 & 12.2164 & 4.7650 & 13.1544 & 4.6537 \\
0.075 & 3.9417 & 4.0129 & 8.9567 & 4.0841 & 9.0795 & 4.0129 \\
0.090 & 3.6154 & 3.6649 & 7.4088 & 3.7143 & 7.2747 & 3.6649 \\
0.110 & 3.3703 & 3.4034 & 6.3552 & 3.4365 & 6.1246 & 3.4034 \\
0.135 & 3.2034 & 3.2254 & 5.6970 & 3.2474 & 5.4563 & 3.2254 \\
0.165 & 3.0944 & 3.1091 & 5.2969 & 3.1239 & 5.0817 & 3.1091 \\
0.200 & 3.0240 & 3.0340 & 5.0532 & 3.0440 & 4.8748 & 3.0340 \\
0.245 & 2.9739 & 2.9805 & 4.8888 & 2.9872 & 4.7512 & 2.9805 \\
0.300 & 2.9405 & 2.9450 & 4.7842 & 2.9494 & 4.6831 & 2.9450 \\
0.370 & 2.9177 & 2.9206 & 4.7149 & 2.9235 & 4.6439 & 2.9206 \\
0.450 & 2.9035 & 2.9054 & 4.6729 & 2.9074 & 4.6230 & 2.9054 \\
0.550 & 2.8936 & 2.8950 & 4.6445 & 2.8963 & 4.6102 & 2.8950 \\
0.670 & 2.8872 & 2.8881 & 4.6260 & 2.8890 & 4.6024 & 2.8881 \\
0.820 & 2.8827 & 2.8833 & 4.6133 & 2.8839 & 4.5975 & 2.8833 \\
1.000 & 2.8798 & 2.8802 & 4.6051 & 2.8806 & 4.5943 & 2.8802 \\
1.200 & 2.8780 & 2.8782 & 4.5999 & 2.8785 & 4.5924 & 2.8782 \\
1.400 & 2.8769 & 2.8771 & 4.5968 & 2.8773 & 4.5913 & 2.8771 \\
1.600 & 2.8761 & 2.8763 & 4.5948 & 2.8764 & 4.5906 & 2.8763 \\
1.800 & 2.8756 & 2.8758 & 4.5935 & 2.8759 & 4.5901 & 2.8758 \\
2.000 & 2.8753 & 2.8754 & 4.5925 & 2.8755 & 4.5898 & 2.8754 \\
2.250 & 2.8750 & 2.8751 & 4.5916 & 2.8751 & 4.5895 & 2.8751 \\
2.500 & 2.8747 & 2.8748 & 4.5910 & 2.8749 & 4.5892 & 2.8748 \\
2.750 & 2.8746 & 2.8746 & 4.5905 & 2.8747 & 4.5891 & 2.8746 \\
3.000 & 2.8745 & 2.8745 & 4.5902 & 2.8745 & 4.5889 & 2.8745 \\
3.250 & 2.8744 & 2.8744 & 4.5899 & 2.8744 & 4.5888 & 2.8744 \\
3.500 & 2.8743 & 2.8743 & 4.5897 & 2.8743 & 4.5888 & 2.8743 \\
4.000 & 2.8742 & 2.8742 & 4.5893 & 2.8742 & 4.5887 & 2.8742 \\
5.000 & 2.8740 & 2.8740 & 4.5890 & 2.8741 & 4.5885 & 2.8740 \\
    \bottomrule
    \end{tabular*}
\end{table}

\FloatBarrier

\paragraph*{Reverse coverage error}
To assess coverage rather than distance of fitted points to the target,
we use repetition zero with
\(N_q^{\mathrm{rev}}=(1024,8192,8192,16384,8192,1024)\) in the same
manifold order. For retained cloud \(\widehat{\mathcal P}\) and fixed
reference set \(\mathcal R=\{y_\ell\}_{\ell=1}^{N_{\mathrm{ref}}}
\subset\mathcal M\), we report
\(\overline E_{\mathrm{rev}}=N_{\mathrm{ref}}^{-1}
\sum_{\ell=1}^{N_{\mathrm{ref}}}d(y_\ell,\widehat{\mathcal P})\) and
\(E_{\mathrm{rev}}=\max_{\ell\le N_{\mathrm{ref}}}
 d(y_\ell,\widehat{\mathcal P})\), where
\(d(y,\widehat{\mathcal P})=\min_{\widehat z\in\widehat{\mathcal P}}
\|y-\widehat z\|_2\).
The reference counts are \(262144\) for all the manifolds. Within each method and manifold,
all estimators share queries and reference points; fitted batches are merged
before distances are computed. Filled markers require at most \(1\%\)
query failures, with no imputation and undefined errors for empty clouds.
These are single-repetition diagnostics, not 20-repetition summaries.

\paragraph*{Query and reference-point construction}
The query coordinates are noiseless base locations on the underlying
manifold; each fitting method subsequently applies its prescribed
initialization and fitting procedure. The forward and reverse query
counts are those specified above.
On \(\mathbb S^1\), queries are equally spaced in angle,
\(\theta_i=2\pi i/N\), \(i=0,\ldots,N-1\), for query count \(N\).
On the Clifford torus, forward queries use a \(16\times16\) tensor grid
in the two angular coordinates, whereas reverse queries use a
\(64\times128\) tensor grid. Each angular grid consists of equally
spaced points on \([0,2\pi)\), without a midpoint shift.

On the three spherical geometries, queries are generated by normalizing
independent standard Gaussian vectors in \(\mathbb R^{d+1}\).
For \(\mathbb S^2\subset\mathbb R^{10}\), the resulting three coordinates are
embedded in the first three ambient coordinates, with the remaining
seven set to zero. The generator is NumPy's \texttt{PCG64}, initialized
with \texttt{SeedSequence([20260912, geometry\_id])}, where the geometry
identifiers are \(2,4,5\) for \(\mathbb S^2\subset\mathbb R^3\),
\(\mathbb S^3\subset\mathbb R^4\), and \(\mathbb S^2\subset\mathbb R^{10}\), respectively.

For the corrugated torus, let \(J\) be the surface-area density specified
in the implementation details and define the normalized meridian area
distribution
\[
F_J(v)=
\frac{\displaystyle\int_0^v J(s)\,ds}
     {\displaystyle\int_0^{2\pi}J(s)\,ds}.
\]
A deterministic \(n_u\times n_v\) query set is constructed as
\[
u_i=\frac{2\pi(i+1/2)}{n_u},
\qquad
v_j=F_J^{-1}\!\left(\frac{j+1/2}{n_v}\right),
\qquad
x_{ij}=X(u_i,v_j),
\]
where \(i=0,\ldots,n_u-1\) and \(j=0,\ldots,n_v-1\).
Forward and reverse queries use \((n_u,n_v)=(32,16)\) and \((64,16)\),
respectively. Numerically, \(F_J\) is constructed by cumulative
trapezoidal integration on \(65536\) equal meridian intervals, and its
inverse is evaluated by linear interpolation.

For the \(N_{\mathrm{ref}}=262144\) reverse reference points, the circle
uses an equally spaced angular grid, and the Clifford torus uses a
\(512\times512\) angular tensor grid. The spherical reference sets use
the same normalized-Gaussian construction, but with the separate seed
\texttt{SeedSequence([20268831, geometry\_id])}.
The corrugated-torus reference set uses the deterministic
\(512\times512\) area-coordinate midpoint grid defined above.

For each geometry, the query design is fixed across scales and
repetitions, and the reverse reference set is shared across STSP, DR,
and YY. The frozen coordinate arrays are included in the public
implementation and are read directly for replay, rather than
regenerated independently for each method.

Figures~\ref{fig:e3-erev-stsp}--\ref{fig:e3-erev-yy} show where the
forward accuracy gains also improve coverage. For STSP on
\(\mathbb S^3\subset\mathbb R^4\), GBC improves on both reliable Single-scale
baselines at tested bandwidths \(h_2=0.45\)--\(1.8\) for mean reverse
error, but only \(1.2\)--\(1.8\) for maximum reverse error. DR GBC is
reliable at 3 of the 28 reverse scales on this manifold. On the corrugated
torus, neither STSP nor DR has two adjacent jointly reliable scales with
lower GBC maximum reverse error than both baselines; GBC-YY likewise has
no such interval on the Clifford torus or \(\mathbb S^2\subset\mathbb R^{10}\)
relative to all three baselines. Reverse error includes point spacing and
missing queries as well as fitting error, so improved forward accuracy
need not produce an equally broad coverage improvement. These
single-repetition, finite-reference diagnostics do not certify continuous
Hausdorff accuracy.

Overall, these experiments support higher fitting accuracy and sustained low error over
a substantial range of reliable scales, most clearly for STSP. This is
the scale-selection benefit illustrated in Section \ref{sec:e3-bandwidth-robustness}: accurate fitting
is not confined to a single bandwidth. The extent of this benefit depends
on the method, geometry, and error metric, with numerical failures
limiting DR and local error amplification plausibly limiting YY's
maximum-error gains.

\begin{figure}[!htbp]
    \centering
    \begin{subfigure}{0.90\linewidth}
        \centering
        \includegraphics[width=\linewidth]{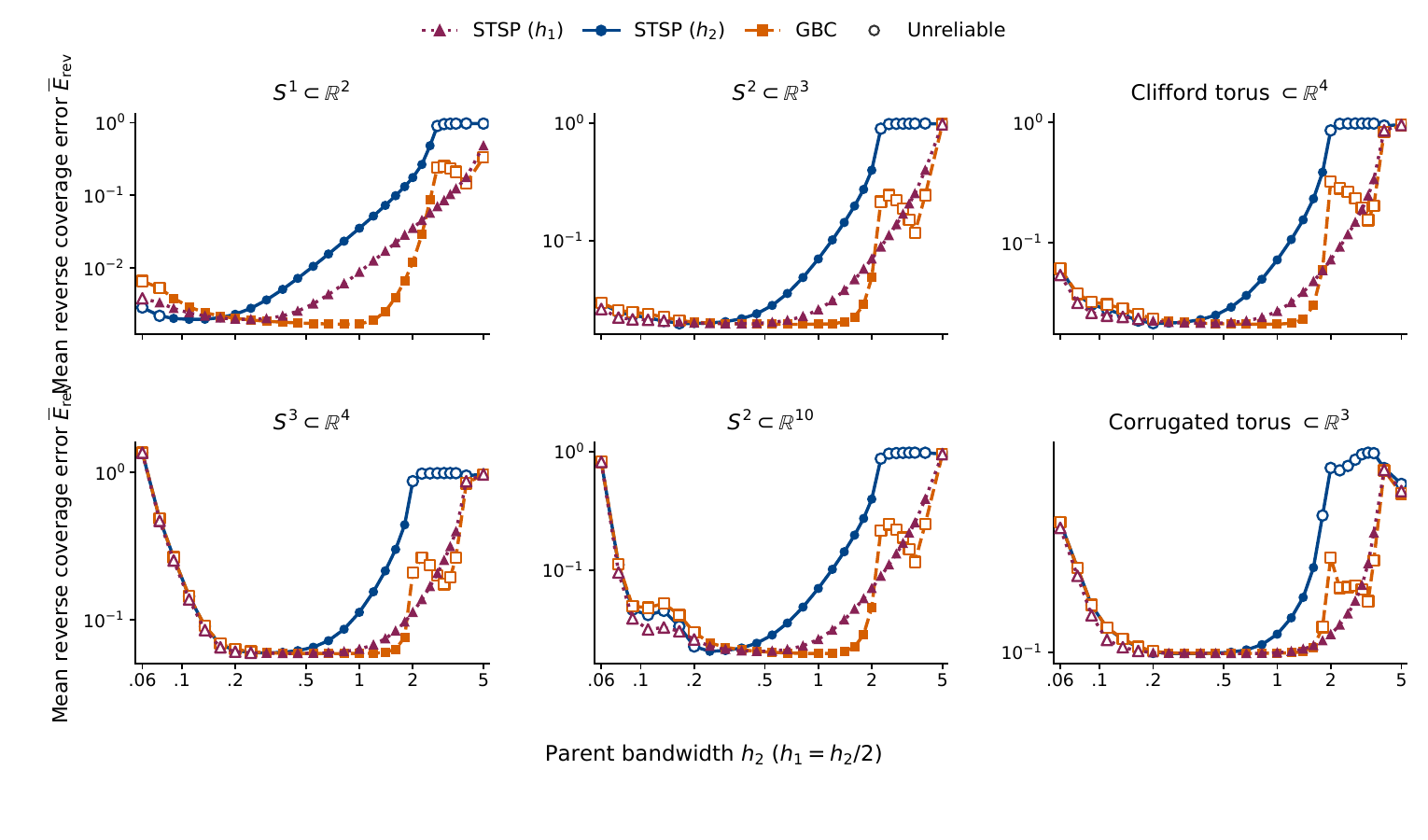}
        \caption{Mean reverse coverage error \(\overline E_{\mathrm{rev}}\).}
        \label{fig:e3-erev-stsp-mean}
    \end{subfigure}
    \par\vspace{3pt}
    \begin{subfigure}{0.90\linewidth}
        \centering
        \includegraphics[width=\linewidth]{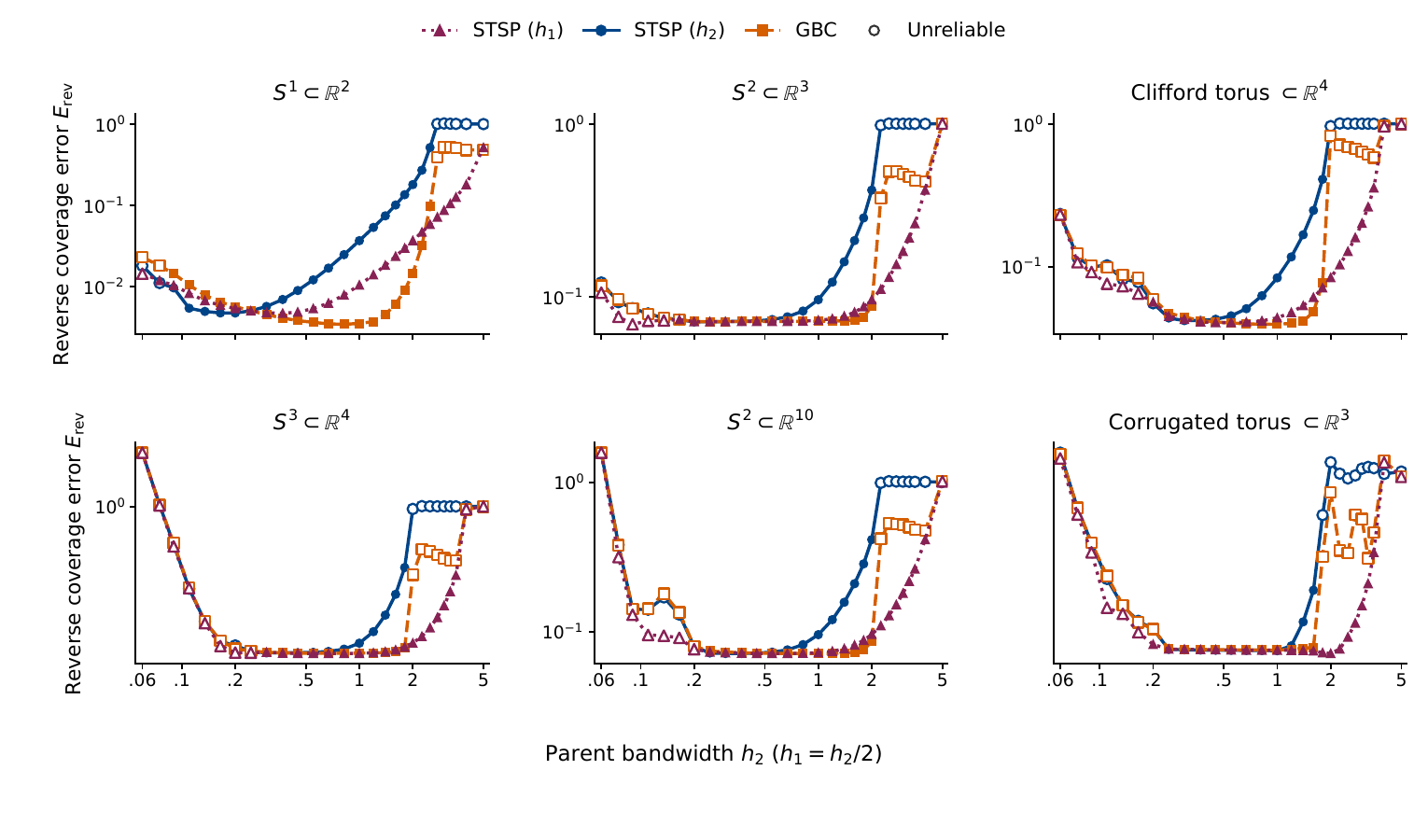}
        \caption{Maximum reverse coverage error \(E_{\mathrm{rev}}\).}
        \label{fig:e3-erev-stsp-max}
    \end{subfigure}
    \caption{STSP reverse coverage in repetition zero: (a) mean and (b) maximum
    distance from the evaluation reference set to the fitted cloud.
    Curves compare Single-scale at \(h_1\), Single-scale at \(h_2\), and GBC,
    with \(h_1=h_2/2\). Hollow markers denote query-failure rates
    above \(1\%\).}
    \label{fig:e3-erev-stsp}
\end{figure}

\begin{figure}[!htbp]
    \centering
    \begin{subfigure}{0.90\linewidth}
        \centering
        \includegraphics[width=\linewidth]{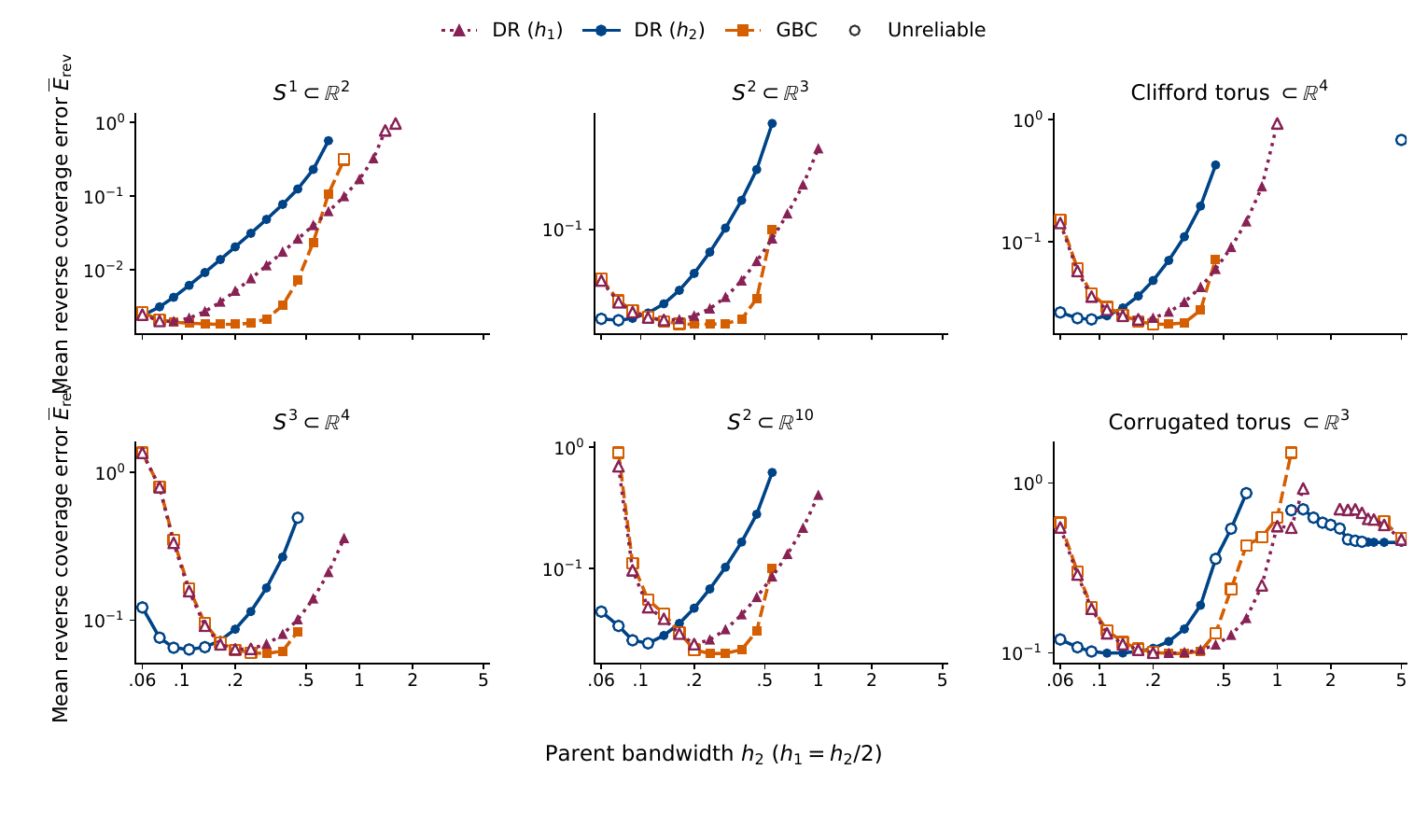}
        \caption{Mean reverse coverage error \(\overline E_{\mathrm{rev}}\).}
        \label{fig:e3-erev-dr-mean}
    \end{subfigure}
    \par\vspace{3pt}
    \begin{subfigure}{0.90\linewidth}
        \centering
        \includegraphics[width=\linewidth]{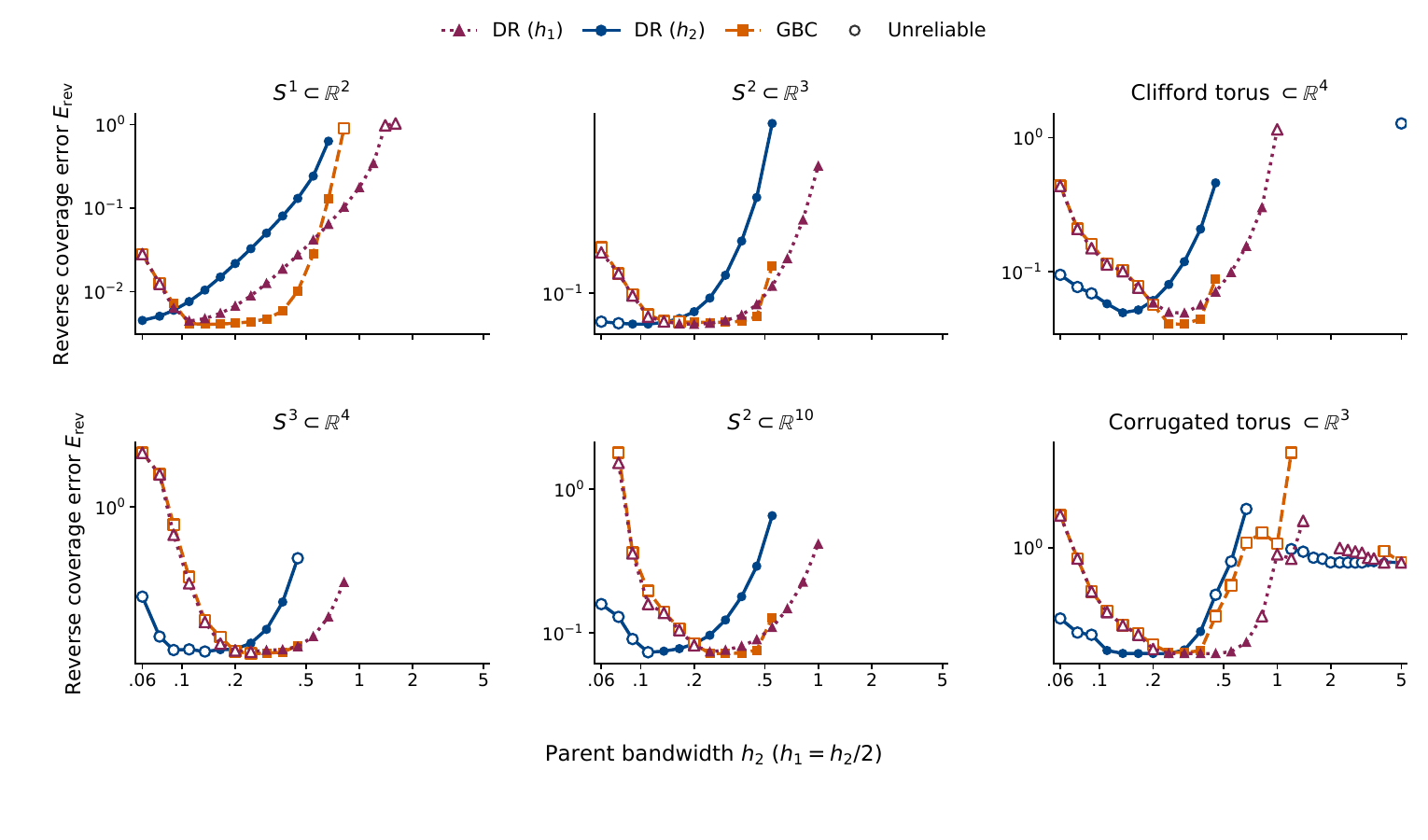}
        \caption{Maximum reverse coverage error \(E_{\mathrm{rev}}\).}
        \label{fig:e3-erev-dr-max}
    \end{subfigure}
    \caption{DR reverse coverage in repetition zero: (a) mean and (b) maximum
    distance from the evaluation reference set to the fitted cloud.
    Curves compare both Single-scale baselines with GBC, with \(h_1=h_2/2\).
    Hollow markers denote query-failure rates above \(1\%\); gaps
    denote undefined values.}
    \label{fig:e3-erev-dr}
\end{figure}

\begin{figure}[!htbp]
    \centering
    \begin{subfigure}{0.90\linewidth}
        \centering
        \includegraphics[width=\linewidth]{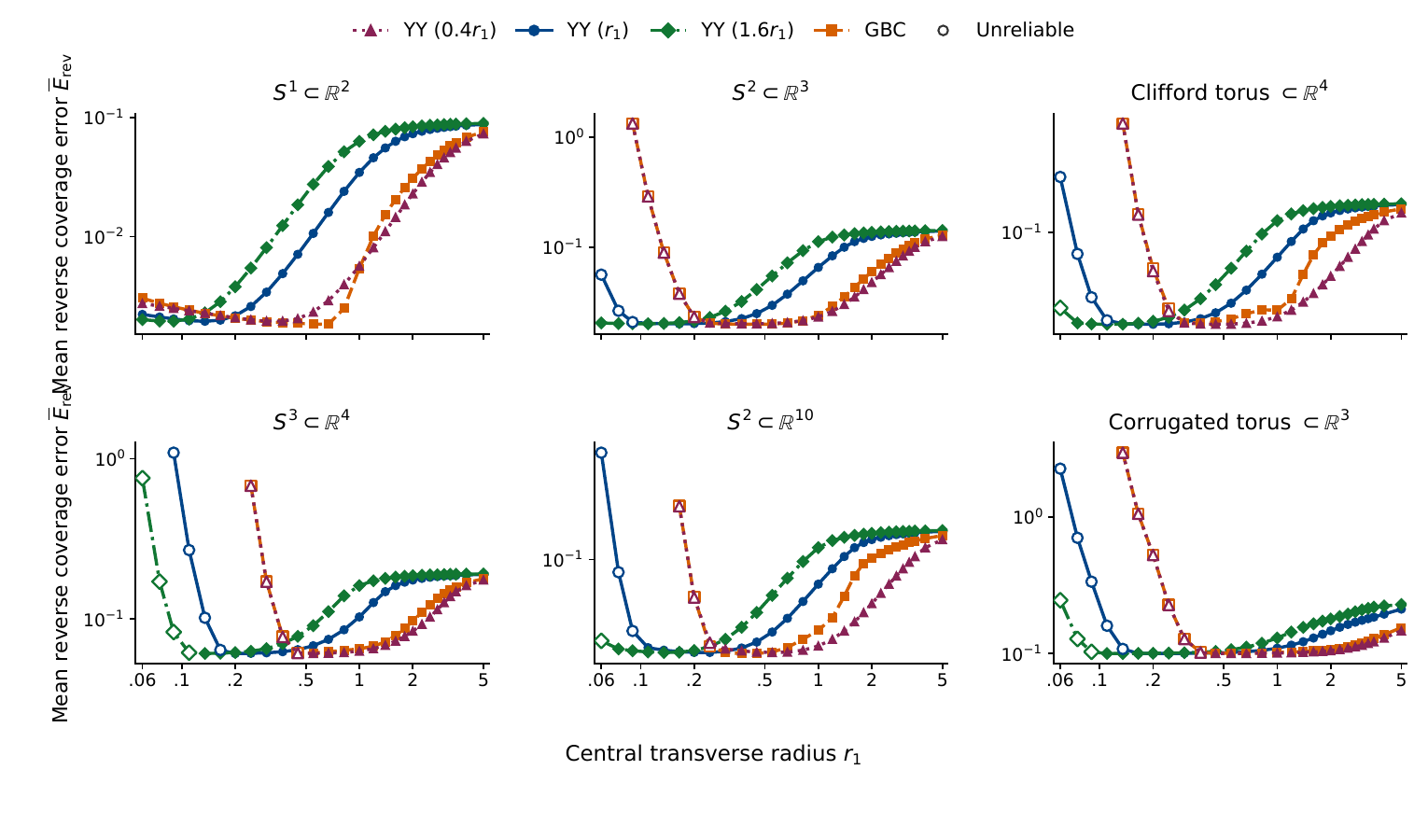}
        \caption{Mean reverse coverage error \(\overline E_{\mathrm{rev}}\).}
        \label{fig:e3-erev-yy-mean}
    \end{subfigure}
    \par\vspace{3pt}
    \begin{subfigure}{0.90\linewidth}
        \centering
        \includegraphics[width=\linewidth]{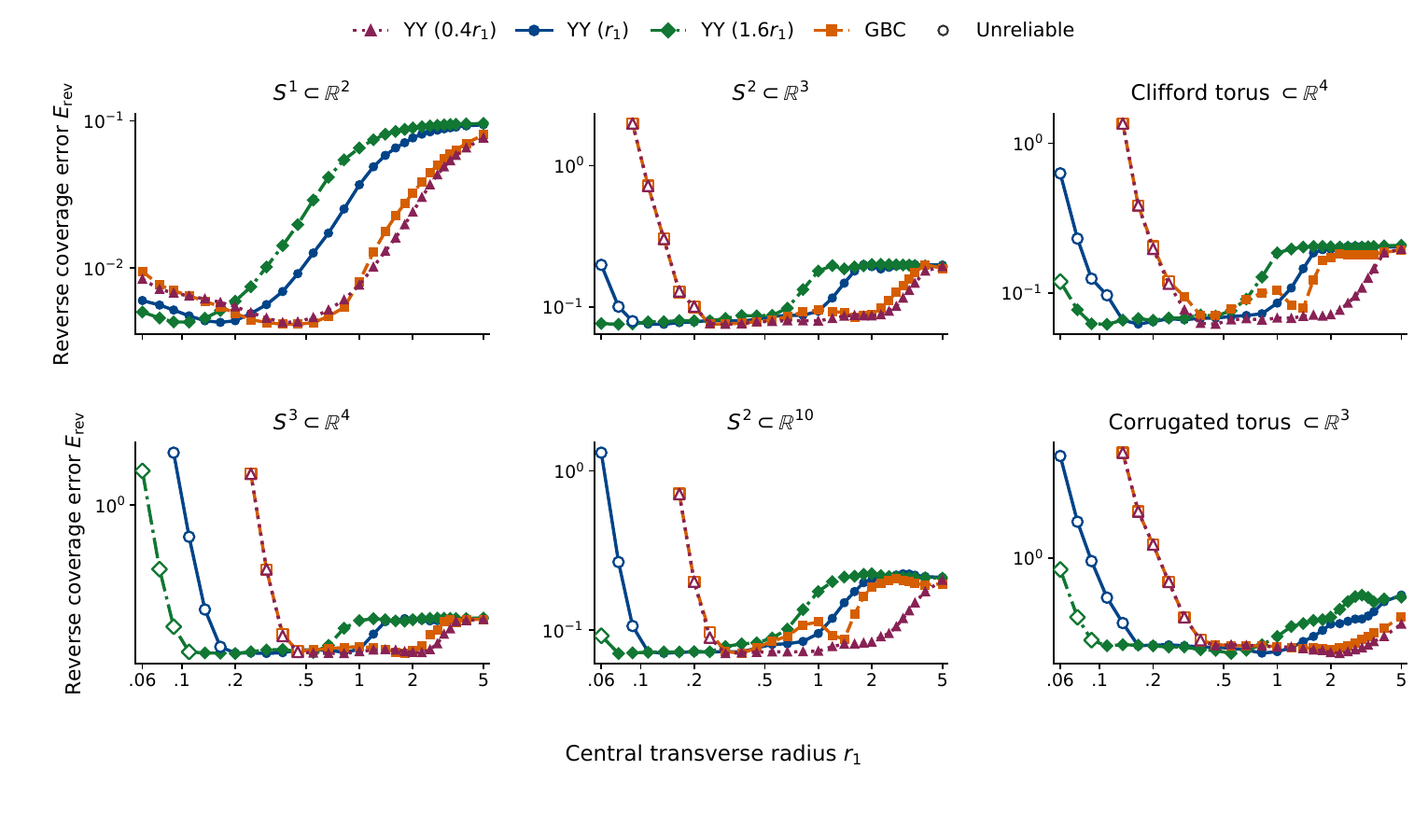}
        \caption{Maximum reverse coverage error \(E_{\mathrm{rev}}\).}
        \label{fig:e3-erev-yy-max}
    \end{subfigure}
    \caption{YY reverse coverage in repetition zero: (a) mean and (b) maximum
    distance from the evaluation reference set to the fitted cloud.
    Curves compare Single-scale at \(0.4r_1\), \(r_1\), and \(1.6r_1\) with
    GBC-YY. Hollow markers denote query-failure rates above \(1\%\);
    gaps denote undefined values.}
    \label{fig:e3-erev-yy}
\end{figure}

\paragraph*{Section \ref{sec:enf_cp_comparison}: Correspondence and Candidate Density}
We use the nonconstant-curvature curve
$\left\{
\bigl(
1.5\cos\theta,\,
\sin\theta,\,
0.3\sin(2\theta),\,
0
\bigr)
:\theta\in[0,2\pi)
\right\} \subset\mathbb R^3$,
$n=200{,}000$, $\sigma=0.06$, $h_1=0.12$, $h_2=0.24$,
$N_q=32$ fixed reference queries, and ten paired repetitions.
All four methods share the same observations, fitted fine-scale
reference roots, and extrapolation weights. The common reference
queries are initialized and validated using true normal charts;
thus the benchmark is not a fully data-driven reconstruction from
arbitrary initial points. True geometry is not supplied to the
ENF or local CP cross-scale matching steps.

ENF intersects the empirical $h_2$-scale fitted branch with
normal planes estimated by $h_1$-scale weighted PCA. Oracle
uses true normal planes as an explicitly labeled reference.
All accepted fitted points satisfy the empirical STSP residual
\[
\widehat R_h(y)=\widehat Q_h(y)
  \{y-\widehat\mu_h(y)\}.
\]
For a local normal basis $B$, numerical roots are accepted only if
\[
\bigl\|B^\top\widehat R_{h_2}(y)\bigr\|_2\le10^{-11},
\qquad
\bigl\|\widehat R_{h_2}(y)\bigr\|_2\le10^{-9},
\]
and if the normalized-mass, eigengap, transversality,
Jacobian, and search-region checks are satisfied.
A small projected residual alone does not constitute a successful
fit when the full residual check fails.

For discrete NN, the coarse-scale candidate counts are
\[
M_{\mathrm{match}}\in\{32,36,40,44,48,52,56,60,64,
80,96,112,128,256,512\}.
\]
Candidate locations are initialized on the true curve in
normalized arclength at
$s_j=((j+0.74)/32)\bmod1$, $0\le j<32$,
and expanded by nested midpoint refinement.
Each candidate is independently fitted at scale $h_2$ using
an empirical STSP residual solve; NN returns the nearest
successful fitted candidate. The candidate construction
uses true geometry for initialization, but matching is
performed on fitted candidates rather than true curve points.

\paragraph*{Implementation of local continuous CP}
For each shared reference point
$z_{1,j}\in\widehat{\mathcal M}_{h_1}$,
we obtain a unit tangent $T_j\in\mathbb R^3$ and an
orthonormal normal basis $B_j\in\mathbb R^{3\times2}$
from empirical $h_2$-scale weighted PCA at $z_{1,j}$.
We construct a continuous local representation of the
empirical target STSP branch by
\begin{equation}\label{eq:supp_local_cp_chart}
\begin{aligned}
z_j(s)&=z_{1,j}+sT_j+B_j u_j(s),\\
B_j^\top\widehat R_{h_2}(z_j(s))&=0.
\end{aligned}
\end{equation}
For each scalar search position $s$, the normal coordinates
$u_j(s)\in\mathbb R^2$ are obtained numerically by solving
the projected residual equation, subject to both residual
tolerances and the same mass, eigengap, transversality,
and Jacobian checks as above.

Let $J_j(s)=D\widehat R_{h_2}(z_j(s))$. Whenever the
projected normal Jacobian is nonsingular, implicit
differentiation gives
\begin{equation}\label{eq:supp_local_cp_derivative}
\begin{aligned}
u_j'(s)
&=-\{B_j^\top J_j(s)B_j\}^{-1}
   B_j^\top J_j(s)T_j,\\
z_j'(s)&=T_j+B_j u_j'(s),\\
\frac{d}{ds}\|z_j(s)-z_{1,j}\|_2^2
&=2\langle z_j(s)-z_{1,j},z_j'(s)\rangle.
\end{aligned}
\end{equation}
The scalar search starts with nine nodes and is refined to
41 or 81 nodes when the sampled distance derivative is
nonmonotone. Bracketed stationary points are refined using
Brent's method; four additional normal-direction initializations
are used to check the selected local sheet.

Local CP thereby optimizes a continuous distance objective on
a constructed empirical STSP sheet, not a finite NN candidate
library. The CP-specific search receives no true curve
parameterization, true normal basis, Oracle correspondence,
or discrete candidate library. The common fine-scale reference
roots nevertheless retain the geometry-assisted initialization
and validation described above. The numerical search does not
certify global optimality over arbitrary additional components
of the empirical STSP residual zero set. Accordingly, we call
this method \emph{local continuous CP}, reserving \emph{global
CP} for the geometric closest-point map in Proposition~\ref{prop:cp_supplement}.

\paragraph*{Error evaluation and timing}
For $A\in\{\mathrm{NN},\mathrm{ENF},\mathrm{locCP},
\mathrm{Or}\}$, the matching discrepancy is
\[
E_{\mathrm{match}}^{(A)}
=\max_{1\le j\le N_q}
\|Z_{2,j}^{(A)}-Z_{2,j}^{\mathrm{Or}}\|,
\]
where Oracle uses the true normal fiber on the same empirical
coarse-scale fitted branch. Final geometric errors are
computed after extrapolation using common weights.
The true manifold is used for post hoc error evaluation.

Local CP timing is measured as single-thread CPU
wall-clock time for batches of $32$ reference queries
across the ten paired sample realizations. The median complete time is $25.562$ seconds; the
separately measured median CP correspondence-stage time is
$24.615$ seconds.
The complete time includes observation-array construction,
both STSP indices, extrapolation weights, shared reference-root
reconstruction and validation, CP-specific empirical PCA,
implicit target-branch residual solves and derivatives,
continuous search, the four additional normal starts, and
the final extrapolation combination. Input-file reading,
true-curve arclength and query-grid generation, truth
projections and error evaluation, output writing, and plotting
are excluded. These measurements describe this specific
implementation; they are not algorithm-independent
complexity estimates.

All $32$ queries are accepted in each of the ten
repetitions for every method and NN candidate-count setting,
so the observed query-matching failure rate is zero.

\subsubsection{Oscillatory normal graph: nonuniqueness of global CP}
\label{sec:supp_oscillatory_cp}

We use a deterministic oscillatory normal graph to illustrate
the distinction between the local well-posedness of ENF and
the global uniqueness required by closest-point matching.
This is a finite-scale illustration of the geometric
counterexample in Appendix~\ref{app:general_enf_proofs}, not a fitted
branch produced by empirical STSP.

Let $\mathcal M=\mathbb S^1\subset\mathbb R^2$ and take
$h_1=0.1$, $h_2=0.2$. Define the reference curve
\[
X_{h_1}(\theta)
=(1+h_1^2)(\cos\theta,\sin\theta)
\]
and the oscillatory target curve
\begin{equation}\label{eq:supp_oscillatory_graph}
\begin{aligned}
X_{h_2}(\theta)&=r(\theta)(\cos\theta,\sin\theta),\\
r(\theta)&=1+h_2^2+\alpha\cos(N\theta),\\
N&=h_2^{-6}=15625,\qquad
\alpha=\frac{h_2^2}{N}=2.56\times10^{-6}.
\end{aligned}
\end{equation}
The oscillation has amplitude $\alpha=O(h_2^8)$ and
first-derivative amplitude $\alpha N=O(h_2^2)$,
whereas $\alpha N^2=O(h_2^{-4})$.
Thus the example has small positional and first-order
perturbations but no scale-uniform second-order
curvature bound.

Consider the common reference point $y_1=(1.01,0)$.
The reference-normal fiber is the horizontal axis;
its unique intersection with the target curve near
$\theta=0$ is
\[
Z^{\mathrm{ENF}}=(1.04+\alpha,0).
\]
By contrast, global CP minimizes the squared distance
\begin{equation}\label{eq:supp_oscillatory_g}
\begin{aligned}
g(\theta)
  &=\|X_{h_2}(\theta)-y_1\|_2^2\\
  &=r(\theta)^2+1.01^2
    -2(1.01)r(\theta)\cos\theta.
\end{aligned}
\end{equation}
The objective is even, and direct differentiation yields
\[
\begin{aligned}
g''(0)
&=2(1.01)r(0)
  -2\alpha N^2\{r(0)-1.01\}\\
&\approx-35.4024<0.
\end{aligned}
\]
Hence $\theta=0$ is a strict local maximum and cannot
be a global minimizer. Also $g(\pi)>g(0)$, so
$\theta=\pi$ cannot be a global minimizer either.
Since $g$ is continuous and attains its minimum on the
circle, any global minimizer is away from the two fixed
points of reflection $\theta\mapsto-\theta$.
By symmetry it has a distinct reflected global minimizer.
Numerically, the two minimizing angles are approximately
\[
\theta_{\mathrm{CP}}\approx\pm1.9035\times10^{-4}.
\]
Consequently, global CP is set-valued at $y_1$, although
the nearby ENF intersection is unique.

\begin{figure}[tbp]
  \centering
  \includegraphics[width=0.96\linewidth]{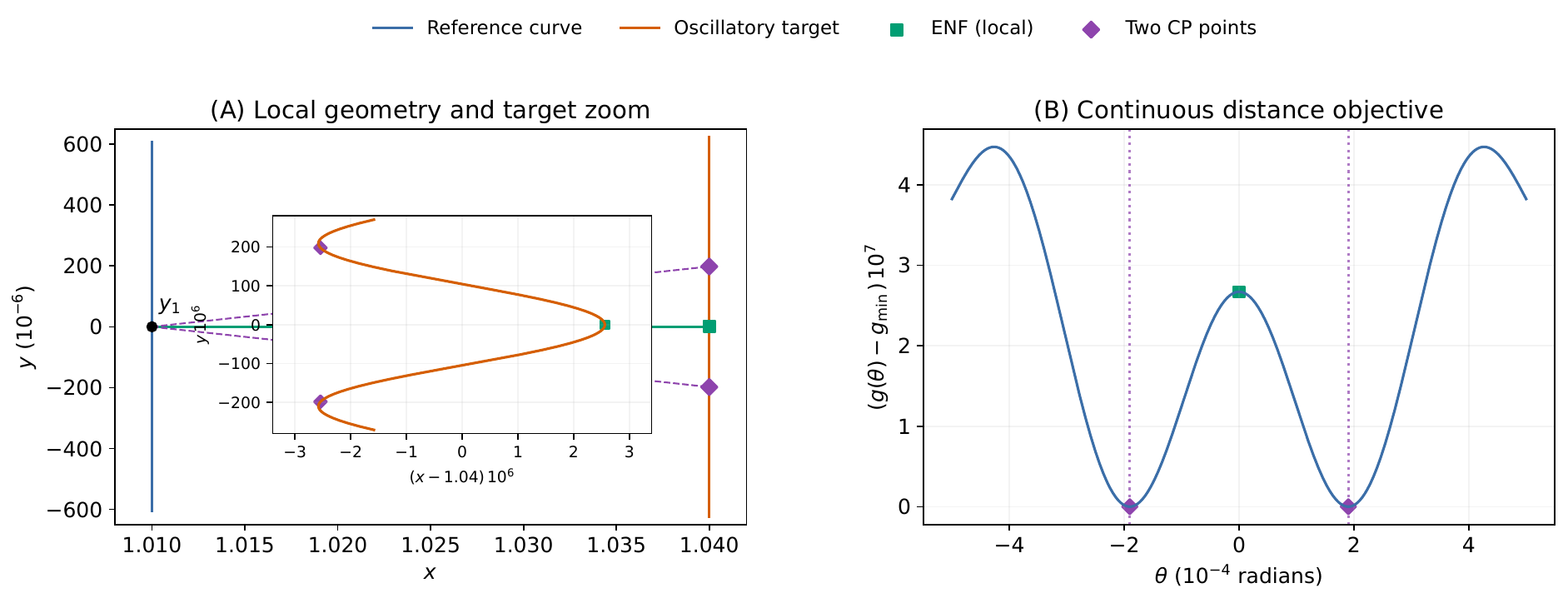}
  \caption{Failure of global closest-point uniqueness on
  an oscillatory normal graph.
  (A) Reference and target curves with a local inset,
  showing the unique nearby ENF intersection and two
  distinct global CP points.
  (B) Rescaled squared-distance objective
  $10^7\{g(\theta)-\min g\}$, with its two symmetric
  global minima.
  This artificial example exhibits global CP nonuniqueness
  without loss of local ENF transversality; it is not an
  empirical STSP result and does not imply that individual
  CP minimizers lose fourth-order localization.}
  \label{fig:supp_oscillatory_cp}
\end{figure}

Figure~\ref{fig:supp_oscillatory_cp} illustrates why
fourth-order positional accuracy and a well-defined,
smooth correspondence map are different questions.
Even when individual closest points remain well localized,
the global CP selection need not be single valued under
weaker curvature regularity. The ENF construction
instead selects a local transverse intersection.

\Needspace{4\baselineskip}
\subsection{Directional Ablation of the YY Estimator}
\label{app:directional-ablation}

Figure~\ref{fig:e5-directional-population} supplements
Table~\ref{tab:e5-population-slopes}. All estimators reuse the same starting
points, first-stage axes, and fitted outputs; only the application of the
weights changes. Scalar-A tracks directional GBC in axial error but retains
second-order complementary error, whereas Scalar-C exhibits the opposite
pattern. Both scalar controls therefore remain approximately second order
in total error. Among the tested estimators, only directional GBC suppresses
both leading components and has approximately fourth-order total error.
This ablation supports the role of component-specific weights in the
specified higher-codimension example, rather than an advantage from merely
combining additional scales.

\begin{figure}[!htbp]
    \centering
    \includegraphics[width=\linewidth]
    {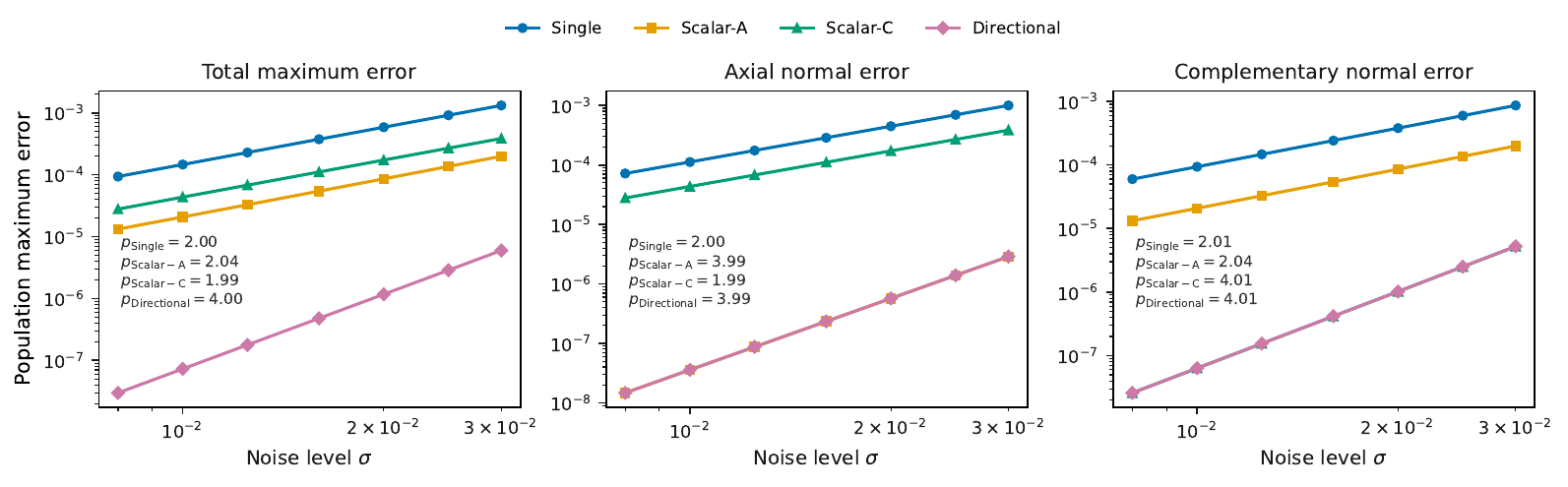}
    \caption{Population directional-weight ablation on
    \(\left\{
\bigl(
1.5\cos\theta,\,
\sin\theta,\,
0.3\sin(2\theta),\,
0
\bigr)
:\theta\in[0,2\pi)
\right\}\).
    Left to right: total maximum point-to-manifold error, axial normal
    error, and complementary normal error versus \(\sigma\).
   The curves show all seven noise levels, while the reported log--log slopes are fitted over \(\sigma\in\{0.008,0.010,0.0125,0.016\}\), as in Table~\ref{tab:e5-population-slopes}.}
    \label{fig:e5-directional-population}
\end{figure}

\FloatBarrier
\Needspace{4\baselineskip}

\subsection{Multiscale Bias Cancellation for STSP under Nonuniform
Sampling}
\label{supp:nonuniform_sampling}

The general multiscale bias-cancellation framework does not
inherently require uniform sampling on the underlying manifold.
Its essential requirement is that the single-scale fitting
procedure admits a suitable cross-scale bias expansion, together
with the geometric regularity and stability conditions needed
for extrapolation. Uniform sampling is imposed in our theoretical
verification of these conditions for specific fitting procedures.

To empirically examine that nonuniform sampling does not
compromise the effectiveness of our multiscale extrapolation,
we consider STSP as a representative example and evaluate
its performance under nonuniform latent sampling. We adopt the population and finite-sample configurations
of Appendices~\ref{app:population-verification} and
\ref{app:asymptotic-behavior}, changing only the latent
sampling distribution.

\paragraph*{Sampling distributions}
We consider the same six manifolds as in the main STSP
experiments. For each manifold $\mathcal M$, the sampling
density $p$ is specified with respect to the intrinsic
volume measure $d\mathrm{Vol}_{\mathcal M}$.

\begin{table}[htbp]
\centering
\caption{Nonuniform sampling densities used in the STSP experiments.}
\label{tab:nonuniform_densities}
\small
\begin{tabular}{@{}ll@{}}
\toprule
Manifold & Sampling density $p$ \\
\midrule
$\mathbb S^1\subset\mathbb R^2$
& $\displaystyle \frac{1+0.5\cos\theta}{2\pi}$ \\[0.5em]
$\mathbb S^2\subset\mathbb R^3$
& $\displaystyle \frac{1+0.5x_1}{4\pi}$ \\[0.5em]
Clifford torus $\subset\mathbb R^4$
& $\displaystyle \frac{1+0.5\cos u}{2\pi^2}$ \\[0.5em]
$\mathbb S^3\subset\mathbb R^4$
& $\displaystyle \frac{1+0.5x_1}{2\pi^2}$ \\[0.5em]
$\mathbb S^2\subset\mathbb R^{10}$
& $\displaystyle \frac{1+0.5x_1}{4\pi}$ \\[0.5em]
Corrugated torus $\subset\mathbb R^3$
& $\displaystyle \frac{1+0.5\cos v}{Z}$ \\
\bottomrule
\end{tabular}
\end{table}

Here $\theta$ denotes the angular coordinate on $\mathbb S^1$,
$(u,v)$ are the angular coordinates on the tori, and $x_1$
denotes the first ambient coordinate. The Clifford torus is
$\mathbb S^1(1/\sqrt{2})\times\mathbb S^1(1/\sqrt{2})$,
while $Z$ normalizes the density on the corrugated torus
with respect to its surface-area measure.
All densities are smooth and strictly positive, with a
maximum-to-minimum ratio of $3$.

\begin{figure}[!htbp]
    \centering
    \includegraphics[width=0.98\linewidth]{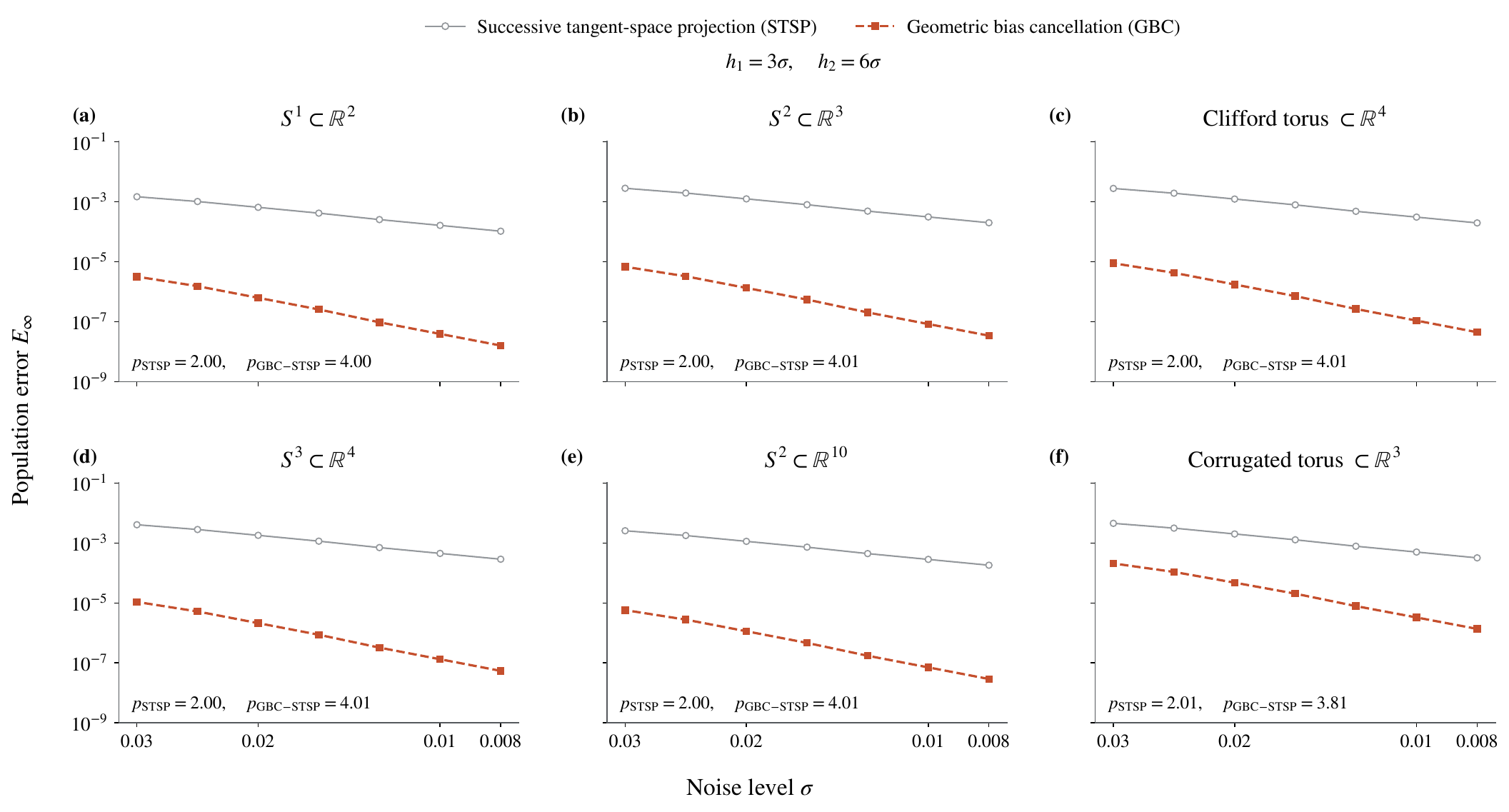}
    \caption{
    Population maximum point-to-manifold error $E_\infty$
    for single-scale and multiscale STSP under nonuniform
    latent sampling. The six panels correspond to the
    manifolds in Table~\ref{tab:nonuniform_densities},
    with $h_1=3\sigma$ and $h_2=6\sigma$.
    Panel annotations give the fitted log--log slopes
    for Single-scale and GBC, respectively.
    }
    \label{fig:supp-nonuniform-pop-max}
\end{figure}

\begin{figure}[!htbp]
    \centering
    \includegraphics[width=0.98\linewidth]{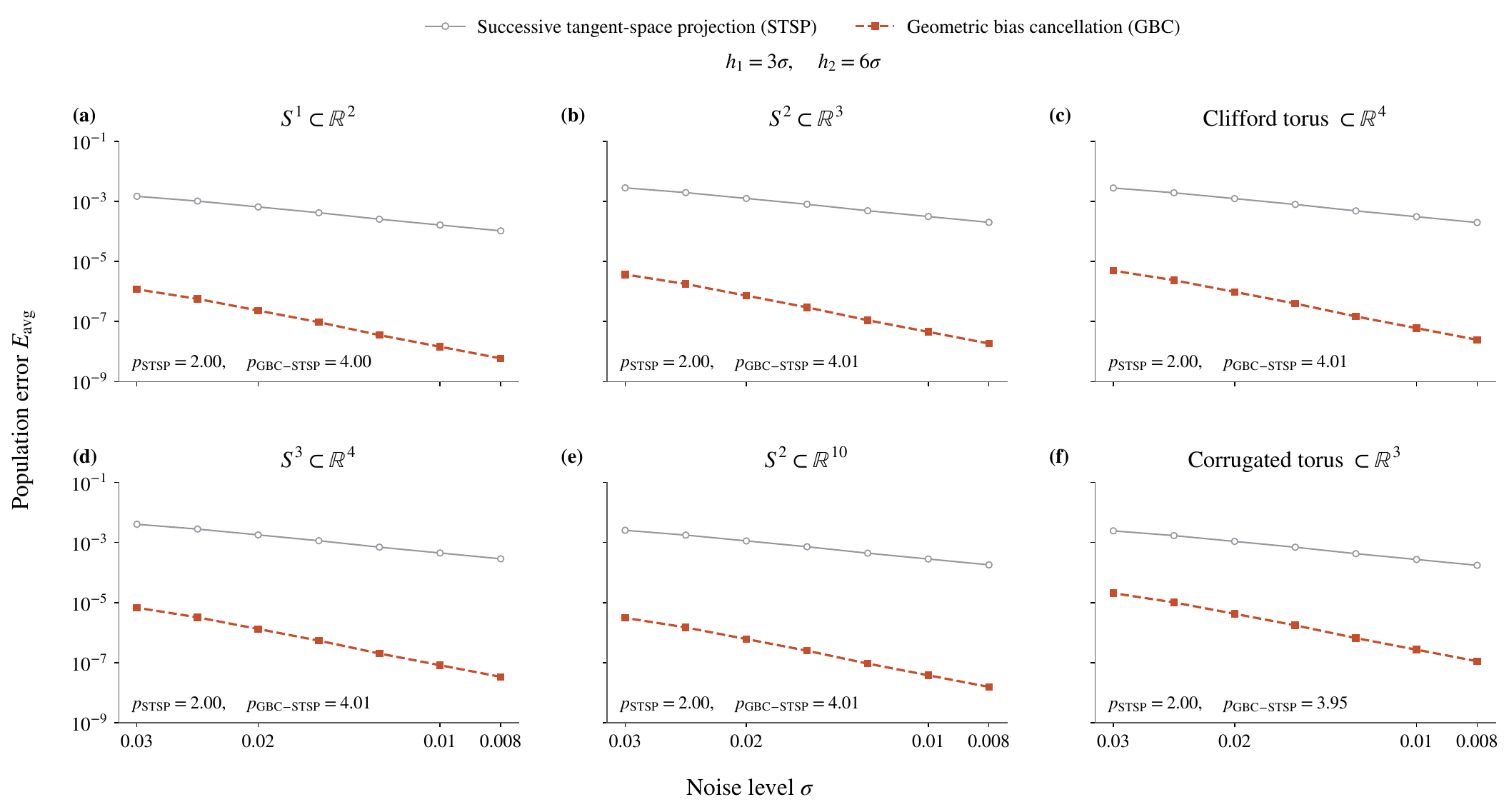}
    \caption{
    Population average point-to-manifold error
    $E_{\mathrm{avg}}$ under the same nonuniform densities
    and scale settings as in
    Figure~\ref{fig:supp-nonuniform-pop-max}.
    Panel annotations give the fitted log--log slopes
    for Single-scale and GBC, respectively.
    }
    \label{fig:supp-nonuniform-pop-avg}
\end{figure}

\begin{figure}[!htbp]
    \centering
    \includegraphics[width=0.98\linewidth]{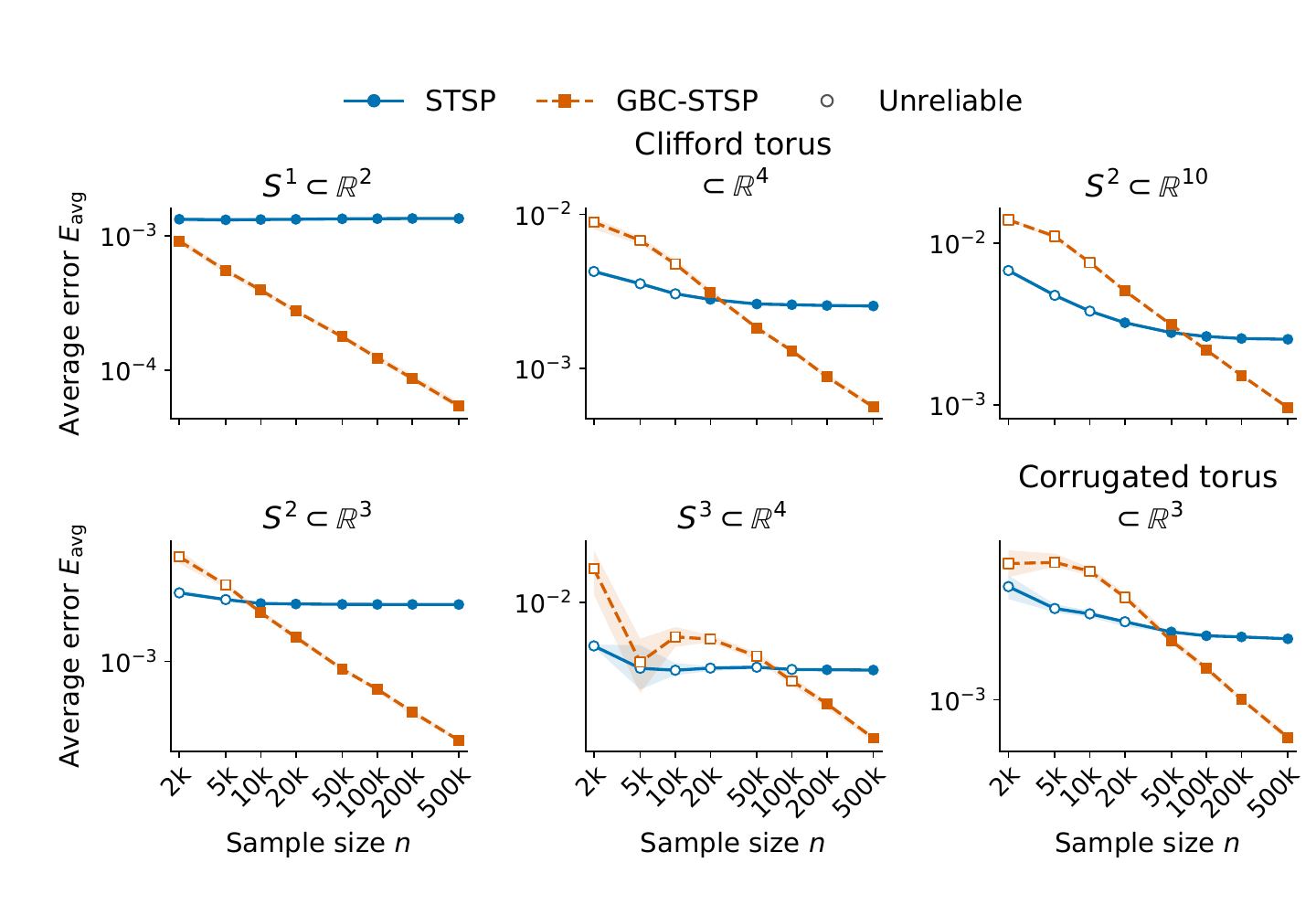}
    \caption{
    Finite-sample average point-to-manifold error
    $E_{\mathrm{avg}}$ for STSP under nonuniform sampling,
    with $\sigma=0.005$, $h_1=0.1$, and $h_2=0.2$.
    Curves show medians and bands show interquartile
    ranges over 20 paired repetitions.
    Hollow markers denote unreliable configurations.
    }
    \label{fig:supp-nonuniform-finite-avg}
\end{figure}

\begin{figure}[!htbp]
    \centering
    \includegraphics[width=0.98\linewidth]{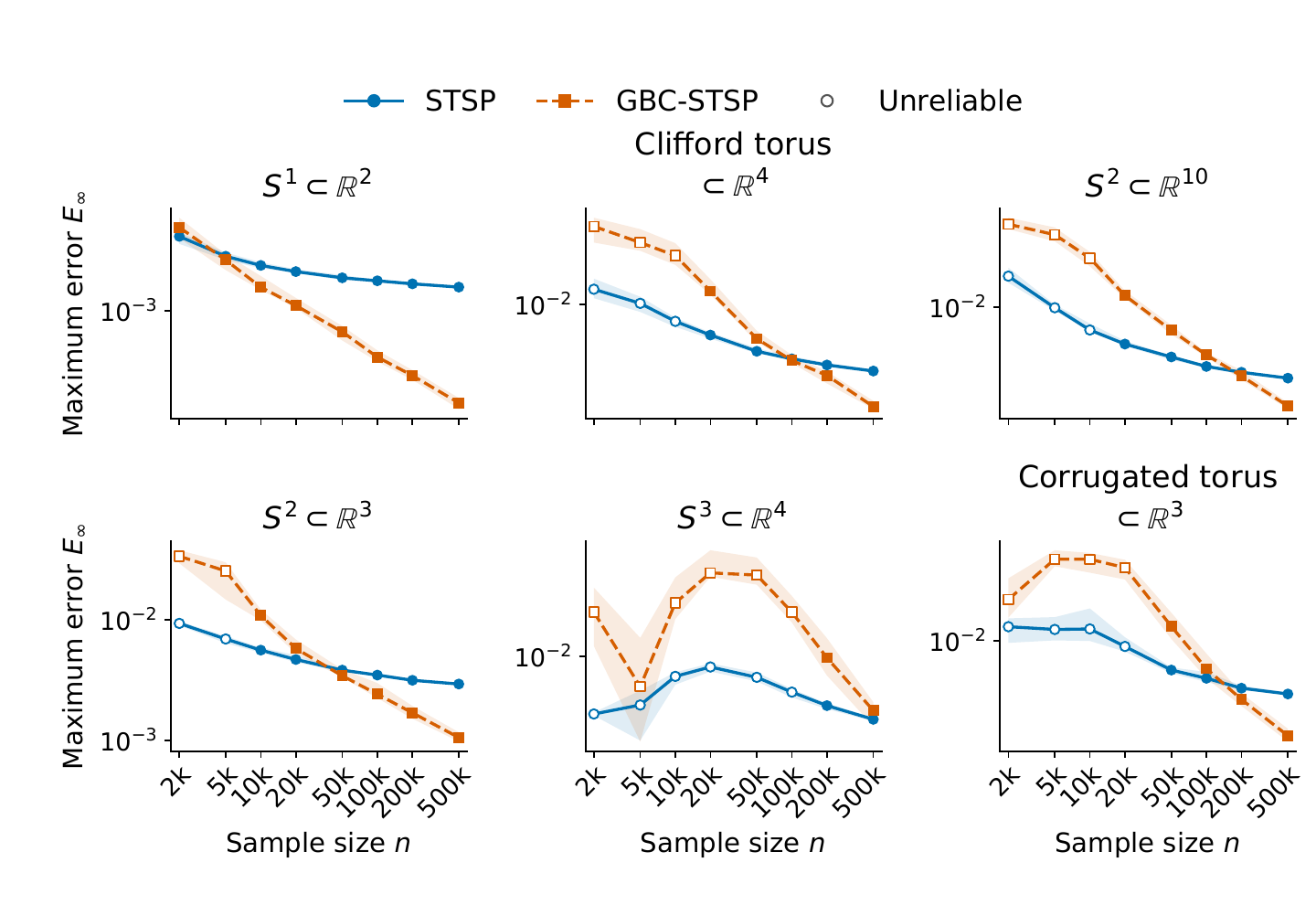}
    \caption{
    Finite-sample maximum point-to-manifold error
    $E_\infty$ under the same nonuniform densities
    and experimental settings as in
    Figure~\ref{fig:supp-nonuniform-finite-avg}.
    Curves show medians and bands show interquartile
    ranges over 20 paired repetitions.
    Hollow markers denote unreliable configurations.
    }
    \label{fig:supp-nonuniform-finite-max}
\end{figure}

\paragraph*{Population results}
Figures~\ref{fig:supp-nonuniform-pop-max} and
\ref{fig:supp-nonuniform-pop-avg} show that the
multiscale bias reduction persists under the
nonuniform sampling densities considered here.

Across the five sphere and Clifford-torus settings,
the fitted population-error slopes are approximately
$2.00$ for Single-scale and $4.00$--$4.01$ for GBC under
both error metrics. On the corrugated torus, the
corresponding GBC slopes are $3.81$ for $E_\infty$
and $3.95$ for $E_{\mathrm{avg}}$, compared with
approximately $2.00$ for Single-scale.

These results indicate that the prescribed STSP
extrapolation weights continue to cancel the leading
second-order population bias under the tested
nonuniform densities, without density-dependent
adjustments to the fitting operator or extrapolation
weights.

\paragraph*{Finite-sample results}
Figures~\ref{fig:supp-nonuniform-finite-avg} and
\ref{fig:supp-nonuniform-finite-max} further examine
the finite-sample behavior under nonuniform sampling.

For $E_{\mathrm{avg}}$, the Single-scale curves approach a
comparatively flat error level as the sample size
increases, whereas GBC continues to improve and
achieves smaller errors at the largest sample sizes
on all six manifolds. The maximum-error comparisons
show a similar overall pattern, although the gains
are less uniform, particularly on
$\mathbb S^3\subset\mathbb R^4$.

Taken together, these experiments provide numerical
evidence that nonuniform latent sampling does not
destroy the population second-to-fourth-order bias
reduction for STSP in the smooth, strictly positive
density settings considered here. The finite-sample
results further indicate that the improvement becomes
more pronounced as sampling variability decreases.
These findings support the distinction between the
conditions required by the general multiscale
framework and the uniform-sampling assumption used
to verify those conditions for specific fitting
procedures.

\subsection{Implementation details}

\paragraph*{Corrugated-torus geometry}
The corrugated torus used throughout the experiments is
\(\mathbb T_{\mathrm{corr}}^2=\{X(u,v):u,v\in[0,2\pi)\}\), where
\(\rho(v)=2+r(v)\cos v\) and
\[
\begin{aligned}
X(u,v)&=(\rho(v)\cos u,\rho(v)\sin u,r(v)\sin v)^\top,\\
r(v)&=0.7[1+0.18\cos(3v)+0.10\sin(2v)].
\end{aligned}
\]
Thus the major radius is \(2\), the baseline minor radius is \(0.7\),
and the two modulation terms have relative amplitudes \(0.18\), \(0.10\)
and frequencies \(3\), \(2\), respectively. With
\(r'(v)=0.7[-0.54\sin(3v)+0.20\cos(2v)]\), the surface-area element is
\(dA=J(v)\,du\,dv\), where
\(J(v)=\rho(v)\sqrt{r(v)^2+r'(v)^2}\). This is the area measure used
for the latent samples, quadrature, and area-based query grids below.

\subsubsection{Population experiments}
\label{app:e1-population}

The queries, quadrature resolutions, and Newton settings
in this subsection apply to the population experiments
of Appendix~\ref{app:population-verification}.
Experiment-specific changes are recorded in
Appendix~\ref{app:experiment-specific}.

\paragraph*{Operators and reference queries}
Population operators are evaluated deterministically, without Monte Carlo
sampling. STSP and the PCA search/matching steps of DR use
\(K_h(y,Y)=(1-\|Y-y\|^2/h^2)_+^6\), with weighted moments
\[
\begin{aligned}
\mu_h(y)&=\frac{\mathbb E[Y K_h(y,Y)]}{\mathbb E[K_h(y,Y)]},\\
\Sigma_h(y)&=\frac{\mathbb E[(Y-\mu_h(y))(Y-\mu_h(y))^\top K_h(y,Y)]}
                   {\mathbb E[K_h(y,Y)]}.
\end{aligned}
\]
The complement of the leading \(d\)-dimensional covariance eigenspace
gives all \(D-d\) estimated normal directions. DR retains its Gaussian
ridge equation; the compact-kernel PCA only guides local search and
matching. For each prescribed \(x\in\mathcal M\), a local fine-scale
solve gives \(z\in\mathcal M_{h_1}\). We recompute PCA at \(z\) using
exactly \(h_1\) and hold its orthonormal normal basis \(E\) fixed for
the coarse-scale solve. With \(t_1=0\), the output is
\(y_{\mathrm{GBC}}=z+E\sum_{k=1}^{2}w_kt_k\), where
\(\sum_kw_k=1\) and \(\sum_kw_kc(h_k,\sigma)=0\).
Single-scale at \(h_2\) has its own local solve. Prescribed base points provide
controlled branch initialization; automatic branch discovery is outside
scope, and subsequent matching uses neither true normals nor true
cross-scale correspondence.

Symmetry reduces STSP/DR queries to \((1,0,\ldots,0)\) on the circle
and spheres and \((1/\sqrt2,0,1/\sqrt2,0)\) on the Clifford torus.
On the corrugated torus, we fix \(u=0\) and use 128 equal-probability
quantiles of the surface-area density \(J(v)\) defined above. Its cumulative
distribution is evaluated on 65,536 intervals and inverted by linear
interpolation.

\paragraph*{YY construction}
YY starts from the prescribed normal graph \(X_0(x)=x+u_0(x)\), not a
preliminary manifold estimated from observations. For codimension-one
circles, spheres, and the corrugated torus, \(u_0=\sigma n\); the torus
normal is the normalized cross product of its parameter derivatives.
For the Clifford torus,
\(u_0=\sigma(1+0.2\cos u)(n_1+2n_2)/\sqrt5\), with
\(n_1=(\cos u,\sin u,0,0)^\top\) and
\(n_2=(0,0,\cos v,\sin v)^\top\).
For \(\mathbb S^2\subset\mathbb R^{10}\),
\(u_0=\sigma(1+0.2x_1)(n_1+n_2)/\sqrt2\), with
\(n_1=(x_1,x_2,x_3,0,\ldots,0)^\top\) and \(n_2=e_4\).
The axis is obtained by normalizing
\(F_{r_0}(X_0(x))-X_0(x)\), and the same center and axis are used at
every transverse scale. Spherical and transverse kernels use
\(K(s)=(1-s^2)_+^6\); the longitudinal profile is
\[
\eta(a)=
\begin{cases}
1,&0\le a\le1/2,\\
P(2a-1),&1/2<a<1,\\
0,&a\ge1,
\end{cases}
\qquad
P(s)=1-35s^4+84s^5-70s^6+20s^7.
\]
Axial weights satisfy \(\sum_kw_k^\parallel=1\) and
\(\sum_k(1+\alpha_k)w_k^\parallel=0\), using the solution with minimum Euclidean norm. In codimension one, \(w_k^\perp=w_k^\parallel\); otherwise,
complementary weights satisfy \(\sum_kw_k^\perp=1\),
\(\sum_k\alpha_kw_k^\perp=0\), and
\(\sum_k(1-\beta_k)w_k^\perp=0\).
The coefficients use the integrals in
Appendix~\ref{sec:ysl_common_moments}.

\paragraph*{YY population queries}
For \(\mathbb S^1\subset\mathbb R^2\), \(\mathbb S^2\subset\mathbb R^3\), and
\(\mathbb S^3\subset\mathbb R^4\), rotational symmetry reduces evaluation to
one representative base point, with average weight one, at both
numerical resolutions. For the Clifford torus, the symmetry-reduced
queries use \(L=32\) production phases and \(L=64\) refined phases,
\[
\phi_j=\frac{2\pi j}{L},\qquad j=0,\ldots,L-1,
\]
with equal average weights \(1/L\); the preliminary-offset amplitude
multiplying \(\sigma\) is \(1+0.2\cos\phi_j\).
For \(\mathbb S^2\subset\mathbb R^{10}\), the reduced coordinate \(x_1\) is
integrated using \(32\)- and \(64\)-point Gauss--Legendre rules on
\([-1,1]\), with weights divided by two and preliminary-offset amplitude
\(1+0.2x_1\). The two endpoints \(x_1=\pm1\) are additionally included
with zero average weight so that they contribute to \(E_\infty\), giving
\(34\) and \(66\) queries, respectively.

For the corrugated torus, rotational symmetry fixes \(u=0\).
Both production and refined calculations use the same \(128\) meridian
queries
\[
v_j=F_J^{-1}\!\left(\frac{j+1/2}{128}\right),
\qquad j=0,\ldots,127,
\]
where \(F_J\) is the normalized meridian area distribution defined in
the query-construction paragraph of Appendix~\ref{app:bandwidth-robustness}.
The base points are \(X(0,v_j)\), the preliminary points are
\(X(0,v_j)+\sigma n(0,v_j)\), and the average weights are \(1/128\).
These are deterministic area-coordinate mid-quantiles, not Monte Carlo
samples.

\paragraph*{Numerical integration}
Table~\ref{tab:e1-quadrature} records the latent-manifold
quadrature rules for the population experiments.
For STSP and matching PCA, Gaussian-noise integration is reduced to
noncentral-\(\chi\) radial moments, using Gauss--Legendre orders 160 and
224 for production and refined calculations, respectively. DR handles
Gaussian-noise integration by Gaussian convolution while retaining the
compact-kernel PCA calculation. YY uses the same two orders for radial
and longitudinal integrations. Compact-kernel calculations use fixed
\(12\sigma\) latent-node padding.

\begin{table}[t]
\centering
\small
\caption{Deterministic integration rules for the population
operators in Appendix~\ref{app:population-verification}.
``Reported resolution'' is used for the reported results,
while ``Refined check'' uses a denser grid to check
numerical resolution.}
\label{tab:e1-quadrature}
\setlength{\tabcolsep}{4pt}
\begin{tabular}{@{}p{0.20\linewidth}p{0.38\linewidth}p{0.16\linewidth}p{0.18\linewidth}@{}}
\toprule
Geometry & Integration scheme & Reported\newline resolution & Refined\newline check\\
\midrule
\(\mathbb S^1\)
& Equally spaced angular points
& \(8192\)
& \(16384\)\\

Clifford torus
& Equally spaced grid in the two periodic angles
& \(512\times512\)
& \(768\times768\)\\

\(\mathbb S^2\subset\mathbb R^3\) or \(\mathbb R^{10}\)
& Weighted latitude quadrature with symmetric tangential directions
& \(768\times4\)
& \(1024\times4\)\\

\(\mathbb S^3\)
& Weighted quadrature with symmetric tangential directions
& \(768\times6\)
& \(1024\times6\)\\

Corrugated torus
& Midpoint grid in \(u\) and equal-area points in \(v\)
& \(2560\times2560\)
& \(3840\times3840\)\\
\bottomrule
\end{tabular}
\end{table}

\paragraph*{Population root solves}
STSP and DR solve \(F(t)=B^\top G_h(z+Bt)=0\) from \(t_0=0\), where
\(B\) is the fixed orthonormal PCA normal basis and \(G_h\) is the
residual returned by the corresponding implementation. All norms are
Euclidean and all stopping tolerances are absolute. Table~\ref{tab:e1-population-newton}
lists the parameters. The Newton correction solves
\(D_tF(t)\Delta t=-F(t)\), with the Jacobian evaluated by central
differences. Success requires the projected-residual, full-residual, and
Newton-correction thresholds to hold simultaneously. At each update,
line search tests \(t+2^{-j}\Delta t\), \(j=0,\ldots,15\), and accepts
only a candidate inside the search region with strictly smaller projected
residual norm. The 16 candidates include the full Newton step, so at most
15 halvings are tried. Failure to find an acceptable candidate, or to
meet all stopping conditions by the update limit, is recorded as failure.

\begin{table}[!htbp]
\centering
\small
\caption{Newton settings for population STSP/DR solves
in Appendix~\ref{app:population-verification}.
The three stopping thresholds must all be satisfied.}
\label{tab:e1-population-newton}
\setlength{\tabcolsep}{5pt}
\begin{tabular}{@{}p{0.46\linewidth}p{0.48\linewidth}@{}}
\toprule
Quantity & Setting\\
\midrule
Search region & \(\|t\|\le4h^2\)\\
Maximum Newton updates & \(24\)\\
Line-search candidates per update & \(16\) (full step plus at most 15 halvings)\\
Projected residual & \(\|F(t)\|\le2\times10^{-12}\)\\
Full residual & \(\|G_h(z+Bt)\|\le2\times10^{-10}\)\\
Newton correction & \(\|\Delta t\|\le2\times10^{-11}\)\\
Central-difference step & \(\varepsilon=\max(10^{-6}h,5\times10^{-10})\)\\
\bottomrule
\end{tabular}
\end{table}

The normalized eigengap must also be at least \(10^{-7}\); for DR,
the largest eigenvalue of the normal Hessian must be strictly below
\(-10^{-8}\). These validity checks do not replace the three stopping
conditions. Nonpositive mass, an insufficient eigengap, a singular Jacobian,
exit from the search region, or nonfinite quantities also cause rejection.
No true-normal or enlarged-neighborhood fallback is used. YY does not use
this Newton solver.

\paragraph*{Distances and numerical checks}
Errors use distance to the continuous target, not to a discretized
reference cloud. For an embedded unit sphere and the Clifford torus,
respectively,
\[
\begin{aligned}
d(y,S^d)
&=\left[(\|y_{1:d+1}\|-1)^2+\|y_{d+2:D}\|^2\right]^{1/2},\\
d(y,\mathbb T_{\mathrm C}^2)
&=\left[
\left(\sqrt{y_1^2+y_2^2}-\frac1{\sqrt2}\right)^2+
\left(\sqrt{y_3^2+y_4^2}-\frac1{\sqrt2}\right)^2
\right]^{1/2}.
\end{aligned}
\]
The sphere's second term is omitted when \(D=d+1\). For the corrugated
torus, an 8192-point grid initializes a continuous closest-point search,
followed by bounded scalar minimization and least-squares refinement at
tolerance \(10^{-14}\). STSP/DR use equal query weights; YY uses its
quadrature weights. The corrugated-torus STSP/DR maximum is over 128
queries, not a certified continuous supremum. For the population STSP/DR experiments, both the
production--refined discrepancy and final Newton
correction must be below $10^{-3}$ times the smallest
GBC error in the six-manifold audit.
These checks assess numerical resolution and convergence, not certified
bounds on total numerical error. The separate YY criterion is defined
below.

\paragraph*{YY numerical-convergence criterion}
Let \(\mathcal G\) denote the six geometries and \(\Sigma\) the prescribed
seven noise levels. Define
\[
m_*=
\min_{\substack{g\in\mathcal G\\ \sigma\in\Sigma}}
E_{\infty,\mathrm{prod}}^{\mathrm{GBC\text{-}YY}}(g,\sigma)
\]
and
\[
\Delta_*=
\max_{\substack{
g\in\mathcal G,\ \sigma\in\Sigma\\
A\in\{\mathrm{Single\text{-}scale},\mathrm{GBC\text{-}YY}\}\\
p\in\{\mathrm{avg},\infty\}
}}
\left|
E_{p,\mathrm{prod}}^A(g,\sigma)
-
E_{p,\mathrm{refined}}^A(g,\sigma)
\right|.
\]
In addition to the original per-package checks, the final six-manifold
validation requires
\(
\frac{\Delta_*}{m_*}<10^{-2}.
\)
Thus, the average and maximum errors of both Single-scale and GBC-YY must
satisfy the criterion under one common denominator: the smallest
production GBC-YY maximum error over all geometry--noise configurations.
This is not a relative difference normalized separately by each
reported metric.

For the archived results,
\[
\begin{aligned}
m_*&=2.3845222241\times10^{-8},\qquad
\Delta_*=4.7628051849\times10^{-11},\\
\frac{\Delta_*}{m_*}&=1.9973834325\times10^{-3}.
\end{aligned}
\]
In particular, the corrugated-torus comparison
refines the integration rule while retaining the same \(128\) query
locations; it does not certify the continuous supremum by query-grid
refinement.

\subsubsection{Finite-sample experiments}
\label{app:e2-finite-sample}

\paragraph*{Data and base locations}
On the same six manifolds, we generate \(Y_i=X_i+\xi_i\), with
\(X_i\) sampled from uniform surface measure and
\(\xi_i\sim N(0,\sigma^2 I_D)\). We fix \(\sigma=0.005\) and use
\(n\in\{2000,\allowbreak 5000,\allowbreak 10000,\allowbreak 20000,\allowbreak 50000,\allowbreak 100000,\allowbreak 200000,\allowbreak 500000\}\), with
20 repetitions. In each repetition, one sample of \(500{,}000\)
observations is generated, and its first \(n\) observations are used
at each smaller sample size. The 144 base locations on the true manifold
are fixed across sample sizes, repetitions, and methods: equally spaced
angles on the circle, a \(12\times12\) angular grid on the Clifford
torus, and normalized Gaussian vectors with a fixed seed on the spheres.
The corrugated torus uses a \(12\times12\) grid with equally spaced
\(u\) and area-based quantiles in \(v\); its clean observations use
the same surface-area distribution. Base locations provide controlled
initializations, not fitted query points or automatic branch discovery.

\paragraph*{Fitting and matching}
STSP and DR use \(h_1=0.1\), \(h_2=0.2\) and empirical versions of
their population operators, with expectations replaced by sample averages.
Each base point initializes a local \(h_1\)-scale solve to obtain the
fitted reference query. The compact-kernel PCA normal space is recomputed
there at \(h_1\) and held fixed for coarse-scale matching; DR retains
its Gaussian ridge residual. GBC follows
Algorithm~\ref{alg:normal_guided_ms} with the theoretical coefficients
\(c(h_k,\sigma)\), not coefficients fitted to observed errors.
Single-scale is computed separately: STSP uses its own fine-reference solve
before the coarse solve, whereas DR starts at the base point with its own
\(h_2\)-PCA. A GBC failure does not discard a valid Single-scale output.
YY maps each base location to the prescribed preliminary manifold in
Appendix~\ref{app:e1-population} and replaces population expectations by
empirical weighted averages. It estimates the axis from
\(\widehat F_{r_0}(X_0(x))-X_0(x)\), using \(r_0=32\sigma\),
\(r_2=12\sigma\sqrt{\log(1/\sigma)}\), and transverse radii
\((r_{1,1},r_{1,2},r_{1,3})=(0.8,1,1.6)40\sigma\).
All scales share the preliminary center and empirical axis. Single-scale uses
the central radius; GBC-YY uses the directional constraints above.
A failed side-radius calculation does not discard a valid central Single-scale.

\paragraph*{Numerical solvers and failures}
Empirical operators are evaluated directly from observations in double
precision, without latent-manifold quadrature. STSP and YY coefficient integrals use adaptive quadrature with absolute and relative tolerances
\(2\times10^{-13}\); weight constraints are checked at \(10^{-10}\).
STSP and DR use damped Newton iteration in fixed PCA coordinates.
STSP allows \(\|t\|\le4h_2^2\) for both reference and coarse solves,
50 Newton steps, and 20 backtracking trials, with projected and full
residual tolerances \(10^{-11}\) and \(10^{-10}\), respectively.
DR uses \(\|t\|\le4h_1^2\) for the fine solve and
\(\|t\|\le4h_2^2\) for the coarse solve, at most 30 Newton steps,
16 backtracking trials, and projected residual tolerance \(10^{-10}\). Insufficient support, spectral separation,
or normal transversality, a Jacobian condition number above \(10^{12}\),
exit from the search region, failed line search, or nonfinite quantities
cause rejection. DR Gaussian sums are truncated at \(8h\), requiring
an omitted-weight upper bound of at most \(10^{-6}\). YY instead requires
\(N_{\mathrm{eff}}=(\sum_i a_i)^2/\sum_i a_i^2\ge8\) for local
nonnegative weights \(a_i\) and first-stage displacement at least
\(0.1\sigma\). Empty neighborhoods, invalid mass, inadequate effective
support, and invalid first-stage signals are failures. No failure is
replaced by true-normal matching, a larger neighborhood, or an alternative
numerical procedure.

\paragraph*{Error aggregation and reliability}
For estimator \(A\) and repetition \(r\), let \(\mathcal S_r^A\)
be the successful query set. For \(\mathcal S_r^A\ne\varnothing\),
we record \(E_{\mathrm{avg},r}^A=|\mathcal S_r^A|^{-1}
\sum_{j\in\mathcal S_r^A}d(\widehat z_{r,j}^A,\mathcal M)\) and
\(E_{\infty,r}^A=\max_{j\in\mathcal S_r^A}
 d(\widehat z_{r,j}^A,\mathcal M)\), using the continuous-target distances above.
All 144 locations have equal weight, including for YY. Each estimator
retains its own successes; failed queries are not imputed, and a repetition
is complete only when all 144 succeed. Medians and the 25th and 75th percentiles
use repetitions with finite errors, including partial successes; the bands
are not confidence intervals. Filled markers require at most \(1\%\)
aggregate query failures and at least 18 complete repetitions; other finite
results use hollow markers. An empty success set leaves the errors
undefined. 

\subsubsection{Experiment-specific modifications}
\label{app:experiment-specific}

The remaining experiments retain the applicable population or empirical
operators above. Only their additional settings and departures from the
shared protocols are listed here.

\paragraph*{Section \ref{sec:e3-bandwidth-robustness}: Scale Sensitivity and Reverse Coverage}
The geometry-specific budgets, fitting radii, forward/reverse query counts,
and reliability rules are given in Appendix~\ref{app:bandwidth-robustness}.
The common grid for \(h_2\) in STSP/DR and \(r_1\) in YY is
\(\mathcal H=\{0.06,\allowbreak 0.075,\allowbreak 0.09,\allowbreak 0.11,\allowbreak 0.135,\allowbreak 0.165,\allowbreak 0.20,\allowbreak 0.245,\allowbreak 0.30,\allowbreak 0.37,\allowbreak 0.45,\allowbreak 0.55,\allowbreak 0.67,\allowbreak 0.82,\allowbreak 1.0,\allowbreak 1.2,\allowbreak 1.4,\allowbreak 1.6,\allowbreak 1.8,\allowbreak 2.0,\allowbreak 2.25,\allowbreak 2.5,\allowbreak 2.75,\allowbreak 3.0,\allowbreak 3.25,\allowbreak 3.5,\allowbreak 4.0,\allowbreak 5.0\}\).
Single-scale is retained at both STSP/DR scales and all three YY radii.
The reverse metrics \(E_{\mathrm{rev,avg}}^A\) and
\(E_{\mathrm{rev},\infty}^A\) are the estimator-specific
\(\overline E_{\mathrm{rev}}\) and \(E_{\mathrm{rev}}\) defined there;
reverse coverage uses repetition zero rather than the 20-repetition
forward summary.

\paragraph*{Section \ref{sec:enf_cp_comparison}: Correspondence and Candidate Density}
This experiment uses the closed curve
\(\gamma(\theta)=(1.5\cos\theta,\sin\theta,
0.3\sin2\theta)^\top\), \(0\le\theta<2\pi\), with
\(n=200000\), \(\sigma=0.06\), \(h_1=0.12\), \(h_2=0.24\),
32 fixed reference queries, and 10 repetitions. Oracle, ENF, and NN share
fitted reference points and weights. ENF uses \(h_1\)-PCA normal planes;
Oracle uses the prescribed true planes. All local solves use the empirical STSP residual
\[
\widehat R_h(y)
=
\widehat Q_h(y)
\bigl(y-\widehat\mu_h(y)\bigr),
\]
as defined in Section~\ref{sec:stsp_validation}.
A numerical root is accepted only if both residual conditions
\[
\bigl\|B^\top\widehat R_h(y)\bigr\|_2\le10^{-11},
\qquad
\bigl\|\widehat R_h(y)\bigr\|_2\le10^{-9}
\]
are satisfied, where \(B\) is the orthonormal basis of the corresponding
search fiber. These are absolute Euclidean residual tolerances.
Acceptance additionally requires passing the support, normalized-mass,
eigengap, transversality, and fiber-Jacobian checks, while remaining
inside the prescribed search tube. Failure of either residual check
causes rejection; in particular, convergence of the projected residual
alone does not constitute a successful root.
NN candidate counts are
\(M_{\mathrm{match}}\in\{32,\allowbreak 36,\allowbreak 40,\allowbreak 44,\allowbreak 48,\allowbreak 52,\allowbreak 56,\allowbreak 60,\allowbreak 64,\allowbreak 80,\allowbreak 96,\allowbreak 112,\allowbreak 128,\allowbreak 256,\allowbreak 512\}\).
Starting from normalized-arclength locations
\(s_j=((j+0.74)/32)\bmod1\), \(j=0,\ldots,31\), nested dyadic
refinement adds candidate locations. An empirical \(h_2\)-solve at each
location gives a fitted candidate, and NN selects the closest successful
candidate in Euclidean distance. True-geometry fibers control the reference
and candidate solves, not the PCA-based ENF matching step; candidate
failures are retained. Total NN time includes reference computation,
candidate fitting, index construction, search, and extrapolation, but
excludes truth-distance evaluation, diagnostics, plotting, file I/O, and
GPU warm-up. The comparison concerns discrete NN, not continuous
closest-point projection.

\paragraph*{Section \ref{sec:exp-directional}: Directional-Weight Ablation}
We use \[\gamma_4(\theta)=(1.5\cos\theta,\sin\theta,
0.3\sin(2\theta),0)^\top,\] \(0\le\theta<2\pi\).
Primes denote derivatives with respect to \(\theta\). Define the unit
tangent \(T_4=\gamma_4'/\|\gamma_4'\|\) and principal normal
\(n_4=(I_4-T_4T_4^\top)\gamma_4''/
\|(I_4-T_4T_4^\top)\gamma_4''\|\). With
\(e_4=(0,0,0,1)^\top\), set
\(a_5=(n_4+e_4)/\sqrt2\) and \(b_5=(n_4-e_4)/\sqrt2\).
Both the population and finite-sample experiments use
\(X_{0,5}(\theta)=\gamma_4(\theta)+u_{0,5}(\theta)\), with
\[
u_{0,5}(\theta)
=\sigma[1+0.1\cos(3\theta)]a_5(\theta)
+\tfrac14\sigma^2\sin(2\theta)b_5(\theta).
\]
The leading offset has amplitude between \(0.9\sigma\) and
\(1.1\sigma\); the signed second-order term in the orthogonal normal
direction is retained. These geometric directions specify a controlled
initialization, not an observation-based preliminary estimate or an imposed
YY axis. At each preliminary point \(z\), the actual first-stage axis is
\(e(z)=(F_{r_0}(z)-z)/\|F_{r_0}(z)-z\|\), with
\(\widehat F_{r_0}\) used in finite-sample runs.

All estimators share the preliminary points, computed first-stage axes,
and three fitted outputs.
For scale-specific displacements \(v_k\) and first-stage axis projector
\(P\), Scalar-A uses \(v_{\mathrm{Scalar-A}}=\sum_kw_k^\parallel v_k\),
Scalar-C uses \(v_{\mathrm{Scalar-C}}=\sum_kw_k^\perp v_k\), and
directional GBC uses
\(v_{\mathrm{GBC\text{-}YY}}=\sum_kw_k^\parallel Pv_k+
\sum_kw_k^\perp(I-P)v_k\). Single-scale uses \(r_1=5\sigma\).
The population experiment uses 64 queries, the seven noise levels in
Appendix~\ref{app:population-verification}, and radii \((1.25,5,8)\sigma\).
The finite-sample experiment uses $\sigma=0.05$, radii
$(5,6,12)\sigma$, $n\in\{240000,\allowbreak 480000,\allowbreak 960000\}$,
20 paired repetitions, and 256 queries. Both experiments use
$r_0=3\sigma$ and $r_2=12\sigma\sqrt{\log(1/\sigma)}$.

The query points are deterministically selected to be equally
spaced in normalized arclength along the true curve.
Define
\[
F_\ell(\theta)
=
\frac{\int_0^\theta \|\gamma_4'(t)\|_2\,dt}
     {\int_0^{2\pi}\|\gamma_4'(t)\|_2\,dt}.
\]
We set
\[
\theta_j
=
F_\ell^{-1}\left(\frac{j+0.37}{N_q}\right),
\qquad j=0,\ldots,N_q-1,
\]
with $N_q=64$ for the population experiment and
$N_q=256$ for the finite-sample experiment.
The true base points are $x_j=\gamma_4(\theta_j)$,
and the preliminary queries are $z_j=X_{0,5}(\theta_j)$.
Within each experiment, all four ablation estimators
use the same fixed query points across noise levels,
sample sizes, and repetitions.

Numerically, \(F_\ell\) is tabulated by cumulative Simpson
integration on \(131073\) equally spaced parameter nodes
over \([0,2\pi]\), including both endpoints.
Its inverse is evaluated by linear interpolation.
The finite-sample radii
\((5,6,12)\sigma\) were chosen based on preliminary experiments and kept
fixed for the 20 independent evaluation repetitions. The corresponding
complementary-moment matrix has a condition number of approximately
\(244.94\). For this finite-sample configuration, let
\(a_k=w_k^\parallel\) and \(b_k=w_k^\perp\) denote the axial and
complementary directional-extrapolation weights, respectively.
Their \(\ell^1\) amplification factors are
\[
W_\parallel:=\sum_{k=1}^{3}|a_k|\approx1.7714,
\qquad
W_\perp:=\sum_{k=1}^{3}|b_k|\approx13.3207,
\]
so the total weight bound is
\[
W:=W_\parallel+W_\perp\approx15.0921.
\]
These quantities describe the potential amplification of perturbations
in the common-start single-scale estimates. They are distinct from
the complementary-moment matrix condition number, which measures the
sensitivity of the weight computation. The displayed values are
rounded, while the implementation uses full-precision weights. Weight-constraint residuals must be at
most \(10^{-10}\); failed queries were recorded, and successful outliers were retained.

\paragraph*{Directional error metrics}
Let \(x_j=\gamma_4(\theta_j)\) be the original noiseless base point,
\(z_j=X_{0,5}(\theta_j)\) the common preliminary point, and
\(y_j^A\) the output of estimator
\begin{center}
\(A\in\{\mathrm{Single\text{-}scale},\mathrm{Scalar\text{-}A},
\mathrm{Scalar\text{-}C},\mathrm{GBC\text{-}YY}\}\).
\end{center}
Define the unit tangent and true normal projector at \(x_j\) by
\[
T_j=\frac{\gamma_4'(\theta_j)}{\|\gamma_4'(\theta_j)\|},
\qquad Q_j=I_4-T_jT_j^\top.
\]
Let \(e_j\) be the actual computed first-stage unit axis at \(z_j\):
\[
e_j=\frac{F_{r_0}(z_j)-z_j}{\|F_{r_0}(z_j)-z_j\|},
\]
with the empirical first-stage map used in finite-sample runs.
For error evaluation, we project this axis onto the true normal
space and normalize it:
\[
\nu_j=\frac{Q_je_j}{\|Q_je_j\|},\qquad
\widetilde P_j=\nu_j\nu_j^\top,\qquad
\widetilde S_j=Q_j-\widetilde P_j.
\]
A vanishing normal component of the axis is treated as a failure.
The axial and complementary normal errors are
\[
\begin{aligned}
E_{\parallel}^{A}
&=\max_{1\le j\le N_q}
  \bigl\|\widetilde P_j(y_j^A-x_j)\bigr\|,\\
E_{\perp}^{A}
&=\max_{1\le j\le N_q}
  \bigl\|\widetilde S_j(y_j^A-x_j)\bigr\|.
\end{aligned}
\]
Table~\ref{tab:e5-population-slopes} and
Figure~\ref{fig:e5-directional-population} here report these finite-query
maxima with \(N_q=64\); all queries succeeded in the archived
population runs.

These diagnostics use the original base points, not the
nearest-point reparameterization of the outputs. The projected
axis and normal projectors above are used only for error evaluation.
In particular, \(\widetilde S_j\) removes tangential components,
whereas the ambient complement \(I_4-U_{0,j}\), with
\(U_{0,j}=e_je_j^\top\), is used in constructing the directional
estimator and may contain tangential directions.
The total geometric error is evaluated separately as
\[
E_{\mathrm{geom}}^{A}
=\max_{1\le j\le N_q}\operatorname{dist}(y_j^A,\mathcal M),
\]
using numerical nearest-point minimization on the continuous curve.
The componentwise maxima do not form a Pythagorean decomposition
of \(E_{\mathrm{geom}}^{A}\).

\subsection{Implementation details for Figure~\ref{fig:intro_stsp_population}}
Figure~\ref{fig:intro_stsp_population} compares our multiscale
STSP with single-scale STSP and the MS-STSP method of
\citet{xia2026stsp}. The geometries, noise levels, query points, numerical
quadrature, and error evaluation follow the population
experiments in Appendix~\ref{app:population-verification}.
Here we describe the implementation of the MS-STSP baseline. 

Our implementation follows the MS-STSP algorithm of
\citet{xia2026stsp}, with only numerical adaptations
to the population-level evaluation, radius selection,
and convergence criteria for consistency with the
population-verification experiments.
We implement population MS-STSP using three radii
$h_\ell\in\{3\sigma,4.5\sigma,6\sigma\}$.
Let $m_h(z)=\mathbb E[Y\mid Y\in B(z,h)]$ denote the
hard-ball conditional mean, and let $Q_h^{\mathrm{MS}}(z)$ be the projector
onto the eigenspace associated with the smallest $D-d$
eigenvalues of the corresponding centered conditional covariance.
For each query $x_j$, all three branches start from
$z_{\ell,j}^{(0)}=x_j$ and follow the updates
\[
z_{\ell,j}^{(k+1)}
=z_{\ell,j}^{(k)}
+Q_{h_\ell}^{\mathrm{MS}}(z_{\ell,j}^{(k)})
\bigl[m_{h_\ell}(z_{\ell,j}^{(k)})
-z_{\ell,j}^{(k)}\bigr].
\]
The noisy observation distribution is held fixed throughout
the iterations. Within each geometry--noise--quadrature configuration,
all branches use a common terminal iteration $K$,
determined by the residual tolerances specified in
Appendix~\ref{app:e1-population}, a maximum of
256 iterations, and an additional eight-step convergence
check against the numerical accuracy criterion in
the same section.

The final MS-STSP estimate is the intercept of an unweighted,
coordinatewise linear regression of $z_{\ell,j}^{(K)}$
on $h_\ell^2$, namely,
\[
y_j^{\mathrm{MS}}
=\sum_{\ell=1}^{3}w_\ell z_{\ell,j}^{(K)},
\qquad
(w_1,w_2,w_3)=\frac{1}{218}(237,92,-111).
\]
These weights satisfy $\sum_\ell w_\ell=1$ and
$\sum_\ell w_\ell h_\ell^2=0$.
Unlike the smooth-kernel STSP used for Single-scale and GBC
in the population-verification experiments, the MS-STSP
baseline uses hard-ball conditional moments.

\FloatBarrier
\endgroup

\FloatBarrier
\bibliographystyle{plainnat}
\bibliography{references}

@article{aamari2019nonasymptotic,
  title         = {Nonasymptotic rates for manifold, tangent space and curvature estimation},
  author        = {Aamari, Eddie and Levrard, Cl{\'e}ment},
  journal       = {The Annals of Statistics},
  volume        = {47},
  number        = {1},
  pages         = {177--204},
  year          = {2019},
  publisher     = {Institute of Mathematical Statistics}
}

@article{aizenbud2025localgeo,
  title         = {Estimation of Local Geometric Structure on Manifolds from Noisy Data},
  author        = {Aizenbud, Yariv and Sober, Barak},
  journal       = {Journal of Machine Learning Research},
  volume        = {26},
  number        = {64},
  pages         = {1--89},
  year          = {2025},
  url           = {https://jmlr.org/papers/v26/25-0183.html}
}

@article{bartlett1963statistical,
  title         = {Statistical estimation of density functions},
  author        = {Bartlett, M. S.},
  journal       = {Sankhy{\=a} Series A},
  volume        = {25},
  pages         = {245--254},
  year          = {1963}
}

@article{belkin2003laplacian,
  title         = {Laplacian eigenmaps for dimensionality reduction and data representation},
  author        = {Belkin, Mikhail and Niyogi, Partha},
  journal       = {Neural Computation},
  volume        = {15},
  number        = {6},
  pages         = {1373--1396},
  year          = {2003},
  publisher     = {MIT Press}
}

@inproceedings{boissonnat2009manifold,
  title         = {Manifold reconstruction in arbitrary dimensions using witness complexes},
  author        = {Boissonnat, Jean-Daniel and Guibas, Leonidas J. and Oudot, Steve Y.},
  booktitle     = {Proceedings of the Twenty-Third Annual Symposium on Computational Geometry},
  pages         = {194--203},
  year          = {2007},
  publisher     = {Association for Computing Machinery},
  doi           = {10.1145/1247069.1247106}
}

@article{chen2015asymptotic,
  title         = {Asymptotic theory for density ridges},
  author        = {Chen, Yen-Chi and Genovese, Christopher R. and Wasserman, Larry},
  journal       = {The Annals of Statistics},
  volume        = {43},
  pages         = {1896--1928},
  year          = {2015}
}

@inproceedings{cheng2005manifold,
  title         = {Manifold Reconstruction from Point Samples},
  author        = {Cheng, Siu-Wing and Dey, Tamal K. and Ramos, Edgar A.},
  booktitle     = {Proceedings of the Sixteenth Annual ACM-SIAM Symposium on Discrete Algorithms},
  pages         = {1018--1027},
  year          = {2005},
  publisher     = {Society for Industrial and Applied Mathematics}
}

@article{coifman2006diffusion,
  title         = {Diffusion maps},
  author        = {Coifman, Ronald R. and Lafon, St{\'e}phane},
  journal       = {Applied and Computational Harmonic Analysis},
  volume        = {21},
  pages         = {5--30},
  year          = {2006}
}

@article{delaigle2009design,
  title         = {A Design-Adaptive Local Polynomial Estimator for the Errors-in-Variables Problem},
  author        = {Delaigle, Aurore and Fan, Jianqing and Carroll, Raymond J.},
  journal       = {Journal of the American Statistical Association},
  volume        = {104},
  number        = {485},
  pages         = {348--359},
  year          = {2009},
  doi           = {10.1198/jasa.2009.0114}
}

@article{fan1993local,
  title         = {Local linear regression smoothers and their minimax efficiencies},
  author        = {Fan, Jianqing},
  journal       = {The Annals of Statistics},
  volume        = {21},
  pages         = {196--216},
  year          = {1993}
}

@article{fan1993nonparametric,
  title         = {Nonparametric Regression with Errors in Variables},
  author        = {Fan, Jianqing and Truong, Young K.},
  journal       = {The Annals of Statistics},
  volume        = {21},
  number        = {4},
  pages         = {1900--1925},
  year          = {1993},
  doi           = {10.1214/aos/1176349402}
}

@article{federer1959curvature,
  title         = {Curvature measures},
  author        = {Federer, Herbert},
  journal       = {Transactions of the American Mathematical Society},
  volume        = {93},
  pages         = {418--491},
  year          = {1959}
}

@article{fefferman2016testing,
  title         = {Testing the manifold hypothesis},
  author        = {Fefferman, Charles and Mitter, Sanjoy and Narayanan, Hariharan},
  journal       = {Journal of the American Mathematical Society},
  volume        = {29},
  number        = {4},
  pages         = {983--1049},
  year          = {2016},
  publisher     = {American Mathematical Society},
  doi           = {10.1090/jams/852}
}

@inproceedings{fefferman2018fitting,
  title         = {Fitting a putative manifold to noisy data},
  author        = {Fefferman, Charles and Ivanov, Sergei and Kurylev, Yaroslav and Lassas, Matti and Narayanan, Hariharan},
  booktitle     = {Proceedings of the 31st Conference on Learning Theory},
  series        = {Proceedings of Machine Learning Research},
  volume        = {75},
  pages         = {688--720},
  year          = {2018},
  publisher     = {PMLR},
  url           = {https://proceedings.mlr.press/v75/fefferman18a.html}
}

@article{fefferman2020reconstruction,
  title         = {Reconstruction and interpolation of manifolds {I}: The geometric {Whitney} problem},
  author        = {Fefferman, Charles and Ivanov, Sergei and Kurylev, Yaroslav and Lassas, Matti and Narayanan, Hariharan},
  journal       = {Foundations of Computational Mathematics},
  volume        = {20},
  pages         = {1035--1133},
  year          = {2020}
}

@article{genovese2012manifold,
  title         = {Manifold estimation and singular deconvolution under {Hausdorff} loss},
  author        = {Genovese, Christopher R. and Perone-Pacifico, Marco and Verdinelli, Isabella and Wasserman, Larry},
  journal       = {The Annals of Statistics},
  volume        = {40},
  number        = {2},
  pages         = {941--963},
  year          = {2012},
  doi           = {10.1214/12-AOS994}
}

@article{genovese2012minimax,
  title         = {Minimax manifold estimation},
  author        = {Genovese, Christopher R. and Perone-Pacifico, Marco and Verdinelli, Isabella and Wasserman, Larry},
  journal       = {Journal of Machine Learning Research},
  volume        = {13},
  number        = {43},
  pages         = {1263--1291},
  year          = {2012},
  url           = {https://jmlr.org/papers/v13/genovese12a.html}
}

@article{genovese2014nonparametric,
  title         = {Nonparametric Ridge Estimation},
  author        = {Genovese, Christopher R. and Perone-Pacifico, Marco and Verdinelli, Isabella and Wasserman, Larry},
  journal       = {The Annals of Statistics},
  volume        = {42},
  number        = {4},
  pages         = {1511--1545},
  year          = {2014},
  doi           = {10.1214/14-AOS1218}
}

@article{hastie1989principal,
  title         = {Principal curves},
  author        = {Hastie, Trevor and Stuetzle, Werner},
  journal       = {Journal of the American Statistical Association},
  volume        = {84},
  number        = {406},
  pages         = {502--516},
  year          = {1989},
  publisher     = {Taylor \& Francis}
}

@article{jones1995simple,
  title         = {A simple bias reduction method for density estimation},
  author        = {Jones, M. C. and Linton, Oliver and Nielsen, Jens Perch},
  journal       = {Biometrika},
  volume        = {82},
  pages         = {327--338},
  year          = {1995}
}

@article{lee1980two,
  title         = {Two algorithms for constructing a {Delaunay} triangulation},
  author        = {Lee, D. T. and Schachter, Bruce J.},
  journal       = {International Journal of Computer and Information Sciences},
  volume        = {9},
  pages         = {219--242},
  year          = {1980}
}

@article{li2026curvature,
  title         = {Curvature-driven manifold fitting under unbounded isotropic noise},
  author        = {Li, Ruowei and Yao, Zhigang},
  journal       = {arXiv preprint arXiv:2601.10133},
  year          = {2026}
}

@article{mcinnes2018umap,
  title         = {{UMAP}: Uniform Manifold Approximation and Projection for Dimension Reduction},
  author        = {McInnes, Leland and Healy, John and Melville, James},
  journal       = {arXiv preprint arXiv:1802.03426v3},
  year          = {2020},
  url           = {https://arxiv.org/abs/1802.03426v3}
}

@article{mohammed2017learning,
  title         = {Manifold learning using kernel density estimation and local principal components analysis},
  author        = {Mohammed, Kitty and Narayanan, Hariharan},
  journal       = {arXiv preprint arXiv:1709.03615},
  year          = {2017},
  url           = {https://arxiv.org/abs/1709.03615}
}

@article{niyogi2008finding,
  title         = {Finding the homology of submanifolds with high confidence from random samples},
  author        = {Niyogi, Partha and Smale, Stephen and Weinberger, Shmuel},
  journal       = {Discrete \& Computational Geometry},
  volume        = {39},
  pages         = {419--441},
  year          = {2008}
}

@article{ozertem2011finding,
  title         = {Locally Defined Principal Curves and Surfaces},
  author        = {Ozertem, Umut and Erdogmus, Deniz},
  journal       = {Journal of Machine Learning Research},
  volume        = {12},
  number        = {34},
  pages         = {1249--1286},
  year          = {2011},
  url           = {https://jmlr.org/papers/v12/ozertem11a.html}
}

@article{panaretos2014principal,
  title         = {Principal flows},
  author        = {Panaretos, Victor M. and Pham, Tung and Yao, Zhigang},
  journal       = {Journal of the American Statistical Association},
  volume        = {109},
  number        = {505},
  pages         = {424--436},
  year          = {2014},
  publisher     = {Taylor \& Francis}
}

@article{qiao2016theoretical,
  title         = {Theoretical analysis of nonparametric filament estimation},
  author        = {Qiao, Wanli and Polonik, Wolfgang},
  journal       = {The Annals of Statistics},
  volume        = {44},
  pages         = {1269--1297},
  year          = {2016},
  publisher     = {Institute of Mathematical Statistics}
}

@article{roweis2000nonlinear,
  title         = {Nonlinear dimensionality reduction by locally linear embedding},
  author        = {Roweis, Sam T. and Saul, Lawrence K.},
  journal       = {Science},
  volume        = {290},
  number        = {5500},
  pages         = {2323--2326},
  year          = {2000},
  publisher     = {American Association for the Advancement of Science}
}

@article{ruppert1994multivariate,
  title         = {Multivariate locally weighted least squares regression},
  author        = {Ruppert, David and Wand, M. P.},
  journal       = {The Annals of Statistics},
  volume        = {22},
  pages         = {1346--1370},
  year          = {1994}
}

@article{schucany1977improvement,
  title         = {Improvement of Kernel Type Density Estimators},
  author        = {Schucany, W. R. and Sommers, John P.},
  journal       = {Journal of the American Statistical Association},
  volume        = {72},
  number        = {358},
  pages         = {420--423},
  year          = {1977},
  doi           = {10.1080/01621459.1977.10481012}
}

@article{sober2020manifold,
  title         = {Manifold approximation by moving least-squares projection},
  author        = {Sober, Barak and Levin, David},
  journal       = {Constructive Approximation},
  volume        = {52},
  pages         = {433--478},
  year          = {2020}
}

@article{su2025principal,
  title   = {Principal Decomposition with Nested Submanifolds},
  author  = {Su, Jiaji and Yao, Zhigang},
  journal = {Biometrika},
  year    = {2026},
  doi     = {10.1093/biomet/asag071},
  url     = {https://doi.org/10.1093/biomet/asag071},
  note    = {Advance online publication}
}

@article{tenenbaum2000global,
  title         = {A global geometric framework for nonlinear dimensionality reduction},
  author        = {Tenenbaum, Joshua B. and de Silva, Vin and Langford, John C.},
  journal       = {Science},
  volume        = {290},
  number        = {5500},
  pages         = {2319--2323},
  year          = {2000},
  publisher     = {American Association for the Advancement of Science}
}

@article{tibshirani1992principal,
  title         = {Principal curves revisited},
  author        = {Tibshirani, Robert},
  journal       = {Statistics and Computing},
  volume        = {2},
  pages         = {183--190},
  year          = {1992}
}

@inproceedings{wang2010manifold,
  title         = {Manifold blurring mean shift algorithms for manifold denoising},
  author        = {Wang, Weiran and Carreira-Perpi{\~n}{\'a}n, Miguel {\'A}.},
  booktitle     = {Proceedings of the IEEE Conference on Computer Vision and Pattern Recognition},
  pages         = {1759--1766},
  year          = {2010}
}

@misc{xia2026stsp,
  title         = {Manifold Fitting by Successive Tangent-Space Projection},
  author        = {Xia, Yuqing and Li, Bingjie and Yao, Zhigang},
  year          = {2026},
  eprint        = {2610.09423},
  archiveprefix = {arXiv},
  primaryclass  = {stat.ME},
  howpublished  = {arXiv preprint arXiv:2610.09423},
  url           = {https://arxiv.org/abs/2610.09423}
}

@article{yao2020principal,
  title   = {Principal Boundary on {Riemannian} Manifolds},
  author  = {Yao, Zhigang and Zhang, Zhenyue},
  journal = {Journal of the American Statistical Association},
  volume  = {115},
  number  = {531},
  pages   = {1435--1448},
  year    = {2020},
  doi     = {10.1080/01621459.2019.1610660}
}

@article{yao2023manifold,
  title         = {Manifold fitting},
  author        = {Yao, Zhigang and Su, Jiaji and Li, Bingjie and Yau, Shing-Tung},
  journal       = {arXiv preprint arXiv:2304.07680},
  year          = {2023}
}

@article{yao2025unbounded,
  title         = {Manifold fitting under unbounded noise},
  author        = {Yao, Zhigang and Xia, Yuqing},
  journal       = {Journal of Machine Learning Research},
  volume        = {26},
  number        = {45},
  pages         = {1--55},
  year          = {2025},
  url           = {https://jmlr.org/papers/v26/21-0039.html}
}

@article{yao2026review,
  title         = {Manifold fitting: A review of methods and applications},
  author        = {Yao, Zhigang and Su, Jiaji},
  journal       = {Harvard Data Science Review},
  volume        = {8},
  number        = {3},
  year          = {2026},
  doi           = {10.1162/99608f92.686fa68c}
}

@article{zhang2004principal,
  title         = {Principal manifolds and nonlinear dimensionality reduction via tangent space alignment},
  author        = {Zhang, Zhenyue and Zha, Hongyuan},
  journal       = {SIAM Journal on Scientific Computing},
  volume        = {26},
  number        = {1},
  pages         = {313--338},
  year          = {2004}
}
\end{document}